\UseRawInputEncoding 

\documentclass{article}

\newcommand{\mS}{\mathbb{S}}
\newcommand{\mI}{\mathbb{I}}
\newcommand{\mT}{\mathbb{T}}
\newcommand{\mU}{\mathbb{U}}

\newcommand{\mW}{\mathbb{W}}

\newcommand{\supp}{{\rm supp}}

\usepackage{proof}
\usepackage{latexsym}
\usepackage{amsfonts}
\usepackage{amssymb}
\usepackage{cite}
\usepackage{bm}

\newcommand{\brem}{\begin{remark}}
\newcommand{\erem}{\end{remark}}
\newcommand{\blem}{\begin{lemma}}
\newcommand{\elem}{\end{lemma}}
\newcommand{\bth}{\begin{theorem}}
\newcommand{\ethm}{\end{theorem}}
\newcommand{\benu}{\begin{enumerate}}
\newcommand{\eenu}{\end{enumerate}}
\newcommand{\bdes}{\begin{description}}
\newcommand{\edes}{\end{description}}
\newcommand{\bdf}{\begin{definition}}
\newcommand{\edf}{\end{definition}}
\newcommand{\bcor}{\begin{cor}}
\newcommand{\ecor}{\end{cor}}
\newcommand{\bprp}{\begin{proposition}}
\newcommand{\eprp}{\end{proposition}}
\newcommand{\bmlem}{\begin{mlemma}}
\newcommand{\emlem}{\end{mlemma}}
\newcommand{\bclm}{\begin{claim}}
\newcommand{\eclm}{\end{claim}}
\newcommand{\bprf}{{\bf Proof}.\hspace{2mm}}

\newcommand{\eprf}{\hspace*{\fill} $\Box$}

\newcommand{\beqn}{\begin{equation}}
\newcommand{\eeqn}{\end{equation}}
\newcommand{\beqnarr}{\begin{eqnarray}}
\newcommand{\eeqnarr}{\end{eqnarray}}
\newcommand{\beqnarrs}{\begin{eqnarray*}}
\newcommand{\eeqnarrs}{\end{eqnarray*}}

\newcommand{\spand}{\,\&\,}

\newcommand{\Natural}{\mathbb{N}}

\newcommand{\restrict}{\!\upharpoonright\!}

\newtheorem{theorem}{Theorem}[section]
\newtheorem{definition}[theorem]{Definition}
\newtheorem{proposition}[theorem]{Proposition}
\newtheorem{lemma}[theorem]{Lemma}
\newtheorem{cor}[theorem]{Corollary}
\newtheorem{remark}[theorem]{Remark}
\newtheorem{mlemma}[theorem]{Main Lemma}
\newtheorem{claim}[theorem]{Claim}

\newcommand{\alp}{\alpha}

\newcommand{\veps}{\varepsilon}
\newcommand{\del}{\delta}
\newcommand{\Del}{\Delta}
\newcommand{\ome}{\omega}
\newcommand{\Ome}{\Omega}
\newcommand{\bet}{\beta}
\newcommand{\gam}{\gamma}
\newcommand{\Gam}{\Gamma}
\newcommand{\kap}{\kappa}
\newcommand{\sig}{\sigma}
\newcommand{\Sig}{\Sigma}
\newcommand{\tht}{\theta}
\newcommand{\Tht}{\Theta}
\newcommand{\lam}{\lambda}
\newcommand{\Lam}{\Lambda}
\newcommand{\vphi}{\varphi}

\newcommand{\fal}{\forall}
\newcommand{\exi}{\exists}
\newcommand{\rarw }{\rightarrow}
\newcommand{\Rarw }{\Rightarrow}

\newcommand{\lrarw}{\leftrightarrow}
\newcommand{\Lrarw}{\Leftrightarrow}

\newcommand{\cala}{{\cal A}}

\newcommand{\calc}{{\cal C}}

\newcommand{\calg}{{\cal G}}
\newcommand{\calh}{{\cal H}}

\newcommand{\calw}{{\cal W}}
\newcommand{\calx}{{\cal X}}
\newcommand{\caly}{{\cal Y}}

\newcommand{\st}{\star}

\newcommand{\la}{\langle}
\newcommand{\ra}{\rangle}
\newcommand{\lc}{\lceil}
\newcommand{\rc}{\rceil}
\newcommand{\rk}{\mbox{{\rm rk}}}

\newcommand{\sfsc}{{\sf SC}}
\newcommand{\sfk}{{\sf k}}

\newcommand{\setm}{\setminus}

\newcommand{\mK}{\mathbb{K}}

\title{
Lectures on Ordinal Analysis
\thanks{This is a lecture notes for a mini-course in Department of Mathematics, Ghent University,
14 Mar.-25 Mar. 2023.
I'd like to thank A. Weiermann and F. Pakhomov for the hospitality during my stay in Gent, Belgium.
}
}
\author{Toshiyasu Arai (University of Tokyo, JAPAN)
}
\date{}

\begin{document}
\maketitle

The lecture rely on the followings, especially on starred ones.

\begin{itemize}

\item
\cite{Buchholz75}
W. Buchholz, 
Normalfunktionen und konstruktive Systeme von Ordinalzahlen.
In: Diller, J., M\"uller, G. H. (eds.)
Proof Theory Symposion Keil 1974, Lect. Notes Math. vol. 500, pp. 4-25, Springer (1975)

\item
$\mbox{\cite{Buchholz}}^{\star}$ W. Buchholz, 
A simplified version of local predicativity, 
in \textit{Proof Theory},
eds. P. H. G. Aczel, H. Simmons and S. S. Wainer
(Cambridge UP,1992), pp. 115--147.

\item
$\mbox{\cite{BuchholzBSL}}^{\st}$
W. Buchholz, Review of the paper: A. Setzer, Well-ordering proofs for Martin-L\"of type theory, Bulletin of Symbolic Logic 6 (2000) 478-479.

\item
$\mbox{\cite{J1}}^{\st}$
G. J\"ager, Zur Beweistheorie der Kripke-Platek Mengenlehre \"uber den nat\"urlichen Zahlen, 
Archiv f. math. Logik u. Grundl., 22(1982), 121-139.

\item
$\mbox{\cite{J2}}^{\st}$
G. J\"ager, 
A well-ordering proof for Feferman's theory $T_{0}$, 
Archiv f. math. Logik u. Grundl., 23(1983), 65-77.

\item
$\mbox{\cite{Rathjen94}}^{\st}$
M. Rathjen, 
Proof theory of reflection, 
Ann. Pure Appl. Logic 68 (1994) 181--224.

\item
\cite{RathjenAFML2}
M. Rathjen, 
An ordinal analysis of parameter free $\Pi^{1}_{2}$-comprehension,
Arch. Math. Logic 44 (2005) 263-362.

\item
(An ordinal analysis of set theory)
$\mbox{\cite{J1}}^{\st}$.

\item
(Operator controlled derivations)
A streamlined technique introduced in $\mbox{\cite{Buchholz}}^{\st}$, and
its extension in
$\mbox{\cite{Rathjen94}}^{\st}$.

\item
(Shrewd cardinals)
\cite{RathjenAFML2}

\item
(Well-foundedness proofs)
Distinguished classes are introduced in \cite{Buchholz75}.
I have learnt it in $\mbox{\cite{J2}}^{\st}$ and its improved version in
$\mbox{\cite{BuchholzBSL}}^{\st}$.

\end{itemize}

\textbf{Plan}

\benu
\item
${\sf KP}\ome$

\item
Rathjen's analysis of $\Pi_{3}$-reflection

Well-foundedness proof in ${\sf KP}\Pi_{3}$ (skipped)

\item
First-order reflection

\item
First-order reflection (contd.)

\item
$\Pi^{1}_{1}$-reflection

\item
$\Pi^{1}_{1}$-reflection (contd.)

\item
$\Pi^{1}_{1}$-reflection (contd.)

\item
$\Pi_{1}$-collection

\item
$\Pi_{1}$-collection (contd.)

\eenu

An ordinal $\alp$ is said to be \textit{recursive} iff there exists a recursive (computable)
well ordering on $\ome$ of type $\alp$.\index{recursive ordinal}
$\ome^{CK}_{1}$ (\textit{Church-Kleene $\ome_{1}$}) denotes the least non-recursive ordinal.
\index{$\ome^{CK}_{1}$, the least non-recursive ordinal}

\bdf\label{df:ti}
{\rm
\benu
\item
$Prg[\prec,U]:\Lrarw \fal x[\fal y\prec x(y\in U)\to x\in U]$ 
\\
($U$ is \textit{progressive}\index{progressive} with respect to $\prec$).
\index{$Prg[\prec,U]$, $U$ is progressive with respect to $\prec$}
\item
${\rm TI}[\prec,A] :\Lrarw Prg[\prec,A]\to \fal x\, A(x)$ for formulas $A(x)$, and
\\
${\rm TI}[\prec,U]\Lrarw Prg[\prec,U]\to \fal x\, U(x)$ (transfinite induction on $\prec$ ).
\index{${\rm TI}[\prec,U]$, transfinite induction}
\item
Let $\prec$ be a computable strict partial order on $\ome$.
If $\prec$ is well-founded, then let
$|n|_{\prec}:=\sup\{|m|_{\prec}+1:m\prec n\}$, and
$|\prec|:=\sup\{|n|_{\prec}+1:n\in\ome\}$ (the order type of $\prec$).
Otherwise let $|\prec|:=\ome_{1}^{CK}$.
\eenu
}
\edf

\bdf\label{df:prfthordinal}
{\rm
For a theory $T$ comprising elementary recursive arithmetic {\sf EA}
the \textit{proof-theoretic ordinal}\index{proof-theoretic ordinal} 
$|T|$
of $T$ is defined by
\beqn\label{eq:prfthordinal}
|T|:=\sup\{|\!\prec \!|: T\vdash {\rm TI}[\prec,U]\, \mbox{for some recursive well order}
\prec\}
\eeqn
where $U$ is a fresh predicate constant.
}
\edf
Now, most brutally speaking, the aim of the ordinal analysis is to compute and/or describe the proof-theoretic ordinals of 
natural theories, thereby measuring the proof-theoretic strengths of theories with respect to
$\Pi^{1}_{1}$-consequences.

\section{Ordinal analysis of ${\sf KP}\ome$}

\subsection{Kripke-Platek set theory}\label{sect:KP}
A fragment {\sf KP} of Zermelo-Fraenkel set theory {\sf ZF}, Kripke-Platek set theory, is
introduced
Let $\mathcal{L}_{set}=\{\in,=\}$ be the set-theoretic language.
In this section we deal only with set-theoretic models $\langle X;\in\restrict(X\times X)\rangle$, 
and the model is identified with the sets $X$.

\bdf\label{df:delfml}{\rm (}$\Del_{0}, \Sig_{1}, \Pi_{2}, \Sig${\rm )}
\benu
\item
{\rm 
A set-theoretic formula is said to be a
\textit{$\Del_{0}$-formula}\index{bounded formula, $\Del_{0}$-formula in set theory}
if every quantifier occurring in it is \textit{bounded} by a set.
\textit{Bounded quantifiers} is of the form $\fal x\in u, \exi x\in u$.
}

\item 
{\rm 
A formula of the form $\exi x A$ with a $\Del_{0}$-matrix $A$ is a
\textit{$\Sig_{1}$-formula}\index{$\Sig_{1}$-formula in set theory}.

Its dual $\fal x A$ is a \textit{$\Pi_{1}$-formula}.
}
\item
{\rm 
The set of $\Sig$-formulas [$\Pi$-formulas] is the smallest class including $\Del_{0}$-formulas, closed under
positive operations $\land, \lor$ , bounded quantifications $\fal x\in u, \exi x\in u$, and
existential (unbounded) quantification $\exi x$ [universal (unbounded) quantification $\fal x$], resp.

For example $\fal x\in u\exi y A\, (A\in\Del_{0})$ is a $\Sig$-formula  but not a $\Sig_{1}$-formula.
}
\item 
{\rm 
A formula of the form $\fal x A$ with a $\Sig_{1}$-matrix $A$ is a
\textit{$\Pi_{2}$-formula}.
}
\eenu
\edf

We see easily that $\Del_{0}$-formulas are \textit{absolute} in the sense that
for any transitive sets $X\subset Y$ ($X$ is transitive iff $\fal y\in X\fal x\in y(x\in X)$), 
$X\models A[\bar{x}] \Lrarw Y\models A[\bar{x}]$ for any $\Del_{0}$-formula $A$ and $\bar{x}=x_{1},\ldots,x_{n}$ with $x_{i}\in X$.

\bdf
{\rm
Axioms of {\sf KP} are {\bf Extensionality} $\fal a,b[\fal x\in a(x\in b)\land\fal x\in b(x\in a)\to a=b]$, 
{\bf Null set}(the empty set $\emptyset$ exists), {\bf Pair} $\fal x,y\exi a(x\in a\land y\in a)$, 
{\bf Union} $\fal a\exi b\fal x\in a\fal y\in x(y\in b)$, and the following three schemata.
\bdes
\item[{\bf $\Del_{0}$-Separation}]
For any set $a$ and any \textit{$\Del_{0}$-formula} $A$, the set $b=\{x\in a: A(x)\}$ exists.
Namely $\exi b\fal x[x\in b \lrarw x\in a\land A(x)]$.

\item[{\bf $\Del_{0}$-Collection}]
$\fal x\in a\exi y \, A(x,y) \rarw \exi b\fal x\in a\exi y\in b\, A(x,y)$ for \textit{$\Del_{0}$-formulas} $A$.

\item[{\bf Foundation} or {\bf $\in$-Induction}]
$\fal x[\fal y\in x F(y) \rarw F(x)] \rarw \fal x F(x)$
\\
 for \textit{arbitrary formula} $F$.
\edes
}
\edf

${\sf KP}\ome$ denotes {\sf KP} 
\textit{plus} Axiom of 
{\bf Infinity} 
$\exi x\neq\emptyset\fal y\in x[y\cup\{y\}\in x]$.

\subsection{Constructible hierarchy and admissible sets}

The constructible hierarchy
$\{L_{\alp}: \alp\in ON\}$.
\benu
\item 
$L_{0}:=\emptyset$.
\item 
$L_{\alp+1}$ is the collection of all definable sets in $(L_{\alp},\in)$.

\item
$L_{\lam}:=\bigcup_{\alp<\lam}L_{\alp}$ for limits $\lam$.
\item
$L:=\bigcup_{\alp\in ON}L_{\alp}$.
\eenu

Note that
$L_{\ome\alp}\models {\sf KP}- \left( \Del_{0}\mbox{{\bf -Collection}}\right)$ for $\alp>0$,
and $\ome\in L_{\ome\alp}$ if $\alp>1$.

\bdf\label{df:admissible}
{\rm
\benu
\item
A transitive set $A$ is \textit{admissible} if $(A;\in)\models{\sf KP}$.

\item 
An ordinal $\alp$ is \textit{admissible} if $L_{\alp}$ is admissible.
\index{admissible ordinal}

\item 
A relation $R$ on an admissible set $A$ is \textit{$A$-recursive} 
[\textit{$A$-recursively enumerable, $A$-r.e.}]
 (\textit{$A$-finite}) if $R$ is $\Del_{1}$ [$\Sig_{1}$] ($R\in A$), resp.
 
\item 
A function on an admissible set $A$ is \textit{$A$-recursive} if its graph is $A$-r.e.

\item 
An ordinal $\alp$ is \textit{recursively regular} iff $L_{\alp}\models {\sf KP}\ome$.
\index{recursively regular ordinal}
\eenu
}
\edf
Observe that an ordinal $\alp$ is recursively regular iff $\alp$ is a multiplicative principal number$>\ome$,
and for any $L_{\alp}$-recursive function $f:\bet \to \alp$ with a $\bet<\alp$,
 $\sup\{f(\gam):\gam<\bet\}<\alp$ holds.

\bth\label{th:pi2rflonL}{\rm (}$\Pi_{2}${\bf -Reflection on} $L${\rm )}\\
For any $\Sig$-predicate $A$
\[
{\sf KP}\ome \vdash \fal x\in L\exi y\in L\, A(x,y) \to \exi z\in L\fal x\in z\exi y\in z\, A(x,y).
\]
In particular for recursively regular ordinals $\Ome$,
\[
\fal \alp<\Ome\exi\bet<\Ome\, A(\alp,\bet) \to \exi \gam<\Ome\fal \alp<\gam\exi\bet<\gam\, A(\alp,\bet).
\]
\end{theorem}

\blem\label{lem:Sigord}
$|{\sf KP}\ome|\leq |{\sf KP}\ome|_{\Sig}
:=\min\{\alp: \fal A\in\Sig({\sf KP}\ome\vdash A\Rarw L_{\alp}\models A)\}$. 
\elem
\bprf
Suppose ${\sf KP}\ome$ proves ${\rm TI}[\prec,U]$ 
for a computable order $\prec$ on $\ome$, where a unary predicate $U$ may occur in Foundation schema,
but not in $\Del_{0}$-Separation nor $\Del_{0}$-Collection.
Then 
 $\fal n\in\ome\exi\alp(\alp=|n|_{\prec}=\sup\{|m|_{\prec}+1:m\prec n\})$ is provable in ${\sf KP}\ome$.
Therefore 
$|{\sf KP}\ome|\leq |{\sf KP}\ome|_{\Sig}$.
\eprf
\\

The \textit{Mostowski collapsing} $\mbox{{\rm clpse}}(b)$
of a set $b$ is defined by
$C_{b}(x)=\{C_{b}(y): y\in x\cap b\}$ and
$\mbox{{\rm clpse}}(b):=C_{b}(b)=\{C_{b}(x): x\in b\}$.

\bdf
{\rm
We say that a class $\calc$ is \textit{$\Pi_{n}$-classes} for $n\geq 2$ if there exists a set-theoretic 
$\Pi_{n}$-formula $F(\bar{a})$ with parameters $\bar{a}$ such that for any transitive set $P$ with 
$\bar{a}\subset P$,
$P\in\calc\Lrarw P\models F(\bar{a})$ holds.
For a whole universe $L$, $L\in\calc$ denotes the formula $F(\bar{a})$.
By a \textit{$\Pi^{1}_{0}$-class} we mean a $\Pi_{n}$-class for some $n\geq 2$.
}
\edf

\subsection{Buchholz' $\psi$-functions}\label{sect:collapsing}

In this section we work in ${\sf KP}\ome$.

We are in a position to introduce a collapsing function $\psi_{\sig}(\alp)<\sig$ (even if $\alp\geq\sig$).
The following definition is due to \cite{Buchholz86}.

\bdf\label{df:HKP}
{\rm 
Let $\Ome=\ome_{1}$ or $\Ome=\ome_{1}^{CK}$.
Define simultaneously by recursion on ordinals $\alp<\Gam_{\Ome+1}$ the classes
 $\calh_{\alp}(X)\, (X\subset \Ome)$ and the ordinals $\psi_{\Ome}(\alp)$ as follows.

$\calh_{\alp}(X)$ is the Skolem hull of $\{0,\Ome\}\cup X$ under the functions
$+, \vphi$, and
$\bet\mapsto\psi_{\Ome}(\bet)\,(\bet<\alp)$.

Let 
\beqn\label{eq:psidf}
\psi_{\Ome}(\alp)=\min(\{\Ome\}\cup\{\bet<\Ome:  \calh_{\alp}(\bet)\cap\Ome\subset \bet\})
\eeqn
}
\edf

Let us interpret $\Ome=\ome_{1}$.
Then we see readily
that 
$\calh_{\alp}(X)$ is countable for any countable $X$.

To see that the ordinal $\psi_{\Ome}(\alp)$ could be defined, it suffices to show the existence of 
an ordinal $\bet<\Ome$ such that 
$\calh_{\alp}(\bet)\cap\Ome\subset\bet$ :
let $\bet=\sup\{\bet_{n}:n\in\ome\}$ with 
$\bet_{n+1}=\min\{\bet<\Ome:\calh_{\alp}(\bet_{n})\cap\Ome\subset\bet\}$
and $\bet_{0}=0<\Ome$.
Then 
$\calh_{\alp}(\bet)\cap\Ome\subset\bet$ since 
$\calh_{\alp}(\bet)=\bigcup_{n}\calh_{\alp}(\bet_{n})$,
and $\bet<\Ome$ since $\Ome>\ome$ is regular.

The ordinal $\psi_{\Ome_{1}}(\varepsilon_{\Ome_{1}+1})$ is called the 
\textit{Bachmann-Howard ordinal}.

\bprp\label{prp:psi}
\benu
\item\label{prp:psi.0}
$\alp_{0}\leq\alp_{1} \land X_{0}\subset X_{1} \Rarw \calh_{\alp_{0}}(X_{0})\subset
 \calh_{\alp_{1}}(X_{1})$.
 
\item\label{prp:psi.1}
$\calh_{\alp}(\psi_{\Ome}(\alp))\cap\Ome=\psi_{\Ome}(\alp)$ and $\psi_{\Ome}(\alp)\not\in \calh_{\alp}(\psi_{\Ome}(\alp))$.

\item\label{prp:psi.1.5}
$\alp_{0}\leq\alp \Rarw \psi_{\Ome}(\alp_{0})\leq\psi_{\Ome}(\alp) \land \calh_{\alp_{0}}(\psi_{\Ome}(\alp_{0}))\subset \calh_{\alp}(\psi_{\Ome}(\alp))$.

\item\label{prp:psi.2}
$\alp_{0}\in \calh_{\alp}(\psi_{\Ome}(\alp))\cap\alp \Rarw \psi_{\Ome}(\alp_{0})<\psi_{\Ome}(\alp)$.
Therefore
\\
$\alp_{0}\in \calh_{\alp_{0}}(\psi_{\Ome}(\alp_{0})) \land \alp\in \calh_{\alp}(\psi_{\Ome}(\alp)) \Rarw 
(\alp_{0}<\alp \lrarw \psi_{\Ome}(\alp_{0})<\psi_{\Ome}(\alp))$.

\item\label{prp:psi.3}
$\psi_{\Ome}(\alp)$ is a strongly critical number such that
$\psi_{\Ome}(\alp)<\Ome$.

\item\label{prp:psi.4}
$\gam\in \calh_{\alp}(\bet) \Lrarw \sfsc(\gam)\subset\calh_{\alp}(\bet)$, where
$\sfsc(0)=\sfsc(\Ome)=\emptyset$, $\sfsc(\gam)=\{\gam\}$ if 
$\gam\neq\Ome$ is strongly critical, and 
$\sfsc(\vphi\gam\del)=\sfsc(\gam+\del)=\sfsc(\gam)\cup\sfsc(\del)$.

\item\label{prp:psi.5}
$\calh_{\alp}(\psi_{\Ome}(\alp))=\calh_{\alp}(0)$ and 
$\psi_{\Ome}(\alp)=\min\{\xi:\xi\not\in \calh_{\alp}(0)\cap\Ome\}$.
\eenu
\eprp

Proposition \ref{prp:psi}.\ref{prp:psi.5} means that
$\psi_{\Ome}(\alp)$ is the Mostowski's collapse of the point $\Ome$ in the iterated
Skolem hull $\calh_{\alp}(0)$ of ordinals $\{0,\Ome\}$ under 
addition $+$ and the binary Veblen function $\vphi$.
This suggests us that the ordinal $\psi_{\Ome}(\alp)$ could be a substitute for $\Ome$
in a restricted situation.
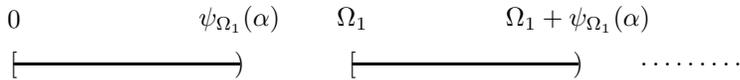
\vspace{-3cm}
\begin{figure}[h]

\unitlength 3mm
\begin{picture}(12,12)
\put(0,2){\makebox(0,0){$0$}}
\put(0,0){\makebox(0,0){$[$}}
\put(0,0){\thicklines\line(1,0){10}}
\put(10,2){\makebox(0,0){$\psi_{\Ome_{1}}(\alp)$}}
\put(10,0){\makebox(0,0){$)$}}
\put(15,2){\makebox(0,0){$\Ome_{1}$}}
\put(15,0){\makebox(0,0){$[$}}
\put(15,0){\thicklines\line(1,0){10}}
\put(25,2){\makebox(0,0){$\Ome_{1}+\psi_{\Ome_{1}}(\alp)$}}
\put(25,0){\makebox(0,0){$)$}}
\put(30,0){\makebox(0,0){$\ldots\ldots\ldots$}}
\end{picture}

\end{figure}

\subsection{Computable notation system $OT(\Ome)$ of ordinals}\label{subsec:TOme}

By Proposition \ref{prp:psi}.\ref{prp:psi.5} we have 
$\calh_{\varepsilon_{\Ome+1}}(0)=\calh_{\varepsilon_{\Ome+1}}(0)=
\calh_{\varepsilon_{\Ome+1}}(\psi_{\Ome}(\varepsilon_{\Ome+1}))$, and hence each ordinal below 
$\psi_{\Ome}(\varepsilon_{\Ome+1})$ can be denoted by terms built up from $0, \Ome, +, \vphi, \psi$.
Although the representation is not uniquely determined from ordinals, e.g., 
$\psi_{\Ome}(\psi_{\Ome}(\Ome))=\psi_{\Ome}(\Ome)$,
$\alp$ can be determined from the ordinal $\psi_{\Ome}(\alp)$ 
if $\alp\in \calh_{\alp}(0)$, 
cf. Propositions \ref{prp:psi}.\ref{prp:psi.2} and \ref{prp:psi}.\ref{prp:psi.5}.
We can devise a recursive notation system $OT(\Ome)$ of ordinals with this restriction
in such a way that the following holds

\bprp
{\sf EA} proves that $(OT(\Ome),<)$ is a linear order.
\eprp

\subsection{Ramified set theory}

\bdf

{\rm
\textit{RS-terms} $t$ and their \textit{levels} $|t|$ are defined recursively as follows.
\benu

\item
For each ordinal $\alp\in OT(\Ome)\cap(\Ome+1)$, 
${\sf L}_{\alp}$ is an RS-term of level $|{\sf L}_{\alp}|=\alp$.

\item
Let $\tht(x,y_{1},\ldots,y_{n})$ be a formula in the set-theoretic language, and
$s_{1},\ldots,s_{n}$ be RS-terms such that $\max\{|s_{i}|:1\leq i\leq n\}<\alp$.
Then the formal expression
$[x\in {\sf L}_{\alp}: \tht^{{\sf L}_{\alp}}(x,s_{1},\ldots,s_{n})]$ is an RS-term
of level
$|[x\in {\sf L}_{\alp}: \tht^{{\sf L}_{\alp}}(x,s_{1},\ldots,s_{n})]|=\alp$.
\eenu
$RS$ denotes the set of all RS-terms.

Let $\tht(x_{1},\ldots,x_{n})$ be a formula 
such that each quantifier is bounded by a variable $y$, $Q x\in y$,
 all free variables occurring in $\tht$ are among the list $x_{1},\ldots,x_{n}$,
and each $x_{i}$ occurs freely in $\tht$.
An \textit{RS-formula} is obtained from such a formula $\tht(x_{1},\ldots,x_{n})$
by substituting RS-terms $t_{i}$ for each $x_{i}$.

Let $\sfk({\sf L}_{\alp}):=\{\alp\}$, 
$\sfk([x\in {\sf L}_{\alp}: \tht^{{\sf L}_{\alp}}(x,s_{1},\ldots,s_{n})])=\{\alp\}\cup\bigcup_{i\leq n}\sfk(s_{i})$
and
\[
 \sfk(\tht(t_{1},\ldots,t_{n})):=\bigcup_{i\leq n}\sfk(t_{i}),\,
|\tht(t_{1},\ldots,t_{n})|:=\max\{|t_{1}|,\ldots,|t_{n}|, 0\}
.\]

The bound ${\sf L}_{\Ome}$ in $\exi x\in {\sf L}_{\Ome}$
and $\fal x\in {\sf L}_{\Ome}$ is the replacements of the unbounded quantifiers $\exi$ and $\fal$, resp.
}
\edf

\bdf
{\rm
Let $s,t$ be RS-terms with $|s|<|t|$.
\[
(s\dot{\in}t) :\equiv\left\{
\begin{array}{ll}
B(s) & t\equiv[x\in {\sf L}_{\alp}: B(x)]
\\
\top & t\equiv{\sf L}_{\alp} 
 \end{array}
 \right.
\]
where $\top$ denotes a true literal, e.g., $\emptyset\not\in\emptyset$.
}
\edf

We assign disjunctions or conjunctions to sentences as follows.
When a disjunction $\bigvee(A_{i})_{i\in J}$ [a conjunction $\bigwedge(A_{i})_{i\in J}$]
is assigned to $A$,
we denote $A\simeq\bigvee(A_{i})_{i\in J}$ [$A\simeq\bigwedge(A_{i})_{i\in J}$], resp.

\bdf\label{df:assigndc}
\begin{enumerate}

\item
$(A_{0}\lor A_{1}):\simeq\bigvee(A_{i})_{i\in J}$
{\rm and}
$(A_{0}\land A_{1}):\simeq\bigwedge(A_{i})_{i\in J}$
{\rm with} $J:=2$.

\item
$(a\in b):\simeq\bigvee(t\dot{\in}b\land t=a)_{t\in J}$
{\rm and}
$(a\not\in b):\simeq\bigwedge(t\dot{\in}b \to t\neq a)_{t\in J}$
{\rm with}
$J:=Tm(|b|):=\{t\in RS: |t|<|b|\}$.

\item
 {\rm Let $a,b$ be set terms.}
 \\
$(a\neq b):\simeq\bigvee(\lnot A_{i})_{i\in J}$
{\rm and}
$(a=b):\simeq\bigwedge(A_{i})_{i\in J}$
{\rm with} $J:=2$ {\rm and}
$A_{0}:\equiv(\fal x\in a(x\in b))$, $A_{1}:\equiv(\fal x\in b(x\in a))$.

\item
$
\exists x\in b\, A(x):\simeq\bigvee(t\dot{\in}b\land A(t))_{t\in J}$
{\rm and}
$
\forall x\in b\, A(x):\simeq\bigwedge(t\dot{\in}b \to A(t))_{t\in J}
$ 
{\rm with}
$
J:=Tm(|b|)$.

\end{enumerate}
\edf

\blem\label{prp:impdegreeKP.simeq.2}
$\fal i\in J\left(\sfk(i)\subset\sfk(A_{i})\subset\sfk(A)\cup\sfk(i)\right)$
for $A\simeq\bigvee(A_{i})_{i\in J}$, where $\sfk(0)=\sfk(1)=\emptyset$.
\elem

The \textit{rank} $\rk(A),\rk(a)<\Ome+\omega$ of RS-formulas $A$ 
and RS-terms $a$ are defined so that the followings hold for any formula $A$.

\bprp\label{prp:degreeKP}
\benu

\item\label{prp:impdegreeKP.-2}
$\rk(\cala)\in\{\ome| \cala|+n:n\in\ome\}$ for RS-terms and RS-formulas $\cala$.

\item\label{prp:degreeKP2}
$\rk(B(t))\in\{\ome|t|+n:n\in\ome\}\cup\{\rk(B({\sf L}_{0}))\}$.

\item\label{prp:impdegreeKP}
Let $A\simeq\bigvee(A_{i})_{i\in J}$. Then
$\fal i\in J(\rk(A_{i})<\rk(A))$.
\eenu
\eprp

\bdf\label{df:sigmaomeRS}
{\rm 
\benu
\item
Let $B(x_{1},\ldots,x_{n})$ be a $\Del_{0}$-formula,
and $a_{1},\ldots,a_{n}\in RS$ be $|a_{i}|<\Ome$.
Then $B(a_{1},\ldots,a_{n})$ is a \textit{$\Del(\Ome)$-formula}.

\item
Let $A(x_{1},\ldots,x_{n})$ be a $\Sig$-formula,
and $a_{1},\ldots,a_{n}\in RS$ be $|a_{i}|<\Ome$.
Then $A^{({\sf L}_{\Ome})}(a_{1},\ldots,a_{n})$ is a \textit{$\Sig(\Ome)$-formula},
where for RS-terms $c$,
$A^{(c)}$ denotes the result of replacing unbounded existential quantifiers
$\exi x(\cdots)$ by $\exi x\in c(\cdots)$.

\item
Let $B\equiv A^{({\sf L}_{\Ome})}$ be a $\Sig(\Ome)$-formula, and 
$\alp\in OT(\Ome)\cap\Ome$.
Then $B^{(\alp,\Ome)}\equiv A^{({\sf L}_{\alp})}$.
For $\Gam\subset\Sig(\Ome)$, $\Gam^{(\alp,\Ome)}:=\{B^{(\alp,\Ome)}: B\in\Gam\}$.

\eenu
}
\edf

Let us define a derivability relation $\calh_{\gam}[\Tht]\vdash^{a}_{b}\Gamma$ for 
finite sets $\Tht$ of ordinals, $\gam,a<\veps_{\Ome+1}$, $b<\Ome+\ome$ and 
RS-sequents, i.e., finite sets of RS-formulas $\Gam$.

\bdf\label{df:controlderreg}
{\rm
$\calh_{\gam}[\Tht]\vdash^{a}_{b}\Gamma$ holds if
\beqn
\label{eq:controlderKP}
\{\gam,a,b\}\cup\sfk(\Gam)\subset\calh_{\gam}[\Tht]
\eeqn
and one of the following
cases holds:

\begin{description}
\item[$(\bigvee)$]
There are $A\in\Gamma$ such that
$A\simeq\bigvee(A_{i})_{i\in J}$, an $i\in J$  with
\beqn
\label{eq:controlder1KP}
|i|<a
\eeqn
and an $a(i)<a$ for which
$\calh_{\gam}[\Tht]\vdash^{a(i)}_{b}\Gamma,A_{i}$ holds.
\[
\infer[(\bigvee)]{\calh_{\gam}[\Tht]\vdash^{a}_{b}\Gamma}
{
\calh_{\gam}[\Tht]\vdash^{a(i)}_{b}\Gamma,A_{i}
}
\,
(|i|<a)
\]

\item[$(\bigwedge)$]
There is an $A\in\Gamma$ such that
$A\simeq\bigwedge(A_{i})_{i\in J}$, and for each $i\in J$, 
there is an $a(i)$ such that $a(i)<a$ for which
$\calh_{\gam}[\Tht\cup\sfk(i)]\vdash^{a(i)}_{b}\Gamma,A_{i}$ holds.
\[
\infer[(\bigwedge)]{\calh_{\gam}[\Tht]\vdash^{a}_{b}\Gamma}
{
\{
\calh_{\gam}[\Tht\cup\sfk(i)]\vdash^{a(i)}_{b}\Gamma,A_{i}
\}_{i\in J}
}
\]

\item[$(cut)$]
There are $C$ and $a_{0}<a$ such that
$\rk(C)<b$, $\calh_{\gam}[\Tht]\vdash^{a_{0}}_{b}\Gamma,\lnot C$ 
and
 $\calh_{\gam}[\Tht]\vdash^{a_{0}}_{b}C,\Gamma$.
\[
\infer[(cut)]{\calh_{\gam}[\Tht]\vdash^{a}_{b}\Gamma}
{
\calh_{\gam}[\Tht]\vdash^{a_{0}}_{b}\Gamma,\lnot C
&
\calh_{\gam}[\Tht]\vdash^{a_{0}}_{b}C, \Gamma
\,(\rk(C)<b)
}
\]

\item[$(\Del_{0}(\Ome)\mbox{{\rm -Coll}})$]
$b\geq\Ome$, and there are a formula
$C\in\Sig(\Ome)$ and an
$a_{0}<a$ such that
$\calh_{\gam}[\Tht]\vdash^{a_{0}}_{b}\Gamma, C$
and
$\calh_{\gam}[\Tht\cup\{\alp\}]\vdash^{a_{0}}_{b}\Gamma, \lnot C^{(\alp,\Ome)}$
for every $\alp<\Ome$.

\[
\infer[(\Del_{0}(\Ome)\mbox{{\rm -Coll}})]{\calh_{\gam}[\Tht]\vdash^{a}_{b}\Gamma}
{
\calh_{\gam}[\Tht]\vdash^{a_{0}}_{b}\Gamma, C
&
\{
\calh_{\gam}[\Tht\cup\{\alp\}]\vdash^{a_{0}}_{b} \lnot C^{(\alp,\Ome)}, \Gam
\}_{\alp<\Ome}
}
\]
\end{description}
}
\edf

\blem\label{lem:tautKP}{\rm (Tautology)}
$
\calh_{0}[\sfk(A)]\vdash_{0}^{2d}\lnot A,A
$
 with $d=\mbox{{\rm rk}}(A)$.
\elem

\blem\label{lem:3.9.7KP} {\rm (Inversion)}\\
$\calh_{\gam}[\Tht]\vdash^{a}_{b}\Gam,A \Rarw \fal i\in J(\calh_{\gam}[\Tht\cup\sfk(i)]\vdash^{a}_{b}\Gam,A_{i})$
for $A\simeq\bigwedge(A_{i})_{i\in J}$.
\elem

\blem\label{lem:boundednessKP}{\rm (Boundedness)}
Let $a\leq \bet\in\calh_{\gam}[\Tht]\cap\Ome$ and $\Lam\subset\Sig(\Ome)$.
Then $
\calh_{\gam}[\Tht]\vdash^{a}_{b}\Gam,\Lam 
\Rarw \calh_{\gam}[\Tht]\vdash^{a}_{b}\Gam,\Lam^{(\bet,\Ome)}
$.
\elem

\blem\label{lem:embedKP}{\rm (Embedding)}\\
Let $\Gam[\vec{x}:=\vec{a}]\, (\vec{a}\subset RS)$ denote a closed instance of a sequent $\Gam$ with restriction of unbounded quantifiers to ${\sf L}_{\Ome}$.
Assume ${\sf KP}\ome\vdash\Gam$. Then
\[
\exi m,l<\ome\fal\vec{a}\subset RS[
\calh_{0}[\sfk(\vec{a})]\vdash^{\Ome+l}_{\Ome+m}
\Gam[\vec{x}:=\vec{a}]]
\]
where $\sfk(\vec{a})=\sfk(a_{1})\cup\cdots\sfk(a_{n})$ for $\vec{a}=(a_{1},\ldots,a_{n})$.
\elem

Let $\tht_{c}(a)$ be the $c$-th iterate of $\tht_{1}(a)=\ome^{a}$.
$\tht_{0}(a)=a$, $\tht_{c\dot{+}d}(a)=\tht_{c}(\tht_{d}(a))$,
and
$\tht_{\ome^{c}}(a)=\vphi_{c}(a)$.

\blem\label{lem:predceKP}{\rm (Predicative Cut-elimination)}\\
$\calh_{\gam}[\Tht]\vdash^{a}_{b+c}\Gam \Rarw \calh_{\gam}[\Tht]\vdash^{\tht_{c}(a)}_{b}\Gam$
if $\lnot(b<\Ome\leq b+c)$.
\elem

\bth\label{th:CollapsingthmKP}{\rm (Collapsing)}\\
Suppose 
\beqn
\label{eq:CollapsingthmassKP}
\Tht\subset\calh_{\gam}(\psi_{\Ome}(\gam))
\eeqn
for a finite set $\Tht$ of ordinals, and $\Gam\subset\Sig(\Ome)$. 
Then for $\hat{a}=\gam+\ome^{a}$ and $\bet=\psi_{\Ome}(\hat{a})$
\[
\calh_{\gam}[\Tht]\vdash^{a}_{\Ome}\Gam \Rarw 
\calh_{\hat{a}+1}[\Tht]\vdash^{\bet}_{\bet}\Gam.
\]
\end{theorem}
\bprf
This is seen by induction on $a$.
Observe that $\sfk(\Gam)\cup\{\bet\}\subset\calh_{\hat{\alp}+1}[\Tht]$ by 
$\gam<\hat{\alp}+1$ and (\ref{eq:controlderKP}).
\\

\noindent
{\bf Case 1}. The last inference is a $(\bigvee)$.

Let $A\in\Gam$ be such that $A\simeq\bigvee(A_{i})_{i\in J}$, and for an $i\in J$ and an 
$a(i)<a$
\[
\infer[(\bigvee)]{\calh_{\gam}[\Tht]\vdash^{a}_{\Ome}\Gam}
{
\calh_{\gam}[\Tht]\vdash^{a(i)}_{\Ome}\Gam,A_{i}
}
\]
By IH it suffices to show $|i|<\psi_{\Ome}(\hat{a})$ for (\ref{eq:controlder1KP}).
We can assume $\sfk(i)\subset\sfk(A_{i})$.
Then $|i|\in \sfk(A_{i})\subset \calh_{\gam}[\Tht]\subset\calh_{\gam}(\psi_{\Ome}(\gam))$ by
(\ref{eq:controlderKP}) and the assumption (\ref{eq:CollapsingthmassKP}).
On the other hand we have $|i|<\Ome$.
Hence $|i|\in\calh_{\hat{a}}(\psi_{\Ome}(\hat{a}))\cap\Ome=\psi_{\Ome}(\hat{a})$.
\\

\noindent
{\bf Case 2}. The last inference is a $(\bigwedge)$.

Let $A\in\Gam$ be such that $A\simeq\bigwedge(A_{i})_{i\in J}$, and for each $i\in J$,
there are 
$a(i)<a$ such that 
\[
\infer[(\bigwedge)]{\calh_{\gam}[\Tht]\vdash^{a}_{\Ome}\Gam}
{
\{
\calh_{\gam}[\Tht\cup\sfk(i)]\vdash^{a(i)}_{\Ome}\Gam,A_{i}
\}_{i\in J}
}
\]
By IH it suffices to show that $\fal i\in J(\sfk(i)\subset\calh_{\gam}(\psi_{\Ome}(\gam)))$.
For example consider the case when $A\equiv(\fal x\in u\, B(x))$ for a set term $u$.
Then $J=\{t\in RS: |t|<|u|\}$.
Since $A$ is a $\Sig(\Ome)$-sentence, we have $|a|<\Ome$.
On the other hand we have $|u|\in\calh_{\gam}[\Tht]$ for $|u|=\max \sfk(u)$, and hence
$\sfk(i)\subset |u|\in\calh_{\gam}(\psi_{\Ome}(\gam))\cap\Ome=\psi_{\Ome}(\gam)$
for any $i\in J$.
\\

\noindent
{\bf Case 3}. The last inference is a $(\Del_{0}(\Ome)\mbox{{\rm -Coll}})$.

There are a sentence $C\in\Sig(\Ome)$ and an
$a_{0}<a$ such that

\[
\infer[(\Del_{0}(\Ome)\mbox{{\rm -Coll}})]{\calh_{\gam}[\Tht]\vdash^{a}_{\Ome}\Gamma}
{
\calh_{\gam}[\Tht]\vdash^{a_{0}}_{\Ome}\Gamma, C
&
\{
\calh_{\gam}[\Tht\cup\{\alp\}]\vdash^{a_{0}}_{\Ome} \lnot C^{(\alp,\Ome)}, \Gam
\}_{\alp<\Ome}
}
\]
Let $\widehat{a_{0}}=\gam+\ome^{a_{0}}$ and $\bet_{0}=\psi_{\Ome}(\widehat{a_{0}})$.
IH yields 
$\calh_{\hat{a}+1}[\Tht]\vdash^{\bet_{0}}_{\bet_{0}}\Gamma,C$.
Boundedness \ref{lem:boundednessKP} yields
$\calh_{\widehat{a_{0}}+1}[\Tht]\vdash^{\bet_{0}}_{\bet_{0}}\Gamma,C^{(\bet_{0},\Ome)}$,
where $\bet_{0}\in\calh_{\widehat{a_{0}}+1}[\Tht]$.
On the other hand we have 
$\calh_{\gam}[\Tht\cup\{\bet_{0}\}]\vdash^{a_{0}}_{\Ome} \lnot C^{(\bet_{0},\Ome)}, \Gam$,
and 
$\calh_{\widehat{a_{0}}+1}[\Tht]\vdash^{a_{0}}_{\Ome} \lnot C^{(\bet_{0},\Ome)}, \Gam$.
IH yields
$\calh_{\widehat{a_{0}}+\ome^{a_{0}}+1}[\Tht]\vdash^{\bet_{1}}_{\bet} \lnot C^{(\bet_{0},\Ome)}, \Gam$,
where
$\bet_{1}=\psi_{\Ome}(\widehat{a_{0}}+\ome^{a_{0}})$ with
$\widehat{a_{0}}+\ome^{a_{0}}=\gam+\ome^{a_{0}}+\ome^{a_{0}}<\hat{a}$.
A $(cut)$ with $\rk(C^{(\bet_{0},\Ome)})<\bet$ yields
$\calh_{\hat{a}+1}[\Tht]\vdash^{\bet}_{\bet}\Gam$.
\\

\noindent
{\bf Case 4}. The last inference is a $(cut)$.

\[
\infer[(cut)]{\calh_{\gam}[\Tht]\vdash^{a}_{\Ome}\Gam}
{
\calh_{\gam}[\Tht]\vdash^{a_{0}}_{\Ome}\Gam,\lnot C 
& 
\calh_{\gam}[\Tht]\vdash^{a_{0}}_{\Ome}C,\Gam
}
\]
We obtain $\rk(C)<\Ome$, and $\rk(C)\in\calh_{\gam}[\Tht]\cap\Ome\subset\psi_{\Ome}(\gam)\leq\bet$.
IH followed by a $(cut)$ yields the lemma.
\eprf

\blem\label{lem:truthKP}{\rm (Truth)}

If
$\calh_{\gam}[\Tht]\vdash^{\alp}_{\Ome}\Gam$ with
$\Gam\subset\Del(\Ome)$,
then
$L_{\Ome}\models\Gam$.
\elem

\bth\label{th:mthKP1}

${\sf KP}\ome\vdash \Gam \mbox{ {\rm and }} \Gam\subset\Sig(\Ome_{1})
\Rarw \exi m<\ome \left[L_{\Ome}\models\Gam^{(\psi_{\Ome}(\ome_{m}(\Ome+1)),\Ome)}\right]$.

\end{theorem}
\bprf
Let ${\sf KP}\ome\vdash \Gam$ for a set $\Gam$ of $\Sig$-sentences. 
By Embedding \ref{lem:embedKP} pick an $m<\ome$ such that
$\calh_{0}[\emptyset]\vdash^{\Ome+m}_{\Ome+m}\Gam$.
Predicative Cut-elimination \ref{lem:predceKP} yields
$\calh_{0}[\emptyset]\vdash^{a}_{\Ome}\Gam$ for $a=\ome_{m}(\Ome+m)$.
Let $\bet=\psi_{\Ome}(\hat{a})$ with $\hat{a}=\ome^{a}=\ome_{m+1}(\Ome+m)$.
We then obtain 
$\calh_{\hat{a}+1}[\emptyset]\vdash^{\bet}_{\bet}\Gam$ by Collapsing \ref{th:CollapsingthmKP},
and $\calh_{\hat{a}+1}[\emptyset]\vdash^{\bet}_{\bet}\Gam^{(\bet,\Ome)}$ by 
Boundedness \ref{lem:boundednessKP}.
We see $L_{\Ome}\models \Gam^{(\bet,\Ome)}$
from Truth \ref{lem:truthKP}.
From $\bet<\psi_{\Ome}(\ome_{m+2}(\Ome+1))$ and the persistency of $\Sig$-formulas, 
we conclude $L_{\Ome}\models \Gam^{(\psi_{\Ome}(\ome_{m+2}(\Ome+1)),\Ome)}$.
\eprf

\subsection{Well-foundedness proof in ${\sf KP}\ome$}\label{sect:wfproofsOme}
In this subsection
$\alp,\bet,\gam,\del,\ldots$ range over ordinal terms in $OT(\Ome)$,
and
$<$ denotes the relation between ordinal terms defined in Definition \ref{df:TOme}.
An ordinal term $\alp$ is identified with the set $\{\bet\in OT(\Ome):\bet<\alp\}$.
For ordinal terms $\alp,\bet$,
ordinal terms $\alp+\bet$ and $\ome^{\alp}$ are defined trivially.

In this subsection we show that the theory ${\sf ID}$ for non-iterated positive elementary inductive definitions on $\mathbb{N}$ proves the fact
that the relation $<$ on $OT(\Ome)$ is well-founded up to each $\alp<\psi_{\Ome}(\varepsilon_{\Ome+1})$.

\bth\label{th:wfID}
For {\rm each} $n<\ome$
\[
{\sf ID}\vdash 
{\rm TI}[< \restrict \psi_{\Ome}(\ome_{n}(\Ome+1)),B]
\]
for any formula $B$ in the language $\mathcal{L}({\sf ID})$.
\end{theorem}

$Acc$ denotes the accessible part of $<$ in $OT(\Ome)$,
which is defined in ${\sf ID}$ as the least fixed point $P_{\cala}$
of the operator $\cala(X,\alp):\Lrarw \alp\subset X\Lrarw(\fal\bet<\alp(\bet\in X))$.
It suffices to show the following, which is equivalent to Theorem \ref{th:wfID}.

\bth\label{th:B3.15}
For each $\alp<\psi_{\Ome}(\varepsilon_{\Ome+1})$, ${\sf ID}\vdash\alp\in Acc$.
\end{theorem}

The least fixed point $Acc$ enjoys $\fal \alp(\alp\subset Acc\to \alp\in Acc)$,
and $\fal \alp(\alp\subset F\to \alp\in F)\to Acc\subset F$.
From these we see easily that $Acc$ is closed under $+,\vphi$ besides $0\in Acc$.
Hence we obtain $\Gam_{0}=\psi_{\Ome}(0)\in Acc$.
Likewise $\Gam_{1}=\psi_{\Ome}(1)\in Acc$ follows.
To prove $\psi_{\Ome}(\Ome)\in Acc$, we need to show $\psi_{\Ome}(\alp)\in Acc$ for any $\alp<\Ome$
such that $\psi_{\Ome}(\alp)$ is an ordinal term, i.e., $G(\alp)<\alp$.
This means that when $\psi_{\Ome}(\bet)$ occurs in $\alp$, then $\bet<\alp$ holds.
Thus we have a chance to prove inductively that $\psi_{\Ome}(\alp)\in Acc$.
The ordinal term $\alp$ is built from $0$, $\Ome$ and some ordinal terms $\psi_{\Ome}(\bet)$ with
$\bet<\alp$ by $+,\vphi$.
Let us assume that each of ordinals $\psi_{\Ome}(\bet)<\Ome$ occurring in $\alp$ is in $W_{0}=Acc\cap\Ome$,
and denote the set of such ordinals $\alp$ by $M_{1}$.
Though we don't have $\Ome\in Acc$ in hand (since this means that $OT(\Ome)\cap\Ome$ is well-founded,
which is the fact we are going to prove), $\Ome$ is in the accessible part $W_{1}$ of the set $M_{1}$.
It turns out that $W_{1}$ is progressive on $M_{1}$, and $\Ome\in W_{1}$.
Moreover $\ome^{\Ome+1}\in W_{1}$ is seen as for the jump set for epsilon numbers.
In this way we see that $\alp\in W_{1}$, i.e., $\psi_{\Ome}(\alp)\in W_{0}$ for \textit{each}
 $\alp<\veps_{\Ome+1}$.
\\

Let $\sfsc(\alp)$ denote the set of strongly critical parts of $\alp$ 
defined in Proposition \ref{prp:psi}.\ref{prp:psi.4}, and let
$\sfsc_{\Ome}(\alp)=\sfsc(\alp)\cap\Ome$.

\bdf
{\rm
 $
M_{1}=\{\alp\in OT(\Ome): \sfsc_{\Ome}(\alp)\subset W_{0}\}
$.
}
\edf

\bprp\label{prp:Ka}

$G(\bet)<\alp \Rarw \sfsc_{\Ome}(\bet)<\psi_{\Ome}(\alp)$
for $\psi_{\Ome}(\alp)\in OT(\Ome)$.
\eprp
\bprf
By induction on the length of ordinal terms $\bet$.
Assume $G(\bet)<\alp$.
By IH we can assume $\bet=\psi_{\Ome}(\gam)$.
Then $\gam\in G(\bet)$ and $\sfsc_{\Ome}(\bet)=\{\bet\}$.
Hence $\gam<\alp$ and $\bet<\psi_{\Ome}(\alp)$.
\eprf

In what follows we work in ${\sf ID}$ except otherwise stated.

\blem\label{lem:endext}
$M_{1}\cap\Ome=W_{0}$.

\elem

\[
\cala(X):=\{\alp\in M_{1}: M_{1}\cap\alp\subset X\}
.\]

\bprp\label{prp:B3.10}
For each formula $F$,
$\cala(F)\subset F \to \Ome\in F$.
\eprp
\bprf
Assuming $\cala(F)\subset F$, we see $\alp\in W_{0} \Rarw \alp\in F$
by induction on $\alp\in W_{0}$.
\eprf

\blem\label{lem:B3.11}
For each formula $F$,
$\cala(F)\subset F \to \cala({\sf j}[F])\subset{\sf j}[F]$, where
${\sf j}[F]:=
\{\bet\in OT(\Ome): \fal\alp
(M_{1}\cap\alp\subset F \to
M_{1}\cap(\alp+\ome^{\bet})\subset F)\}$.
\elem

\blem\label{lem:B3.12}
For each formula $F$ and {\rm each} $n<\ome$,
$\cala(F)\subset F\to \ome_{n}(\Ome+1)\in F$.
\elem

\[
\alp\in W :\Lrarw 
\left(
\psi_{\Ome}(\alp)\in OT(\Ome) \to 
\psi_{\Ome}(\alp)\in W_{0}
\right)
.
\]

\blem\label{lem:B3.14}
$\cala(W)\subset W$.
\elem
\bprf
Assume $\alp\in\cala(W)$ and 
$\psi_{\Ome}(\alp)\in OT(\Ome)$.
Then $\alp\in M_{1}$ and $M_{1}\cap\alp\subset W$.
We show 
\[
\gam<\psi_{\Ome}(\alp) \to\gam\in W_{0}
\]
by induction on the length of ordinal terms $\gam$.
We can assume that $\gam=\psi_{\Ome}(\bet)$.
Then $\bet<\alp$.
We see $\bet\in M_{1}$ from IH.
Therefore $\bet\in M_{1}\cap\alp\subset W$, which yields
$\gam=\psi_{\Ome}(\bet)\in W_{0}$.
Therefore
$\psi_{\Ome}(\alp)\subset W_{0}$.
\eprf
\\

Let us show Theorem \ref{th:B3.15}.
We show that 
${\sf ID}$ proves $\psi_{\Ome}(\ome_{n}(\Ome+1))\in W_{0}$ for each 
$n<\ome$.
By Lemmas \ref{lem:B3.12} and \ref{lem:B3.14} we obtain
$\ome_{n}(\Ome+1)\in W$.
Thus $\psi_{\Ome}(\ome_{n}(\Ome+1))\in W_{0}$ by the definition of $W$.

\section{Rathjen's analysis of $\Pi_{3}$-reflection}\label{sec:pi3}

Given an analysis of ${\sf KP}\ome$ for a single recursively regular ordinal,
it is not hard to extend it to an analysis of theories of recursively regular ordinals of
a given order type, e.g., to ${\sf KP}\ell$, or equivalently to $\Pi^{1}_{1}\mbox{-CA+BI}$.
Or to an iteration of recursively regularities in another manner.
Specifically an ordinal analysis of ${\sf KP}M$ for recursively Mahlo ordinals
is not an obstacle.

Let us introduce a $\Pi_{i}$-recursively Mahlo operation $RM_{i}$ and its iterations.
A $\Pi_{i}$-recursively Mahlo operation $RM_{i}$ for $2\leq i<\omega$, 
is defined through a universal $\Pi_{i}$-formula 
$\Pi_{i}(a)$ such that for each $\Pi_{i}$-formula $\vphi(x)$ there exists a natural number $n$ such that
${\sf KP}\vdash \fal x[\vphi(x)\lrarw \Pi_{i}(\la n,x\ra)]$.
Let $\mathcal{X}$ be a collection of sets.
\begin{eqnarray*}
P\in RM_{i}(\mathcal{X}) & :\Leftrightarrow & \forall b\in P
\left[
P\models\Pi_{i}(b) \to \exists Q\in \mathcal{X}\cap P(b\in Q\models\Pi_{i}(b))
\right]
\\
&&\mbox{(read: } P \mbox{ is } \Pi_{i}\mbox{-reflecting on } \mathcal{X}\mbox{.)}
\end{eqnarray*}
Let $RM_{i}=RM_{i}(V)$, and $V$ is \textit{$\Pi_{i}$-reflecting} if $V\in RM_{i}$.
Under the axiom $V=L$ of constructibility,
$V\in RM_{2}$ iff $V\models{\sf KP}\ome$,
and
$V\in RM_{2}(RM_{2})$ iff $V$ is \textit{recursively Mahlo} universe.
When $V=L_{\sig}$, the ordinal $\sig$ is recursively Mahlo ordinal.

Let ${\sf KP}M$ denote a set theory for recursively Mahlo universes.
For an ordinal analysis of ${\sf KP}M$, it suffices for us to have two step collapsings
$\alp\mapsto\sig=\psi_{M}(\alp)\in RM_{2}$ and
$(\sig,\bet)\mapsto\psi_{\sig}(\bet)$.

Assume that $P\in\mathcal{X}$ is given by a $\Del_{0}$-formula.
Then there exists a $\Pi_{i+1}$-formula $rm_{i}$ such that
for any non-empty transitive sets $P\in V\cup\{V\}$,
$P\in RM_{i}(\mathcal{X})\lrarw rm_{i}^{P}$,
where $rm_{i}^{P}$ denotes the result of restricting unbounded quantifiers in $rm_{i}$ to $P$.

An iteration of $RM_{i}$ along a definable relation $\prec$ is defined as follows.
\[
P\in RM_{i}(a;\prec) :\Lrarw a\in P\in\bigcap\{ RM_{i}(RM_{i}(b;\prec)): b\in P\models b\prec a\}.
\]
Assume that $b\prec a$ is given by a $\Sig_{1}$-formula.
Then there exists a $\Pi_{i+1}$-formula $rm_{i}(a,\prec)$ such that
for any non-empty transitive sets $P\in V\cup\{V\}$ and $a\in P$,
$P\in RM_{i}(a;\prec)\lrarw rm_{i}^{P}(a,\prec)$.

For $2\leq N<\ome$, ${\sf KP}\Pi_{N}$ denotes a set theory for $\Pi_{N}$-reflecting universes $V$,
which is obtained from ${\sf KP}\ome$ by adding an axiom $V\in RM_{N}$
(the axiom for $\Pi_{N}$-reflection) stating that
its universe is $\Pi_{N}$-reflecting.
This means that for each $\Pi_{N}$-formula $\vphi$,
$\vphi(a) \to \exi c[ad_{N}^{c}\land a\in c \land \vphi^{c}(a)]$ is an axiom,
where $ad_{2}^{c}:\equiv(\fal x\in c\fal y\in x(y\in c))$, i.e., $c$ is transitive,
and for $N>2$, $ad\equiv ad_{N}$ denotes 
a $\Pi_{3}$-sentence such that $P\models ad\Lrarw P\models\mbox{{\sf KP}}\ome$
for any transitive and well-founded sets $P$.
${\sf KP}\Pi_{2}$ is a subtheory of ${\sf KP}\ome+(V=L)$, which is interpreted in ${\sf KP}\ome$:
${\sf KP}\ome+(V=L)\vdash\vphi \Rarw {\sf KP}\ome\vdash\vphi^{L}$, cf.\,Theorem \ref{th:pi2rflonL}.

${\sf KP}\Pi_{N+1}$ is much stronger than ${\sf KP}\Pi_{N}$ since $\Pi_{N}$-recursively Mahlo  operation
$RM_{N}$ can be iterated in ${\sf KP}\Pi_{N+1}$.
For example, ${\sf KP}\Pi_{N+1}$ proves $\fal \alp\in ON[V\in RM_{N}(\alp;<)]$ by induction on
ordinals $\alp$.
Suppose $\fal\bet<\alp[V\in RM_{N}(\bet;<)]$.
Let $\vphi$ be a $\Pi_{N}$-formula such that
$V\models\vphi$, and $\bet<\alp$.
We can reflect a $\Pi_{N+1}$-formula $V\in RM_{N}(\bet;<) \land\vphi$, and obtain
a set $P$ such that $P\in RM_{N}(\bet;<) \land P\models\vphi$.
Hence $V\in RM_{N}(\alp;<)$.
This means that $V$ is in the diagonal intersection 
$\triangle_{\alp}RM_{N}(\alp;<)$, i.e., $V\in\bigcap\{RM_{N}(\alp;<): \alp\in ON\cap V\}$.
Since this is a $\Pi_{N+1}$-formula, the $\Pi_{N+1}$-reflecting universe $V$ reflects it:
there exists a set $P\in V$ such that $P$ is in the diagonal intersection,
i.e., $P\in\bigcap\{RM_{N}(\alp;<): \alp\in ON\cap P\}$, and so forth.

Let $ON\subset V$ denote the class of ordinals, $ON^{\veps}\subset V$ and $<^{\veps}$ be $\Del$-predicates such that
for any transitive and well-founded model $V$ of $\mbox{{\sf KP}}\ome$,
$<^{\veps}$ is a well order of type $\veps_{\mathbb{K}+1}$ on $ON^{\veps}$
for the order type $\mathbb{K}$ of the class $ON$ in $V$.
$\lceil\ome_{n}(\mathbb{K}+1)\rceil\in ON^{\veps}$ denotes the code of the `ordinal' 
$\ome_{n}(\mathbb{K}+1)$, which is assumed to be a closed `term' built from the code 
$\lceil \mathbb{K}\rceil$ and $n$,
e.g., $\lceil\alp\rceil=\la 0,\alp\ra$ for $\alp\in ON$, 
$\lceil \mathbb{K}\rceil=\la 1,0\ra$ and 
$\lceil\ome_{n}(\mathbb{K}+1)\rceil=\la 2,\la 2,\cdots\la 2, \la 3,\lceil \mathbb{K}\rceil, \la 0,1\ra\ra\ra\cdots\ra\ra$.

$<^{\veps}$ is assumed to be a standard epsilon order with base $\mathbb{K}$ (not on $\Natural$, but on $V$) such that 
$\mbox{{\sf KP}}\ome$ proves the fact that $<^{\veps}$ is a linear ordering, and for any formula $\vphi$
and each $n<\ome$,
\beqn\label{eq:trindveps}
\mbox{{\sf KP}}\ome\vdash\fal x(\fal y<^{\veps}x\,\vphi(y)\to\vphi(x)) \to \fal x<^{\veps}\lceil\ome_{n}(\mathbb{K}+1)\rceil\vphi(x)
\eeqn

\bth\label{th:resolutionpi3}{\rm (\cite{consv})}

For each $N\geq 2$,
{\sf KP}$\Pi_{N+1}$ is $\Pi_{N+1}$-conservative over
the theory 
\[
\mbox{{\sf KP}}\ome+
\{
V\in RM_{N}(\lceil\ome_{n}(\mathbb{K}+1)\rceil;<^{\veps}): 
n\in\ome\}
.\]
\end{theorem}

From (\ref{eq:trindveps}) we see that {\sf KP}$\Pi_{N+1}$ proves $V\in RM_{N}(\lceil\ome_{n}(\mathbb{K}+1)\rceil;<^{\veps})$
for each $n\in\ome$.

Let us consider the simplest case $N=3$, i.e., an ordinal analysis of 
set theory ${\sf KP}\Pi_{3}$ for $\Pi_{3}$-reflecting universe.
It turns out that ${\sf KP}\Pi_{3}$ is proof-theoretically reducible to
iterations of recursively Mahlo operations
$V\in RM_{2}(\lceil\ome_{n}(\mathbb{K}+1)\rceil;<^{\veps})\, (n\in\ome)$,
but how to analyze it proof-theoretically?
Here we need a breakthrough done by \cite{Rathjen94}.

\subsection{Ordinals for ${\sf KP}\Pi_{3}$}

In this subsection we define collapsing functions $\psi_{\sig}^{\xi}(a)$ for ${\sf KP}\Pi_{3}$.
It is much easier for us to justify the definitions with an existence of a small large cardinal.
Let $\mathbb{K}$ be the least weakly compact cardinal, i.e., $\Pi^{1}_{1}$-indescribable
cardinal, and $\Ome=\ome_{1}$.
In general for $n\geq 0$, $A\subset ON$ is \textit{$\Pi^{1}_{n}$-indescribable} in an ordinal $\pi$
iff for every $\Pi^{1}_{n}(P)$-formula $\vphi(P)$ with a predicate $P$ and $C\subset\pi$,
if $(L_{\pi},C)\models\vphi(P)$, then
$(L_{\alp},C\cap\alp)\models\vphi(P)$ for an $\alp\in A\cap\pi$.
First let us introduce the Mahlo operation.
Let $A\subset\mathbb{K}$ be a set, and $\alpha\leq\mathbb{K}$ a limit ordinal.
$\alpha\in M_{2}(A)$ iff $A\cap\alpha$ is $\Pi^{1}_{0}$-indescribable in $\alpha$.

As in Definition \ref{df:HKP} we define the Skolem hull $\calh_{a}(X)$ and
simultaneously classes $Mh_{2}^{a}(\xi)$ as follows.

\bdf\label{df:HKPpi3}
{\rm
Define simultaneously by recursion on ordinals $a<\veps_{\mK+1}$ the classes
 $\calh_{a}(X)\, (X\subset \Gam_{\mK+1})$, $Mh_{2}^{a}(\xi)\,(\xi<\veps_{\mK+1})$ and the ordinals $\psi_{\sig}^{\xi}(a)$ as follows.

\benu
\item
$\calh_{a}(X)$ denotes the Skolem hull of $\{0,\Ome,\mK\}\cup X$ under the functions
$+, \vphi$, and
$(\sig,\nu,b)\mapsto\psi_{\sig}^{\nu}(b)\, (b<a)$.

\item
Let for $\xi>0$,
\beqn\label{eq:Mh2}
\pi\in Mh_{2}^{a}(\xi):\Lrarw
\{a,\xi\}\subset\calh_{a}(\pi) \spand
\fal\nu\in\calh_{a}(\pi)\cap\xi\left(
\pi\in M_{2}(Mh_{2}^{a}(\nu))
\right)
\eeqn

$\pi\in Mh_{2}^{a}(0)$ iff $\pi$ is a limit ordinal.

\item
For $0\leq\xi<\veps_{\mK+1}$,
\beqn\label{eq:psiPi3}
\psi_{\pi}^{\xi}(a)=
\min\left(\{\pi\}\cup
\{\kap\in Mh_{2}^{a}(\xi): \{\xi,\pi,a\}\subset\calh_{a}(\kap) \spand
\calh_{a}(\kap)\cap\pi\subset\kap\}
\right)
\eeqn
and
$\psi_{\Ome}(\alp)=\min\{\bet<\Ome:  \calh_{\alp}(\bet)\cap\Ome\subset \bet\}$.
\eenu
}
\edf

We see that each of 
$x=\mathcal{H}_{a}(y)$,
$x=\psi^{\xi}_{\kappa}a$ and
$x\in Mh^{a}_{2}(\xi)$, 
is a $\Sig_{1}$-predicate as fixed points in $\mbox{{\sf ZFL}}$

Since 
the cardinality of the set $\mathcal{H}_{\veps_{\mK+1}}(\pi)$ is $\pi$ for any infinite
cardinal $\pi\leq\mK$,
pick an injection $f :\mathcal{H}_{\veps_{\mK+1}}(\mK)\to \mK$ so that
$f "\mathcal{H}_{\veps_{\mK+1}}(\pi)\subset\pi$ for any weakly inaccessibles $\pi\leq\mK$.

\begin{lemma}\label{lem:welldefinedness}
{\rm (Cf.\,Theorem 4.12 in \cite{Rathjen94}.)}
\benu

\item\label{lem:welldefinedness.3}
There exists a $\Pi^{1}_{1}$-formula $mh^{a}_{2}(x)$ such that
$\pi\in Mh^{a}_{2}(\xi)$ iff $L_{\pi}\models mh^{a}_{2}(\xi)$ 
for any weakly inaccessible cardinals $\pi\leq\mK$ with $f"(\{a,\xi\})\subset L_{\pi}$.

\item\label{lem:welldefinedness.1pi3}
$\mK\in Mh^{a}_{2}(\veps_{\mK+1})\cap M_{2}(Mh^{a}_{2}(\veps_{\mK+1}))$ for every $a<\veps_{\mK+1}$.

\eenu
\end{lemma}
\bprf
\ref{lem:welldefinedness}.\ref{lem:welldefinedness.3}.
Let $\pi$ be a weakly inaccessible cardinal and
$f$ an injection such that $f"\calh_{\veps_{\mK+1}}(\pi)\subset L_{\pi}$.
Assume that $f"(\{a,\xi\})\subset L_{\pi}$.
Then for $f(\xi)\in f" \calh_{a}(\pi)$,
$\pi\in Mh^{a}_{2}(\xi)$ iff
for any $f(\nu)\in L_{\pi}$,
if $f(\nu)\in f"\calh_{a}(\pi)$ and 
$\nu<\xi$, then $\pi\in M_{2}(Mh^{a}_{2}(\nu))$, where
$f"\calh_{a}(\pi)\subset L_{\pi}$ is a class in $L_{\pi}$.
\\

\noindent
\ref{lem:welldefinedness}.\ref{lem:welldefinedness.1pi3}.
We show the following $B(\xi)$ is progressive in $\xi<\veps_{\mK+1}$:
\[
B(\xi)  :\Lrarw 
 \mK\in Mh^{a}_{2}(\xi)\cap M_{2}(Mh^{a}_{2}(\xi))
\]
Note that $\xi\in \mathcal{H}_{a}(\mK)$ holds for any $\xi<\veps_{\mK+1}$.

Suppose $\forall \nu<\xi\, B(\nu)$.
We have to show that $Mh^{a}_{2}(\xi)$ is $\Pi^{1}_{0}$-indescribable in $\mK$.
It is easy to see that if
$\pi\in M_{2}(Mh^{a}_{2}(\xi))$, then $\pi\in Mh^{a}_{2}(\xi)$
by induction on $\pi$.
Let $\tht(P)$ be a first-order formula with a predicate $P$ such that
$(L_{\mK},C)\models\tht(P)$ for $C\subset\mK$.

By IH we have $\forall \nu<\xi[\mK\in M_{2}(Mh^{a}_{2}(\nu))]$.
In other words, $\mK\in Mh^{a}_{2}(\xi)$, i.e., $(L_{\mK},C)\models mh^{a}_{2}(\xi)\land \tht(P)$.
Since the universe $L_{\mK}$ is $\Pi^{1}_{1}$-indescribable, pick a $\pi<\mK$ such that
$(L_{\pi},C\cap\pi)$ enjoys the $\Pi^{1}_{1}$-sentence $mh^{a}_{2}(\xi)\land \tht(P)$, 
and $\{f(a),f(\xi)\}\subset L_{\pi}$.
Therefore  $\pi\in Mh^{a}_{2}(\xi)$.
Thus $\mK\in M_{2}(Mh^{a}_{2}(\xi))$.
\eprf

\blem\label{lem:psiK}
For every $\{a,\xi\}\subset\veps_{\mK+1}$,
$\psi_{\mK}^{\xi}(a)<\mK$ for the $\Pi^{1}_{1}$-indescribable cardinal $\mK$.
\elem
\bprf
Let $\{a,\xi\}\subset\veps_{\mK+1}$.
By Lemma \ref{lem:welldefinedness}.\ref{lem:welldefinedness.1pi3} we obtain
$\mK\in M_{2}(Mh^{a}_{2}(\xi))$.
On the other, $\{\kap<\mK: \{\xi,a\}\subset\calh_{a}(\kap),
\calh_{a}(\kap)\cap\mK\subset\kap\}$ is a club subset of $\mK$.
Hence $\psi_{\mK}^{\xi}(a)<\mK$ by
the definition (\ref{eq:psiPi3}).
\eprf
\\

From the definition (\ref{eq:psiPi3}) we see
\[
\pi\in Mh_{2}^{a}(\mu)\cap \calh_{a}(\pi) \spand \xi\in \calh_{a}(\pi)\cap\mu
\Rarw \pi\in M_{2}(Mh_{2}^{a}(\xi))\spand \psi_{\pi}^{\xi}(a)<\pi
\]
In what follows $M_{2}$ denote the $\Pi_{2}$-recursively Mahlo operation $RM_{2}$.

\subsection{Operator controlled derivations for ${\sf KP}\Pi_{3}$}
$OT(\Pi_{3})$ denotes a computable notation system of ordinals with
collapsing functions $\psi_{\sig}^{\nu}(b)$.
$\kap=\psi_{\sig}^{\nu}(b)\in OT(\Pi_{3})$ 
if $\{\sig,\nu,b\}\subset OT(\Pi_{3})\cap\calh_{b}(\kap)$, 
$\nu=m_{2}(\kap)<m_{2}(\sig)$
and
\beqn\label{eq:OTpi3m2}
SC_{\mK}(\nu)\subset\kap
\spand \nu\leq b
\eeqn
where $m_{2}(\Ome)=1$ and
$m_{2}(\mK)=\veps_{\mK+1}$.
We need the condition (\ref{eq:OTpi3m2}) in our well-foundedness proof of $OT(\Pi_{3})$, 
cf.\,Proposition \ref{prp:m2} and Lemma \ref{th:id5wf21}.

Operator controlled derivations for ${\sf KP}\Pi_{3}$ are defined as in 
Definition \ref{df:controlderreg} for ${\sf KP}\ome$ together with the following inference rules.
For ordinals $\pi=\psi_{\sig}^{\xi}(a)$, let
$m_{2}(\pi)=\xi$.

\bdes
\item[$({\rm rfl}_{\Pi_{3}}(\mK))$]
$b\geq\mK$.
There exist 
an ordinal $a_{0}\in\calh_{\gam}[\Tht]\cap a$,
and a $\Sig_{3}(\mK)$-sentence $A$ enjoying the following conditions:

\[
\infer[({\rm rfl}_{\Pi_{3}}(\mK))]{
\calh_{\gam}[\Tht]\vdash^{a}_{b}\Gam}
{
\mathcal{H}_{\gam}[\Tht]\vdash^{a_{0}}_{b}\Gamma, \lnot A
&
\{
\mathcal{H}_{\gam}[\Tht\cup \{\rho\}]\vdash^{a_{0}}_{b}\Gamma, 
A^{(\rho,\mK)}: \rho<\mK\}
}
\]
The inference says that $\mK\in RM_{3}$.

\item[$({\rm rfl}_{\Pi_{2}}(\alp,\pi,\nu))$]
There exist ordinals 
$\alp<\pi\leq b<\mK$, 
$\nu< 
m_{2}(\pi)$ such that $SC_{\mK}(\nu)\subset\pi$ and 
$\nu\leq\gam$, cf.\,(\ref{eq:OTpi3m2}), 
$a_{0}<a$,
and a finite set $\Del$ of $\Sig_{2}(\pi)$-sentences enjoying the following conditions:

\benu

\item
$\{\alp,\pi,
\nu\}\subset\calh_{\gam}[\Tht]$.

 \item
For each $\del\in\Del$,
$
\mathcal{H}_{\gam}[\Tht]\vdash^{a_{0}}_{b}\Gamma, \lnot\del
$.

\item
For each
$\alp<\rho\in Mh_{2}(\nu)\cap\pi$,
$\mathcal{H}_{\gam}[\Tht\cup\{\rho\}]\vdash^{a_{0}}_{b}\Gamma, 
\Del^{(\rho,\pi)}$ holds.

By $\rho\in Mh_{2}(\nu)$ we mean
$\nu\leq m_{2}(\rho)$.

\eenu
{\small
\[
\hspace{-15mm}
\infer[({\rm rfl}_{\Pi_{2}}(\alp,\pi,\nu))]{
\calh_{\gam}[\Tht]\vdash^{a}_{b}\Gam}
{
\{
\mathcal{H}_{\gam}[\Tht]\vdash^{a_{0}}_{b}\Gamma, \lnot\del
\}_{\del\in\Del}
&
\{
\mathcal{H}_{\gam}[\Tht\cup\{\rho\}]\vdash^{a_{0}}_{b}\Gamma, 
\Del^{(\rho,\pi)}: \alp<\rho\in Mh_{2}(\nu)\cap \pi
\}
}
\]
}
The inference says that 
$\pi\in M_{2}(Mh_{2}^{\gam}(\nu))$
provided that $\{m_{2}(\pi),\gam, \nu\}\subset\calh_{\gam}(\pi)$.
\edes

The axiom for $\Pi_{3}$-reflection follows from the inference $({\rm rfl}_{\Pi_{3}}(\mK))$ as follows.
Let $A\in\Sig_{3}(\mK)$ with $d=\rk(A)<\mK+\ome$, and 
$d_{\rho}=\rk(A^{(\rho,\mK)})$ for
$\rho<\mK$.
\[
\infer[({\rm rfl}_{\Pi_{3}}(\mK))]{\calh_{0}[\sfk(A)]\vdash^{\mK+\ome}_{\mK}\lnot A,\exi z\, A^{(z,\mK)}}
{
 \calh_{0}[\sfk(A)]\vdash^{2d}_{0} A,\lnot A
 &
 \infer{\calh_{0}[\sfk(A)\cup\{\rho\}]\vdash^{\mK}_{0}\exi z\, A^{(z,\mK)},\lnot A^{(\rho,\mK)}}
 {
 \calh_{0}[\sfk(A)\cup\{\rho\}]\vdash^{2d_{\rho}}_{0} A^{(\rho,\mK)},\lnot A^{(\rho,\mK)}
 }
}
\]

\begin{quote}
An appropriate name for this collapsing technique would be stationary collapsing since in order for
this procedure to work, a single derivation has to be collapsed into a ``stationary" family of derivations.
\cite{Rathjen94}
\end{quote}

We see from the following proof that $\alp=\psi_{\mK}(\gam+\mK)$ holds
in
every inference $({\rm rfl}_{\Pi_{2}}(\alp,\kap,a_{0}))$ 
occurring
in a witnessed derivation of
$\calh_{\hat{a}+1}[\Tht\cup\{\kap\}]\vdash^{\bet}_{\bet}\Gam^{(\kap,\mK)}$.
Let us call the unique ordinal $\alp$ a \textit{base}.

\blem\label{lem:lowerPi3}
Assume $\Gam\subset\Sig_{2}(\mK)$, $\Tht\subset\calh_{\gam}(\psi_{\mK}(\gam))$,
and
$\calh_{\gam}[\Tht]\vdash^{a}_{\mK}\Gam$ with $a\leq\gam$.
Then 
$\calh_{\hat{a}+1}[\Tht\cup\{\kap\}]\vdash^{\bet}_{\bet}\Gam^{(\kap,\mK)}$ holds for any 
$\kap\in Mh_{2}(a)\cap\psi_{\mK}(\gam+\mK\cdot\ome)$ such that
$\psi_{\mK}(\gam+\mK)<\kap$,
where 
$\hat{a}=\gam+\ome^{\mK+a}$ and $\bet=\psi_{\mK}(\hat{a})$.
\elem
\bprf
By induction on $a$.
Note that there exists a $\kap\in OT(\Pi_{3})$ such that 
$\psi_{\mK}(\gam+\mK)<\kap\in Mh_{2}(a)\cap\psi_{\mK}(\gam+\mK\cdot\ome)$.
F.e. $\kap=\psi_{\mK}^{a}(\gam+\mK+1)$.
\\
\textbf{Case 1}.
Consider the case when the last inference is a $({\rm rfl}_{\Pi_{3}}(\mK))$.
For $\Sig_{3}\ni A\simeq\bigvee(A_{i})_{i\in J}$,
\[
\infer[({\rm rfl}_{\Pi_{3}}(\mK))]{
\calh_{\gam}[\Tht]\vdash^{a}_{\mK}\Gam}
{
\mathcal{H}_{\gam}[\Tht]\vdash^{a_{0}}_{\mK}\Gamma, \lnot A
&
\{
\mathcal{H}_{\gam}[\Tht\cup \{\rho\}]\vdash^{a_{0}}_{\mK}\Gamma, 
A^{(\rho,\mK)}: \rho<\mK\}
}
\]
Let 
\[
\psi_{\mK}(\gam+\mK)\leq\sig\in Mh_{2}(a_{0})\cap\kap 
.\]
Let $i\in Tm(\sig)$, i.e., $\sfk(i)\subset\sig$.
For each $i\in Tm(\sig)$ Inversion yields
$\mathcal{H}_{\gam+|i|}[\Tht\cup\sfk(i)]\vdash^{a_{0}}_{\mK}\Gamma, \lnot A_{i}$
with $\sfk(i)<\psi_{\mK}(\gam+|i|)$.
By IH we obtain
$\mathcal{H}_{\hat{a}+1}[\Tht\cup \{\sig\}\cup\sfk(i)]\vdash^{\bet_{0}}_{\bet}\Gamma^{(\sig,\mK)}, \lnot A_{i}^{(\sig,\mK)}$ for every $i\in Tm(\sig)$, where
$\bet_{0}=\psi_{\mK}(\widehat{a_{0}})$ with $\widehat{a_{0}}=\gam+\ome^{\mK+a_{0}}=\gam+|i|+\ome^{\mK+a_{0}}$.
A $(\bigwedge)$ yields
\[
\mathcal{H}_{\hat{a}+1}[\Tht\cup\{\sig\}]\vdash^{\bet_{0}+1}_{\bet}\Gamma^{(\sig,\mK)}, \lnot A^{(\sig,\mK)}
\]
On the other hand we have
$\mathcal{H}_{\gam+\sig}[\Tht\cup \{\sig\}]\vdash^{a_{0}}_{\mK}\Gamma, 
A^{(\sig,\mK)}$ with $\sig\in\calh_{\gam+\sig}(\psi_{\mK}(\gam+\sig))$, but
$\sig\not\in\calh_{\gam}(\psi_{\mK}(\gam+\mK))$.
We obtain 
$\kap\in Mh_{2}(a_{0})$
by $a_{0}<a$, and $\gam+\sig+\mK=\gam+\mK$.
IH yields

\[
\mathcal{H}_{\hat{a}+1}[\Tht\cup \{\kap,\sig\}]\vdash^{\bet_{0}}_{\bet}\Gamma^{(\kap,\mK)}, 
A^{(\sig,\mK)}
\]
A $(cut)$ of the cut formula $A^{(\sig,\mK)}$ with $\rk(A^{(\sig,\mK)})<\kap<\psi_{\mK}(\gam+\mK\cdot\ome)\leq\bet$ yields
\[
\mathcal{H}_{\hat{a}+1}[\Tht\cup \{\kap,\sig\}]\vdash^{\bet_{0}+2}_{\bet}\Gamma^{(\kap,\mK)},\Gam^{(\sig,\mK)}
\]

On the other side
\[
\mathcal{H}_{\gam}[\Tht\cup\{\kap\}]\vdash^{2d}_{0}\lnot\tht^{(\kap,\mK)},\Gam^{(\kap,\mK)}
\]
holds
for each $\tht\in\Gam\subset\Sig_{2}(\mK)$, where
$d=\max\{\rk(\tht^{(\kap,\mK)}):\tht\in\Gam\}<\kap+\ome<\bet$.

Moreover we have $a_{0}<\hat{a}$,
$SC_{\mK}(a_{0})\subset\calh_{\gam}[\Tht]\cap\mK\subset\calh_{\gam}(\psi_{\mK}(\gam))\cap\mK\subset\kap$.
A $({\rm rfl}_{\Pi_{2}}(\del,\kap,a_{0}))$ with $\del=\psi_{\mK}(\gam+\mK)$,
$\{\del,\kap,a_{0}\}\subset\calh_{\hat{a}+1}[\Tht\cup\{\kap\}]$
yields
$\mathcal{H}_{\hat{a}+1}[\Tht\cup\{\kap\}]\vdash^{\bet}_{\bet}\Gam^{(\kap,\mK)}$.
{
\small
\[
\hspace{-13mm}
\infer[({\rm rfl}_{\Pi_{2}}(\del,\kap,a_{0}))]{\mathcal{H}_{\hat{a}+1}[\Tht\cup\{\kap\}]\vdash^{\bet}_{\bet}\Gam^{(\kap,\mK)}}
{
\{\mathcal{H}_{\gam}[\Tht\cup\{\kap\}]\vdash^{2d}_{0}\lnot\tht^{(\kap,\mK)},\Gam^{(\kap,\mK)}\}_{\tht\in\Gam}
&
\hspace{-2mm}
\infer{
\{\mathcal{H}_{\hat{a}+1}[\Tht\cup\{\kap,\sig\}]\vdash^{\bet_{0}+2}_{\bet}\Gamma^{(\kap,\mK)},\Gam^{(\sig,\mK)}\}_{
\del<\sig\in Mh_{2}(a_{0})\cap\kap
 }
}
 {
 \mathcal{H}_{\hat{a}+1}[\Tht\cup\{\sig\}]\vdash^{\bet_{0}+1}_{\bet}\Gamma^{(\sig,\mK)}, \lnot A^{(\sig,\mK)}
 &
 \mathcal{H}_{\hat{a}+1}[\Tht\cup\{\kap,\sig\}]\vdash^{\bet_{0}}_{\bet}\Gamma^{(\kap,\mK)}, 
A^{(\sig,\mK)}
 }
}
\]
}
\\
\textbf{Case 2}.
The last inference is a $(cut)$ of a cut formula $C$ with $\rk(C)<\mK$.
Then
$\rk(C)\in\calh_{\gam}[\Tht]\cap\mK\subset\psi_{\mK}(\gam)<\bet$ by (\ref{eq:controlderKP}),
Proposition \ref{prp:impdegreeKP}.\ref{prp:impdegreeKP.-2} and the assumption 
$\Tht\subset\calh_{\gam}(\psi_{\mK}(\gam))$.
\\
\textbf{Case 3}.
The last inference is a $(\bigwedge)$ with a main formula $\Pi_{1}(\mK)\ni A\simeq\bigwedge(A_{\iota})_{\iota\in J}$.
We may assume $J=Tm(\mK)$.
Then $A^{(\kap,\mK)}\simeq\bigwedge(A_{\iota})_{\iota\in Tm(\kap)}$, and
we obtain the lemma by pruning the branches for $\iota\not\in Tm(\kap)$.
\\
\textbf{Case 4}.
The last inference is a $(\bigvee)$ with a main formula 
$\Sig_{2}(\mK)\ni A\simeq\bigvee(A_{\iota})_{\iota\in J}$.
We may assume $J=Tm(\mK)$.
Then $A^{(\kap,\mK)}\simeq\bigvee(A_{\iota})_{\iota\in Tm(\kap)}$.
\[
\infer[(\bigvee)]{\calh_{\gam}[\Tht]\vdash^{a}_{\mK}\Gam}
{
\mathcal{H}_{\gam}[\Tht]\vdash^{a_{0}}_{\mK}\Gamma, A_{\iota}
}
\]
We may assume that $\sfk(\iota)\subset\sfk(A_{\iota})$.
Then by (\ref{eq:controlderKP}) and $\Tht\subset\calh_{\gam}(\psi_{\mK}(\gam))$ we obtain
$\sfk(\iota)\subset\calh_{\gam}[\Tht]\cap\mK\subset\calh_{\gam}(\psi_{\mK}(\gam))\cap\mK\subset\kap$,
and $\iota\in Tm(\kap)$.
\eprf
\\

An ordinal term $\alp$ in $OT(\Pi_{3})$ is said to be \textit{regular} if
either $\alp\in\{\Ome,\mK\}$ or $\alp=\psi_{\sig}^{\nu}(a)$ for some
$\sig,a$ and $\nu>0$.

\blem\label{lem:lowerPi2}
Let $\lam$ be regular, $\Gam\subset\Sig_{1}(\lam)$ and
$\calh_{\gam}[\Tht]\vdash^{a}_{b}\Gam$, where $a<\mK$, 
$\calh_{\gam}[\Tht]\ni\lam\leq b<\mK$, and
$\fal\kap\in[\lam,b)(\Tht\subset\calh_{\gam}(\psi_{\kap}(\gam)))$.
Let $\hat{a}=\gam+\tht_{b}(a)$ and $\bet=\psi^{\eta}_{\lam}(\hat{a})$ such that
$0\leq\eta\in\calh_{\gam}[\Tht]$, $\eta<m_{2}(\lam)$, $SC_{\mK}(\eta)\subset\bet$ and $\eta\leq\gam$.
Then
$\calh_{\hat{a}+1}[\Tht]\vdash^{\bet}_{\bet}\Gam$ holds.
\elem
\bprf 
By main induction on $b$ with subsidiary induction on $a$ as in Theorem \ref{th:CollapsingthmKP}.
\\
\textbf{Case 1}.
Consider first the case when the last inference is a 
$({\rm rfl}_{\Pi_{2}}(\alp,\sig,\nu))$ with $b\geq\sig>\alp$.
{\small
\[
\hspace{-5mm}
\infer[({\rm rfl}_{\Pi_{2}}(\alp,\sig,\nu))]{
\calh_{\gam}[\Tht]\vdash^{a}_{b}\Gam}
{
\{
\mathcal{H}_{\gam}[\Tht]\vdash^{a_{0}}_{b}\Gamma, \lnot\del
\}_{\del\in\Del}
&
\{
\mathcal{H}_{\gam}[\Tht\cup\{\rho\}]\vdash^{a_{0}}_{b}\Gamma, 
\Del^{(\rho,\sig)}: \alp<\rho\in Mh_{2}(\nu)\cap \sig\}
}
\]
}
where $\Del\subset\Sig_{2}(\sig)$, 
$\{\alp,\sig,\nu\}\subset\calh_{\gam}[\Tht]$,
$\nu<
m_{2}(\sig)$, $\nu\leq\gam$ and
$SC_{\mK}(\nu)\subset\sig$.
\\
\textbf{Case 1.1}. $\sig<\lam$:
Then $\{\lnot\del\}\cup\Del^{(\rho,\sig)}\subset\Del_{0}(\lam)$ for each $\del\in\Del$.
For any $\lam\leq\kap<b$, we obtain $\rho<\sig\in\mathcal{H}_{\gam}[\Tht]\cap\kap\subset\psi_{\kap}(\gam)$.
SIH yields the lemma.
\\
\textbf{Case 1.2}. $\sig\geq\lam$:
For each $\del\in\Del$, let $\del\simeq\bigvee(\del_{i})_{i\in J}$. We may assume
$J=Tm(\sig)$.
Inversion yields $\mathcal{H}_{\gam+|i|}[\Tht\cup\sfk(i)]\vdash^{a_{0}}_{b}\Gamma, \lnot\del_{i}$.
Let $\widehat{{a}_{0}}=\gam+\tht_{b}(a_{0})$ and
$\rho=\psi_{\sig}^{\nu}(\widehat{{a}_{0}}+\alp)$,
where $\Tht\subset\calh_{\gam}(\rho)$ by the assumption,
$\{\alp,\sig,\nu,\widehat{{a}_{0}}\}\subset\calh_{\gam}[\Tht]$ with $\nu<m_{2}(\sig)$.
Hence 
$\{\alp,\sig,\nu,\widehat{{a}_{0}}\}\subset\calh_{\gam}(\rho)$ and 
$\alp<\rho$ by $\alp<\sig$.
Therefore, cf.\,(\ref{eq:OTpi3m2}), $SC_{\mK}(\nu)\subset\rho\in Mh_{2}(\nu)\cap\sig\cap\calh_{\widehat{a_{0}}+\alp+1}[\Tht]$.

For each $\sfk(i)\subset\rho$ and $\lnot\del_{i}\in \Sig_{1}(\sig)$, we obtain
$\gam+|i|+\tht_{b}(a_{0})=\widehat{a_{0}}$ by $|i|<\rho<\sig\leq b$, and
$\mathcal{H}_{\widehat{a_{0}}+1}[\Tht\cup\sfk(i)]\vdash^{\rho_{0}}_{\rho_{0}}\Gamma, \lnot\del_{i}$ by SIH
for $\rho_{0}=\psi_{|sig}(\widehat{a_{0}})\leq\rho$.
Hence
$\mathcal{H}_{\widehat{a_{0}}+\alp+1}[\Tht\cup\sfk(i)]\vdash^{\rho}_{\rho}\Gamma, \lnot\del_{i}$
By Boundedness we obtain
$\mathcal{H}_{\widehat{a_{0}}+\alp+1}[\Tht\cup\sfk(i)]\vdash^{\rho}_{\rho}\Gamma, \lnot\del_{i}^{(\rho,\sig)}$.
A $(\bigwedge)$ yields
\[
\mathcal{H}_{\widehat{a_{0}}+\alp+1}[\Tht]\vdash^{\rho+1}_{\rho}\Gamma, \lnot\del^{(\rho,\sig)}
.\]

On the other hand we have
$\mathcal{H}_{\gam}[\Tht\cup\{\rho\}]\vdash^{a_{0}}_{b}\Gamma, 
\Del^{(\rho,\sig)}$, and 
$\mathcal{H}_{\widehat{a_{0}}+\alp+1}[\Tht]\vdash^{a_{0}}_{b}\Gamma, 
\Del^{(\rho,\sig)}$.
By SIH we obtain
\[
\mathcal{H}_{\widehat{a_{1}}+1}[\Tht]\vdash^{\bet_{1}}_{\bet_{1}}\Gamma, 
\Del^{(\rho,\sig)}
\]
 for 
$\bet_{1}=\psi_{\sig}(\widehat{a}_{1})>\rho$, with
$\widehat{a_{1}}=\widehat{a_{0}}+\alp+\tht_{b}(a_{0})\leq\gam+\tht_{b}(a_{0})\cdot 3<\hat{a}$.
Therefore we obtain
$\mathcal{H}_{\widehat{a_{1}}+1}[\Tht]\vdash^{\bet_{1}+\ome}_{\bet_{1}}\Gamma$
by several $(cut)$'s of $\rk(\del^{(\rho,\sig)})<\rho+\ome<\bet_{1}$.

If $\sig=\lam$, then we are done.
Let $\lam<\sig\leq b$. Then $\lam\in\calh_{\gam}[\Tht]\cap\sig\subset\bet_{1}$.
MIH yields
$\mathcal{H}_{\widehat{a_{2}}+1}[\Tht]\vdash^{\bet_{2}}_{\bet_{2}}\Gamma$,
where
$\widehat{a_{2}}=\widehat{a_{1}}+\tht_{\bet_{1}}(\bet_{1}+\ome)<\hat{a}$
by $\bet_{1}<\sig\leq b$,
and
$\bet_{2}=\psi_{\lam}(\widehat{a_{2}})<\psi_{\lam}(\hat{a})\leq\bet$.
\\
\textbf{Case 2}.
Next the last inference is a $(cut)$ of a cut formula $C$ with $d=\rk(C)<b$.
\[
\infer[(cut)]{\calh_{\gam}[\Tht]\vdash^{a}_{b}\Gam}
{
\calh_{\gam}[\Tht]\vdash^{a_{0}}_{b}\Gam,\lnot C
&
\calh_{\gam}[\Tht]\vdash^{a_{0}}_{b}\Gam, C
}
\]
If $d<\lam$, then SIH yields the lemma. Let $\lam\leq d$ and
$\widehat{a_{0}}=\gam+\tht_{b}(a_{0})$.
\\
\textbf{Case 2.1}.
There exists a regular $\sig\in\calh_{\gam}[\Tht]$ such that
$d<\sig\leq b$:
For $\{\lnot C,C\}\subset\Del_{0}(\sig)$, we obtain
$\calh_{\widehat{a_{0}}+1}[\Tht]\vdash^{\bet_{0}}_{\bet_{0}}\Gam,C$
and
$\calh_{\widehat{a_{0}}+1}[\Tht]\vdash^{\bet_{0}}_{\bet_{0}}\Gam,\lnot C$ for $\bet_{0}=\psi_{\sig}(\widehat{a_{0}})$ by SIH.
A $(cut)$ yields
$\calh_{\widehat{a_{0}}+1}[\Tht]\vdash^{\bet_{0}+1}_{\bet_{0}}\Gam$.
MIH yields
$\calh_{\widehat{a_{1}}+1}[\Tht]\vdash^{\bet_{1}}_{\bet_{1}}\Gam$,
where $\widehat{a_{1}}=\widehat{a_{0}}+\tht_{\bet_{0}}(\bet_{0}+1)<\hat{a}$
and $\bet_{1}=\psi_{\lam}(\widehat{a_{1}})<\psi_{\lam}(\hat{a})\leq\bet$.
\\
\textbf{Case 2.2}. Otherwise:
Then there is no regular $\sig\in\calh_{\gam}[\Tht]$ such that $d<\sig\leq b$.
Let $d+c=b$.
Then by Cut-elimination we obtain
$\calh_{\gam}[\Tht]\vdash^{\tht_{c}(a)}_{d}\Gam$.
MIH yields
$\calh_{\hat{a}+1}[\Tht]\vdash^{\psi_{\lam}(\hat{a})}_{\psi_{\lam}(\hat{a})}\Gam$, where
$\gam+\tht_{d}(\tht_{c}(a))=\gam+\tht_{b}(a)=\hat{a}$.
\eprf

\bth\label{th:Pi3}
Assume ${\sf KP}\Pi_{3}\vdash\tht^{L_{\Ome}}$ for $\tht\in\Sig$.
Then there exists an $n<\ome$ such that
$L_{\alp}\models\tht$ for $\alp=\psi_{\Ome}(\ome_{n}(\mK+1))$ in $OT(\Pi_{3})$.
\end{theorem}
\bprf
By Embedding there exists an $m>0$ such that
$\calh_{0}[\emptyset]\vdash^{\mK+m}_{\mK+m}\tht^{L_{\Ome}}$.
By Cut-elimination,
$\calh_{0}[\emptyset]\vdash^{a}_{\mK}\tht^{L_{\Ome}}$
and
$\calh_{a}[\emptyset]\vdash^{a}_{\mK}\tht^{L_{\Ome}}$
for $a=\ome_{m}(\mK+m)$.
By Lemma \ref{lem:lowerPi3} we obtain
$\calh_{\ome^{a}+1}[\{\kap\}]\vdash^{\bet}_{\bet}\tht^{L_{\Ome}}$,
where $\bet=\psi_{\mK}(\ome^{a})$, 
$a+\ome^{\mK+a}=\ome^{a}$, $(\tht^{L_{\Ome}})^{(\kap,\mK)}\equiv \tht^{L_{\Ome}}$ and
$\psi_{\mK}(a+\mK)<\kap\in Mh_{2}(a)\cap\psi_{\mK}(a+\mK\cdot\ome)$.
F.e. $\kap=\psi_{\mK}^{a}(a+\mK+1)\in \calh_{a+\mK+2}[\emptyset]$.
Hence
$\calh_{\ome^{a}+\mK+2}[\emptyset]\vdash^{\bet}_{\bet}\tht^{L_{\Ome}}$.
Lemma \ref{lem:lowerPi2} then yields
$\calh_{\gam+1}[\emptyset]\vdash^{\bet_{1}}_{\bet_{1}}\tht^{L_{\Ome}}$
for $\gam=\ome^{a}+\mK+\tht_{\bet}(\bet)$ and
$\bet_{1}=\psi_{\Ome}(\gam)<\psi_{\Ome}(\ome^{a}+\mK\cdot 2)<\psi_{\Ome}(\ome_{m+2}(\mK+1))=\alp$.
Therefore
$L_{\alp}\models\tht$.
\eprf

\section{Well-foundedness proof in ${\sf KP}\Pi_{3}$}\label{sec:wfpi3}
$OT(\Pi_{3})$ denotes the computable notation system in section \ref{sec:pi3}.
$\kap=\psi_{\sig}^{\nu}(b)\in OT(\Pi_{3})$ only if $\nu=m_{2}(\kap)<m_{2}(\sig)$, 
$SC_{\mK}(\nu)\subset\kap$ and $\nu\leq b$, cf.\,(\ref{eq:OTpi3m2}).
In this section we show the
\bth\label{th:wfpi3}
${\sf KP}\Pi_{3}$ proves
 the well-foundedness of $OT(\Pi_{3})$
up to \textit{each} $\alp<\psi_{\Ome}(\veps_{\mK+1})$.
\end{theorem}
We assume a standard encoding $OT(\Pi_{3})\ni\alp\mapsto\lc\alp\rc\in\ome$,
and identify ordinal terms $\alp$ with its code $\lc\alp\rc$.

\subsection{Distinguished sets}\label{subsec:CX}

In this subsection we work in ${\sf KP}\ell$.

\bdf\label{df:CX}{\rm \cite{BuchholzBSL}.}\\
{\rm For} $\alp\in OT(\Pi_{3}), X\subset OT(\Pi_{3})${\rm , let} 
\beqnarr
\calc^{\alp}(X) & := & 
\mbox{{\rm closure of }} \{0,\Ome,\mK\}\cup(X\cap\alp) \mbox{ {\rm under} } +,\vphi
\nonumber \\
&& \mbox{ {\rm and} } (\sig,\alp,\nu)\mapsto \psi_{\sig}^{\nu}(\alp) \mbox{ {\rm for} } \sig>\alp \mbox{ {\rm in }} OT(\Pi_{3})
\label{eq:CX}
\eeqnarr
\edf

$\alpha^{+}=\Omega_{a+1}$ denotes the least regular term above $\alpha$ if such a term exists. 
Otherwise $\alpha^{+}:=\infty$.

\begin{proposition}\label{lem:CX2} 
Assume $\forall\gamma\in X[ \gamma\in\mathcal{C}^{\gamma}(X)]$ for a set $X\subset OT(\Pi_{3})$.

\benu
\item\label{lem:CX2.3} 
$\alpha\leq\beta \Rarw \mathcal{C}^{\beta}(X)\subset \mathcal{C}^{\alpha}(X)$.

\item\label{lem:CX2.4} 
$\alpha<\beta<\alpha^{+} \Rarw \mathcal{C}^{\beta}(X)=\mathcal{C}^{\alpha}(X)$.
\eenu
\end{proposition}
\bprf
\ref{lem:CX2}.\ref{lem:CX2.3}.
We see by induction on $\ell\gamma\,(\gamma\in OT(\Pi_{3}))$ that
\begin{equation}\label{eq:CX2.3}
\forall\beta\geq\alpha[\gamma\in \mathcal{C}^{\beta}(X) \Rarw \gamma\in \mathcal{C}^{\alpha}(X)\cup(X\cap\beta)]
\end{equation}
For example, if $\psi_{\pi}^{\nu}(\delta)\in \mathcal{C}^{\beta}(X)$ with $\pi>\beta\geq\alpha$
and $\{\pi,\delta,\nu\}\subset\mathcal{C}^{\alpha}(X)\cup(X\cap\beta)$, then 
$\pi\in\mathcal{C}^{\alpha}(X)$, and 
for any $\gamma\in\{\delta,\nu\}$, either $\gamma\in\mathcal{C}^{\alpha}(X)$ or
$\gamma\in X\cap\beta$. If $\gamma<\alpha$, then $\gamma\in X\cap\alpha\subset\mathcal{C}^{\alpha}(X)$.
If $\alpha\leq\gamma\in X\cap\beta$, then $\gamma\in\mathcal{C}^{\gamma}(X)$ by the assumption, and
by IH we have $\gamma\in\mathcal{C}^{\alpha}(X)\cup(X\cap\gamma)$, i.e., $\gamma\in\mathcal{C}^{\alpha}(X)$.
Therefore $\{\pi,\delta,\nu\}\subset\mathcal{C}^{\alpha}(X)$, and 
$\psi_{\pi}^{\nu}(\delta)\in\mathcal{C}^{\alpha}(X)$.

Using (\ref{eq:CX2.3}) we see from the assumption that
$
\forall\beta\geq\alpha[ \gamma\in\mathcal{C}^{\beta}(X) \Rarw \gamma\in\mathcal{C}^{\alpha}(X)]
$.
\\
\ref{lem:CX2}.\ref{lem:CX2.4}.
Assume $\alpha<\beta<\alpha^{+}$. Then by Proposition \ref{lem:CX2}.\ref{lem:CX2.3} we have
$\mathcal{C}^{\beta}(X)\subset\mathcal{C}^{\alpha}(X)$.
$\mathcal{C}^{\alpha}(X) \subset\mathcal{C}^{\beta}(X)$ is easily seen from $\beta<\alpha^{+}$.
\eprf

\begin{definition}\label{df:wftg}
{\rm
\benu
\item 
$Prg[X,Y] :\Lrarw \forall\alpha\in X(X\cap\alpha\subset Y \to \alpha\in Y)$.

\item 
For a definable class $\mathcal{X}$, $TI[\mathcal{X}]$ denotes the schema:\\
$TI[\mathcal{X}] :\Lrarw Prg[\mathcal{X},\mathcal{Y}]\to \mathcal{X}\subset\mathcal{Y} \mbox{ {\rm holds for} any definable classes } \mathcal{Y}$.
\item
For $X\subset OT(\Pi_{3})$, $W(X)$ denotes the \textit{well-founded part} of $X$. 
\item 
$Wo[X] : \Lrarw X\subset W(X)$.
\eenu
}
\end{definition}
Note that for $\alpha\in OT(\Pi_{3})$,
$W(X)\cap\alpha=W(X\cap\alpha)$.

\begin{definition} \label{df:3wfdtg32}
{\rm 
For $X\subset OT(\Pi_{3})$ and
$\alpha\in OT(\Pi_{3})$,
\benu
\item\label{df:3wfdtg.832}
\begin{equation}\label{eq:distinguishedclass}
D[X] :\Lrarw 
\forall\alpha(\alpha\leq X\to W(\mathcal{C}^{\alpha}(X))\cap\alpha^{+}= X\cap\alpha^{+})
\end{equation}

A set $X$ is said to be a \textit{distinguished set} if $D[X]$.

\item\label{df:3wfdtg.932}
$\mathcal{W}:=\bigcup\{X :D[X]\}$.

\eenu
}
\end{definition}

Let $\alpha\in X$ for a distinguished set $X$. 
Then $W(\mathcal{C}^{\alpha}(X))\cap\alpha^{+}= X\cap\alpha^{+}$.
Hence $X$ is a well order.

\begin{proposition}\label{lem:3.11.632}
Let $X$ be a distinguished set. Then
$\alpha\in X \Rarw \forall\beta[\alpha\in \mathcal{C}^{\beta}(X)]$.
\end{proposition}
\bprf
Let $D[X]$ and $\alpha\in X$. 
Then $\alpha\in X\cap\alpha^{+}=W(\mathcal{C}^{\alpha}(X))\cap\alpha^{+}\subset \mathcal{C}^{\alpha}(X)$.
Hence $\forall\gamma\in X(\gamma\in\mathcal{C}^{\gamma}(X))$, and
$\alpha\in\mathcal{C}^{\beta}(X)$ for any $\beta\leq\alpha$ by Proposition \ref{lem:CX2}.\ref{lem:CX2.3}. 
Moreover for $\beta>\alpha$ we have $\alpha\in X\cap\beta\subset\mathcal{C}^{\beta}(X)$.
\eprf

\begin{proposition}\label{lem:5uv.232general}
$X\cap\alpha=Y\cap\alpha \Rarw \forall\beta<\alpha^{+}
\left[\mathcal{C}^{\beta}(X)=\mathcal{C}^{\beta}(Y)\right]$
if $\fal\gam\in X(\gam\in\mathcal{C}^{\gam}(X))$ and $\fal\gam\in Y(\gam\in\mathcal{C}^{\gam}(Y))$.
\end{proposition}
\bprf
Assume that $X\cap\alpha=Y\cap\alpha$ and $\alpha\leq\beta<\alpha^{+}$.
We obtain $\mathcal{C}^{\alpha}(X)=\mathcal{C}^{\alpha}(Y)$.
On the other hand we have $\mathcal{C}^{\beta}(X)=\mathcal{C}^{\alpha}(X)$ and similarly for 
$\mathcal{C}^{\beta}(Y)$
by Proposition \ref{lem:CX2}.\ref{lem:CX2.4}.
Hence $\mathcal{C}^{\beta}(X)=\mathcal{C}^{\beta}(Y)$.
\eprf

\begin{proposition}\label{lem:3wf5}
$\alpha\leq X \spand \alpha\leq Y \Rarw X\cap\alpha^{+}=Y\cap\alpha^{+}$
if $D[X]$ and $D[Y]$.
\end{proposition}
\bprf
For distinguished set $X$, $\alpha\leq X \Rarw X\cap\alpha^{+}=W(\mathcal{C}^{\alp}(X))\cap\alp^{+}$.
Hence the proposition follows from
Propositions \ref{lem:3.11.632} and \ref{lem:5uv.232general}.
\eprf

\begin{proposition}\label{lem:3wf6max}
$\mathcal{W}$ is the maximal distinguished class.
\end{proposition}
\bprf
First we show $\fal\gam\in\calw(\gam\in\mathcal{C}^{\gam}(\calw))$.
Let $\gam\in\calw$, and pick a distinguished set $X$ such that $\gam\in X$.
Then $\gam\in\mathcal{C}^{\gam}(X)\subset\mathcal{C}^{\gam}(\calw)$ by $X\subset\calw$.

Let $\alp\leq\calw$. Pick a distinguished set $X$ such that $\alp\leq X$.
We claim that $\calw\cap\alp^{+}=X\cap\alp^{+}$.
Let $Y$ be a distinguished set and $\bet\in Y\cap\alp^{+}$.
Then $\bet\in Y\cap\bet^{+}=X\cap\bet^{+}$ by Proposition \ref{lem:3wf5}.
The claim yields
$W(\mathcal{C}^{\alp}(\calw))\cap\alp^{+}=W(\mathcal{C}^{\alp}(X))\cap\alp^{+}=
X\cap\alp^{+}=\calw\cap\alp^{+}$.
Hence $D[\mathcal{W}]$.
\eprf

\begin{definition}\label{df:calg}
$\mathcal{G}(X):=\{\alpha\in OT(\Pi_{3}) :\alpha\in \mathcal{C}^{\alpha}(X) \spand
 \mathcal{C}^{\alpha}(X)\cap\alpha\subset X\}$.
\end{definition}

\begin{lemma}\label{lem:wf5.332}
For $D[X]$,
$X\subset\mathcal{G}(X)$.
\end{lemma}
\bprf
Let $\gam\in X$.
We have $\gamma\in W(\mathcal{C}^{\gamma}(X))\cap\gamma^{+}=X\cap\gamma^{+}$.
Hence $\gamma\in\mathcal{C}^{\gamma}(X)$.
Assume $\alpha\in\mathcal{C}^{\gamma}(X)\cap\gamma$. 
Then $\alpha\in W(\mathcal{C}^{\gamma}(X))\cap\gamma^{+}\subset X$. 
Therefore
$\mathcal{C}^{\gamma}(X)\cap\gamma\subset X$.
\eprf

\begin{definition}\label{df:EGFk}

{\rm 
For ordinal terms $\alpha,\delta\in OT(\Pi_{3})$, finite sets 
$G_{\del}(\alpha)\subset OT(\Pi_{3})$ 
 are defined recursively as follows.
\benu

\item
$G_{\del}(\alpha)=\emptyset$ for $\alpha\in\{0,\Omega,\mK\}$.
$G_{\del}(\alpha_{m}+\cdots+\alpha_{0})=\bigcup_{i\leq m}G_{\del}(\alpha_{i})$.
$G_{\del}(\varphi\beta\gamma)=G_{\del}(\beta)\cup G_{\del}(\gamma)$.

\item
$\displaystyle{
G_{\delta}(\psi_{\pi}^{\nu}(a))=\left\{
\begin{array}{ll}
G_{\delta}(\{\pi,a,\nu\}) & \delta<\pi
\\
\{\psi_{\pi}^{\nu}(a)\} & \pi\leq\delta
\end{array}
\right.
}$.

\eenu

}
\end{definition}

\begin{proposition}\label{prp:G}
For $\{\alpha,\delta, a,b,\rho\}\subset OT(\Pi_{3})$,
\benu
\item\label{prp:G1}
$G_{\delta}(\alpha)\leq\alpha$.

\item\label{prp:G2}
$\alpha\in\mathcal{H}_{a}(b)
\Rarw G_{\delta}(\alpha)\subset\mathcal{H}_{a}(b)$.

\eenu
\end{proposition}
\bprf
These are shown simultaneously by induction on the lengths $\ell\alpha$ of ordinal terms $\alp$.
It is easy to see that
\begin{equation}\label{eq:G}
G_{\delta}(\alpha)\ni\beta \Rarw \beta<\delta \spand \ell\beta\leq \ell\alpha
\end{equation}
\ref{prp:G}.\ref{prp:G1}.
Consider the case $\alpha=\psi_{\pi}^{\nu}(a)$ with $\delta<\pi$.
Then $G_{\delta}(\alpha)=G_{\delta}(\{\pi,a,\nu\})$.
On the other hand we have  
$\{\pi,a,\nu\}\subset\mathcal{H}_{a}(\alpha)$.
Proposition \ref{prp:G}.\ref{prp:G2} with (\ref{eq:G}) 
yields 
$G_{\delta}(\{\pi,a,\nu\})\subset\mathcal{H}_{a}(\alpha)\cap\pi\subset\alpha$.
Hence $G_{\delta}(\alpha)<\alpha$.
\\
\ref{prp:G}.\ref{prp:G2}.
Since $G_{\delta}(\alpha)\leq\alpha$ by Proposition \ref{prp:G}.\ref{prp:G1}, we can assume $\alpha\geq b$.

Consider the case $\alpha=\psi_{\pi}^{\nu}(a)$ with $\delta<\pi$.
Then $\{\pi,a,\nu\}\subset\mathcal{H}_{a}(b)$ and 
$G_{\delta}(\alpha)= G_{\delta}(\{\pi,a,\nu\})$.
IH yields the proposition.
\eprf

\begin{proposition}\label{lem:CX3} 
Let $\gamma<\beta$.
Assume $\alpha\in\mathcal{C}^{\gamma}(X)$ and $\forall\kappa\leq\beta[G_{\kappa}(\alpha)<\gamma]$.
Moreover assume 
$\forall\delta[\ell\delta\leq\ell\alpha\spand\delta\in\mathcal{C}^{\gamma}(X)\cap\gamma
\Rarw \delta\in \mathcal{C}^{\beta}(X)]$.
Then $\alpha\in\mathcal{C}^{\beta}(X)$.
\end{proposition}
\bprf
By induction on $\ell\alpha$. 
If $\alpha<\gamma$, then $\alpha\in \mathcal{C}^{\gamma}(X)\cap\gamma$.
The third assumption yields $\alpha\in \mathcal{C}^{\beta}(X)$. 
Assume $\alpha\geq\gamma$. 
Except the case $\alpha=\psi_{\pi}^{\nu}(a)$ for some $\pi,a,\nu$, 
IH yields $\alpha\in\mathcal{C}^{\beta}(X)$. 
Suppose $\alpha=\psi_{\pi}^{\nu}(a)$ for some $\{\pi,a,\nu\}\subset\mathcal{C}^{\gamma}(X)$ and 
$\pi>\gamma$.
If $\pi\leq\beta$, then $\{\alpha\}=G_{\pi}(\alpha)<\gamma$ by the second assumption. 
Hence this is not the case, and we obtain 
$\pi>\beta$.
 Then $G_{\kappa}(\{\pi,a,\nu\})=G_{\kappa}(\alpha)<\gamma$ for any $\kappa\leq\beta<\pi$. 
IH yields $\{\pi,a,\nu\}\subset\mathcal{C}^{\beta}(X)$. 
We conclude $\alpha\in\mathcal{C}^{\beta}(X)$ from $\pi>\beta$.
\eprf

\blem\label{lem:6.21}
Suppose $D[Y]$ and $\alp\in\mathcal{G}(Y)$.
Let 
$X=W(\mathcal{C}^{\alp}(Y))\cap\alp^{+}$.
Assume that the following condition (\ref{eq:6.21.55})
is fulfilled.
Then $\alp\in X$ and $D[X]$.

\beqn
\fal\bet\left(
Y\cap\alp^{+}<\bet \spand \bet^{+}<\alp^{+} \to
W(\mathcal{C}^{\bet}(Y))\cap\bet^{+}\subset Y
\right)
\label{eq:6.21.55}
\eeqn

\elem
\bprf
Let $\alp\in\mathcal{G}(Y)$.
By $\mathcal{C}^{\alp}(Y)\cap\alp\subset Y$ and $Wo[Y]$
we obtain by Proposition \ref{lem:3.11.632}
\beqn\label{eq:6.21.57}
X\cap\alp=Y\cap\alp=\mathcal{C}^{\alp}(Y)\cap\alp
\eeqn
Hence $\alp\in X$.

\bclm\label{clm:6.22-1}
$\alp^{+}=\gam^{+} \spand \gam\in X \Rarw \gam\in\mathcal{C}^{\gam}(X)$.
\eclm
\textbf{Proof} of Claim \ref{clm:6.22-1}.
Let $\alp^{+}=\gam^{+}$ and
$\gam\in X=W(\mathcal{C}^{\alp}(Y))\cap\alp^{+}$. 
We obtain $\gam\in \mathcal{C}^{\alp}(Y)=\mathcal{C}^{\gam}(Y)$ by 
Propositions \ref{lem:3.11.632} and \ref{lem:CX2}.
Hence
 $Y\cap\gam\subset \mathcal{C}^{\gam}(Y)\cap\gam=\mathcal{C}^{\alp}(Y)\cap\gam$.
$\gam\in W(\mathcal{C}^{\alp}(Y))$ yields
$Y\cap\gam\subset X$.
Therefore
we obtain
$\gam\in\mathcal{C}^{\gam}(Y)\subset\mathcal{C}^{\gam}(X)$.
\hspace*{\fill} $\Box$ of Claim \ref{clm:6.22-1}.

\bclm\label{clm:6.22}
$D[X]$.
\eclm
\textbf{Proof} of Claim \ref{clm:6.22}.
 We have $X\cap\alp=Y\cap\alp$ by (\ref{eq:6.21.57}).
Let $\bet\leq X$. 
We show $W(\mathcal{C}^{\bet}(X))\cap\bet^{+}=X\cap\bet^{+}$.
\\
\textbf{Case 1}. $\bet^{+}=\alp^{+}$:
We obtain $\mathcal{C}^{\bet}(X)=\mathcal{C}^{\alp}(X)=\mathcal{C}^{\alp}(Y)$ by
Proposition \ref{lem:CX2}, Claim \ref{clm:6.22-1} and (\ref{eq:6.21.57}).
\\
\textbf{Case 2}. $\bet^{+}<\alp^{+}$: Then $\bet^{+}\leq\alp$.

First let $Y\cap\alp^{+}<\bet$. Then the assumption (\ref{eq:6.21.55}) 
yields
$W(\mathcal{C}^{\bet}(Y))\cap\bet^{+}\subset Y$.
We obtain
$W(\mathcal{C}^{\bet}(X))\cap\bet^{+}\subset Y\cap\bet^{+}=X\cap\bet^{+}$
by (\ref{eq:6.21.57}).
It remains to show $Y\cap\bet^{+}\subset W(\mathcal{C}^{\bet}(Y))$.
Let $\gam\in Y\cap\bet^{+}$.
We obtain $\gam\in W(\mathcal{C}^{\gam}(Y))$ by $D[Y]$.
On the other hand we have $\mathcal{C}^{\bet}(Y)\subset\mathcal{C}^{\gam}(Y)$
by Proposition \ref{lem:CX2}.
Moreover Proposition \ref{lem:3.11.632} yields 
$\gam\in \mathcal{C}^{\bet}(Y)$.
Hence $\gam\in W(\mathcal{C}^{\bet}(Y))$.

Next let $\bet\leq Y\cap\alp^{+}$.
We obtain $Y\cap\bet^{+}=\mathcal{W}(\mathcal{C}^{\bet}(Y))\cap\bet^{+}$, and
$X\cap\bet^{+}=\mathcal{W}(\mathcal{C}^{\bet}(X))\cap\bet^{+}$ by (\ref{eq:6.21.57}).
\hspace*{\fill} $\Box$ of Claim \ref{clm:6.22}.

This completes a proof of Lemma \ref{lem:6.21}.
\eprf

\bprp\label{prp:noncritical}
Let $D[X]$.
\benu
\item\label{prp:noncritical.1}
Let $\{\alp,\bet\}\subset X$ with $\alp+\bet=\alp\#\bet$ and $\alp>0$.
Then $\gam=\alp+\bet\in X$.

\item\label{prp:noncritical.2}
If $\{\alp,\bet\}\subset X$, then $\vphi_{\alp}(\bet)\in X$.
\eenu
\eprp
\bprf
Proposition \ref{prp:noncritical}.\ref{prp:noncritical.2} is seen 
by main induction on $\alp\in X$ with
subsidiary induction on $\bet\in X$ using 
Proposition \ref{prp:noncritical}.\ref{prp:noncritical.1}.
We show Proposition \ref{prp:noncritical}.\ref{prp:noncritical.1}.
We obtain
$\alp\in X\cap\gam^{+}=W(\mathcal{C}^{\gam}(X))\cap\gam^{+}$ with 
$\gam^{+}=\alp^{+}$.
We see that $\alp+\bet\in W(\mathcal{C}^{\gam}(X))$ by induction on 
$\bet\in X\cap\alp\subset \mathcal{C}^{\gam}(X)$.
\eprf

\bprp\label{lem:6.29}
Let $X_{0}=W(\mathcal{C}^{0}(\emptyset))\cap 0^{+}$ with $0^{+}=\Ome$, and
$X_{1}=W(\mathcal{C}^{\Ome}(X_{0}))\cap\Ome^{+}$. Then
$0\in X_{0}$, $\Ome\in X_{1}$ and $D[X_{i}]$ for $i=0,1$.
\eprp
\bprf
For each $\alp\in\{0,\Ome\}$ and any set $Y\subset OT(\Pi_{3})$ we have
$\alp\in \mathcal{C}^{\alp}(Y)$.
First 
we obtain 
$0\in\mathcal{G}(\emptyset)$ and $D[\emptyset]$.
Also there is no $\bet$ such that 
$\bet^{+}<0^{+}$.
Hence the condition (\ref{eq:6.21.55}) is fulfilled, and we obtain
$0\in X_{0}$ and $D[X_{0}]$ by 
Lemma \ref{lem:6.21}.

Next 
let $\gam\in\mathcal{C}^{\Ome}(X_{0})\cap\Ome$.
We show $\gam\in X_{0}$ by induction on the lengths $\ell\gam$ of ordinal terms $\gam$ as follows.
We see that
each strongly critical number $\gam\in\mathcal{C}^{\Ome}(X_{0})\cap\Ome$ is in $X_{0}$
since 
if $\psi_{\sig}^{\nu}(\bet)<\Ome$, then $\sig=\Ome$.
Otherwise $\gam\in X_{0}$ is seen from IH using Proposition \ref{prp:noncritical} and 
$0\in X_{0}$.
Therefore we obtain
$\alp\in\mathcal{G}(X_{0})$.
Let $\bet^{+}<\alp^{+}$. Then $\bet^{+}=\Ome$ and $\bet<\Ome$.
Then $W(\mathcal{C}^{\bet}(X_{0}))\cap\Ome=
W(\mathcal{C}^{0}(X_{0}))\cap\Ome=X_{0}$ by Proposition \ref{lem:CX2}.
Hence the condition (\ref{eq:6.21.55}) is fulfilled, and we obtain
$\Ome\in X_{1}$ and $D[X_{1}]$ by
Lemma \ref{lem:6.21}.
\eprf

\bdf
{\rm
$\bet\prec\alp$ iff there exists a sequence $\{\sig_{i}\}_{i\leq n}(n>0)$ such that
$\alp=\sig_{0}$, $\bet=\sig_{n}$ and 
for each $i<n$, there are some $\nu_{i},a_{i}$ such that $\sig_{i+1}=\psi_{\sig_{i}}^{\nu_{i}}(a_{i})$.
}
\edf

Note that $\bet\prec\alp \Rarw m_{2}(\bet)<m_{2}(\alp)$.

\begin{lemma}\label{th:3wf16}
Suppose $D[Y]$ 
with $\{0,\Ome\}\subset Y$, and for $\eta\in OT(\Pi_{3})$
\begin{equation}\label{eq:3wf16hyp.132}
\eta\in\mathcal{G}(Y)
\end{equation}
and
\begin{equation}\label{eq:3wf16hyp.232}
\fal\gam\prec\eta (\gam\in\mathcal{G}(Y)
\Rarw \gam\in Y)
\end{equation}
Let $X=W(\mathcal{C}^{\eta}(Y))\cap\eta^{+}$.
Then
$\eta\in X$ and $D[X]$.
\end{lemma}
\bprf
By Lemma \ref{lem:6.21} 
and the hypothesis (\ref{eq:3wf16hyp.132}) it suffices to show 
(\ref{eq:6.21.55}), i.e., 
\[
\fal\bet\left(
Y\cap\eta^{+}<\bet \spand \bet^{+}<\eta^{+} \to
W(\mathcal{C}^{\bet}(Y))\cap\bet^{+}\subset Y
\right)
.\]
Assume $Y\cap\eta^{+}<\bet$ and $\bet^{+}<\eta^{+}$. 
We have to show  
$W(\mathcal{C}^{\beta}(Y))\cap\bet^{+}\subset Y$. 
We prove this by induction on $\gamma\in W(\mathcal{C}^{\beta}(Y))\cap\bet^{+}$. 
Suppose $\gamma\in \mathcal{C}^{\beta}(Y)\cap\bet^{+}$ and 
\[
\mbox{{\rm MIH : }} \mathcal{C}^{\beta}(Y)\cap\gamma\subset Y.
\]
We show $\gamma\in Y$. 
We can assume that
\begin{equation}
\label{eq:3wf9hyp.232X}
Y\cap\eta^{+}<\gamma
\end{equation}
since if $\gamma\leq \delta$ for some $\delta\in Y\cap\eta^{+}$, then by 
$Y\cap\eta^{+}<\beta$ and $\gamma\in \mathcal{C}^{\beta}(Y)$ 
we obtain $\delta<\beta$, 
$\gamma\in \mathcal{C}^{\delta}(Y)$ and 
$\delta\in W(\mathcal{C}^{\delta}(Y))\cap\delta^{+}=Y\cap\delta^{+}$.
Hence $\gamma\in W(\mathcal{C}^{\delta}(Y))\cap\delta^{+}\subset Y$.

We show first 
\begin{equation}
\label{eq:3wf9hyp.232}
\gamma\in\mathcal{G}(Y)
\end{equation}
First $\gamma\in \mathcal{C}^{\gamma}(Y)$ by 
$\gamma\in \mathcal{C}^{\beta}(Y)\cap\beta^{+}$ and Proposition \ref{lem:CX2}. 
Second we show the following claim by induction on $\ell\alpha$:

\begin{equation}\label{clm:3wf1632}
\alpha\in\mathcal{C}^{\gamma}(Y)\cap\gamma \Rarw  \alpha\in Y
\end{equation}
\textbf{Proof} of (\ref{clm:3wf1632}). 
Assume $\alpha\in\mathcal{C}^{\gamma}(Y)$.
We can assume $\gamma^{+}\leq\beta$ for otherwise we have 
$\alpha\in \mathcal{C}^{\gamma}(Y)\cap\gamma=\mathcal{C}^{\beta}(Y)\cap\gamma\subset Y$ by MIH.

By induction hypothesis on lengths, $\alp<\gam<\bet^{+}<\eta^{+}$,
Proposition \ref{prp:noncritical}, and $\{0,\Ome\}\subset Y$, 
 we can assume that
$\alpha=\psi_{\pi}^{\nu}(a)$ for some $\pi>\gamma$ such that
 $\{\pi,a,\nu\}\subset\mathcal{C}^{\gamma}(Y)$.
\\
\textbf{Case 1}. $\beta<\pi$: 
Then $G_{\bet}(\{\pi,a,\nu\})=G_{\bet}(\alpha)<\alpha<\gamma$ by Proposition \ref{prp:G}.\ref{prp:G1}.
 Proposition \ref{lem:CX3} with induction hypothesis on lengths yields 
$\{\pi,a,\nu\}\subset\mathcal{C}^{\beta}(Y)$.
Hence $\alpha\in \mathcal{C}^{\beta}(Y)\cap\gamma$ by $\pi>\beta$.
MIH yields $\alpha\in Y$.
\\
\textbf{Case 2}. $\beta\geq\pi$: 
We have $\alpha<\gamma<\pi\leq\beta$. 
It suffices to show that $\alpha\leq Y\cap\eta^{+}$.
Then by (\ref{eq:3wf9hyp.232X}) we have $\alpha\leq\delta\in Y\cap\eta^{+}$ for some $\delta<\gamma$.
$\mathcal{C}^{\delta}(Y)\ni\alpha\leq\delta\in Y\cap\delta^{+}=W(\mathcal{C}^{\delta}(Y))\cap\delta^{+}$
yields $\alpha\in W(\mathcal{C}^{\delta}(Y))\cap\delta^{+}\subset Y$.

Assume first that $\gam$ is not a strongly critical number.
By $\alpha=\psi_{\pi}^{\nu}(a)<\gamma$,
we can assume that 
$\gamma\neq 0$.
Let $\delta$ denote the largest immediate subterm of $\gamma$.
We obtain $\delta\in \mathcal{C}^{\beta}(Y)\cap\gamma$
by (\ref{eq:3wf9hyp.232X}), 
$Y\cap\eta^{+}<\gamma\in\mathcal{C}^{\beta}(Y)$.
Hence $\delta\in Y$ by MIH.
Also by $\alpha<\gamma$, we obtain $\alpha\leq\delta$, i.e., $\alpha\leq Y$, and we are done.

Next let $\gamma=\psi_{\kappa}^{\xi}(b)$ for some $b,\xi$ and $\kappa>\beta$ by (\ref{eq:3wf9hyp.232X}) and $\gamma\in \mathcal{C}^{\beta}(Y)$.
We have $\alpha<\gamma<\pi\leq\beta<\kappa$.
We obtain
$\pi\not\in\mathcal{H}_{b}(\gamma)$ since otherwise by $\pi<\kappa$ we would have $\pi<\gamma$.
Therefore $\alp=\psi_{\pi}^{\nu}(a) <\psi_{\kappa}^{\xi}(b)=\gam<\pi<\kap$ with 
$\pi\in\mathcal{H}_{a}(\alp)$ and $\pi\not\in\mathcal{H}_{b}(\gamma)$. This
yields $a> b$ and 
$\{\kappa,b,\xi\}\not\subset\mathcal{H}_{a}(\alpha)$.

On the other hand we have 
$\{\kappa,b,\xi\}\subset\mathcal{H}_{a}(\gamma)$.
This means that there exists a subterm
$\delta<\gam$ of one of $\kappa,b,\xi$
such that $\del\not\in\mathcal{H}_{a}(\alpha)$.
Also we have $\{\kappa,b,\xi\}\subset\mathcal{C}^{\beta}(Y)$.
Then $\delta\in\mathcal{C}^{\beta}(Y)\cap\gamma$.
By MIH
we obtain 
$\alpha\leq\delta\in \mathcal{C}^{\beta}(Y)\cap\gamma\subset Y$.

\hspace*{\fill} $\Box$ of (\ref{clm:3wf1632}) and (\ref{eq:3wf9hyp.232}).
\\

\noindent
Hence we obtain $\gamma\in\mathcal{G}(Y)$.
We have $\gamma<\beta^{+}\leq\eta$ and 
$\gamma\in\mathcal{C}^{\gamma}(Y)$.
If $\gamma\prec\eta$, then the hypothesis 
(\ref{eq:3wf16hyp.232}) yields $\gamma\in Y$.
In what follows assume $\gamma\not\prec\eta$.

If $G_{\eta}(\gamma)<\gamma$, then  Proposition \ref{lem:CX3} yields 
$\gamma\in\mathcal{C}^{\eta}(Y)\cap\eta\subset Y$ by $\eta\in\mathcal{G}(Y)$.

Suppose $G_{\eta}(\gamma)=\{\gamma\}$. 
This means, by $\gamma\not\prec\eta$, that
$\gamma\prec\tau$ for a $\tau<\eta$.
Let $\tau$ denote the maximal such one.
We have $\gamma<\tau<\eta$.
From $\gamma\in\mathcal{C}^{\gamma}(Y)$ we see 
$\tau\in\mathcal{C}^{\gamma}(Y)$.
Next we show that
\begin{equation}\label{eq:3wf1632last}
G_{\eta}(\tau)<\gamma
\end{equation}
Let $\tau=\psi_{\kap}^{\mu}(b)$.
Then $\eta<\kap$ by the maximality of $\tau$, and
$G_{\eta}(\tau)=G_{\eta}(\{\kap,b,\mu\})<\tau$ by Proposition \ref{prp:G}.\ref{prp:G1}.
On the other hand we have $\tau\in\mathcal{H}_{a}(\gam)$.
Proposition \ref{prp:G}.\ref{prp:G2} yields
$G_{\eta}(\tau)\subset\mathcal{H}_{a}(\gam)$.
We see $G_{\eta}(\tau)<\gam$ inductively.

Proposition \ref{lem:CX3} with (\ref{eq:3wf1632last}) yields 
$\tau\in\mathcal{C}^{\eta}(Y)$, and
$\tau\in\mathcal{C}^{\eta}(Y)\cap\eta\subset Y$ by $\eta\in\mathcal{G}(Y)$.
Therefore $Y\cap\eta^{+}<\gamma<\tau\in Y$.
This is not the case by (\ref{eq:3wf9hyp.232X}).
We are done.
\hspace*{\fill} $\Box$

\begin{proposition}\label{prp:updis}
$\alpha\leq \mathcal{W}\cap\beta^{+} \spand \alpha\in\mathcal{C}^{\beta}(\mathcal{W}) \Rarw \alpha\in \mathcal{W}$.
\end{proposition}
\bprf
This is seen from Propositions \ref{lem:CX2}, \ref{lem:3.11.632} an \ref{lem:3wf6}.
\eprf

\subsection{Mahlo universes}

In Proposition \ref{lem:3wf6max}, we saw that $\mathcal{W}$ is the maximal distinguished class,
which is $\Sig^{1-}_{2}$-definable and a proper class in ${\sf KP}\Pi_{3}$.
$\calw^{P}$ in Definition \ref{df:3awp} denotes the maximal distinguished class \textit{inside} a set
$P$. $\calw^{P}$ exists as a set.

Let $ad$ denote a $\Pi^{-}_{3}$-sentence such that a transitive set $z$ is admissible 
iff $(z;\in)\models ad$.
Let $lmtad :\Lrarw \fal x\exi y(x\in y \land ad^{y})$.
Observe that $lmtad$ is a $\Pi^{-}_{2}$-sentence.

\bdf\label{df:3auni}
{\rm
$L$ denotes a whole universe, which is a model of ${\sf KP}\Pi_{3}$.
\benu
\item 
By a \textit{universe} 
we mean either the whole universe $L$ or a transitive set $Q\in L$ with $\ome\in Q$. 
Universes are denoted by $P,Q,\ldots$

\item For a universe $P$ and a set-theoretic sentence $\vphi$, 
$P\models\vphi :\Lrarw (P;\in)\models\vphi$.

\item 
A universe $P$ is said to be a \textit{limit universe} if $lmtad^{P}$ holds, i.e., $P$
 is a limit of admissible sets. The class of limit universes is denoted by $Lmtad$.
 
\eenu
}
\edf

\blem\label{lem:3ahier} 
$W(\mathcal{C}^{\alp}(X))$ as well as $D[X]$ are absolute for limit universes $P$.

\elem
\bprf
Let $P$ be a limit universe and $X\in \mathcal{P}(\ome)\cap P$. 
Then $W(X)$ is $\Del_{1}$ in $P$, and so is $W(\mathcal{C}^{\alp}(X))$.
Hence $W(\mathcal{C}^{\alp}(X))=\{\bet\in OT(\Pi_{3}):P\models \bet\in W(\mathcal{C}^{\alp}(X))\}$,
and $D[X]\Lrarw P\models D[X]$.
\eprf

\bdf\label{df:3awp}
{\rm 
For a universe $P$, let
$
\calw^{P}:=\bigcup\{X\in P:D[X]\}
$.
}
\edf

\blem\label{lem:3afin}
Let $P$ be a universe closed under finite unions, and $\alp\in OT(\Pi_{3})$.
\benu
\item\label{lem:3afin.1}
There is a finite set $K(\alp)\subset OT(\Pi_{3})$ such that
$
\fal Y\in P\fal\gam
[
K(\alp)\cap Y=K(\alp)\cap \calw^{P} 
\Rarw 
\left(
\alp\in \mathcal{C}^{\gam}(\calw^{P})
\Lrarw \alp\in \mathcal{C}^{\gam}(Y)
\right)
]
$.

\item\label{lem:3afin.2}
There exists a distinguished set $X\in P$ such that
$
\fal Y\in P\fal\gam
[
X\subset Y\spand D[Y]
\Rarw (\alp\in \mathcal{C}^{\gam}(\calw^{P})
\Lrarw \alp\in \mathcal{C}^{\gam}(Y)
)
]
$.
\eenu
\elem
\bprf
\ref{lem:3afin}.\ref{lem:3afin.1}.
F.e. the set of subterms of $\alp$ enjoys the condition for $K(\alp)$.
\\
\ref{lem:3afin}.\ref{lem:3afin.2}.
By $X,Y\in P\Rarw X\cup Y\in P$, pick a distinguished set $X\in P$ such that
$K(\alp)\cap\calw^{P}\subset X$.
\eprf

\begin{proposition}\label{lem:3wf6}
For each limit universe $P$,
$D[\mathcal{W}^{P}]$ holds, and
 $\exists X(X=\mathcal{W}^{P})$ if $P$ is a set.
\end{proposition}
\bprf
$D[\mathcal{W}^{P}]$ is seen as in Proposition \ref{lem:3wf6max}.
\eprf
\\

For a universal $\Pi_{n}$-formula $\Pi_{n}(a)\, (n>0)$ uniformly on admissibles, let
\[
P\in M_{2}(\calc) :\Lrarw P\in Lmtad \spand 
\fal b\in P[P\models\Pi_{2}(b)\rarw \exi Q\in \calc\cap P(Q\models\Pi_{2}(b))]
.\]

\blem\label{lem:4acalg} 
Let $\calc$ be a $\Pi^{1}_{0}$-class such that $\calc\subset Lmtad$.
Suppose $P\in M_{2}(\calc)$ and $\alp\in\calg(\calw^{P})$. 
Then there exists a universe $Q\in \calc$ such that $\alp\in\calg(\calw^{Q})$. 
\elem
\bprf 
Suppose $P\in M_{2}(\calc)$ and $\alp\in\calg(\calw^{P})$. 
First by $\alp\in \mathcal{C}^{\alp}(\calw^{P})$ and Lemma \ref{lem:3afin} pick a distinguished set
 $X_{0}\in P$ such that $\alp\in \mathcal{C}^{\alp}(X_{0})$ and
$K(\alp)\cap\calw^{P}\subset X_{0}$.
Next writing $\mathcal{C}^{\alp}(\calw^{P})\cap\alp\subset\calw^{P}$ analytically we have

\[
\fal\bet<\alp[\bet\in \mathcal{C}^{\alp}(\calw^{P}) \Rarw 
\exi Y\in P(D[Y]\spand \bet\in Y)]
\]
By Lemma \ref{lem:3afin} we obtain
$\bet\in \mathcal{C}^{\alp}(\calw^{P})\Lrarw
\exi X\in P\{D[X] \spand K(\bet)\cap\calw^{P}\subset X \spand \bet\in \mathcal{C}^{\alp}(X) \}$.
Hence
for any $\bet<\alp$ and any distinguished set $X\in P$,
 there are $\gam\in K(\bet)$, $Z\in P$ and a distinguished set $Y\in P$ such that 
 if $\gam\in Z \spand D[Z] \to \gam\in X$ and $\bet\in \mathcal{C}^{\alp}(X)$, then 
 $\bet\in Y$.
By Lemma \ref{lem:3ahier} $D[X]$ is absolute for limit universes. 
Hence the following $\Pi_{2}$-predicate holds in the universe $P\in M_{2}(\calc)$:
\beqnarr
&& \fal\bet<\alp\fal X\exi\gam\in K(\bet)\exi Z\exi Y
[
\{
D[X]\spand 
 (\gam\in Z \spand D[Z] \to \gam\in X)
  \spand \bet\in \mathcal{C}^{\alp}(X)
\}
\nonumber
\\
&&
 \Rarw 
 \left(
 D[Y]\spand \bet\in Y
 \right)
 ]
 \label{eq:3acalg}
\eeqnarr

Now pick a universe $Q\in \calc\cap P$ with $X_{0}\in Q$ and $Q\models (\ref{eq:3acalg})$. 
Tracing the above argument backwards in the limit universe $Q$ we obtain 
$\mathcal{C}^{\alp}(\calw^{Q})\cap\alp\subset\calw^{Q}$ and 
$X_{0}\subset\calw^{Q}=\bigcup\{X\in Q: Q\models D[X]\}\in P$. 
Thus Lemma \ref{lem:3afin} yields $\alp\in \mathcal{C}^{\alp}(\calw^{Q})$. 
We obtain $\alp\in\calg(\calw^{Q})$.
\eprf

\bdf\label{df:UVM}
{\rm 
We define the class $M_{2}(\alp)$ of $\alp$-recursively Mahlo universes for 
$\alp\in OT(\Pi_{3})$ as follows:
\beqn\label{eqarr:4a(1)}
 P\in M_{2}(\alp) \Lrarw 
 P\in Lmtad \spand 
 \fal\bet\prec\alp
 [
SC_{\mK}(m_{2}(\bet))\subset\calw^{P}
 \Rarw P\in M_{2}(M_{2}(\bet))] 
\eeqn
$M_{2}(\alp)$ is a $\Pi_{3}$-class.

}
\edf

\bprp\label{prp:m2}
If $\eta\in\calg(Y)$, then $SC_{\mK}(m_{2}(\eta))\subset Y$.
\eprp
\bprf
Let $\nu=m_{2}(\eta)$. 
Then $SC_{\mK}(\nu)\subset\eta$ by (\ref{eq:OTpi3m2}).
From $\eta\in\mathcal{C}^{\eta}(Y)$ we see $SC_{\mK}(\nu)\subset\mathcal{C}^{\eta}(Y)$.
Hence $SC_{\mK}(\nu)\subset\mathcal{C}^{\eta}(Y)\cap\eta\subset Y$ by $\eta\in\calg(Y)$.
\eprf

\blem\label{lem:3awf16.1} 
If $\eta\in\calg(\calw^{P})$
and $P\in M_{2}(M_{2}(\eta))$, 
then $\eta\in\calw^{P}$.
\elem
\bprf 
We show this by induction on $\in$. Suppose, as IH, the lemma holds for any $Q\in P$.
By Lemma \ref{lem:4acalg} pick a $Q\in P$ such that $Q\in M_{2}(\eta)$, and for $Y=\calw^{Q}\in P$,
$\{0,\Ome\}\subset Y$ and
\beqn\renewcommand{\theequation}{\ref{eq:3wf16hyp.132}} 
\eta\in\calg(Y)
\eeqn
\addtocounter{equation}{-1}
On the other the definition (\ref{eqarr:4a(1)}) yields 
$\fal\gam\prec\eta[
SC_{\mK}(m_{2}(\gam))\subset\calw^{Q}
 \Rarw Q\in M_{2}(M_{2}(\gam))]$.
Hence by Proposition \ref{prp:m2}
$\fal\gam\prec\eta[\gam\in \calg(\calw^{Q}) \Rarw  
Q\in M_{2}(M_{2}(\gam))]$.

IH yields with $Y=\calw^{Q}$
\beqn\renewcommand{\theequation}{\ref{eq:3wf16hyp.232}}
\fal\gam\prec\eta(\gam\in\calg(Y) \Rarw  \gam\in Y)
\eeqn
\addtocounter{equation}{-1}
Therefore by Lemma \ref{th:3wf16} we conclude 
$\eta\in X$ and $D[X]$ for $X=W(\mathcal{C}^{\eta}(Y))\cap\eta^{+}$.

$X\in P$ follows from $Y\in P\in Lmtad$. 
Consequently $\eta\in\calw^{P}$.
\eprf

\blem\label{lem:3wf19-1}
\benu

\item 
$\mathcal{C}^{\mK}(\calw)\cap\mK=\calw\cap\mK$.

\label{lem:3wf19-1.1}
\item 
$\mK\in\mathcal{C}^{\mK}(\calw)$.
\label{lem:3wf19-1.2}

\item 
For {\rm each} $n\in\ome,$ 
$TI[\mathcal{C}^{\mK}(\calw)\cap\ome_{n}(\mK+1)]$.
\label{lem:3wf19-1.4}
\eenu
\elem
\bprf
We show Lemma \ref{lem:3wf19-1}.\ref{lem:3wf19-1.4}.
It suffices to show $TI[\calw]$.
Assume $Prg[\calw,A]$ for a formula $A$, and $\alp\in\calw$.
Pick a distinguished set $X$ such that $\alp\in X$. Then $X\cap\alp^{+}=\calw\cap\alp^{+}$, and
hence $Prg[X\cap(\alp+1),A]$. $Wo[X]$ yields $A(\alp)$.
\eprf

\blem\label{lem:4aro}
 $\fal\eta[m_{2}(\eta)\in\mathcal{C}^{\mK}(\calw)\cap\ome_{n}(\mK+1) \Rarw L\in M_{2}(M_{2}(\eta))]$
 holds for {\rm each} $n\in\ome$.
\elem
\bprf 
We show the lemma by induction on $\nu=m_{2}(\eta)\in\mathcal{C}^{\mK}(\calw)$ up to each 
$\ome_{n}(\mK+1)$. Suppose $\nu\in\mathcal{C}^{\mK}(\calw)$ and $L\models\Pi_{2}(b)$ for a $b\in L$. We have to find a universe $Q\in L$ such that $b\in Q$, $Q\in M_{2}(\eta)$ and $Q\models\Pi_{2}(b)$.

By the definition (\ref{eqarr:4a(1)}) $L\in M_{2}(\eta)$ is equivalent to
$\fal\gam\prec\eta [ m_{2}(\gam)\in\mathcal{C}^{\mK}(\calw) \Rarw L\in M_{2}(M_{2}(\gam))]$.
We obtain $\gam\prec\eta \Rarw m_{2}(\gam)<m_{2}(\eta)=\nu$. 
Thus IH yields $L\in M_{2}(\eta)$.
Let $g$ be a primitive recursive function in the sense of set theory such that 
$L\in M_{2}(\eta) \Lrarw P\models\Pi_{3}(g(\eta))$.
Then $L\models\Pi_{2}(b) \land \Pi_{3}(g(\eta))$. 
Since this is a $\Pi_{3}$-formula which holds in a $\Pi_{3}$-reflecting universe $L$, 
we conclude for some $Q\in L$, $Q\models\Pi_{2}(b)\land \Pi_{3}(g(\eta))$
 and hence $Q\in M_{2}(\eta)$. We are done.
\eprf

\brem
{\rm
Only here we need $\Pi_{3}$-reflection. Therefore it sufffices for a whole universe $L$ to admit iterations of $\Pi_{2}$-recursively Mahlo operations along a well founded relation $\prec$ which is $\Sig$ on $L$: $L\in M_{2}^{\prec}(\mu)=\bigcap\{M_{2}(M_{2}^{\prec}(\nu)):L\models\nu\prec\mu\}$. Hence our wellfoundednes proof is formalizable in a set theory axiomatizing such universes $L$.
}
\erem

\blem\label{lem:3awf16.2}
For {\rm each} $n\in\ome$,
$m_{2}(\eta)<\ome_{n}(\mK+1)\spand \eta\in\calg(\calw) \Rarw \eta\in\calw$.
\elem
\bprf 
Assume $\nu=m_{2}(\eta)<\ome_{n}(\mK+1)$ and $\eta\in\calg(\calw)$. 
By Proposition \ref{prp:m2} we obtain $\nu\in\mathcal{C}^{\mK}(\calw)$.
Lemma \ref{lem:4aro} yields $L\in M_{2}(M_{2}(\eta))$. 
From this we see $L\in M_{2}(\calc)$ with $\calc=M_{2}(M_{2}(\eta))$ 
as in the proof of Lemma \ref{lem:4aro} using $\Pi_{3}$-reflection of the whole universe $L$ once again. 
Then by Lemma \ref{lem:4acalg} pick a set $P\in L$ such that 
$\eta\in\calg(\calw^{P})$ and $P\in\calc=M_{2}(M_{2}(\eta))$. 
Lemma \ref{lem:3awf16.1} yields $\eta\in\calw^{P}\subset\calw$.
\eprf

\subsection{Well-foundedness proof (concluded)}
\begin{definition}\label{df:EKk}

{\rm 
For terms $\alpha,\kappa,\delta\in OT(\Pi_{3})$, finite sets 
$\mathcal{E}(\alpha), K_{\delta}(\alpha),  k_{\delta}(\alpha)\subset OT(\Pi_{3})$ 
 are defined recursively as follows.
\benu

\item
$\mathcal{E}(\alpha)=\emptyset$ for $\alpha\in\{0,\Omega,\mK\}$.
$\mathcal{E}(\alpha_{m}+\cdots+\alpha_{0})=\bigcup_{i\leq m}\mathcal{E}(\alpha_{i})$.
$\mathcal{E}(\varphi\beta\gamma)=\mathcal{E}(\beta)\cup\mathcal{E}(\gamma)$.
$\mathcal{E}(\psi_{\pi}^{\nu}(a))=\{\psi_{\pi}^{\nu}(a)\}$.

\item
$\mathcal{A}(\alpha)=\bigcup\{\mathcal{A}(\beta): \beta\in\mathcal{E}(\alpha)\}$
 for $\mathcal{A}\in\{K_{\delta},k_{\delta}\}$.

\item

$
K_{\del}(\psi_{\pi}^{\nu}(a))=\left\{
\begin{array}{ll}
\{a\}\cup K_{\del}(\{\pi,a\}\cup SC_{\mK}(\nu)) & \psi_{\pi}^{\nu}(a)\geq\del
\\
\emptyset & \psi_{\pi}^{\nu}(a)<\del
\end{array}
\right.
$.

\item

$k_{\del}(\psi_{\pi}^{\nu}(a))=\left\{
\begin{array}{ll}
\{\psi_{\pi}^{\nu}(a)\}\cup k_{\del}(\{\pi,a\}\cup SC_{\mK}(\nu)) & \psi_{\pi}^{\nu}(a)\geq\del
\\
\emptyset & \psi_{\pi}^{\nu}(a)<\del
\end{array}
\right.
$.

\eenu

}
\end{definition}

Note that $K_{\del}(\alp)<a \Lrarw \alp\in\calh_{a}(\del)$.

\begin{definition}\label{df:id4wfA}
{\rm 
For $a,\nu\in OT(\Pi_{3})$, define:
\begin{eqnarray}
 A(a,\nu) & :\Lrarw &
 \forall\sigma\in\mathcal{C}^{\mK}(\mathcal{W})[\psi_{\sigma}^{\nu}(a)\in OT(\Pi_{3}) \Rarw \psi_{\sigma}^{\nu}(a)\in\mathcal{W}].
\\
\mbox{{\rm MIH}}(a) & :\Lrarw &
 \forall b\in\mathcal{C}^{\mK}(\calw)\cap a\forall \nu\in\mathcal{C}^{\mK}(\calw) \, A(b,\nu).
\\
\mbox{{\rm SIH}}(a,\nu) & :\Lrarw &
 \forall \xi\in \mathcal{C}^{\mK}(\calw) [\xi<\nu  \Rarw A(a,\xi)].
\end{eqnarray}
}
\end{definition}

\begin{lemma}\label{th:id5wf21}
For {\rm each} $n$ the following holds:
Assume $\{a,\nu\}\subset\mathcal{C}^{\mK}(\calw)\cap\ome_{n}(\mK+1)$, 
$\mbox{{\rm MIH}}(a)$, and 
$\mbox{{\rm SIH}}(a,\nu)$ in Definition \ref{df:id4wfA}.
Then
\[
 \forall\kappa\in\mathcal{C}^{\mK}(\calw)[\psi_{\kappa}^{\nu}(a)\in OT(\Pi_{3}) \Rarw 
 \psi_{\kappa}^{\nu}(a)\in\mathcal{W}].
\]
\end{lemma}
\bprf
Let $\alpha_{1}=\psi_{\kappa}^{\nu}(a)\in OT(\Pi_{3})$ with 
$\{a,\kap,\nu\}\subset\mathcal{C}^{\mK}(\calw)$ and $\nu\leq a<\ome_{n}(\mK+1)$, cf.\,(\ref{eq:OTpi3m2}). 
By Lemma \ref{lem:3awf16.2} it suffices to show $\alpha_{1}\in\mathcal{G}(\calw)$.

By Proposition \ref{lem:3.11.632} we have
$\{\kappa,a, \nu\}\subset\mathcal{C}^{\alpha_{1}}(\mathcal{W})$, and hence $\alpha_{1}\in\mathcal{C}^{\alpha_{1}}(\mathcal{W})$.
It suffices to show the following claim.
\begin{equation}\label{clm:id5wf21.1}
\forall\beta_{1}\in\mathcal{C}^{\alpha_{1}}(\mathcal{W})\cap\alpha_{1}[\beta_1\in\mathcal{W}].
\end{equation}
\textbf{Proof} of (\ref{clm:id5wf21.1}) by induction on $\ell\beta_1$. 
Assume $\beta_{1}\in\mathcal{C}^{\alpha_{1}}(\mathcal{W})\cap\alpha_{1}$ and let
\[
\mbox{LIH} :\Lrarw
\forall\gamma\in\mathcal{C}^{\alpha_{1}}(\mathcal{W})\cap\alpha_{1}[\ell\gamma<\ell\beta_{1} \Rarw \gamma\in\mathcal{W}].
\]

We show $\beta_1\in\mathcal{W}$. 
By Propositions \ref{prp:noncritical}, \ref{lem:6.29} and LIH, we may assume that
$\beta_{1}=\psi_{\pi}^{\xi}(b)$ for some $\pi,b,\xi$ such that
$\{\pi,b,\xi\}\subset\mathcal{C}^{\alpha_{1}}(\mathcal{W})$.

$\beta_{1}=\psi_{\pi}^{\xi}(b)<\psi_{\kap}^{\nu}(a)=\alp_{1}$ holds iff
one of the following holds:
(1) $\pi\leq\alp_{1}$. (2) $b<a$, $\beta_{1}<\kappa$ and $\{\pi,b,\xi\}\subset\calh_{a}(\alp_{1})$.
(3) $b=a$, $\pi=\kappa$, $\xi\in\calh_{a}(\alpha_{1})$ and $\xi<\nu$.
(4) $a\leq b$ and $\{\kap,a,\nu\}\not\subset\calh_{b}(\bet_{1})$.
\\
\textbf{Case 1}. $\pi\leq\alpha_{1}$: 
Then $\beta_{1}\in\mathcal{W}$ by $\beta_{1}\in\mathcal{C}^{\alpha_{1}}(\mathcal{W})$.
\\
\textbf{Case 2}.
$b<a$, $\beta_{1}<\kappa$ and $\{\pi,b,\xi\}\subset\calh_{a}(\alp_{1})$:
Let $B$ denote a set of subterms of $\beta_{1}$ defined recursively as follows.
First $\{\pi,b\}\cup SC_{\mK}(\xi)\subset B$.
Let $\alpha_{1}\leq\beta\in B$. 
If $\beta=_{NF}\gamma_{m}+\cdots+\gamma_{0}$, then $\{\gamma_{i}:i\leq m\}\subset B$.
If $\beta=_{NF}\varphi\gamma\delta$, then $\{\gamma,\delta\}\subset B$.
If $\beta=\psi_{\sigma}^{\mu}(c)$, then $\{\sigma,c\}\cup SC_{\mK}(\mu)\subset B$.

Then from $\{\pi,b,\xi\}\subset\mathcal{C}^{\alpha_{1}}(\mathcal{W})$ we see inductively that
$B\subset\mathcal{C}^{\alpha_{1}}(\mathcal{W})$.
Hence by LIH we obtain $B\cap\alpha_{1}\subset\mathcal{W}$.
Moreover if $\alpha_{1}\leq\psi_{\sigma}^{\mu}(c)\in B$, then 
we see $c<a$ 
from $\{\pi,b,\xi\}\subset\calh_{a}(\alp_{1})$.
We claim that
\begin{equation}\label{eq:case2A}
\forall\beta\in B(\beta\in\mathcal{C}^{\mK}(\calw))
\end{equation}
\textbf{Proof} of (\ref{eq:case2A}) by induction on $\ell\beta$.
Let $\beta\in B$. We can assume that $\alpha_{1}\leq\beta=\psi_{\sigma}^{\mu}(c)$ by induction hypothesis on the lengths.
Then by induction hypothesis we have $\{\sigma,c\}\cup SC_{\mK}(\mu)\subset\mathcal{C}^{\mK}(\calw)$.
On the other hand we have $\mu\leq c<a$ by (\ref{eq:OTpi3m2}).
$\mbox{MIH}(a)$ yields $\beta\in\mathcal{W}$.
Thus (\ref{eq:case2A}) is shown.
\hspace*{\fill} $\Box$
\\

In particular we obtain $\{\pi,b\}\cup SC_{\mK}(\xi)\subset \mathcal{C}^{\mK}(\calw)$.
Moreover we have $\xi\leq b<a$ by (\ref{eq:OTpi3m2}).
Therefore once again $\mbox{MIH}(a)$ yields $\beta_{1}\in\mathcal{W}$.
\\
\textbf{Case 3}. 
$b=a$, $\pi=\kappa$, $\xi\in\calh_{a}(\alpha_{1})$ and $\xi<\nu\leq a$: 
As in (\ref{eq:case2A}) we see that $SC_{\mK}(\xi)\subset\mathcal{W}$ from $\mbox{MIH}(a)$.
$\mbox{SIH}(a,\nu)$ yields $\beta_{1}\in\mathcal{W}$.
\\
\textbf{Case 4}.
$a\leq b$ and $\{\kap,a,\nu\}\not\subset\calh_{b}(\bet_{1})$:
It suffices to find a $\gamma$ such that $\beta_{1}\leq\gamma\in\mathcal{W}\cap\alpha_{1}$.
Then $\beta_{1}\in\mathcal{W}$ follows from $\beta_{1}\in\mathcal{C}^{\alpha_{1}}(\mathcal{W})$ and Proposition \ref{prp:updis}.

$k_{\delta}(\alpha)$ denotes the set in Definition \ref{df:EKk}.
In general we see that $a\in K_{\delta}(\alpha)$ iff $\psi_{\sigma}^{h}(a)\in k_{\delta}(\alpha)$ for some $\sigma,h$,
and for each $\psi_{\sigma}^{h}(a)\in k_{\delta}(\psi_{\sigma_{0}}^{h_{0}}(a_{0}))$ there exists a sequence
$\{\alpha_{i}\}_{i\leq m}$ of subterms of $\alpha_{0}=\psi_{\sigma_{0}}^{h_{0}}(a_{0})$ such that 
$\alpha_{m}=\psi_{\sigma}^{h}(a)$, 
$\alpha_{i}=\psi_{\sigma_{i}}^{h_{i}}(a_{i})$ for some $\sigma_{i},a_{i},h_{i}$,
and for each $i<m$,
$\delta\leq\alpha_{i+1}\in\mathcal{E}(C_{i})$ for $C_{i}=\{\sigma_{i},a_{i}\}\cup SC_{\mK}(h_{i})$.

Let $\delta\in SC_{\mK}(f)\cup\{\kappa,a\}$ such that $b\leq \gamma$
for a $\gamma\in K_{\beta_{1}}(\delta)$.
Pick an $\alpha_{2}=\psi_{\sigma_{2}}^{h_{2}}(a_{2})\in \mathcal{E}(\delta)$ 
such that $\gamma\in K_{\beta_{1}}(\alpha_{2})$, and 
an $\alpha_{m}=\psi_{\sigma_{m}}^{h_{m}}(a_{m})\in k_{\beta_{1}}(\alpha_{2})$ for some $\sigma_{m},h_{m}$ and
$a_{m}\geq b\geq a$.
We have $\alpha_{2}\in\mathcal{W}$ by $\delta\in\mathcal{W}$.
If $\alpha_{2}<\alpha_{1}$, then
$\beta_{1}\leq\alpha_{2}\in\mathcal{W}\cap\alpha_{1}$, and we are done.
Assume $\alpha_{2}\geq\alpha_{1}$.
Then $a_{2}\in K_{\alpha_{1}}(\alpha_{2})<a\leq b$, and $m>2$.

Let $\{\alpha_{i}\}_{2\leq i\leq m}$ be the sequence of subterms of $\alpha_{2}$ such that
$\alpha_{i}=\psi_{\sigma_{i}}^{h_{i}}(a_{i})$ for some $\sigma_{i},a_{i},h_{i}$,
and for each $i<m$,
$\beta_{1}\leq\alpha_{i+1}\in\mathcal{E}(C_{i})$ for $C_{i}=\{\sigma_{i},a_{i}\}\cup SC_{\mK}(h_{i})$.

Let $\{n_{j}\}_{0\leq j\leq k}\, (0<k\leq m-2)$ be the increasing sequence $n_{0}<n_{1}<\cdots<n_{k}\leq m$ 
defined recursively by $n_{0}=2$, and assuming $n_{j}$ has been defined so that
$n_{j}<m$ and $\alpha_{n_{j}}\geq\alpha_{1}$, $n_{j+1}$ is defined by
$n_{j+1}=\min(\{i: n_{j}< i<m, \alpha_{i}<\alpha_{n_{j}}\}\cup\{m\})$.
If either $n_{j}=m$ or $\alpha_{n_{j}}<\alpha_{1}$, then $k=j$ and $n_{j+1}$ is undefined.
Then we claim that
\begin{equation}\label{eq:case4A}
\forall j\leq k(\alpha_{n_{j}}\in\mathcal{W}) \spand \alpha_{n_{k}}<\alpha_{1}
\end{equation}
\textbf{Proof} of (\ref{eq:case4A}).
By induction on $j\leq k$ we show first that $\forall j\leq k(\alpha_{n_{j}}\in\mathcal{W})$. 
We have $\alpha_{n_{0}}=\alpha_{2}\in\mathcal{W}$.
Assume $\alpha_{n_{j}}\in\mathcal{W}$ and $j<k$.
Then $n_{j}<m$, i.e., $\alpha_{n_{j+1}}<\alpha_{n_{j}}$, and 
by $\alpha_{n_{j}}\in\mathcal{C}^{\alpha_{n_{j}}}(\mathcal{W})$, we have $C_{n_{j}}\subset\mathcal{C}^{\alpha_{n_{j}}}(\mathcal{W})$,
and hence $\alpha_{n_{j}+1}\in\mathcal{E}(C_{n_{j}})\subset\mathcal{C}^{\alpha_{n_{j}}}(\mathcal{W})$.
We see inductively that
$\alpha_{i}\in \mathcal{C}^{\alpha_{n_{j}}}(\mathcal{W})$ for any $i$ with $n_{j}\leq i\leq n_{j+1}$.
Therefore $\alpha_{n_{j+1}}\in \mathcal{C}^{\alpha_{n_{j}}}(\mathcal{W})\cap\alpha_{n_{j}}\subset\mathcal{W}$ by Proposition \ref{prp:updis}.

Next we show that $\alpha_{n_{k}}<\alpha_{1}$.
We can assume that $n_{k}=m$.
This means that $\forall i(n_{k-1}\leq i<m \Rarw \alpha_{i}\geq\alpha_{n_{k-1}})$.
We have
$\alpha_{2}=\alpha_{n_{0}}>\alpha_{n_{1}}>\cdots>\alpha_{n_{k-1}}\geq\alpha_{1}$, and
$\forall i<m(\alpha_{i}\geq\alpha_{1})$.
Therefore $\alpha_{m}\in k_{\alpha_{1}}(\alpha_{2})\subset k_{\alpha_{1}}(\{\kappa,a\}\cup SC_{\mK}(h))$, i.e.,
$a_{m}\in K_{\alpha_{1}}(\{\kappa,a\}\cup SC_{\mK}(h))$ for $\alpha_{m}=\psi_{\sigma_{m}}^{h_{m}}(a_{m})$.
On the other hand we have $K_{\alpha_{1}}(\{\kappa,a\}\cup SC_{\mK}(h))<a$ for 
$\alpha_{1}=\psi_{\sigma}^{h}(a)$.
Thus $a\leq a_{m}<a$, a contradiction.

(\ref{eq:case4A}) is shown, and we obtain $\beta_{1}\leq\alpha_{n_{k}}\in\mathcal{W}\cap\alpha_{1}$.

This completes a proof of (\ref{clm:id5wf21.1}) and of the lemma.
\hspace*{\fill} $\Box$

\begin{lemma}\label{lem:psiw}
For {\rm each} $\alp\in OT(\Pi_{3})$,
$\alp\in\mathcal{C}^{\mK}(\calw)$.
\end{lemma}
\bprf
This is seen by meta-induction on $\ell\alp$.
By Propositions \ref{prp:noncritical}, \ref{lem:6.29}, and Lemma \ref{lem:3wf19-1},
we may assume $\alp=\psi_{\kap}^{\nu}(a)$.
By IH pick an $n<\ome$ such that $\{\kap,\nu,a\}\subset\mathcal{C}^{\mK}(\calw)\cap\ome_{n+1}(\mK+1)$.
Lemma \ref{th:id5wf21} yields $\alp\in\calw$.
\hspace*{\fill} $\Box$
\\

Theorem \ref{th:wfpi3} follows from Lemma \ref{lem:psiw} and
the fact $\calw\cap\Ome=W(\mathcal{C}^{0}(\emptyset))\cap\Ome=W(OT(\Pi_{3}))\cap\Ome$.

\section{$\Pi_{4}$-reflection}\label{sec:pi4}
\begin{quote}
In this paper we focus on the ordinal analysis of $\Pi_{3}$ reflection.
This means no genuine loss of generality, as the removal of $\Pi_{3}$ reflection rules in derivations 
already exhibits the pattern of cut elimination that applies for arbitrary $\Pi_{n}$ reflection rules as well.
(\cite{Rathjen94})
\end{quote}

In this section $\mK$ denotes either a $\Pi^{1}_{2}$-indescribable cardinal or
a $\Pi_{4}$-reflecting ordinal.
Skolem hull $\calh_{a}(X)$ and a Mahlo class $Mh_{3}^{a}(\xi)$ are defined as in
Definition \ref{df:HKPpi3}:
Let for $\xi>0$,
\[
\pi\in Mh_{3}^{a}(\xi):\Lrarw
\left[
\{a,\xi\}\subset\calh_{a}(\pi)  \spand
\fal\nu\in\calh_{a}(\pi)\cap\xi\left(
\pi\in M_{3}(Mh_{3}^{a}(\nu))
\right)
\right]
\]
where
$\alp\in M_{3}(A)$ iff $A$ is $\Pi^{1}_{1}$-indescribable in $\alp$ or
$\alp$ is $\Pi_{3}$-reflecting on $A$.

Then as in (\ref{eq:psiPi3})
\[
\psi_{\pi}^{\xi}(a)=
\min\left(\{\pi\}\cup
\{\kap\in Mh_{3}^{a}(\xi): \{\xi,\pi,a\}\subset\calh_{a}(\kap) \spand
\calh_{a}(\kap)\cap\pi\subset\kap\}
\right)
\]
where $\xi=m_{3}(\psi_{\pi}^{\xi}(a))$.

As in Lemmas \ref{lem:welldefinedness} and \ref{lem:psiK} we see the following
for $\Pi^{1}_{2}$-indescribable cardinal $\mK$.

\blem\label{lem:psiK3}
Let $a\in\calh_{a}(\mK)\cap\veps_{\mK+1}$. Then
$\mK\in M_{3}(Mh^{a}_{3}(\veps_{\mK+1}))$.
For every $\xi\in\calh_{a}(\mK)\cap\veps_{\mK+1}$,
$\psi_{\mK}^{\xi}(a)<\mK$.
\elem

Operator controlled derivations for ${\sf KP}\Pi_{4}$ are closed under the following inference rules.
For convenience let us attach an assignment 
$\bar{m}:\pi\mapsto\bar{m}(\pi)=(\bar{m}_{2}(\pi),\bar{m}_{3}(\pi))$ to the derivations,
where $\bar{m}_{i}(\pi)\leq m_{i}(\pi)$ for $i=2,3$.
Although our derivability relation should be written as
$(\calh_{\gam}[\Tht],\bar{m})\vdash^{a}_{b}\Gam$,
let us write $\calh_{\gam}[\Tht]\vdash^{a}_{b}\Gam$.

\bdes
\item[$({\rm rfl}_{\Pi_{4}}(\mK))$]
$b\geq\mK$.
There exist 
an ordinal $a_{0}\in\calh_{\gam}[\Tht]\cap a$,
and a $\Sig_{4}(\mK)$-sentence $A$ enjoying the following conditions:

\[
\infer[({\rm rfl}_{\Pi_{4}}(\mK))]{
\calh_{\gam}[\Tht]\vdash^{a}_{b}\Gam}
{
\mathcal{H}_{\gam}[\Tht]\vdash^{a_{0}}_{b}\Gamma, \lnot A
&
\{
\mathcal{H}_{\gam}[\Tht\cup \{\rho\}]\vdash^{a_{0}}_{b}\Gamma, 
A^{(\rho,\mK)}: \rho<\mK\}
}
\]

\item[$({\rm rfl}_{\Pi_{3}}(\alp,\pi,\nu))$]

There exist ordinals 
$\alp<\pi\leq b<\mK$, 
$\nu<\bar{m}_{3}(\pi)\leq m_{3}(\pi)$ with $SC_{\mK}(\nu)\subset\pi$ and $\nu\leq\gam$, 
$a_{0}<a$,
and a finite set $\Del$ of $\Sig_{3}(\pi)$-sentences enjoying the following conditions:

\benu

\item 
$\{\alp,\pi,\nu\}\cup\bar{m}(\pi)\subset\calh_{\gam}[\Tht]$.

 \item
For each $\del\in\Del$,
$
\mathcal{H}_{\gam}[\Tht]\vdash^{a_{0}}_{b}\Gamma, \lnot\del
$.

\item

Let
\[
\rho\in Mh_{3}(\nu):\Lrarw
\nu\leq m_{3}(\rho)
.\]

Then for each $\alp<\rho\in Mh_{3}(\nu)\cap\pi$,
$\mathcal{H}_{\gam}[\Tht\cup \{\rho\}]\vdash^{a_{0}}_{b}\Gamma, 
\Del^{(\rho)}$.
\eenu

{\small
\[
\hspace{-10mm}
\infer[({\rm rfl}_{\Pi_{3}}(\alp,\pi,\nu))]{
\calh_{\gam}[\Tht]\vdash^{a}_{b}\Gam}
{
\{
\mathcal{H}_{\gam}[\Tht]\vdash^{a_{0}}_{b}\Gamma, \lnot\del
\}_{\del\in\Del}
&
\{
\mathcal{H}_{\gam}[\Tht\cup \{\rho\}]\vdash^{a_{0}}_{b}\Gamma, 
\Del^{(\rho,\pi)}\}_{\alp<\rho\in Mh_{3}(\nu)\cap \pi}
}
\]
}

\edes

Finite proofs in ${\sf KP}\Pi_{4}$ are embedded to controlled derivations with 
inferences $({\rm rfl}_{\Pi_{4}}(\mK))$, and then $({\rm rfl}_{\Pi_{4}}(\mK))$ is replaced by
inferences $({\rm rfl}_{\Pi_{3}}(\alp,\pi,\nu))$ as in Lemma \ref{lem:lowerPi3}.

\blem\label{lem:lowerPi4}
Assume $\Gam\subset\Sig_{3}(\mK)$, $\Tht\subset\calh_{\gam}(\psi_{\mK}(\gam))$,
and
$\calh_{\gam}[\Tht]\vdash^{a}_{\mK}\Gam$ with $a\leq\gam$.
Then ,
$\calh_{\hat{a}+1}[\Tht\cup \{\kap\}]\vdash^{\bet}_{\bet}\Gam^{(\kap,\mK)}$ holds for every
$\kap\in Mh_{3}(a)\cap \psi_{\mK}(\gam+\mK\cdot\ome)$ such that $\psi_{\mK}(\gam+\mK)<\kap$,
where $\hat{a}=\gam+\ome^{\mK+a}$ and $\bet=\psi_{\mK}(\hat{a})$.
\elem

Let us try to eliminate inferences $({\rm rfl}_{\Pi_{3}}(\alp,\pi,\nu))$ from the resulting derivations
following the proof of Lemma \ref{lem:lowerPi3}.
Let $Mh_{2}(\xi;a)$ be a Mahlo class for which the following holds.

\blem\label{lem:lowerPi43}
Let $\Gam\subset\Sig_{2}(\pi)$ with $\xi=m_{3}(\pi)$, 
 and
$\calh_{\gam}[\Tht]\vdash^{a}_{\pi}\Gam$.
Then for any $\kap\in Mh_{2}(\xi;a)\cap \pi$,
$\calh_{\gam}[\Tht\cup\{\kap\}]\vdash^{\kap+\ome a}_{\pi}\Gam^{(\kap,\pi)}$ holds\footnote{Here we don't need to collapse derivations and cut ranks$<\pi$.}.
\elem
Consider the crucial case.
Let $\Del\subset\Sig_{3}(\pi)$, $\pi\in Mh_{3}(\xi)$ and $\nu<\xi$.

\[
\infer[({\rm rfl}_{\Pi_{3}}(\alp,\pi,\nu))]{
\calh_{\gam}[\Tht]\vdash^{a}_{\pi}\Gam}
{
\{
\mathcal{H}_{\gam}[\Tht]\vdash^{a_{0}}_{\pi}\Gamma, \lnot\del
\}_{\del\in\Del}
&
\{
\mathcal{H}_{\gam}[\Tht\cup\{\rho\}]\vdash^{a_{0}}_{\pi}\Gamma, 
\Del^{(\rho,\pi)}: \alp<\rho\in Mh_{3}(\nu)\cap \pi\}
}
\]
Let $\sig\in Mh_{2}(\xi;a_{0})\cap \kap$.
By IH with Inversion we obtain
$\mathcal{H}_{\gam}[\Tht\cup\{\sig\}]\vdash^{\kap+\ome a_{0}+1}_{\pi}\Gamma^{(\sig,\pi)}, \lnot \del^{(\sig,\pi)}$
for each $\del\in\Del$.

On the other hand we have
$\mathcal{H}_{\gam}[\Tht\cup\{\sig\}]\vdash^{a_{0}}_{\pi}\Gamma, 
\Del^{(\sig,\pi)}$
for $\alp<\sig\in Mh_{3}(\nu)\cap \pi$.
Assume 
$Mh_{2}(\xi;a)\subset Mh_{2}(\xi;a_{0})$.
IH yields
$\mathcal{H}_{\gam}[\Tht\cup\{\kap,\sig\}]\vdash^{\kap+\ome a_{0}}_{\pi}\Gamma^{(\kap,\pi)}, 
\Del^{(\sig,\pi)}$.

Let
$
\alp<\sig\in Mh_{2}(\xi;a_{0})\cap Mh_{3}(\nu)\cap \kap$.
A $(cut)$ of the cut formulas $\del^{(\sig,\pi)}$ then yields
$\mathcal{H}_{\gam}[\Tht\cup\{\kap,\sig\}]\vdash^{\kap+\ome a_{0}+p}_{\pi}\Gamma^{(\kap,\pi)},\Gam^{(\sig,\pi)}$ for a $p<\ome$.

On the other hand we have 
$\mathcal{H}_{\gam}[\Tht\cup\{\kap\}]\vdash^{2d}_{0}\lnot\tht^{(\kap,\pi)},\Gam^{(\kap,\pi)}$
for each $\tht\in\Gam\subset\Sig_{2}(\pi)$, where
$d=\max\{\rk(\tht^{(\kap,\pi)}):\tht\in\Gam\}<\kap+\ome<\pi$.

Now $\kap\in Mh_{2}(\xi;a)\cap \pi$
needs to reflect $\Pi_{2}(\kap)$-formulas
$\lnot\tht^{(\kap,\pi)}$ down to some $\alp<\sig\in Mh_{2}(\xi;a_{0})\cap Mh_{3}(\nu)\cap \kap$.

\[
a_{0}<a \spand \nu<\xi \Rarw
Mh_{2}(\xi;a)\subset M_{2}(Mh_{2}(\xi;a_{0})\cap Mh_{3}(\nu))
\]

Thus we arrive at the following definition of the Mahlo classes $Mh_{2}^{\gam}(\xi;a)$,
which is a $\Pi_{3}$-class in the sense that
there is a $\Pi_{3}$-formula $\tht(\gam,\xi,a)$ such that
$\alp\in Mh_{2}^{\gam}(\xi;a)$  iff $L_{\alp}\models \tht(\gam,\xi,a)$,
while $Mh_{3}^{\gam}(\nu)$ is a $\Pi_{4}$-class.

$\pi\in Mh_{2}^{\gam}(\xi;a)$ iff $\{\gam,\xi, a\}\subset\calh_{\gam}(\pi)$ and 
\[
\fal \{\nu,b\}\subset\calh_{\gam}(\pi)\left[\nu<\xi \spand b<a  \Rarw
\pi\in M_{2}\left(
Mh_{2}^{\gam}(\xi;b)\cap Mh_{3}^{\gam}(\nu)
\right)
\right]
.\]

It turns out that we need Mahlo classes $Mh_{2}^{\gam}(\bar{\xi};\bar{a})$
for finite sequences $\bar{\xi}$ and $\bar{a}$ in our proof-theoretic study, 
cf.\,Lemma \ref{lem:lowerPi43b}.
Let us explain the classes intuitively in the next subsection.

\subsection{Mahlo classes}\label{subsec:Mclasspi3}
Let $M_{i}=RM_{i}$ and
$P,Q,\ldots$ denote transitive classes in $L\cup\{L\}$ for a $\Pi_{4}$-reflecting universe $L$.
For classes $\calx,\caly$ and $i=2,3$ let
\[
\calx\prec_{i}\caly :\Lrarw
\fal P\in\caly(P\in M_{i}(\calx))
\]

\bdf
{\rm Let
\[
M_{2}(\xi;a):= 
\bigcap\{M_{2}\left(M_{2}(\xi;b)\cap M_{3}(\nu)\right):
\nu<\xi, b<a\}
.\]
In general for classes $\caly$ let
\[
M_{2}^{\caly}(\xi;a):= \caly\cap
\bigcap\{M_{2}\left(M_{2}^{\caly}(\xi;b)\cap M_{3}(\nu)\right):
\nu<\xi, b<a\}
.\]
}
\edf

\bprp
\label{prp:5aprl1}
For a $\Pi_{3}$-class $\caly$ and $\mu<\xi$,
$
M_{2}^{\caly}(\xi;a)\cap M_{3}(\mu) \prec_{2}
\caly\cap M_{3}(\xi)
$ and
$M_{2}^{\caly}(\xi;a)\supset
\caly\cap M_{3}(\xi)$.
\eprp
\bprf
By induction on $a$, we show
$P\in\caly\cap M_{3}(\xi) \Rarw 
P\in M_{2}^{\caly}(\xi;a)$.

Let $P\in\caly\cap M_{3}(\xi)$, $\nu<\xi$ and $b<a$.
By IH we obtain $P\in M_{2}^{\caly}(\xi;b)$. 
Since $M_{2}^{\caly}(\xi;b)$ is a $\Pi_{3}$-class, we obtain
$P\in M_{2}\left(M_{2}^{\caly}(\xi;b)\cap M_{3}(\nu)\right)$ by $P\in M_{3}(\xi)$.
Therefore $P\in M_{2}^{\caly}(\xi;a)$.

Since $M_{2}^{\caly}(\xi;a)$ is a $\Pi_{3}$-class and $P\in M_{3}(\xi)\subset M_{3}(M_{3}(\mu))$,
we obtain $P\in M_{2}(M_{2}^{\caly}(\xi;a)\cap M_{3}(\mu))$.
\eprf
\\

Let $\nu<\mu<\xi$. From Proposition \ref{prp:5aprl1} we see
$
M_{2}(\xi;a)\cap M_{3}(\mu) \prec_{2}M_{3}(\xi)
$, and 
$M_{2}^{\caly}(\mu;b)\cap M_{3}(\nu) \prec_{2} \caly\cap M_{3}(\mu)$ for $\caly=M_{2}(\xi;a)$.

Let us write $M_{2}((\xi,\mu);(a,b))$ for $M_{2}^{\caly}(\mu;b)$, where $\xi>\mu$.
Let $\nu<\mu<\xi$.
We obtain
$M_{2}((\xi,\mu);(a,b))\cap M_{3}(\nu) \prec_{2} M_{2}(\xi;a)\cap M_{3}(\mu)\prec_{2}M_{3}(\xi)$.

\bprp\label{prp:pi4pair}
Let $\xi_{1},\zeta<\xi$, $c<b$ and $d<a$. Then
$M_{2}((\xi,\mu);(a,c))\cap M_{3}(\nu)\prec_{2}M_{2}((\xi,\mu);(a,b))$ and
$M_{2}((\xi,\xi_{1});(d,e))\cap M_{3}(\zeta) \prec_{2} M_{2}((\xi,\mu);(a,b))$.
\eprp
\bprf
Let $\caly=M_{2}(\xi;a)$. Then
$M_{2}((\xi,\mu);(a,c))\cap M_{3}(\nu)=M_{2}^{\caly}(\mu;c)\cap M_{3}(\nu)\prec_{2}
M_{2}^{\caly}(\mu;b)=M_{2}((\xi,\mu);(a,b))$ by $c<b$ and $\nu<\mu$.

Next we show $M_{2}^{\calx}(\xi_{1};e)\cap M_{3}(\zeta) \prec_{2}\caly\supset M_{2}^{\caly}(\mu;b)$,
where
$\calx=M_{2}(\xi;d)$ and $M_{2}((\xi,\xi_{1});(d,e))=M_{2}^{\calx}(\xi_{1};e)$.
We have 
$\calx\cap M_{3}(\xi_{1})\cap M_{3}(\zeta)=M_{2}(\xi;d)\cap M_{3}(\xi_{1})\cap M_{3}(\zeta)\prec_{2}M_{2}(\xi;a)=\caly$
by $d<a$ and $\xi_{1},\zeta<\xi$.
On the other hand we have
$M_{2}^{\calx}(\xi_{1};e)\supset\calx\cap M_{3}(\xi_{1})$
by Proposition \ref{prp:5aprl1}.
Hence $M_{2}^{\calx}(\xi_{1};e)\cap M_{3}(\zeta) \prec_{2}\caly$.
\eprf

The same argument applies not only to pairs $(\xi>\mu)$, $(a,b)$, but also to
triples, and so forth.

Let $\bar{\xi}=(\xi_{0}>\xi_{1}>\cdots>\xi_{n})$ and
$\bar{a}=(a_{0},a_{1},\ldots,a_{n})$ be sequences in the same lengths.
By iterating the process $\caly\mapsto \{M_{2}^{\caly}(\xi;a)\}_{a}$
with $M_{3}(\xi)$, we now define 
classes $M_{2}(\bar{\xi};\bar{a})$ 
by induction on the length $n$ of the sequences $\bar{\xi}, \bar{a}$
as follows. 

$M_{2}(\langle\rangle;\langle\rangle)$ denotes the class of transitive sets in $\mbox{L}\cup\{\mbox{L}\}$.

For  $\bar{\xi}*(\xi)=(\xi_{0}>\cdots>\xi_{n}>\xi)$ and $\bar{a}*(a)=(a_{0},\ldots,a_{n},a)$ define 
for  the $\Pi_{3}$-class
 $\caly=M_{2}(\bar{\xi};\bar{a})$
\[
M_{2}(\bar{\xi}*(\xi);\bar{a}*(a)) =M_{2}^{\caly}(\xi;a)
\]
Namely
\[
M_{2}(\bar{\xi}*(\xi);\bar{a}*(a)) = M_{2}(\bar{\xi};\bar{a})\cap 
\bigcap\{M_{2}\left(M_{2}(\bar{\xi}*(\xi);\bar{a}*(b))\cap M_{3}(\nu)\right):  \nu<\xi, b<a\}  
\]

Proposition \ref{prp:pi4pair} is extended to finite sequences.
To state an extension, let us redefine classes $M_{2}(\bar{\xi};\bar{a})$
through ordinals
$\alp=\Lam^{\xi_{0}}a_{0}+\cdots+\Lam^{\xi_{n}}a_{n}$ as follows, where
$\Lam$ is a big enough ordinal such that $\Lam>a_{0}$.

Let $\alp=\Lam^{\xi_{0}}a_{0}+\cdots+\Lam^{\xi_{n}}a_{n}$, where $\xi_{0}>\cdots>\xi_{n}$ and
$a_{0},\ldots,a_{n}\neq 0$.
\[
M_{2}(\alp):=\bigcap\{
M_{2}(M_{2}(\bet)\cap M_{3}(\nu)):(\bet,\nu)<\alp\}
\]
where
 for segments $\alp_{i}=\Lam^{\xi_{0}}a_{0}+\cdots+\Lam^{\xi_{i}}a_{i}$ of 
$\alp=\Lam^{\xi_{0}}a_{0}+\cdots+\Lam^{\xi_{n}}a_{n}$
\[
(\bet,\nu)< \alp :\Lrarw \exi i\leq n\left[
\bet<\alp_{i} \spand \nu<\xi_{i}
\right]
.\]
F.e. in Proposition \ref{prp:pi4pair} we have
$(\Lam^{\xi}a+\Lam^{\mu}c,\nu)<\Lam^{\xi}a+\Lam^{\mu}b$ and $(\Lam^{\xi}d+\Lam^{\xi_{1}}e, \zeta)<\Lam^{\xi}a+\Lam^{\mu}b$, but
$(\Lam^{\xi}a+\Lam^{\mu}c,\mu)\not<\Lam^{\xi}a+\Lam^{\mu}b$,
where
$\nu<\mu<\xi$, $\xi_{1},\zeta<\xi$, $c<b$ and $d<a$.

\bprp\label{prp:less1}
$(\bet,\nu)<\alp<\gam \Rarw (\bet,\nu)<\gam$.
\eprp

$\alp\dot{+}\bet$ designates that $\alp+\bet=\alp\#\bet$.

\blem\label{lem:stepdownint}{\rm (Cf.\,Lemma 3.2 in \cite{LMPS}.)}\\
If $\xi>0$ and $\bet<\Lam^{\xi+1}$, then
$M_{2}(\alp\dot{+}\bet)\prec_{2} M_{2}(\alp,\xi):=M_{2}(\alp)\cap M_{3}(\xi)$.
\elem
\bprf
Suppose $P\in M_{2}(\alp,\xi)=M_{2}(\alp)\cap M_{3}(\xi)$ and $\bet<\Lam^{\xi+1}$.

We show $P\in M_{2}(\alp\dot{+}\bet)$ 
by induction on ordinals $\bet$.
Let $(\gam,\nu)<\alp\dot{+}\bet$.
We need to show that $P \in M_{2}(M_{2}(\gam,\nu))$.

Let $\del$ be a segment of $\alp\dot{+}\bet$ such that
$\gam<\del$ and $\nu<\mu$ where
$\del=\cdots+\Lam^{\mu}b$.
If $\del$ is a segment of $\alp$, then $P\in M_{2}(M_{2}(\gam,\nu))$
by $P\in M_{2}(\alp)$.

Let
$\del=\alp\dot{+}\bet_{0}$, where $\bet_{0}$ is a segment of $\bet$.
Then $\nu<\mu\leq\xi$.
We claim that $P\in M_{2}(\gam)$.
If $\gam<\alp$, then Proposition \ref{prp:less1} 
yields $P\in M_{2}(\alp)\subset M_{2}(\gam)$.
Let $\gam=\alp\dot{+}\gam_{0}<\alp\dot{+}\bet_{0}$.
IH yields $P\in M_{2}(\gam)$.
Thus the claim is shown.
On the other hand we have $P\in M_{3}(\xi)$ and $\nu<\xi$.
Since $M_{2}(\gam)$ is a $\Pi_{3}$-class, we obtain
$P\in M_{3}(M_{2}(\gam,\nu))\subset M_{2}(M_{2}(\gam,\nu))$.
$P\in M_{2}(\alp\dot{+}\bet)$ is shown.

By $P\in M_{2}(\alp\dot{+}\bet)$ and $P\in M_{3}(\xi)\subset M_{3}$ with $\xi>0$, 
we obtain $P\in M_{3}(M_{2}(\alp\dot{+}\bet))\subset M_{2}(M_{2}(\alp\dot{+}\bet))$.
\eprf

\subsection{Skolem hulls and collapsing functions}

We can assume $\xi<\veps_{\mK+1}$ and $a<\Lam=\mK$.
For $\alp<\Lam^{\veps_{\mK+1}}$, let us define $Mh_{2}^{\gam}(\alp)$ as follows.
$(\bet,\nu)$ denotes pairs of ordinals $\bet<\Lam^{\veps_{\mK+1}}$ and $\nu<\veps_{\mK+1}$
such that $\bet+\Lam^{\nu}=\bet\#\Lam^{\nu}$.
Let $\alp=\Lam^{\bet_{0}}a_{0}+\cdots+\Lam^{\bet_{n}}a_{n}$, where
$\veps_{\mK+1}>\bet_{0}>\cdots>\bet_{n}$ and $0<a_{0},\ldots,a_{n}<\Lam$.
Then $\pi\in Mh_{2}^{\gam}(\alp)$ iff 
$\{\gam,\alp\}\subset\calh_{\gam}(\pi)$ and 
\[
\fal \{\nu,\bet\}\subset\calh_{\gam}(\pi)\left[(\bet,\nu)< \alp  \Rarw
\pi\in M_{2}\left(
Mh_{2}^{\gam}(\bet)\cap Mh_{3}^{\gam}(\nu)
\right)
\right]
\]
where for segments $\alp_{i}=\Lam^{\bet_{0}}a_{0}+\cdots+\Lam^{\bet_{i}}a_{i}$ of 
$\alp=\Lam^{\bet_{0}}a_{0}+\cdots+\Lam^{\bet_{n}}a_{n}$
\[
(\bet,\nu)< \alp :\Lrarw \exi i\leq n\left[
\bet<\alp_{i} \spand \nu<\bet_{i}
\right]
.\]
For example, if $\nu<\xi$ and $a_{0}<a$, then $(\Lam^{\xi}a_{0},\nu)<\Lam^{\xi}a$.
The exponents $\bet_{i}$ of $\alp$ designate `$\Pi_{3}$-Mahlo degrees'.

\bprp\label{prp:less}
$(\bet,\nu)<\alp<\gam \Rarw (\bet,\nu)<\gam$.
\eprp

\bdf\label{df:HKPpi4}
{\rm
Define simultaneously by recursion on ordinals $a<\veps_{\mK+1}$ the classes
 $\calh_{a}(X)\, (X\subset \Gam_{\mK+1})$, $Mh_{2}^{a}(\alp)\,(\xi<\veps_{\mK+1})$,
  the ordinals $\psi_{\sig}^{(\alp,\xi)}(a)$ as follows.

\benu
\item
$\calh_{a}(X)$ denotes the Skolem hull of $\{0,\Ome,\mK\}\cup X$ under the functions
$+, \vphi$, 
and the following.

 Let $\{\sig,b,\alp,\xi\}\subset\calh_{a}(X)$, $\alp\in\{0\}\cup[\Lam,\Lam^{\veps_{\mK+1}})$, 
 $\xi\in[0,\veps_{\mK+1})$  and 
 $b<a$.
 Then $\psi_{\sig}^{(\alp,\xi)}(b)\in\calh_{a}(X)$.

\item
$
\pi\in Mh_{3}^{a}(\xi):\Lrarw
\{a,\xi\}\subset\calh_{a}(\pi) \spand
\fal\nu\in\calh_{a}(\pi)\cap\xi\left(
\pi\in M_{3}(Mh_{3}^{a}(\nu))
\right)
$,
where $\alp\in Mh_{3}^{a}(0)$ iff $\alp$ is a limit ordinal.

\item
For $\alp<\Lam^{\veps_{\mK+1}}$ and $a<\veps_{\mK+1}$,
$\pi\in Mh_{2}^{a}(\alp)$ iff $\{a,\alp\}\subset\calh_{a}(\pi)$ and
\[
\fal\{\bet,\nu\}\subset\calh_{a}(\pi)
\left[
(\bet,\nu)<\alp \to
\pi\in M_{2}\left(
Mh_{2}^{a}(\bet,\nu)
\right)
\right]
\]
where
\[
Mh_{2}^{a}(\bet,\nu)=Mh_{2}^{a}(\bet)\cap Mh_{3}^{a}(\nu)
\]
and $\alp\in Mh_{2}^{a}(0)$ iff $\alp$ is a limit ordinal.
Note that $Mh_{2}^{a}(\alp)$ is a $\Pi_{3}$-class.

\item
Let $m_{2}(\mK)=0$, $m_{3}(\mK)=\veps_{\mK+1}$, $m_{2}(\Ome)=1$ and $m_{3}(\Ome)=0$.
 \benu
\item 
For $\{\xi,a\}\subset\calh_{a}(\mK)\cap\veps_{\mK+1}$ with $0<\xi\leq a$, let
\\
$\psi_{\mK}^{(0,\xi)}(a)=
\min\left(\{\mK\}\cup
\{\kap\in Mh_{3}^{a}(\xi): \{\xi,a\}\subset\calh_{a}(\kap) \spand
\calh_{a}(\kap)\cap\mK\subset\kap\}
\right)$.
$m_{2}(\psi_{\mK}^{(0,\xi)}(a))=0$ and $m_{3}(\psi_{\mK}^{(0,\xi)}(a))=\xi$.

 \item
Let
$0\leq\alp< \Lam^{\veps_{\mK+1}}$ and $0<\xi<\veps_{\mK+1}$
be ordinals, 
$0<c\leq a<\Lam=\mK$ with $c\in\calh_{a}(\sig)$
and
$\sig\in Mh_{2}^{a}(\alp,\xi)$.
Then for $\bet=\alp\dot{+}\Lam^{\xi}c$

\[
\psi_{\sig}^{(\bet,0)}(a)=\min\left(\{\sig\}\cup
\{\kap\in Mh_{2}^{a}(\bet): \{\sig,\alp,\xi,c,a\}\subset\calh_{a}(\kap) \spand
\calh_{a}(\kap)\cap\sig\subset\kap\}
\right).
\]
$m_{2}(\psi_{\sig}^{(\bet,0)}(a))=\bet$ and $m_{3}(\psi_{\sig}^{(\bet,0)}(a))=0$.

\item
Let $0<\bet,\alp<\Lam^{\veps_{\mK+1}}$ and $0< \nu<\veps_{\mK+1}$ be such that
$\{\bet,\nu\}\subset\calh_{a}(\sig)$, $SC_{\mK}(\bet,\nu)\subset (a+1)<\mK$ and
$(\bet,\nu)<\alp$.
Then for
$\sig\in Mh_{2}^{a}(\alp)$ with $m_{3}(\sig)=0$
\[
\psi_{\sig}^{(\bet,\nu)}(a)=\min\left(\{\sig\}\cup
\{\kap\in Mh_{2}^{a}(\bet,\nu): \{\sig,\bet,\nu,a\}\subset\calh_{a}(\kap) \spand
\calh_{a}(\kap)\cap\sig\subset\kap\}
\right).
\]
$m_{2}(\psi_{\sig}^{(\bet,\nu)}(a))=\bet$ and $m_{3}(\psi_{\sig}^{(\bet,\nu)}(a))=\nu$.

\item
\[
\psi_{\sig}(a)=\min
\{\kap\leq\sig : \{\sig,a\}\subset\calh_{a}(\kap) \spand
\calh_{a}(\kap)\cap\sig\subset\kap\}
.
\]

 \eenu

We write $\psi_{\sig}(a)$ for $\psi_{\sig}^{(0,0)}(a)$.
\eenu
}
\edf

Let $\mK$ be a $\Pi^{1}_{2}$-indescribable cardinal. As in Lemmas \ref{lem:welldefinedness} and \ref{lem:psiK} we see that
$\psi_{\mK}^{(0,\xi)}(a)<\mK$ for every $\{a,\xi\}\subset\calh_{a}(\mK)\cap\veps_{\mK+1}$.

It is easy to see that $\psi_{\sig}^{(\bet,\nu)}(a)<\sig$
if $(\bet,\nu)<\alp$, $\sig\in Mh_{2}^{a}(\alp)$ and
$\{\bet,\nu\}\subset\calh_{a}(\sig)$.

\blem\label{lem:stepdown}{\rm (Cf.\,Lemma 3.2 in \cite{LMPS}.)}
Assume $\mK\geq\sig\in Mh^{a}_{2}(\alp,\xi)$ with
$0<\xi<\veps_{\mK+1}$,
$\bet<\Lam^{\xi+1}$ and
$\bet\in\calh_{a}(\sig)$.
Then
$\sig\in M_{3}(Mh^{a}_{2}(\alp\dot{+}\bet))$ holds, a fortiori $\sig\in M_{2}(Mh^{a}_{2}(\alp\dot{+}\bet))$.
\elem
\bprf
Suppose $\sig\in Mh^{a}_{2}(\alp,\xi)=Mh^{a}_{2}(\alp)\cap Mh^{a}_{3}(\xi)$ and $\bet\in\calh_{a}(\sig)$ with $\bet<\Lam^{\xi+1}$.
We show $\sig\in Mh^{a}_{2}(\alp\dot{+}\bet)$ 
by induction on ordinals $\bet$.
Let $\{\gam,\nu\}\subset\calh_{a}(\sig)$ and
$(\gam,\nu)<\alp\dot{+}\bet$.
We need to show that $\sig\in M_{2}(Mh^{a}_{2}(\gam,\nu))$.

Let $\del$ be a segment of $\alp\dot{+}\bet$ such that
$\gam<\del$ and $\nu<\mu$ where
$\del=\cdots+\Lam^{\mu}b$.
If $\del$ is a segment of $\alp$, then $\sig\in M_{2}(Mh^{a}_{2}(\gam,\nu))$
by $\sig\in Mh^{a}_{2}(\alp)$.

Let
$\del=\alp\dot{+}\bet_{0}$, where $\bet_{0}$ is a segment of $\bet$.
Then $\nu<\mu\leq\xi$.
We claim that $\sig\in Mh_{2}^{a}(\gam)$.
If $\gam<\alp$, then Proposition \ref{prp:less} with $\gam\in\calh_{a}(\sig)$
yields $\sig\in Mh^{a}_{2}(\alp)\subset Mh^{a}_{2}(\gam)$.
Let $\gam=\alp\dot{+}\gam_{0}<\alp\dot{+}\bet_{0}$ with $\gam_{0}\in\calh_{a}(\sig)$.
IH yields $\sig\in Mh_{2}^{a}(\gam)$.
Thus the claim is shown.
On the other hand we have $\sig\in Mh_{3}^{a}(\xi)$ and $\nu\in\calh_{a}(\sig)\cap\xi$.
Since $Mh_{2}^{a}(\gam)$ is a $\Pi_{3}$-class, we obtain
$\sig\in M_{3}(Mh_{2}^{a}(\gam,\nu))\subset M_{2}(Mh_{2}^{a}(\gam,\nu))$ 
with $Mh_{2}^{a}(\gam,\nu)=Mh_{2}^{a}(\gam)\cap Mh_{3}^{a}(\nu)$.
$\sig\in Mh^{a}_{2}(\alp\dot{+}\bet)$ is shown.

By $\sig\in Mh^{a}_{2}(\alp\dot{+}\bet)$ and $\sig\in Mh^{a}_{3}(\xi)\subset M_{3}$ with $\xi>0$, 
we obtain $\sig\in M_{3}(Mh^{a}_{2}(\alp\dot{+}\bet))$.
\eprf

\bcor\label{cor:stepdownpi4}
If $\sig\in Mh^{a}_{2}(\alp,\xi)$
 and $c\in\calh_{a}(\sig)\cap\Lam$ with $\xi>0$, then
$\psi_{\sig}^{(\bet,0)}(a)<\sig$ for $\bet=\alp\dot{+}\Lam^{\xi}c$.
\ecor
\bprf
We obtain $\sig\in M_{2}(Mh_{2}^{a}(\bet))$ by Lemma \ref{lem:stepdown}.
Since $\{\kap<\sig: \{\bet,a,\sig\}\subset\calh_{a}(\kap), \calh_{a}(\kap)\cap\sig\subset\kap\}$
is a club subset of $\sig$,
we obtain $\psi_{\sig}^{(\bet,0)}(a)<\sig$.
\eprf
\\

$OT(\Pi_{4})$ denotes a computable notation system of ordinals with
collapsing functions $\psi_{\sig}^{(\alp,\xi)}(a)$.
Although in our well-foundedness proof in ${\sf KP}\Pi_{4}$,
ordinal terms $\psi_{\sig}^{(\bet,\nu)}(a)$ has to obey some restrictions such as
(\ref{eq:OTpi3m2}) for $OT(\Pi_{3})$,
it is cumbersome to verify the conditions, and let us skip it.

Operator controlled derivations for ${\sf KP}\Pi_{4}$ are closed under the inference rules
$({\rm rfl}_{\Pi_{4}}(\mK))$, $({\rm rfl}_{\Pi_{3}}(\alp,\pi,\nu))$ and the following.

\bdes
\item[$({\rm rfl}_{\Pi_{2}}(\alp,\pi,\bet,\nu))$]

There exist ordinals 
$\alp<\pi\leq b<\mK$,
$(\bet,\nu)< \bar{m}_{2}(\pi)\leq m_{2}(\pi)$, 
$a_{0}<a$,
and a finite set $\Del$ of $\Sig_{2}(\pi)$-sentences enjoying the following conditions:

\benu

\item
$\{\alp,\pi,\bet,\nu\}\cup\bar{m}(\pi)\subset\calh_{\gam}[\Tht]$.

 \item
For each $\del\in\Del$,
$
\mathcal{H}_{\gam}[\Tht]\vdash^{a_{0}}_{b}\Gamma, \lnot\del
$.

\item
For each $\alp<\rho\in Mh_{2}(\bet,\nu)\cap \pi$,
$\mathcal{H}_{\gam}[\Tht\cup \{\rho\}]\vdash^{a_{0}}_{b}\Gamma, 
\Del^{(\rho,\pi)}$.

\eenu

\[
\hspace{-10mm}
\infer[({\rm rfl}_{\Pi_{2}}(\alp,\pi,\bet,\nu))]{
\calh_{\gam}[\Tht]\vdash^{a}_{b}\Gam}
{
\{
\mathcal{H}_{\gam}[\Tht]\vdash^{a_{0}}_{b}\Gamma, \lnot\del
\}_{\del\in\Del}
&
\{
\mathcal{H}_{\gam}[\Tht\cup \{\rho\}]\vdash^{a_{0}}_{b}\Gamma, 
\Del^{(\rho,\pi)}\}_{\alp<\rho\in Mh_{2}(\bet,\nu)\cap \pi}
}
\]
This inference says that
$\pi\in M_{2}(Mh_{2}^{\gam}(\bet)\cap Mh_{3}^{\gam}(\nu))$.

\edes

\blem\label{lem:lowerPi43b}
Let $\Gam\subset\Sig_{2}(\pi)$.
Assume
$\calh_{\gam}[\Tht]\vdash^{a}_{\pi}\Gam$ for a $\pi<\mK$, and $\{\xi,\alp\}\subset\calh_{\gam}[\Tht]$ 
for $\alp=\bar{m}_{2}(\pi)$, $\xi=\bar{m}_{3}(\pi)$.
Let $\eta$ be the base for $({\rm rfl}_{\Pi_{3}}(\eta,\pi,\nu))$ in $\calh_{\gam}[\Tht]\vdash^{a}_{\pi}\Gam$.
Then for any $\eta<\kap\in Mh_{2}(\alp\dot{+}\Lam^{\xi}(1+a))\cap \pi$,
$\calh_{\gam}[\Tht\cup \{\kap\}]\vdash^{\kap+\ome a}_{\pi}\Gam^{(\kap,\pi)}$ holds, where
$\alp\dot{+}\Lam^{\xi}(1+a)\leq \bar{m}_{2}(\kap)\in\calh_{\gam}[\Tht]$.
Moreover when $\Tht\subset\calh_{\gam}(\kap)$,
$\calh_{\gam}[\Tht\cup \{\kap\}]\vdash^{\kap+\ome a}_{\kap}\Gam^{(\kap,\pi)}$ holds.
\elem
\bprf 
By induction on $a$.
Let $\pi^{\prime}=\kap$ if $\Tht\subset\calh_{\gam}(\kap)$. Otherwise $\pi^{\prime}=\pi$.
Note that there exists a $\kap$ such that $\kap\in Mh_{2}(\alp\dot{+}\Lam^{\xi}(1+a))\cap \pi$
if $\Tht\cup\{\pi\}\subset\calh_{\gam}(\pi)$.
F.e. $\kap=\psi_{\pi}^{(\alp+\Lam^{\xi}(1+a),0)}(\gam+\max\Tht)$.

Let $\eta$ be the base for $({\rm rfl}_{\Pi_{3}}(\eta,\pi,\nu))$ in $\calh_{\gam}[\Tht]\vdash^{a}_{\pi}\Gam$.
\\
\textbf{Case 1}. $({\rm rfl}_{\Pi_{3}}(\eta,\pi,\nu))$:
Then $\eta<\pi$, $\{\eta,\pi,\nu\}\cup\bar{m}(\pi)\subset\calh_{\gam}[\Tht]$, $SC_{\mK}(\nu)\subset\pi$,
and $\nu<\bar{m}_{3}(\pi)\leq m_{3}(\pi)$.
Let $\Del\subset\Sig_{3}(\pi)$.
{\small
\[
\infer[({\rm rfl}_{\Pi_{3}}(\eta,\pi,\nu))]{
\calh_{\gam}[\Tht]\vdash^{a}_{\pi}\Gam}
{
\{
\mathcal{H}_{\gam}[\Tht]\vdash^{a_{0}}_{\pi}\Gamma, \lnot\del
\}_{\del\in\Del}
&
\{
\mathcal{H}_{\gam}[\Tht\cup \{\rho\}]\vdash^{a_{0}}_{\pi}\Gamma, 
\Del^{(\rho,\pi)}\}_{\eta<\rho\in Mh_{3}(\nu)\cap \pi}
}
\]
}
Let $\alp_{0}=\alp\dot{+}\Lam^{\xi}(1+a_{0})$. Then $(\alp_{0},\nu)<\alp_{1}=\alp\dot{+}\Lam^{\xi}(1+a)$.
We obtain $\{\kap,\alp_{1},\nu,\alp_{0}\}\subset\calh_{\gam}[\Tht\cup\{\kap\}]$.
In the following derivation $\alp_{1}\leq\bar{m}_{2}(\kap)$ with $\bar{m}(\kap)\subset\calh_{\gam}[\Tht]$.
{\footnotesize
\[
\hspace{-20mm}
\infer[({\rm rfl}_{\Pi_{2}}(\eta,\kap,\alp_{0},\nu))]{\mathcal{H}_{\gam}[\Tht\cup \{\kap\}]\vdash^{\kap+\ome a}_{\pi^{\prime}}\Gam^{(\kap,\pi)}}
{
\{\mathcal{H}_{\gam}[\Tht\cup \{\kap\}]\vdash^{2d}_{0}\lnot\tht^{(\kap,\pi)},\Gam^{(\kap,\pi)}\}_{\tht\in\Gam}
&
\hspace{-15mm}
\infer{
\{\mathcal{H}_{\gam}[\Tht\cup \{\kap,\sig\}]\vdash^{\kap+\ome a_{0}+p}_{\pi^{\prime}}\Gamma^{(\kap,\pi)},\Gam^{(\sig,\pi)}\}_{\eta<\sig\in Mh_{2}(\alp_{0},\nu)\cap \pi}
}
 {
 \{
 \mathcal{H}_{\gam}[\Tht\cup \{\sig\}]\vdash^{\sig+\ome a_{0}+1}_{\pi^{\prime}}\Gamma^{(\sig,\pi)}, \lnot \del^{(\sig,\pi)}
 \}_{\del\in\Del}
 &
 \mathcal{H}_{\gam}[\Tht\cup \{\kap,\sig\}]\vdash^{\kap+\ome a_{0}}_{\pi^{\prime}}\Gamma^{(\kap,\pi)}, 
\Del^{(\sig,\pi)}
 }
}
\]
}
\textbf{Case 2}. $({\rm rfl}_{\Pi_{2}}(\mu,\pi,\bet,\nu))$:
$(\bet,\nu)<\alp=\bar{m}_{2}(\pi)\leq m_{2}(\pi)$, $\mu<\pi$, 
$\{\mu,\pi,\alp,\bet,\nu\}\subset\calh_{\gam}[\Tht]$ and $\Del\subset\Sig_{2}(\pi)$.
\[
\infer[({\rm rfl}_{\Pi_{2}}(\pi,\bet,\nu))]{
\calh_{\gam}[\Tht]\vdash^{a}_{\pi}\Gam}
{
\{
\mathcal{H}_{\gam}[\Tht]\vdash^{a_{0}}_{\pi}\Gamma, \lnot\del
\}_{\del\in\Del}
&
\{
\mathcal{H}_{\gam}[\Tht\cup \{\rho\}]\vdash^{a_{0}}_{\pi}\Gamma, 
\Del^{(\rho,\pi)}\}_{\mu<\rho\in Mh_{2}(\bet,\nu)\cap \pi}
}
\]
Then $(\bet,\nu)<\alp_{1}=\alp\dot{+}\Lam^{\xi}(1+a)\leq \bar{m}_{2}(\kap)$
with the segment $\alp$ of $\alp\dot{+}\Lam^{\xi}(1+a)$.
We have $\Del^{(\rho,\pi)}=\left(\Del^{(\kap,\pi)}\right)^{(\rho,\kap)}$ and
$\{\kap,\alp_{1},\bet,\nu\}\subset\calh_{\gam}[\Tht\cup\{\kap\}]$.

{\footnotesize
\[
\hspace{-15mm}
\infer[({\rm rfl}_{\Pi_{2}}(\mu,\kap,\bet,\nu))]{
\calh_{\gam}[\Tht\cup \{\kap\}]\vdash^{\kap+\ome a}_{\pi^{\prime}}\Gam^{(\kap,\pi)}}
{
\{
\mathcal{H}_{\gam}[\Tht\cup \{\kap\}]\vdash^{\kap+\ome a_{0}+1}_{\pi^{\prime}}\Gam^{(\kap,\pi)}, 
\lnot\del^{(\kap,\pi)}
\}_{\del\in\Del}
&
\hspace{-2mm}
\{
\mathcal{H}_{\gam}[\Tht\cup \{\kap,\rho\}]\vdash^{\kap+\ome a_{0}}_{\pi^{\prime}}\Gam^{(\kap,\pi)}, 
\Del^{(\rho,\pi)}\}_{\mu<\rho\in Mh_{2}(\bet,\nu)\cap \kap}
}
\]
}
\textbf{Case 3}. The last inference is a $(cut)$ of a cut formula $C$:
Then $\rk(C)\in\calh_{\gam}[\Tht]\cap\pi$ and $C\in\Del_{0}(\pi)$.
If $\Tht\subset\calh_{\gam}(\kap)$, then $\rk(C)<\kap$.
\\
\textbf{Case 4}. The last inference is either a $({\rm rfl}_{\Pi_{3}}(\sig,\nu))$ or a 
$({\rm rfl}_{\Pi_{2}}(\sig,\del,\nu))$ with $\sig\in\calh_{\gam}[\Tht]\cap\pi$:
IH yields the lemma.
If $\Tht\subset\calh_{\gam}(\kap)$, then $\sig<\kap$.
\eprf
\\

We see from the above proof,
if there is a base $\eta$ for inferences $({\rm rfl}_{\Pi_{3}}(\mu_{3},\sig,\nu))$ and simultaneously for 
$({\rm rfl}_{\Pi_{2}}(\mu_{2},\sig,\del,\nu))$ in $\calh_{\gam}[\Tht]\vdash^{a}_{\pi}\Gam$
(in the sense that
$\eta=\mu_{3}=\mu_{2}$),
then the same $\eta$ is a base for
inferences $({\rm rfl}_{\Pi_{3}}(\mu_{3},\sig,\nu))$ and simultaneously for 
$({\rm rfl}_{\Pi_{2}}(\mu_{2},\sig,\del,\nu))$ in
$\calh_{\gam}[\Tht\cup \{\kap\}]\vdash^{\kap+\ome a}_{\pi^{\prime}}\Gam^{(\kap,\pi)}$.

\blem\label{lem:lowerPi24}
Let $\Gam\subset\Sig_{1}(\lam)$ and
$\calh_{\gam}[\Tht]\vdash^{a}_{b}\Gam$ with $a<\Lam=\mK$,
$\calh_{\gam}[\Tht]\ni\lam\leq b<\mK$ and
$\lam$ regular, and assume
$\fal\kap\in[\lam,b)(\Tht\subset\calh_{\gam}(\psi_{\kap}(\gam)))$.

Let $\hat{a}=\gam+\tht_{b}(a)$ and 
$\del=\psi_{\lam}^{(\bet,\nu)}(\hat{a})$ when $\lam\in Mh_{2}^{\gam}(\alp)$,
$m_{3}(\lam)=0$ and $(\bet,\nu)<\alp$
with
$\{\bet,\nu\}\subset\calh_{\gam}[\Tht]$.
Then
$\calh_{\hat{a}+1}[\Tht]\vdash^{\del}_{\del}\Gam$ holds.
\elem
\bprf
By main induction on $b$ with subsidiary induction on $a$ as in Lemma \ref{lem:lowerPi2}.
Let $\eta$ be a base for reflection inferences in $\calh_{\gam}[\Tht]\vdash^{a}_{b}\Gam$.
\\
\textbf{Case 1}. 
Consider the case when the last inference is a 
$({\rm rfl}_{\Pi_{3}}(\eta,\sig,\nu))$ with $b\geq\sig$.
{\small
\[
\infer[({\rm rfl}_{\Pi_{3}}(\eta,\sig,\nu))]{
\calh_{\gam}[\Tht]\vdash^{a}_{b}\Gam}
{
\{
\mathcal{H}_{\gam}[\Tht]\vdash^{a_{0}}_{b}\Gamma, \lnot\del
\}_{\del\in\Del}
&
\{
\mathcal{H}_{\gam}[\Tht\cup \{\rho\}]\vdash^{a_{0}}_{b}\Gamma, 
\Del^{(\rho,\sig)}\}_{\eta<\rho\in Mh_{3}(\nu)\cap \sig}
}
\]
}
where 
$\Del\subset\Sig_{3}(\sig)$, $SC_{\mK}(\nu)\subset\sig$, $\nu<\xi=\bar{m}_{3}(\sig)\leq m_{3}(\sig)$, 
$\alp=\bar{m}_{2}(\sig)\leq m_{2}(\sig)$, $\eta<\sig$
and $\{\eta,\sig,\xi,\alp,\nu\}\subset\calh_{\gam}[\Tht]$.
We may assume that $\sig\geq\lam$.
\\
\textbf{Case 1.1}. There exists a regular $\pi\in\calh_{\gam}[\Tht]$ such that $\sig<\pi\leq b$:
Then $\Del\subset\Del_{0}(\pi)$ and $\sig<b_{0}=\psi_{\pi}(\widehat{a_{0}})$ for
$\widehat{a_{0}}=\gam+\tht_{b}(a_{0})$.
SIH yields
$\mathcal{H}_{\widehat{a_{0}}+1}[\Tht]\vdash^{b_{0}}_{b_{0}}\Gamma, \lnot\del$ for each $\del\in\Del$,
and
$\mathcal{H}_{\widehat{a_{0}}+1}[\Tht\cup \{\rho\}]\vdash^{b_{0}}_{b_{0}}\Gamma, 
\Del^{(\rho,\sig)}$ for each 
$\eta<\rho\in Mh_{3}(\nu)\cap \sig$.
A $({\rm rfl}_{\Pi_{3}}(\eta,\sig,\nu))$ yields
$\mathcal{H}_{\widehat{a_{0}}+1}[\Tht]\vdash^{b_{0}+1}_{b_{0}}\Gamma$, where $b_{0}<b$.
Let
$\del_{0}=\psi_{\lam}(\widehat{a_{1}})$ with 
$\widehat{a_{1}}=\widehat{a_{0}}+\tht_{b_{0}}(b_{0}+1)=\gam+\tht_{b}(a_{0})+\tht_{b_{0}}(b_{0}+1)<
\gam+\tht_{b}(a)=\hat{a}$.
We obtain
$\mathcal{H}_{\widehat{a_{1}}+1}[\Tht]\vdash^{\del_{0}}_{\del_{0}}\Gamma$
by MIH, and the lemma follows.
\\
\textbf{Case 1.2}. Otherwise: 
By Cut-elimination we obtain
$\mathcal{H}_{\gam}[\Tht]\vdash^{\tht_{b}(a_{0})}_{\sig}\Gamma, \lnot\del$ for each $\del\in\Del$, and
$\mathcal{H}_{\gam}[\Tht\cup \{\rho\}]\vdash^{\tht_{b}(a_{0})}_{\sig}\Gamma, 
\Del^{(\rho,\sig)}$ for each $\eta<\rho\in Mh_{3}(\nu)\cap \sig$.
A $({\rm rfl}_{\Pi_{3}}(\eta,\sig,\nu))$ yields
$\mathcal{H}_{\gam}[\Tht]\vdash^{a_{1}}_{\sig}\Gamma$ for $a_{1}=\tht_{b}(a_{0})+1$.
Let $\bet=\alp+\Lam^{\xi}(1+a_{1})$ for 
$\alp=\bar{m}_{2}(\sig)\leq m_{2}(\sig)$ and $\xi=\bar{m}_{3}(\pi)\leq m_{3}(\sig)$, and 
$\kap=\psi_{\sig}^{(\bet,0)}(\gam)$.
We obtain $\Tht\subset\calh_{\gam}(\kap)$ by the assumption.
Hence $\{\gam,\sig, \bet\}\subset\calh_{\gam}(\kap)$, and $\eta<\kap\in Mh_{2}(\bet)\cap \sig$,
cf.\,Corollary \ref{cor:stepdownpi4}.
Moreover we have $\kap\in\calh_{\gam+1}[\Tht]$.

Lemma \ref{lem:lowerPi43b} yields
$\mathcal{H}_{\gam}[\Tht\cup \{\kap\}]\vdash^{\kap+\ome a_{1}}_{\kap}\Gamma^{(\kap,\sig)}$
and 
$\mathcal{H}_{\gam+1}[\Tht]\vdash^{\kap+\ome a_{1}}_{\kap}\Gamma^{(\kap,\sig)}$,
where $\bet\leq\bar{m}_{2}(\kap)$ 
with
$\bar{m}(\kap)\subset\calh_{\gam}[\Tht]$, and
$\Gamma^{(\kap,\sig)}=\Gam$ if $\lam<\sig$, and $\Gamma^{(\kap,\sig)}=\Gamma^{(\kap,\lam)}$
otherwise.
In each case we obtain 
$\mathcal{H}_{\gam+1}[\Tht]\vdash^{\kap+\ome a_{1}}_{\kap}\Gamma$.
MIH then yields
$\mathcal{H}_{\widehat{a_{1}}+1}[\Tht]\vdash^{\del_{1}}_{\del_{1}}\Gamma$,
where
$\del_{1}=\psi_{\lam}(\widehat{a_{1}})$ with
$\widehat{a_{1}}=\gam+\tht_{\kap}(\kap+\ome a_{1})<\gam+\tht_{b}(a)=\hat{a}$
by $\kap<\sig\leq b$ and $a_{1}<\tht_{b}(a)$.
\\
\textbf{Case 2}. 
Consider the case when the last inference is a 
$({\rm rfl}_{\Pi_{2}}(\eta,\sig,\bet,\nu))$ with $b\geq\sig$.
{\small
\[
\infer[({\rm rfl}_{\Pi_{2}}(\eta,\sig,\bet,\nu))]{
\calh_{\gam}[\Tht]\vdash^{a}_{b}\Gam}
{
\{
\mathcal{H}_{\gam}[\Tht]\vdash^{a_{0}}_{b}\Gamma, \lnot\del
\}_{\del\in\Del}
&
\{
\mathcal{H}_{\gam}[\Tht\cup \{\rho\}]\vdash^{a_{0}}_{b}\Gamma, 
\Del^{(\rho,\sig)}\}_{\eta<\rho\in Mh_{2}(\bet,\nu)\cap \sig}
}
\]
}
where 
$\Del\subset\Sig_{2}(\sig)$, 
$(\bet,\nu)<\alp=\bar{m}_{2}(\sig)\leq m_{2}(\sig)$, 
$\xi=\bar{m}_{3}(\sig)\leq m_{3}(\sig)$, $\eta<\sig$
and $\{\eta,\sig,\alp,\xi,\bet,\nu\}\subset\calh_{\gam}[\Tht]$.

We may assume that $\sig\geq\lam$.
For each $\del\in\Del$, let $\del\simeq\bigvee(\del_{i})_{i\in J}$. We may assume
$J=Tm(\sig)$.
Inversion yields $\mathcal{H}_{\gam+|i|}[\Tht\cup\sfk(i)]\vdash^{a_{0}}_{b}\Gamma, \lnot\del_{i}$,
where $\Gam\cup\{\lnot\del_{i}\}\subset\Sig_{1}(\sig)$.
Let $\widehat{{a}_{0}}=\gam+\tht_{b}(a_{0})$ and
$\rho=\psi_{\sig}^{(\bet,\nu)}(\widehat{{a}_{0}})$,
where $\Tht\subset\calh_{\gam}(\rho)$ by the assumption,
$\{\eta,\sig,\bet,\nu,\widehat{{a}_{0}}\}\subset\calh_{\gam}[\Tht]$ with $(\bet,\nu)<m_{2}(\sig)$.
Hence 
$\{\eta,\sig,\bet,\nu,\widehat{{a}_{0}}\}\subset\calh_{\gam}(\rho)$ and 
$\calh_{\gam}(\rho)\cap\sig\subset\rho$.
Therefore $<\eta<\rho\in Mh_{2}(\bet,\nu)\cap\sig\cap\calh_{\widehat{a_{0}}+1}[\Tht]$.

We see the lemma as in Lemma \ref{lem:lowerPi2} by Inversion, picking the $\rho$-th branch from
the right upper seqeunts, and then introducing several $(cut)$'s
instead of $({\rm rfl}_{\Pi_{2}}(\eta,\sig,\bet,\nu))$.
Use MIH when $\lam<\sig$.
\\
\textbf{Case 3}.
As in Lemma \ref{lem:lowerPi2} we see the case when the last inference is a $(cut)$ of a cut formula $C$ with $d=\rk(C)<b$.
\eprf

\bth\label{th:Pi4}
Assume ${\sf KP}\Pi_{4}\vdash\tht^{L_{\Ome}}$ for $\tht\in\Sig$.
Then there exists an $n<\ome$ such that
$L_{\alp}\models\tht$ for $\alp=\psi_{\Ome}(\ome_{n}(\mK+1))$ in $OT(\Pi_{4})$.
\end{theorem}
\bprf
By Embedding there exists an $m>0$ such that
$\calh_{0}[\emptyset]\vdash^{\mK+m}_{\mK+m}\tht^{L_{\Ome}}$.
By Cut-elimination,
$\calh_{0}[\emptyset]\vdash^{a}_{\mK}\tht^{L_{\Ome}}$
for $a=\ome_{m}(\mK+m)$.
By Lemma \ref{lem:lowerPi4} we obtain
$\calh_{\ome^{a}+1}[\{\kap\}]\vdash^{\bet}_{\bet}\tht^{L_{\Ome}}$,
where $\bet=\psi_{\mK}(\ome^{a})$, $\mK+a=a$, $(\tht^{L_{\Ome}})^{(\kap,\mK)}\equiv \tht^{L_{\Ome}}$ and
$\kap\in Mh_{2}(a)\cap\psi_{\mK}(\mK)$.
F.e. $\kap=\psi_{\mK}^{(0,a)}(0)\in \calh_{1}[\emptyset]$.
Hence
$\calh_{\ome^{a}+1}[\emptyset]\vdash^{\bet}_{\bet}\tht^{L_{\Ome}}$.
Lemma \ref{lem:lowerPi24} then yields
$\calh_{\gam+1}[\emptyset]\vdash^{\bet_{1}}_{\bet_{1}}\tht^{L_{\Ome}}$
for $\gam=\ome^{a}+\tht_{\bet}(\bet)$ and
$\bet_{1}=\psi_{\Ome}(\gam)<\psi_{\Ome}(\ome^{a}+\mK)<\psi_{\Ome}(\ome_{m+2}(\mK+1))=\alp$.
Therefore
$L_{\alp}\models\tht$.
\eprf

\section{First order reflection}\label{sec:pi10}

Having established an ordinal analysis for $\Pi_{4}$-reflection in section \ref{sec:pi4},
it is not hard to extend it to first-order reflection.
As expected, an exponential ordinal structure emerges in resolving higher Mahlo classes.

Let $\mK=\Lam$ be either a $\Pi^{1}_{N-2}$-indescribable cardinal or a $\Pi_{N}$-reflecting ordinal
for an integer $N\geq 3$.
Let for $k>0$, $\alp\in M_{k+2}(A)$ iff $A$ is $\Pi^{1}_{k}$-indescribable in $\alp$ or
$\alp$ is $\Pi_{k+2}$-reflecting on $A$.
Let $(\nu_{k},\nu_{k+1},\ldots,\nu_{N-1})$ be a sequence of ordinals 
$\nu_{i}<\veps_{\Lam+1}$, and $\veps_{\Lam+1}>\alp=\Lam^{\bet_{0}}a_{0}+\cdots+\Lam^{\bet_{n}}a_{n}$
with
$\bet_{0}>\cdots>\bet_{n}$ and $0<a_{0},\ldots,a_{n}<\Lam$.
Then $(\nu_{k},\nu_{k+1},\ldots,\nu_{N-1})<\alp$ iff there exists a segment $\alp_{i}=\Lam^{\bet_{0}}a_{0}+\cdots+\Lam^{\bet_{i}}a_{i}$
of $\alp$ such that $\nu_{k}<\alp_{i}$ and $(\nu_{k+1},\ldots,\nu_{N-1})<\bet_{i}$.

\bprp\label{prp:lesspiN}
$\vec{\nu}<\alp<\gam \Rarw \vec{\nu}<\gam$.
\eprp

\subsection{Mahlo classes for $\Pi_{N}$-reflection}
As in subsection \ref{subsec:Mclasspi3}
$P\in M_{i}(\calx)$ designates that $P$ is $\Pi_{i}$-reflecting on $\calx$.
Let
\[
M_{k}(\alp):=\bigcap\{M_{k}(M_{k}(\bar{\nu})) : \bar{\nu}=(\nu_{k},\nu_{k+1},\ldots,\nu_{N-1})<\alp\}
\]
where
\[
M_{k}((\nu_{k},\nu_{k+1},\ldots,\nu_{N-1})):=\bigcap_{i\geq k}M_{i}(\nu_{i}).
\]
By Proposition \ref{prp:lesspiN} we obtain
$\alp_{0}>\alp \Rarw M_{k}(\alp_{0})\subset M_{k}(\alp)$.
Hence for $(\max\{\bar{\nu},\bar{\mu}\})_{i}=\max\{\nu_{i},\mu_{i}\}$, cf.\,\textbf{Case 1} in Lemma \ref{lem:lowerPiN3}, 
\[
M_{2}(\bar{\nu})\cap M_{2}(\bar{\mu})=M_{2}(\max\{\bar{\nu},\bar{\mu}\}).
\]

Let $\bar{\nu}=(\nu_{2},\ldots,\nu_{N-1})$ and $\bar{\mu}=(\mu_{2},\ldots,\mu_{N-1})$.
Then let
\[
\bar{\nu}\prec_{k}\bar{\mu}:\Lrarw M_{2}(\bar{\nu})\prec_{k} M_{2}(\bar{\mu}).
\]

\bprp\label{prp:stepdownpiN}
Let $\bar{\mu}=(\mu_{2},\ldots,\mu_{k-1})$, $\bar{\nu}=(\nu_{k+1},\ldots,\nu_{N-1})$, and
$\bar{\xi}=(\xi_{k+1},\ldots,\xi_{N-1})$.
\benu
\item\label{prp:stepdownpiN.1}
If $(\nu_{k})*\bar{\nu}<\xi_{k}$, then
$\bar{\mu}*(\nu_{k})*\bar{\nu}\prec_{k}\bar{\mu}*(\xi_{k})*\bar{\xi}$.

\item\label{prp:stepdownpiN.2}
(Cf.\,Lemma \ref{lem:stepdownint})
If $\xi_{k+1},a>0$, then
$\bar{\mu}*(\xi_{k}\dot{+}\Lam^{\xi_{k+1}}a)*\bar{0}\prec_{k}\bar{\mu}*(\xi_{k})*\bar{\xi}$.
\eenu
\eprp
\bprf
\ref{prp:stepdownpiN}.\ref{prp:stepdownpiN.1}.
Let $P\in M_{2}(\bar{\mu}*(\xi_{k})*\bar{\xi})\subset M_{2}(\bar{\mu}*\bar{0})\cap M_{k}(\xi_{k})$.
By $(\nu_{k})*\bar{\nu}<\xi_{k}$ we obtain $P\in M_{k}(M_{k}((\nu_{k})*\bar{\nu}))$.
Since $P\in M_{2}(\bar{\mu}*\bar{0})$ is $\Pi_{k}$ on $P$,
we conclude $P\in M_{k}(M_{2}(\bar{\mu}*\bar{0}) \cap M_{k}((\nu_{k})*\bar{\nu}))=M_{k}(M_{k}(\bar{\mu}*(\nu_{k})*\bar{\nu}))$.
\\
\ref{prp:stepdownpiN}.\ref{prp:stepdownpiN.2}.
It suffices to show that
$M_{k}(\xi_{k}\dot{+}\Lam^{\xi_{k+1}}a)\prec_{k}M_{k}(\xi_{k})\cap M_{k+1}(\xi_{k+1})$, and this follows from
$M_{k}(\xi_{k})\cap M_{k+1}(\xi_{k+1}) \subset M_{k}(\xi_{k}\dot{+}\Lam^{\xi_{k+1}}a)$.
The latter is shown by induction on $a$ as in Lemma \ref{lem:stepdownint}
using the fact that $P\in M_{k}(\gam)\cap M_{k+1}(\xi_{k+1})\Rarw P\in M_{k}(M_{k}(\gam)\cap M_{k+1}(\nu))$ for $\nu<\xi_{k+1}$.
\eprf

\subsection{Ordinals for first order reflection}

\bdf\label{df:HKPpiN}
{\rm
Define simultaneously by recursion on ordinals $a<\veps_{\mK+1}$ the classes
 $\calh_{a}(X)\, (X\subset \Gam_{\mK+1})$, $Mh_{k}^{a}(\vec{\nu})\,(lh(\vec{\nu})=N-k)$,
  the ordinals $\psi_{\sig}^{\vec{\nu}}(a)$ as follows.

\benu
\item
$\calh_{a}(X)$ denotes the Skolem hull of $\{0,\Ome,\mK\}\cup X$ under the functions
$+, \vphi$, and the following.

 Let $\vec{\nu}=(\nu_{2}, \ldots,\nu_{N-1})$,
 $\{\sig,b\}\cup\vec{\nu}\subset\calh_{a}(X)$ and 
 $b<a$.
 Then $\psi_{\sig}^{\vec{\nu}}(b)\in\calh_{a}(X)$.

\item

For $2\leq k<N$, 
$\pi\in Mh_{k}^{a}(\alp)$ iff $\{a,\alp\}\subset\calh_{a}(\pi)$ and
\[
\fal\vec{\nu}=(\nu_{k},\ldots,\nu_{N-1})\subset\calh_{a}(\pi)
\left[
\vec{\nu}<\alp \to
\pi\in M_{k}\left(
Mh_{k}^{a}(\vec{\nu})
\right)
\right]
\]
where
\[
Mh_{k}^{a}(\vec{\nu})=\bigcap_{i\geq k}Mh_{i}^{a}(\nu_{i})
.\]
Note that $Mh_{k}^{a}(\alp)$ is a $\Pi_{k+1}$-class.

\item
$\psi_{\sig}(a)=\min\left(\{\sig\}\cup
\{\kap<\sig : \{a,\sig\}\subset\calh_{a}(\kap) \spand
\calh_{a}(\kap)\cap\sig\subset\kap\}\right)$.

$m_{i}(\psi_{\sig}(a))=0$ for $i<N$.

\item
Let $\sig\in Mh_{2}^{a}(\vec{\xi})$ for $\vec{\xi}=(\xi_{2},\ldots,\xi_{N-1})$ with $\xi_{k+1}>0$,
and $0<c<\Lam=\mK$ with $c\in\calh_{a}(\sig)$.
Let $\vec{\nu}=(\xi_{2},\ldots,\xi_{k-1},\xi_{k}\dot{+}\Lam^{\xi_{k+1}}c,0,\ldots,0)$.
Then
\[
\psi_{\sig}^{\vec{\nu}}(a)=\min\left(\{\sig\}\cup
\{\kap\in Mh_{2}^{a}(\vec{\nu})\cap\sig: \{a\}\cup\vec{\nu}\subset\calh_{a}(\kap) \spand
\calh_{a}(\kap)\cap\sig\subset\kap\}
\right).
\]
$m_{i}(\psi_{\sig}^{\vec{\nu}}(a))=\nu_{i}$ for $i<N$,
cf.\,Proposition \ref{prp:stepdownpiN}.\ref{prp:stepdownpiN.2}.

\item
Let $\sig\in Mh_{2}^{a}(\vec{\mu}*\vec{\xi})$ with 
$\vec{\mu}=(\mu_{2},\ldots,\mu_{k-1})$ and $\vec{\xi}=(\xi_{k},\ldots,\xi_{N-1})$, and
$\vec{\nu}=(\nu_{k},\ldots,\nu_{N-1})<\xi_{k}$, cf.\,Proposition \ref{prp:stepdownpiN}.\ref{prp:stepdownpiN.1}.
\[
\psi_{\sig}^{\vec{\mu}*\vec{\nu}}(a)=\min\left(\{\sig\}\cup
\{\kap\in Mh_{2}^{a}(\vec{\mu}*\vec{\nu})\cap\sig: \{a\}\cup\vec{\mu}\cup\vec{\nu}\subset\calh_{a}(\kap) \spand
\calh_{a}(\kap)\cap\sig\subset\kap\}
\right).
\]

$m_{i}(\psi_{\sig}^{\vec{\mu}*\vec{\nu}}(a))=\mu_{i}$ for $i< k$, and
$m_{i}(\psi_{\sig}^{\vec{\mu}*\vec{\nu}}(a))=\nu_{i}$ for $i\geq k$.

\eenu
}
\edf

As in section \ref{sec:pi4} for $\Pi_{4}$-reflection we see the following lemmas
for $\Pi^{1}_{N-2}$-indescribable cardinal $\mK$.

\blem\label{lem:psiKN}
Let $a\in\calh_{a}(\mK)\cap\veps_{\mK+1}$. Then
$\mK\in M_{N-1}(Mh^{a}_{N-1}(\veps_{\mK+1}))$, where
$\veps_{\mK+1}$ denotes the sequence $\vec{\nu}=\vec{0}*(\nu_{N-1})$ with $\nu_{N-1}=\veps_{\mK+1}$.
For every $\xi\in\calh_{a}(\mK)\cap\veps_{\mK+1}$,
$\psi_{\mK}^{\vec{0}*(\xi)}(a)<\mK$.
\elem

\blem\label{lem:stepdownpiN}
Let $\vec{\nu}=(\xi_{2},\ldots,\xi_{k-1},\xi_{k}\dot{+}\Lam^{\xi_{k+1}}c,0,\ldots,0)$, where
$\vec{\xi}=(\xi_{2},\ldots,\xi_{N-1})$ with $\xi_{k+1}>0$,
and $0<c<\Lam$ with $c\in\calh_{a}(\sig)$.

Assume $\sig\in Mh_{2}^{a}(\vec{\xi})$.
Then
$\sig\in M_{2}(Mh^{a}_{2}(\vec{\nu}))$ and
$\psi_{\sig}^{\vec{\nu}}(a)<\sig$, cf.\,Proposition \ref{prp:stepdownpiN}.\ref{prp:stepdownpiN.2}.
\elem

\blem\label{lem:piNdown}
Let $\vec{\mu}=(\mu_{2},\ldots,\mu_{k-1})$ and
$\vec{\nu}=(\nu_{k},\ldots,\nu_{N-1})<\xi$.
Assume $\vec{\nu}\subset\calh_{a}(\sig)$ and 
$\sig\in Mh_{2}^{a}(\vec{\mu}*(\xi))$.
Then $\psi_{\sig}^{\vec{\mu}*\vec{\nu}}(a)<\sig$, cf.\,Proposition \ref{prp:stepdownpiN}.\ref{prp:stepdownpiN.1}.
\elem

\subsection{Operator controlled derivations for first order reflection}

Operator controlled derivations for ${\sf KP}\Pi_{N}$ are closed under the following inference rules.
$\bar{m}:\pi\mapsto\bar{m}(\pi)=(\bar{m}_{2}(\pi),\ldots,\bar{m}_{N-1}(\pi))$ is an
 additional data for the derivations,
where $\bar{m}_{i}(\pi)\leq m_{i}(\pi)$ for $2\leq i\leq N-1$.

\bdes
\item[$({\rm rfl}_{\Pi_{N}}(\mK))$]
$b\geq\mK$.
There exist 
an ordinal $a_{0}\in\calh_{\gam}[\Tht]\cap a$,
and a $\Sig_{N}(\mK)$-sentence $A$ enjoying the following conditions:

\[
\infer[({\rm rfl}_{\Pi_{N}}(\mK))]{
\calh_{\gam}[\Tht]\vdash^{a}_{b}\Gam}
{
\mathcal{H}_{\gam}[\Tht]\vdash^{a_{0}}_{b}\Gamma, \lnot A
&
\{
\mathcal{H}_{\gam}[\Tht\cup \{\rho\}]\vdash^{a_{0}}_{b}\Gamma, 
A^{(\rho,\mK)}: \rho<\mK\}
}
\]

\item[$({\rm rfl}_{\Pi_{k}}(\eta,\pi,\vec{\nu}))$]
for each $2\leq k\leq N-1$, cf.\,Proposition \ref{prp:stepdownpiN}.\ref{prp:stepdownpiN.1}.

There exist ordinals 
$\eta<\pi\leq b<\mK$,
$\vec{\nu}=(\nu_{k},\ldots,\nu_{N-1})< \bar{m}_{k}(\pi)\leq m_{k}(\pi)$, 
$a_{0}<a$,
and a finite set $\Del$ of $\Sig_{k}(\pi)$-sentences enjoying the following conditions:

\benu

\item
$\{\eta,\pi\}\cup\vec{\nu}\cup\bar{m}(\pi)\subset\calh_{\gam}[\Tht]$.

 \item
For each $\del\in\Del$,
$
\mathcal{H}_{\gam}[\Tht]\vdash^{a_{0}}_{b}\Gamma, \lnot\del
$.

\item
For any $\eta<\rho\in Mh_{2}(\bar{m}_{<k}(\pi)*\vec{\nu})\cap \pi$,
$\mathcal{H}_{\gam}[\Tht\cup \{\rho\}]\vdash^{a_{0}}_{b}\Gamma, 
\Del^{(\rho,\pi)}$, where
$\bar{m}_{<k}(\pi)=(\bar{m}_{2}(\pi),\ldots,\bar{m}_{k-1}(\pi))$ and
$\rho\in Mh_{k}(\vec{\nu})$ iff $\nu_{i}\leq m_{i}(\rho)$ for every $k\leq i\leq N-1$.

\eenu

\[
\hspace{-10mm}
\infer[({\rm rfl}_{\Pi_{k}}(\eta,\pi,\vec{\nu}))]{
\calh_{\gam}[\Tht]\vdash^{a}_{b}\Gam}
{
\{
\mathcal{H}_{\gam}[\Tht]\vdash^{a_{0}}_{b}\Gamma, \lnot\del
\}_{\del\in\Del}
&
\{
\mathcal{H}_{\gam}[\Tht\cup \{\rho\}]\vdash^{a_{0}}_{b}\Gamma, 
\Del^{(\rho,\pi)}\}_{\eta<\rho\in Mh_{2}(\bar{m}_{<k}(\pi)*\vec{\nu})\cap \pi}
}
\]

\edes

\blem\label{lem:lowerPiN}
Assume $\Gam\subset\Sig_{N-1}(\mK)$, $\Tht\subset\calh_{\gam}(\psi_{\mK}(\gam))$,
and
$\calh_{\gam}[\Tht]\vdash^{a}_{\mK}\Gam$.
Then 
$\calh_{\hat{a}+1}[\Tht\cup \{\kap\}]\vdash^{\bet}_{\bet}\Gam^{(\kap,\mK)}$ holds for any 
$\eta=\psi_{\mK}(\gam+\mK)<\kap\in Mh_{N-1}(a)\cap \psi_{\mK}(\gam+\mK\cdot\ome)$,
where $\hat{a}=\gam+\ome^{\mK+a}$ and $\bet=\psi_{\mK}(\hat{a})$.
\elem

\blem\label{lem:lowerPiN3}
Assume $\bar{m}(\pi)\subset\calh_{\gam}[\Tht]$, and
there exists a $2\leq k<N-1$ such that
$\bar{m}_{k+1}(\pi)>0$, and let
$k=\max\{k: \bar{m}_{k+1}(\pi)>0\}$ and
$\alp=\bar{m}_{k}(\pi)$, $\xi=\bar{m}_{k+1}(\pi)$.
Moreover assume $\calh_{\gam}[\Tht]\vdash^{a}_{\pi}\Gam$ for $a,\pi<\mK$
and $\Gam\subset\Sig_{k}(\pi)$.

Then for any $\eta<\kap\in Mh_{2}(\bar{m}_{<k}(\pi))\cap Mh_{k}(\alp\dot{+}\Lam^{\xi}(1+a))\cap \pi$,
$\calh_{\gam}[\Tht\cup \{\kap\}]\vdash^{\kap+\ome a}_{\pi}\Gam^{(\kap,\pi)}$ holds, where
$\eta$ is a base,
$\alp\dot{+}\Lam^{\xi}(1+a)\leq \bar{m}_{k}(\kap)\in\calh_{\gam}[\Tht]$ and $\bar{m}_{<k}(\kap)=\bar{m}_{<k}(\pi)$.
Moreover when $\Tht\subset\calh_{\gam}(\kap)$,
$\calh_{\gam}[\Tht\cup \{\kap\}]\vdash^{\kap+\ome a}_{\kap}\Gam^{(\kap,\pi)}$ holds.
\elem
\bprf
This is seen as in Lemma \ref{lem:lowerPi43b} by induction on $a$.
Let $\pi^{\prime}=\kap$ if $\Tht\subset\calh_{\gam}(\kap)$. Otherwise $\pi^{\prime}=\pi$.
Consider the cases when the last inference is a $({\rm rfl}_{\Pi_{n}}(\eta,\pi,\vec{\nu}))$.
We have $n\leq k+1$, $\eta<\pi$,
$\{\eta,\pi\}\cup\vec{\nu}\cup\bar{m}(\pi)\subset\calh_{\gam}[\Tht]$,
$\vec{\nu}=(\nu_{n},\ldots,\nu_{N-1})<\bar{m}_{n}(\pi)\leq m_{n}(\pi)$ and
$\Del\subset\Sig_{n}(\pi)$.
{\small
\[
\infer[({\rm rfl}_{\Pi_{n}}(\eta,\pi,\vec{\nu}))]{
\calh_{\gam}[\Tht]\vdash^{a}_{\pi}\Gam}
{
\{
\mathcal{H}_{\gam}[\Tht]\vdash^{a_{0}}_{\pi}\Gamma, \lnot\del
\}_{\del\in\Del}
&
\{
\mathcal{H}_{\gam}[\Tht\cup \{\rho\}]\vdash^{a_{0}}_{\pi}\Gamma, 
\Del^{(\rho,\pi)}\}_{\eta<\rho\in Mh_{2}(\bar{m}_{<n}(\pi)*\vec{\nu})\cap \pi}
}
\]
}
\textbf{Case 1}. $n=k+1$:
Let $\alp_{0}=\alp\dot{+}\Lam^{\xi}(1+a_{0})$. Then $\vec{\mu}=(\alp_{0})*\vec{\nu}<\alp_{1}=\alp\dot{+}\Lam^{\xi}(1+a)$
by $\vec{\nu}<\xi=\bar{m}_{k+1}(\pi)$.
We obtain $\eta<\kap$, $\{\eta,\kap,\alp_{0}\}\cup\bar{m}(\kap)\cup\vec{\nu}\subset\calh_{\gam}[\Tht\cup\{\kap\}]$.
In the following derivation $\alp_{1}\leq\bar{m}_{k}(\kap)$ with $\bar{m}(\kap)\subset\calh_{\gam}[\Tht]$.
Note that
$\bar{m}_{<k}(\kap)*\vec{\mu}=\bar{m}_{<k}(\pi)*(\alp_{0})*\vec{\nu}=\max\{(\bar{m}_{<k}(\pi)*(\alp_{0})*\bar{0}), (\bar{m}_{<k}(\pi)*(\alp)*\bar{\nu})\}$.
{\footnotesize
\[
\hspace{-20mm}
\infer[({\rm rfl}_{\Pi_{k}}(\eta,\kap,\vec{\mu}))]{\mathcal{H}_{\gam}[\Tht\cup \{\kap\}]\vdash^{\kap+\ome a}_{\pi^{\prime}}\Gam^{(\kap,\pi)}}
{
\{\mathcal{H}_{\gam}[\Tht\cup \{\kap\}]\vdash^{2d}_{0}\lnot\tht^{(\kap,\pi)},\Gam^{(\kap,\pi)}\}_{\tht\in\Gam}
&
\hspace{-15mm}
\infer{
\{\mathcal{H}_{\gam}[\Tht\cup \{\kap,\sig\}]\vdash^{\kap+\ome a_{0}+p}_{\pi^{\prime}}\Gamma^{(\kap,\pi)},\Gam^{(\sig,\pi)}\}_{\eta<\sig\in Mh_{2}(\bar{m}_{<k}(\kap)*\vec{\mu})\cap \kap}
}
 {
 \{
 \mathcal{H}_{\gam}[\Tht\cup \{\sig\}]\vdash^{\sig+\ome a_{0}+1}_{\pi^{\prime}}\Gamma^{(\sig,\pi)}, \lnot \del^{(\sig,\pi)}
 \}_{\del\in\Del}
 &
 \mathcal{H}_{\gam}[\Tht\cup \{\kap,\sig\}]\vdash^{\kap+\ome a_{0}}_{\pi^{\prime}}\Gamma^{(\kap,\pi)}, 
\Del^{(\sig,\pi)}
 }
}
\]
}
\textbf{Case 2}. $n\leq k$:
If $n<k$, then $\vec{\nu}<\bar{m}_{n}(\pi)=\bar{m}_{n}(\kap)\leq m_{n}(\kap)$.
If $n=k$, then $\vec{\nu}<\alp\dot{+}\Lam^{\xi}(1+a)\leq \bar{m}_{k}(\kap)$
with the segment $\alp$ of $\alp\dot{+}\Lam^{\xi}(1+a)$.

{\footnotesize
\[
\hspace{-15mm}
\infer[({\rm rfl}_{\Pi_{n}}(\eta,\kap,\vec{\nu}))]{
\calh_{\gam}[\Tht\cup \{\kap\}]\vdash^{\kap+\ome a}_{\pi^{\prime}}\Gam^{(\kap,\pi)}}
{
\{
\mathcal{H}_{\gam}[\Tht\cup \{\kap\}]\vdash^{\kap+\ome a_{0}+1}_{\pi^{\prime}}\Gam^{(\kap,\pi)}, 
\lnot\del^{(\kap,\pi)}
\}_{\del\in\Del}
&
\hspace{-2mm}
\{
\mathcal{H}_{\gam}[\Tht\cup \{\kap,\rho\}]\vdash^{\kap+\ome a_{0}}_{\pi^{\prime}}\Gam^{(\kap,\pi)}, 
\Del^{(\rho,\pi)}\}_{\eta<\rho\in Mh_{2}(\bar{m}_{<n}(\kap)*\vec{\nu})\cap \kap}
}
\]
}
\eprf

\blem\label{lem:lowerPi2N}
Let $\Gam\subset\Sig_{1}(\lam)$ and
$\calh_{\gam}[\Tht]\vdash^{a}_{b}\Gam$ with $a<\mK$,
$\calh_{\gam}[\Tht]\ni\lam\leq b<\mK$ and
$\lam$ regular.
Assume
$\fal\kap\in[\lam,b)(\Tht\subset\calh_{\gam}(\psi_{\kap}(\gam)))$.

Let $\hat{a}=\gam+\tht_{b}(a)$ and 
$\del=\psi_{\lam}^{\vec{\nu}}(\hat{a})$ when $\lam\in Mh_{k}^{\gam}(\alp)$ and $\vec{\nu}<\alp$
with
$\vec{\nu}\subset\calh_{\gam}[\Tht]$.
Then
$\calh_{\hat{a}+1}[\Tht]\vdash^{\del}_{\del}\Gam$ holds.
\elem
\bprf
This is seen as in Lemma \ref{lem:lowerPi24} by main induction on $b$ with subsidiary induction on $a$.
Let $\eta$ be a base.
\\
\textbf{Case 1}.
Consider the case when the last inference is a 
$({\rm rfl}_{\Pi_{k+1}}(\eta,\sig,\vec{\nu}))$ with $2\leq k< N-1$ and $b\geq\sig$.
{\small
\[
\infer[({\rm rfl}_{\Pi_{k+1}}(\eta,\sig,\vec{\nu}))]{
\calh_{\gam}[\Tht]\vdash^{a}_{b}\Gam}
{
\{
\mathcal{H}_{\gam}[\Tht]\vdash^{a_{0}}_{b}\Gamma, \lnot\del
\}_{\del\in\Del}
&
\{
\mathcal{H}_{\gam}[\Tht\cup \{\rho\}]\vdash^{a_{0}}_{b}\Gamma, 
\Del^{(\rho,\sig)}\}_{\eta<\rho\in Mh_{2}(\bar{m}_{\leq k}(\sig)*\vec{\nu})\cap \sig}
}
\]
}
where 
$\Del\subset\Sig_{k+1}(\sig)$, $\vec{\nu}<\xi=\bar{m}_{k+1}(\sig)\leq m_{k+1}(\sig)$, $\eta<\sig$
and $\{\eta,\sig\}\cup\bar{m}(\sig)\cup\vec{\nu}\subset\calh_{\gam}[\Tht]$.
We may assume that $\sig\geq\lam$ and
there is no regular $\pi\in\calh_{\gam}[\Tht]$ such that $\sig<\pi\leq b$.

We obtain the lemma by Cut-elimination, Lemma \ref{lem:lowerPiN3} for
$\kap=\psi_{\sig}^{\bar{m}_{<k}(\sig)*(\bet)*\vec{0}}(\gam)$ with
$\bet=\bar{m}_{k}(\sig)\dot{+}\Lam^{\bar{m}_{k+1}(\sig)}(1+a_{1})$ and
$a_{1}=\tht_{b}(a_{0})+1$,
and MIH.
\\
\textbf{Case 2}.
Next consider the case when the last inference is a 
$({\rm rfl}_{\Pi_{2}}(\eta,\sig,\vec{\nu}))$ with $b\geq\sig$.
{\small
\[
\infer[({\rm rfl}_{\Pi_{2}}(\eta,\sig,\vec{\nu}))]{
\calh_{\gam}[\Tht]\vdash^{a}_{b}\Gam}
{
\{
\mathcal{H}_{\gam}[\Tht]\vdash^{a_{0}}_{b}\Gamma, \lnot\del
\}_{\del\in\Del}
&
\{
\mathcal{H}_{\gam}[\Tht\cup \{\rho\}]\vdash^{a_{0}}_{b}\Gamma, 
\Del^{(\rho,\sig)}\}_{\eta<\rho\in Mh_{2}(\vec{\nu})\cap \sig}
}
\]
}
where 
$\Del\subset\Sig_{2}(\sig)$, $\vec{\nu}<\xi=\bar{m}_{2}(\sig)\leq m_{2}(\sig)$, $\eta<\sig$
and $\{\eta,\sig\}\cup\bar{m}(\sig)\cup\vec{\nu}\subset\calh_{\gam}[\Tht]$.
We may assume that $\sig\geq\lam$.
Let $\rho=\psi_{\sig}^{\vec{\nu}}(\widehat{a_{0}})$. 
We see $\eta<\rho\in Mh_{2}(\vec{\nu})\cap\sig\cap\calh_{\widehat{a_{0}}+1}[\Tht]$ from
the assumption $\Tht\subset\calh_{\gam}(\rho)$.

We see the lemma as in Lemma \ref{lem:lowerPi2} by Inversion, picking the $\rho$-th branch from
the right upper seqeunts, and then introducing several $(cut)$'s
instead of $({\rm rfl}_{\Pi_{2}}(\sig,\vec{\nu}))$.
Use MIH when $\lam<\sig$.
\eprf
\\

$OT(\Pi_{N})$ denotes a computable notation system of ordinals with
collapsing functions $\psi_{\sig}^{\vec{\nu}}(a)$.

\bth\label{th:PiN}
Assume ${\sf KP}\Pi_{N}\vdash\tht^{L_{\Ome}}$ for $\tht\in\Sig$.
Then there exists an $n<\ome$ such that
$L_{\alp}\models\tht$ for $\alp=\psi_{\Ome}(\ome_{n}(\mK+1))$ in $OT(\Pi_{N})$.
\end{theorem}
\bprf
This is seen from Lemmas \ref{lem:lowerPiN} and \ref{lem:lowerPi2N}.
\eprf

\section{$\Pi^{1}_{1}$-reflection}\label{sec:pi11}

\bdf
{\rm 
$\sig$ is said to be \textit{$\alp$-stable} for $\alp>\sig$ if
$L_{\sig}\prec_{\Sig_{1}}L_{\alp}$.
}
\edf
It is known that
$\sig$ is $(\sig+1)$-stable iff $\sig$ is $\Pi^{1}_{0}$-reflecting, and
$\sig$ is $\sig^{+}$-stable iff $\sig$ is $\Pi^{1}_{1}$-reflecting, where
$\sig^{+}$ denotes the next admissible ordinal above $\sig$, cf.\cite{Richter-Aczel74}.

Let $S_{1}$ denote the theory obtained from ${\sf KP}\ome+(V=L)$ by adding the following axioms
for an individual constant $\mS$: $\mS$ is a limit ordinal and
\[
L_{\mS}\prec_{\Sig_{1}}L
.\]
The latter denotes a schema
\[
\exi x\, B(x,v) \land v\in L_{\mS} \to \exi x\in L_{\mS}\, B(x,v)
\]
for each $\Del_{0}$-formula $B$.
Let $L=L_{\mS^{+}}\models S_{1}$.

An exponential structure emerges in iterating (recursively) Mahlo operations
to resolve first-order reflections $M_{N}$ in terms of Mahlo classes $Mh_{k}^{a}(\alp)$ and
$Mh_{k}^{a}(\vec{\nu})$.
Viewing the vector $\vec{\nu}=(\nu_{2},\nu_{3},\ldots,\nu_{N-1})$ as a function
$\{2,3,\ldots,N-1\}\ni k\mapsto \nu_{k}$,
each $k$ in its domain designates the class of $\Pi_{k}$-formulas or
the Mahlo operation $M_{k}$, while its value $\nu_{k}$ corresponds to the height of
derivations, cf.\,\textbf{Case 1} in the proof of Lemma \ref{lem:lowerPiN3}.

On the other side, the axiom $L_{\mS}\prec_{\Sig_{1}}L_{\mS^{+}}$ says that
$\mS$ `reflects' $\Pi_{\mS^{+}}$-formulas in transfinite levels.
In place of vectors in finite lengths, we need functions
$f:\mS^{+}\to ON$.
Each $c$ in the domain of the function $f$ corresponds to
formulas of ranks$<c$ in inference rules for higher reflections.
Its support ${\rm supp}(f)=\{c<\mS^{+}: f(c)\neq 0\}$ may be assumed to be \textit{finite},
while its value $f(c)<\veps_{\mS^{+}+1}$.
A Veblen function $\tilde{\tht}_{b}(\xi)$ is used to denote ordinals
instead of the exponential function $\tilde{\tht}_{1}(\xi)=(\mS^{+})^{\xi}$.
The relation $\vec{\nu}<\alp$ in section \ref{sec:pi10} is replaced by a relation
$f<^{c}\xi$ for ordinals $c,\xi$ and finite function $f$.
$f<^{c}\xi$ holds if $f(c)<\mu$ for a segment $\mu=\cdots+\tilde{\tht}_{b}(\nu)$ of $\xi$, and
$f(c+d)<\tilde{\tht}_{-d}(\tilde{\tht}_{b}(\nu))$ for $d=\min\{d>0: c+d\in{\rm supp}(f)\}$,
and so forth, where $\tilde{\tht}_{-d}(\xi)$ denotes an inverse of the function
$\xi\mapsto\tilde{\tht}_{d}(\xi)$.

Mahlo classes $Mh_{c}^{a}(\xi)$ introduced in (\ref{eq:dfMhkh})
reflects every fact 
$\pi\in Mh_{0}^{a}(g_{c})=\bigcap\{Mh^{a}_{d}(g(d)): c>d\in {\rm supp}(g)\}$ 
on the ordinals $\pi\in Mh_{c}^{a}(\xi)$ in lower level,
down to
`smaller' Mahlo classes $Mh_{c}^{a}(f)=\bigcap\{Mh^{a}_{d}(f(d)): c\leq d\in {\rm supp}(f)\}$,
where
$f<^{c}\xi$.

This apparatus would suffice to analyze reflections in transfinite levels.
We need another for the axiom $L_{\mS}\prec_{\Sig_{1}}L_{\mS^{+}}$ of $\Pi^{1}_{1}$-reflection, i.e.,
a (formal) \textit{Mostowski collapsing}:
Assume that $B(u,v)$ with $v\in L_{\mS}$ for a $\Del_{0}$-formula $B$.
We need to find a substitute $u^{\prime}\in L_{\mS}$ for $u\in L_{\mS^{+}}$, i.e.,
$B(u^{\prime},v)$.
For simplicity let us assume that $v=\bet<\mS$ and $u=\alp<\mS^{+}$ are ordinals.
We may assume that $\alp\geq\mS$.
Let $\rho<\mS$ be an ordinal, which is bigger than every ordinal$<\mS$
occurring in the `context' of $B(\alp,\bet)$.
This means that if an ordinal $\del<\mS$ occurs in a `relevant' branch of a
derivation of $B(\alp,\bet)$, $\del<\rho$ holds.
Then we can define a Mostwosiki collapsing $\alp\mapsto\alp[\rho/\mS]$ for
ordinal terms $\alp$ such that
$\bet[\rho/\mS]=\bet$ for each relevant $\bet<\mS$, $\mS[\rho/\mS]=\rho$ and
$\alp[\rho/\mS]<(\mS^{+})[\rho/\mS]=\rho^{+}<\mS$, cf.\,Definition \ref{df:Mostwskicollaps}.
Then we see that $B(\alp[\rho/\mS],\bet)$ holds.

Although the above scheme would seem to work, how to implement the plan?
Let $E^{\mS}_{\rho}$ denote the set of ordinal terms $\alp$ such that
every subterm $\bet<\mS$ of $\alp$ is smaller than $\rho$.
It turns out that 
$\mathcal{H}_{\gamma}(E^{\mathbb{S}}_{\rho})\subset E^{\mathbb{S}}_{\rho}$
if $\mathcal{H}_{\gamma}(\rho)\cap\mathbb{S}\subset\rho$.
Let $\calh_{\gam}[\Tht]\vdash^{a}_{b}\Gamma$, and
assume that (\ref{eq:controlderKP}), $\{\gam,a,b\}\cup\sfk(\Gam)\subset\calh_{\gam}[\Tht]$
holds in Definition \ref{df:controlderreg}.
Moreover let us assume that $\Tht\subset E^{\mS}_{\rho}$ holds.
Then we obtain
$\{\gam,a,b\}\cup\sfk(\Gam)\subset\calh_{\gam}[\Tht]\subset\calh_{\gam}(E^{\mS}_{\rho})\subset E^{\mS}_{\rho}$.
This means that $\sfk(\Gam)\subset E^{\mS}_{\rho}$ holds
as long as $\Tht\subset E^{\mS}_{\rho}$ holds, i.e., as long as
we are concerned with branches for $\sfk(\iota)\subset E^{\mS}_{\rho}$
in, e.g., inferences $(\bigwedge)$: $A\simeq\bigwedge(A_{\iota})_{\iota\in J}$
\[
\infer[(\bigwedge)]{\calh_{\gam}[\Tht]\vdash^{a}_{b}\Gam,A}
{
\{
\calh_{\gam}[\Tht]\vdash^{a_{0}}_{b}\Gam,A, A_{\iota}
\}_{\iota\in J}
}
\leadsto
\infer[(\bigwedge)]{\calh_{\gam}[\Tht]\vdash^{a}_{b}\Gam,A}
{
\{
\calh_{\gam}[\Tht]\vdash^{a_{0}}_{b}\Gam,A, A_{\iota}
\}_{\iota\in J, \sfk(\iota)\subset E^{\mS}_{\rho}}
}
\]
and dually $\sfk(\iota)\subset E^{\mS}_{\rho}$ for a minor formula $A_{\iota}$
of a $(\bigvee)$ with the main formula $A\simeq\bigvee(A_{\iota})_{\iota\in J}$,
\textit{provided that} 
$\mathcal{H}_{\gamma}(\rho)\cap\mathbb{S}\subset\rho$.
The proviso means that $\gam_{1}\geq\gam$ when $\rho=\psi_{\mS}^{f}(\gam_{1})$.
Such a $\rho\in\calh_{\gam}[\Tht]$ only when $\rho\in\Tht$.
Let us try to replace the inferences for the stability of $\mS$
\[
\infer[({\rm stbl})]{(\mathcal{H}_{\gamma},\Theta
)
\vdash
\Gamma
}
{
(\mathcal{H}_{\gamma},\Theta
)\vdash
\Gamma, B(u)
&
\{
(\mathcal{H}_{\gamma},\Theta\cup\{\sigma\}
)
\vdash
\Gamma, \lnot B(u)^{[\sigma/\mathbb{S}]}
\}_{\Theta\subset E^{\mathbb{S}}_{\sigma}}
}
\]
by inferences for reflection of $\rho$ with $\Tht\subset E^{\mS}_{\rho}$:
If $B(u)^{[\rho/\mS]}$ holds, then $B(u)^{[\sig/\mS]}$ holds for some $\sig<\rho$.
{\small
\[
\hspace{-5mm}
\infer[({\rm rfl})]{(\mathcal{H}_{\gamma},\Theta\cup\{\rho\}
)
\vdash
\Gamma^{[\rho/\mS]}
}
{
(\mathcal{H}_{\gamma},\Theta\cup\{\rho\}
)\vdash
\Gamma^{[\rho/\mS]}, B(u)^{[\rho/\mS]}
&
\{
(\mathcal{H}_{\gamma},\Theta\cup\{\rho,\sigma\}
)
\vdash
\Gamma^{[\rho/\mS]}, \lnot B(u)^{[\sigma/\mathbb{S}]}
\}_{\Theta\subset E^{\mathbb{S}}_{\sigma},\sig<\rho}
}
\]
}

However we need to eliminate the inferences for reflections in transfinite levels.
In view of analysis in section \ref{sec:pi10} for first-order reflection,
$\Gamma^{[\rho/\mS]}, B(u)^{[\rho/\mS]}$ is replaced by $\Gamma^{[\sig/\mS]}, B(u)^{[\sig/\mS]}$,
and
$\Gamma^{[\rho/\mS]}, \lnot B(u)^{[\sigma/\mathbb{S}]}$ by $\Gamma^{[\kap/\mS]}, \lnot B(u)^{[\sigma/\mathbb{S}]}$ with $\sig<\kap<\rho$.
{\small
\[
\hspace{-10mm}
\infer[({\rm rfl})]{(\mathcal{H}_{\gamma},\Theta\cup\{\kap,\rho\})\vdash\Gam^{[\kap/\mS]}
}
{
\{
(\mathcal{H}_{\gamma},\Theta\cup\{\kap\})\vdash\Gam^{[\kap/\mS]},\lnot\tht^{[\kap/\mS]}
\}_{\tht\in\Gam}
&
\hspace{-20mm}
\infer[(cut)]{
\{(\mathcal{H}_{\gamma},\Theta\cup\{\kap,\rho,\sig\}
)
\vdash
\Gamma^{[\kap/\mS]},\Gamma^{[\sig/\mS]}
\}_{\sig}
}
{
(\mathcal{H}_{\gamma},\Theta\cup\{\rho,\sig\}
)\vdash
\Gamma^{[\sig/\mS]}, B(u)^{[\sig/\mS]}
&
(\mathcal{H}_{\gamma},\Theta\cup\{\kap,\rho,\sigma\}
)
\vdash
\Gamma^{[\kap/\mS]}, \lnot B(u)^{[\sigma/\mathbb{S}]}
}
}
\]
}
We are replacing formulas $\Gam^{[\rho/\mS]}$ by $\Gam^{[\sig/\mS]}$ or by
$\Gam^{[\kap/\mS]}$.
This means that $\alp[\sig/\mS]$ is substituted for each $\alp[\rho/\mS]$.
Namely a composition of uncollapsing and collapsing $\alp[\rho/\mS]\mapsto\alp\mapsto\alp[\sig/\mS]$
arises.
Hence we need $\alp\in E^{\mS}_{\sig}\subsetneq E^{\mS}_{\rho}$ for $\sig<\rho$.
However we have $\Tht\cup\{\rho\}\not\subset E^{\mS}_{\sig}$, and the schema seems to be broken.
Moreover the finite sets $\Tht\cup\{\rho\}$ becomes bigger to $\Tht\cup\{\kap,\rho\}$.
Is it remain finite in eliminating inferences of reflections in transfinite level?

Looking back at the proof of Lemma \ref{lem:lowerPi43b},
for $\Gam\subset\Sig_{2}$ and $\Del\subset\Pi_{2}$
{\small
\[
\infer[({\rm rfl}_{\Pi_{3}})]{
\calh_{\gam}[\Tht]\vdash\Gam^{(\pi,\mK)}
}
{
\{
\mathcal{H}_{\gam}[\Tht]\vdash\Gamma^{(\pi,\mK)}, \lnot\del^{(\pi,\mK)}
\}_{\del\in\Del}
&
\{
\mathcal{H}_{\gam}[\Tht\cup \{\rho\}]\vdash\Gamma^{(\pi,\mK)}, 
\Del^{(\rho,\mK)}
\}_{\rho}
}
\]
}
 is rewritten to
{\small
\[
\hspace{-10mm}
\infer[({\rm rfl}_{\Pi_{2}})]{\mathcal{H}_{\gam}[\Tht\cup \{\kap\}]\vdash\Gam^{(\kap,\mK)}}
{
\{\mathcal{H}_{\gam}[\Tht\cup \{\kap\}]\vdash\lnot\tht^{(\kap,\mK)},\Gam^{(\kap,\mK)}\}_{\tht\in\Gam}
&
\hspace{-15mm}
\infer{
\{\mathcal{H}_{\gam}[\Tht\cup \{\kap,\sig\}]\vdash
\Gamma^{(\kap,\mK)},\Gam^{(\sig,\mK)}\}_{\sig}
}
 {
 \{
 \mathcal{H}_{\gam}[\Tht\cup \{\sig\}]\vdash\Gamma^{(\sig,\mK)}, \lnot \del^{(\sig,\mK)}
 \}_{\del\in\Del}
 &
 \mathcal{H}_{\gam}[\Tht\cup \{\kap,\sig\}]\vdash\Gamma^{(\kap,\mK)}, 
\Del^{(\sig,\mK)}
 }
}
\]
}
This is done by replacing the restriction ${}^{(\pi,\mK)}$ by ${}^{(\sig,\mK)}$ or ${}^{(\kap,\mK)}$, and
ordinals $\pi,\sig,\kap$ enter derivations, but do we need to control these ordinals?
Instead of the restriction ${}^{(\pi,\mK)}$, formulas could put on \textit{caps} $\pi,\sig,\kap$
in such a way that $\sfk(A^{(\sig)})=\sfk(A)$.
This means that the cap $\sig$ does not `occur' in a capped formula $A^{(\sig)}$.
If we choose an ordinal $\gam_{0}$ big enough (depending on a given finite proof figure),
 every ordinal `occurring' in derivations (including the subscript $\gam\leq\gam_{0}$ in
the operators $\calh_{\gam}$)
 is in $\calh_{\gam_{0}}=\calh_{\gam_{0}}(\emptyset)$
for the ordinal $\gam_{0}$, while
each cap $\rho$ exceeds the \textit{threshold} $\gam_{0}$ in the sense that
$\rho\not\in\calh_{\gam_{0}}(\rho)\cap\mS\subset\rho$.
Then every ordinal `occurring' in derivations is in the domain $E^{\mS}_{\rho}$
of the Mostowski collapsing $\alp\mapsto\alp[\rho/\mS]$.
Now details follow.

\subsection{Ordinals for one stable ordinal}\label{sect:ordinalnotation}

For a while, $\mS$ denotes a weakly inaccessible cardinal.

\bdf\label{df:Lam}
{\rm
Let $\Lambda=\omega_{\mathbb{S}+1}$
or $\Lambda=\mS^{+}$.
$\varphi_{b}(\xi)$ denotes the binary Veblen function on 
$\Lambda^{+}$ with $\varphi_{0}(\xi)=\omega^{\xi}$, and
$\tilde{\varphi}_{b}(\xi):=\varphi_{b}(\Lambda\cdot \xi)$ for the epsilon number 
$\Lambda$.

Let $b,\xi<\Lambda^{+}$.
$\theta_{b}(\xi)$ [$\tilde{\theta}_{b}(\xi)$] denotes
a $b$-th iterate of $\varphi_{0}(\xi)=\omega^{\xi}$ [of $\tilde{\varphi}_{0}(\xi)=\Lambda^{\xi}$], resp.
}
\edf

\bdf\label{df:Lam2}
{\rm
Let $\xi<\varphi_{\Lambda}(0)$ be a non-zero ordinal with its normal form:
\begin{equation}\label{eq:CantornfLam}
\xi=\sum_{i\leq m}\tilde{\theta}_{b_{i}}(\xi_{i})\cdot a_{i}=_{NF}
\tilde{\theta}_{b_{m}}(\xi_{m})\cdot a_{m}+\cdots+\tilde{\theta}_{b_{0}}(\xi_{0})\cdot a_{0}
\end{equation}
where
$\tilde{\theta}_{b_{i}}(\xi_{i})>\xi_{i}$,
$\tilde{\theta}_{b_{m}}(\xi_{m})>\cdots>\tilde{\theta}_{b_{0}}(\xi_{0})$, 
$b_{i}=\omega^{c_{i}}<\Lambda$, and
$0<a_{0},\ldots,a_{m}<\Lambda$.
$SC_{\Lam}(\xi)=\bigcup_{i\leq m}(\{a_{i}\}\cup SC_{\Lam}(\xi_{i}))$.

$\tilde{\theta}_{b_{0}}(\xi_{0})$ is said to be the \textit{tail} of $\xi$, denoted 
$\tilde{\theta}_{b_{0}}(\xi_{0})=tl(\xi)$, and
$\tilde{\theta}_{b_{m}}(\xi_{m})$ the \textit{head} of $\xi$, denoted 
$\tilde{\theta}_{b_{m}}(\xi_{m})=hd(\xi)$.

\begin{enumerate}
\item\label{df:Exp2.3}
 $\zeta$ is a \textit{segment} of $\xi$
 iff there exists an $n\, (0\leq n\leq m+1)$
 such that
 $\zeta=_{NF}\sum_{i\geq n}\tilde{\theta}_{b_{i}}(\xi_{i})\cdot a_{i}=
 \tilde{\theta}_{b_{m}}(\xi_{m})\cdot a_{m}+\cdots+\tilde{\theta}_{b_{n}}(\xi_{n})\cdot a_{n}$
 for $\xi$ in (\ref{eq:CantornfLam}).

\item\label{df:thtminus}
Let $\zeta=_{NF}\tilde{\theta}_{b}(\xi)$ with $\tilde{\theta}_{b}(\xi)>\xi$ and $b=\omega^{b_{0}}$,
and $c$ be ordinals.
An ordinal $\tilde{\theta}_{-c}(\zeta)$ is defined recursively as follows.
If $b\geq c$, then $\tilde{\theta}_{-c}(\zeta)=\tilde{\theta}_{b-c}(\xi)$.
Let $c>b$.
If $\xi>0$, then
$\tilde{\theta}_{-c}(\zeta)=\tilde{\theta}_{-(c-b)}(\tilde{\theta}_{b_{m}}(\xi_{m}))$ for the head term 
$hd(\xi)=\tilde{\theta}_{b_{m}}(\xi_{m})$ of 
$\xi$ in (\ref{eq:CantornfLam}).
If $\xi=0$, then let $\tilde{\theta}_{-c}(\zeta)=0$.

\end{enumerate}
}
\edf

\bdf\label{df:Lam3}
{\rm

  \begin{enumerate}
  \item
A function $f:\Lam\to \varphi_{\Lam}(0)$ with a \textit{finite} support
${\rm supp}(f)=\{c<\Lam: f(c)\neq 0\}\subset \Lam$ is said to be a \textit{finite function}
if
$\forall i>0(a_{i}=1)$ and $a_{0}=1$ when $b_{0}>1$
in
$f(c)=_{NF}\tilde{\theta}_{b_{m}}(\xi_{m})\cdot a_{m}+\cdots+\tilde{\theta}_{b_{0}}(\xi_{0})\cdot a_{0}$
for any $c\in{\rm supp}(f)$.

It is identified with the finite function $f\!\upharpoonright\! {\rm supp}(f)$.
When $c\not\in {\rm supp}(f)$, let $f(c):=0$.
$SC_{\Lam}(f):=\bigcup\{\{c\}\cup SC_{\Lam}(f(c))\}: c\in {\rm supp}(f)\}$.
$f,g,h,\ldots$ range over finite functions.

For an ordinal $c$, $f_{c}$ and $f^{c}$ are restrictions of $f$ to the domains
${\rm supp}(f_{c})=\{d\in{\rm supp}(f): d< c\}$ and ${\rm supp}(f^{c})=\{d\in{\rm supp}(f): d\geq c\}$.
$g_{c}*f^{c}$ denotes the concatenated function such that
${\rm supp}(g_{c}*f^{c})={\rm supp}(g_{c})\cup {\rm supp}(f^{c})$, 
$(g_{c}*f^{c})(a)=g(a)$ for $a<c$, and
$(g_{c}*f^{c})(a)=f(a)$ for $a\geq c$.

\item\label{df:Exp2.5}
Let $f$ be a finite function and $c,\xi$ ordinals.
A relation $f<^{c}\xi$ is defined by induction on the
cardinality of the finite set $\{d\in {\rm supp}(f): d>c\}$ as follows.
If $f^{c}=\emptyset$, then $f<^{c}\xi$ holds.
For $f^{c}\neq\emptyset$,
$f<^{c}\xi$ iff
there exists a segment $\mu$ of $\xi$ such that
$f(c)< \mu$
and 
$f<^{c+d} \tilde{\theta}_{-d}(tl(\mu))$ 
for $d=\min\{c+d\in {\rm supp}(f): d>0\}$.

\end{enumerate}

}
\edf

\bprp\label{prp:idless}
$f<^{c}\xi\leq\zeta \Rightarrow f<^{c}\zeta$.
\eprp

\subsection{Mahlo classes for $\Pi^{1}_{1}$-reflection}
In Lemma \ref{lem:stepdownint}  and Proposition \ref{prp:stepdownpiN}.\ref{prp:stepdownpiN.2},
it is crucial the fact that
$P\in M_{k}(\gam) \Rarw P\in M_{k}(M_{k}(\gam)\cap M_{k+1}(\nu))$ if
$P\in M_{k+1}(\xi_{k+1})$ and  $\nu<\xi_{k+1}$.
This means that if $P$ is in a higher Mahlo class, then $P$ reflects
a fact on $P$ in lower Mahlo classes.

$P\in M_{c}(\xi)$ is defined by main induction on $c$ with subsidiary induction on $P$.

\beqn\label{eq:PMc}
P\in M_{c}(\xi) :\Lrarw
 \forall f<^{c}\xi 
 \forall g \left[
P\in M_{0}(g_{c})
 \Rightarrow P\in M_{2}(M_{0}(g_{c}*f^{c}))
 \right]
\eeqn
where $f,g$ range over finite functions and
\[
M_{c}(f)  :=  \bigcap\{M_{d}(f(d)): d\in {\rm supp}(f^{c})\}
=
\bigcap\{M_{d}(f(d)): c\leq d\in {\rm supp}(f)\}.
\]
From Proposition \ref{prp:idless} we see
$\xi<\zeta\Rarw M_{c}(\xi)\supset M_{c}(\zeta)$.

For classes $\calx$ let
\[
P\in M_{c}(\calx) :\Lrarw \fal g\left[P\in M_{0}(g_{c}) \Rarw P\in M_{2}(M_{0}(g_{c})\cap\calx)\right].
\]
Then by $M_{0}(g_{c}*f^{c})=M_{0}(g_{c})\cap M_{c}(f^{c})$,
$P\in M_{c}(\xi)\Lrarw \fal f<^{c}\xi\left[P\in M_{c}(M_{c}(f^{c}))\right]$, i.e.,
$M_{c}(\xi)=\bigcap_{f<^{c}\xi}M_{c}(M_{c}(f^{c}))$.

\bprp\label{prp:MMh}
Suppose $P\in M_{c}(\xi)$.
\begin{enumerate}
\item\label{prp:MMh.1}
Let $f<^{c}\xi$.
Then
$P\in M_{c}(M_{c}(f^{c}))$.

\item\label{prp:MMh.2}
Let $P\in M_{d}(\calx)$ for $d>c$.
Then $P\in M_{c}(M_{c}(\xi)\cap \calx)$.
\end{enumerate}
\eprp
\bprf
\ref{prp:MMh}.\ref{prp:MMh.1}.
Let $g$ be a function such that $P\in M_{0}(g_{c})$.
By the definition (\ref{eq:PMc}) of $P\in M_{c}(\xi)$ we obtain
$P\in M_{2}\left( M_{0}(g_{c}) \cap M_{c}(f^{c}) \right)$.
\\
\ref{prp:MMh}.\ref{prp:MMh.2}.
Let $P\in M_{d}(\calx)$ for $d>c$.
Let $g$ be a function such that $P\in M_{0}(g_{c})$.
We obtain by $d>c$ with the function $g_{c}*h$,
$P\in M_{2}\left( M_{0}(g_{c}) \cap M_{c}(\xi)\cap \calx \right)$, where
${\rm supp}(h)=\{c\}$ and $h(c)=\xi$.
\eprf

\blem\label{lem:stepdownpi11int}
Assume $P\in M_{d}(\xi)\cap M_{c}(\xi_{0})$, $\xi_{0}\neq 0$,
and $d<c$.
Moreover let 
$\xi_{1}\leq\tilde{\theta}_{c-d}(\xi_{0})$.
Then
$P\in M_{d}(\xi\dot{+}\xi_{1})\cap M_{d}(M_{d}(\xi\dot{+}\xi_{1}))$.
\elem
\bprf
This is seen as in Lemma \ref{lem:stepdown}.
 
 We obtain $P\in M_{c}(\xi_{0})\subset M_{c}(M_{c}(\emptyset))$
 by Proposition \ref{prp:MMh}.\ref{prp:MMh.1}.
 Let $P\in M_{d}(\xi\dot{+}\xi_{1})\cap M_{0}(g_{d})$ for a function $g$.
 We show $P\in M_{2}\left(M_{0}(g_{d})\cap M_{d}(\xi\dot{+}\xi_{1}) \right)$.
 Let $h=g_{d}\cup\{(d,\xi\dot{+}\xi_{1})\}$.
 Then $P\in M_{0}(h_{c})$ by $d<c$.
 $P\in M_{c}(M_{c}(\emptyset))$ yields
$P\in M_{2}\left(M_{0}(h_{c})\cap M_{c}(\emptyset) \right)$, and hence
$P\in M_{2}\left(M_{0}(g_{d})\cap M_{d}(\xi\dot{+}\xi_{1}) \right)$.
Therefore $P\in M_{d}(M_{d}(\xi\dot{+}\xi_{1}))$.

Let $f$ be a finite function such that 
$f<^{d}\xi+\xi_{1}$.
We show $P\in M_{d}(M_{d}(f^{d}))$ by main induction on 
the cardinality of the finite set $\{e\in {\rm supp}(f): e>d\}$
 with subsidiary induction on $\xi_{1}$.

First let $f<^{d}\mu$ for a segment $\mu$ of $\xi$. We obtain
$P\in M_{d}(\mu)$ and $P\in M_{d}(M_{d}(f^{d}))$.

In what follows let $f(d)=\xi\dot{+}\zeta$ with $\zeta<\xi_{1}$.
By SIH we obtain $P\in M_{d}(f(d))\cap M_{d}(M_{d}(f(d)))$.
If $\{e\in {\rm supp}(f): e>d\}=\emptyset$, then $M_{d}(f^{d})=M_{d}(f(d))$, and we are done.
Otherwise let $e=\min\{e\in {\rm supp}(f): e>d\}$.

By SIH we can assume $f<^{e}\tilde{\theta}_{-(e-d)}(tl(\xi_{1}))$.
By $\xi_{1}\leq\tilde{\theta}_{c-d}(\xi_{0})$, 
we obtain $f<^{e}\tilde{\theta}_{-(e-d)}(\tilde{\theta}_{c-d}(\xi_{0}))=\tilde{\theta}_{-e}(\tilde{\theta}_{c}(\xi_{0}))$.
We claim that $P\in M_{c_{0}}(M_{c_{0}}(f^{c_{0}}))$ for $c_{0}=\min\{c,e\}$.
If $c=e$, then the claim follows from the assumption $P\in M_{c}(\xi_{0})$ and $f<^{e}\xi_{0}$.
Let $e=c+e_{0}>c$. Then $\tilde{\theta}_{-e}(\tilde{\theta}_{c}(\xi_{0}))=\tilde{\theta}_{-e_{0}}(hd(\xi_{0}))$, and
$f<^{c}\xi_{0}$ with $f(c)=0$ yields the claim.
Let $c=e+c_{1}>e$. Then $\tilde{\theta}_{-e}(\tilde{\theta}_{c}(\xi_{0}))=\tilde{\theta}_{c_{1}}(\xi_{0})$.
MIH yields the claim.

On the other hand we have
$M_{d}(f^{d})=M_{d}(f(d))\cap M_{c_{0}}(f^{c_{0}})$.
$P\in M_{d}(f(d))\cap M_{c_{0}}(M_{c_{0}}(f^{c_{0}}))$ with $d<c_{0}$ yields by 
Proposition \ref{prp:MMh}.\ref{prp:MMh.2},
$P\in M_{d}\left(M_{d}(f(d))\cap M_{c_{0}}(f^{c_{0}})\right)$, i.e.,
$P\in M_{d}(M_{d}(f^{d}))$.
\eprf

For finite functions $f$ and $g$,
\[
 M_{0}(g)\prec M_{0}(f)
 :\Leftrightarrow 
\forall P\in M_{0}(f)
\left(
P\in M_{2}(M_{0}(g))
\right)
.
\]

\bcor\label{cor:stepdownpi11int}
Let $f,g$ be finite functions and $c\in{\rm supp}(f)$.
Assume  that 
there exists an ordinal
$d<c$ 
such that
$(d,c)\cap {\rm supp}(f)=(d,c)\cap {\rm supp}(g)=\emptyset$, 
$g_{d}=f_{d}$, 
$g(d)<f(d)\dot{+}\tilde{\theta}_{c-d}(f(c))\cdot\omega$,
and
$g<^{c}f(c)$.
Then
$M_{0}(g)\prec M_{0}(f)$ holds.
\ecor
\bprf
By Lemma \ref{lem:stepdownpi11int}.
\eprf

\bdf\label{df:irreducible}
{\rm
An \textit{irreducibility} of finite functions $f$
 is defined by induction on the cardinality
$n$ of the finite set ${\rm supp}(f)$.
If $n\leq 1$, $f$ is defined to be irreducible.
Let $n\geq 2$ and $c<c+d$ be the largest two elements in ${\rm supp}(f)$, and let $g$ be 
a finite function
such that ${\rm supp}(g)={\rm supp}(f_{c})\cup\{c\}$, $g_{c}=f_{c}$ and
$g(c)=f(c)+\tilde{\theta}_{d}(f(c+d))$.

Then $f$ is irreducible iff 
$tl(f(c))>\tilde{\theta}_{d}(f(c+d))$ and
$g$ is irreducible.

}
\edf

\bdf\label{df:lx}
 {\rm 
 Let  $f,g$
 be irreducible finite functions, and $b$ an ordinal.
Let us define a relation $f<^{b}_{lx}g$
by induction on the cardinality $\#\{e\in{\rm supp}(f)\cup{\rm supp}(g): e\geq b\}$ as follows.
$f<^{b}_{lx}g$ holds iff $f^{b}\neq g^{b}$ and
for the ordinal $c=\min\{c\geq b : f(c)\neq g(c)\}$,
one of the following conditions is met:

\begin{enumerate}

\item\label{df:lx.23}
$f(c)<g(c)$ and let $\mu$ be the shortest part of $g(c)$ such that $f(c)<\mu$.
Then for any $c<c+d\in{\rm supp}(f)$,  
if $tl(\mu)\leq\tilde{\theta}_{d}(f(c+d))$, then 
$f<_{lx}^{c+d}g$ holds.

\item\label{df:lx.24}
$f(c)>g(c)$ and let $\nu$ be the shortest part of $f(c)$ such that $\nu>g(c)$.
Then there exist a $c<c+d\in {\rm supp}(g)$ such that
$f<_{lx}^{c+d}g$ and
$tl(\nu)\leq \tilde{\theta}_{d}(g(c+d))$.

\end{enumerate}

}
\edf

\bprp\label{lem:psinucomparison}
If $f<^{0}_{lx}g$, then
$M_{0}(f)\prec M_{0}(g)$.

\eprp
\bprf
This is seen from Corollary \ref{cor:stepdownpi11int}.
\eprf

\subsection{Skolem hulls and collapsing functions}\label{subsec:Skolemh}

\bdf\label{df:Cpsiregularsm}
{\rm
Let
$\mathbb{K}=\omega_{\mathbb{S}+1}$,
$a<\varepsilon_{\mathbb{K}+1}$ and $X\subset\Gamma_{\mathbb{K}+1}$.

\begin{enumerate}
\item\label{df:Cpsiregularsm.1}
$\mathcal{H}_{a}(X)$ denotes the Skolem hull of $\{0,\Omega,\mathbb{S},\mK\}\cup X$
under the functions
$+, \vphi$, 
$\bet\mapsto\psi_{\Ome}(\bet)\,(\bet<a)$,
$\mS>\alp\mapsto\alp^{+}$ and
$(\pi,b,f)\mapsto \psi_{\pi}^{f}(b)$, where $b<a$ and $f$ is a finite function such that
$f\in \mathcal{H}_{a}(X):\Lrarw SC_{\mK}(f)\subset\mathcal{H}_{a}(X)$.

\item\label{df:Cpsiregularsm.2}
Let $c<\mathbb{K}$, $a<\varepsilon_{\mathbb{K}+1}$ and $\xi<\varphi_{\mathbb{K}}(0)$.
$\pi\in Mh^{a}_{c}(\xi)$ iff 
$\{a,c,\xi\}\subset\mathcal{H}_{a}(\pi)$ and
{\small
\begin{equation}\label{eq:dfMhkh}
 \forall f<^{c}\xi 
 \forall g \left(
 SC_{\mK}(f)\cup SC_{\mK}(g) \subset\mathcal{H}_{a}(\pi) 
 \,\&\, 
\pi\in Mh^{a}_{0}(g_{c})
 \Rightarrow \pi\in M_{2}(Mh^{a}_{0}(g_{c}*f^{c}))
 \right)
\end{equation}
}
where 
\[
Mh^{a}_{c}(f)  :=  \bigcap\{Mh^{a}_{d}(f(d)): d\in {\rm supp}(f^{c})\}
=
\bigcap\{Mh^{a}_{d}(f(d)): c\leq d\in {\rm supp}(f)\}.
\]

\item\label{df:Cpsiregularsm.3}
{\small
\begin{equation}\label{eq:Psivec}
\psi_{\pi}^{f}(a)
 :=  \min(\{\pi\}\cup\{\kappa\in Mh^{a}_{0}(f)\cap\pi:   \mathcal{H}_{a}(\kappa)\cap\pi\subset\kappa ,
   \{\pi,a\}\cup SC_{\mK}(f)\subset\mathcal{H}_{a}(\kappa)
\})
\end{equation}
}

\end{enumerate}
}

\edf

\textit{Shrewd cardinals} are introduced by \cite{RathjenAFML2}.
A cardinal $\kappa$ is \textit{shrewd} iff for any $\eta>0$, $P\subset V_{\kappa}$,
and formula $\varphi(x,y)$, if
$V_{\kappa+\eta}\models\varphi[P,\kappa]$, then there are $0<\kappa_{0},\eta_{0}<\kappa$ such that
$V_{\kappa_{0}+\eta_{0}}\models\varphi[P\cap V_{\kappa_{0}},\kappa_{0}]$.
$\tilde{T}$ denotes the extension of ${\sf ZFC}$ by the axiom stating that
$\mathbb{S}$ is a shrewd cardinal.

\blem\label{lem:welldefinedness.1pi11}
$\tilde{T}$ proves that
$\mathbb{S}\in Mh^{a}_{c}(\xi)\cap M_{2}(Mh^{a}_{c}(\xi))$
for every $a<\varepsilon_{\mathbb{K}+1}$, $c<\mathbb{K}$,
$\xi<\varphi_{\mathbb{K}}(0)$ such that
$\{a,c,\xi\}\subset\mathcal{H}_{a}(\mathbb{S})$.
\elem 
\bprf
We show the lemma by induction on $\xi<\varphi_{\mathbb{K}}(0)$.

Let $\{a,c,\xi\}\cup SC_{\mK}(f)\subset\mathcal{H}_{a}(\mathbb{S})$ and $f<^{c}\xi$.
We show $\mathbb{S}\in Mh^{a}_{c}(f^{c})$, 
and $\mS\in M_{2}\left( Mh_{0}^{a}(g_{c}) \cap Mh^{a}_{c}(f^{c})\right)$
assuming $\mS\in Mh_{0}^{a}(g_{c})$ and $SC_{\mK}(g_{c})\subset\mathcal{H}_{a}(\mS)$.

For each $d\in{\rm supp}(f^{c})$ we obtain $f(d)<\xi$ by 
$\tilde{\theta}_{-e}(\zeta)\leq\zeta$.
IH yields $\mathbb{S}\in Mh^{a}_{c}(f^{c})$.

We have to show $\mathbb{S}\in M_{2}(A\cap B)$ for $A=Mh_{0}^{a}(g_{c})\cap\mathbb{S}$ and 
$B=Mh^{a}_{c}(f^{c})\cap\mathbb{S}$.
Let $C$ be a club subset of $\mathbb{S}$.

We have $\mathbb{S}\in Mh_{0}^{a}(g_{c})\cap Mh^{a}_{c}(f^{c})$,
and
$\{a,c\}\cup SC_{\mK}(g_{c},f^{c})\subset\mathcal{H}_{a}(\mathbb{S})$.
Pick a $b<\mathbb{S}$ so that $\{a,c\}\cup SC_{\mK}(g_{c},f^{c})\subset\mathcal{H}_{a}(b)$, and a bijection 
$F:\mathbb{S}\to \mathcal{H}_{a}(\mathbb{S})$.
Each $\alpha
\in\mathcal{H}_{a}(\mathbb{S})\cap\Gamma_{\mathbb{K}+1}$ is identified with its code, denoted
by $F^{-1}(\alpha)$.
Let $P$ be the class
$P=\{(\pi,d,\alpha)\in\mathbb{S}^{3} : \pi\in Mh^{a}_{F(d)}(F(\alpha))\}$,
where $F(d)<\mathbb{K}$ and $F(\alpha)<\varphi_{\mathbb{K}}(0)$ with $\{F(d),F(\alpha)\}\subset\mathcal{H}_{a}(\pi)$.
For fixed $a$, the set
$\{(d,\eta)\in \mathbb{K}\times\varphi_{\mathbb{K}}(0)   : \mathbb{S}\in Mh^{a}_{d}(\eta)\}$
is defined from the class $P$ by recursion on ordinals $d<\mathbb{K}$.
Let $\varphi$ be a formula such that $V_{\mathbb{S}+\mathbb{K}}\models\varphi[P, C,\mathbb{S},b]$ iff
$\mathbb{S}\in Mh_{0}^{a}(g_{c})\cap Mh^{a}_{c}(f^{c})$ and $C$ is a club subset of $\mathbb{S}$.
Since $\mathbb{S}$ is shrewd, pick $b<\mathbb{S}_{0}<\mathbb{K}_{0}<\mathbb{S}$ such that
$V_{\mathbb{S}_{0}+\mathbb{K}_{0}}\models\varphi[P\cap \mathbb{S}_{0},C\cap\mathbb{S}_{0},\mathbb{S}_{0},b]$.
We obtain $\mathbb{S}_{0}\in A\cap B\cap C$.
Therefore $\mathbb{S}\in Mh^{a}_{c}(\xi)$ is shown.
$\mathbb{S}\in M_{2}(Mh^{a}_{c}(\xi))$ is seen from the shrewdness of 
$\mathbb{S}$.
\eprf

\bcor\label{cor:welldefinedness.1}
$\tilde{T}$ proves that
$\forall a<\varepsilon_{\mathbb{K}+1}\forall c<\mathbb{K}[\{a,c,\xi\}\subset\mathcal{H}_{a}(\mathbb{S}) \to \psi_{\mathbb{S}}^{f}(a)<\mathbb{S})]$
for every $\xi<\varphi_{\mathbb{K}}(0)$ and finite functions 
$f$ such that ${\rm supp}(f)=\{c\}$, $c<\mathbb{K}$ and
$f(c)=\xi$.
\ecor

\blem\label{lem:stepdownpi11}
Assume $\mathbb{S}\geq\pi\in Mh^{a}_{d}(\xi)\cap Mh^{a}_{c}(\xi_{0})$, $\xi_{0}\neq 0$,
and $d<c$.
Moreover let $\xi_{1}\in\mathcal{H}_{a}(\pi)$ for 
$\xi_{1}\leq\tilde{\theta}_{c-d}(\xi_{0})$.
Then
$\pi\in Mh^{a}_{d}(\xi\dot{+}\xi_{1})\cap M^{a}_{d}(Mh^{a}_{d}(\xi\dot{+}\xi_{1}))$.
\elem
\bprf
As in Lemma \ref{lem:stepdownpi11int}.
\eprf

\bdf\label{df:nfform2}
{\rm
For finite functions $f$ and $g$,
\[
 Mh^{a}_{0}(g)\prec Mh^{a}_{0}(f)
 :\Leftrightarrow 
\forall\pi\in Mh^{a}_{0}(f)
\left(
SC_{\mK}(g)\subset\mathcal{H}_{a}(\pi) \Rightarrow \pi\in M_{2}(Mh^{a}_{0}(g))
\right)
.
\]
}
\edf

\bcor\label{cor:stepdownpi11}
Let $f,g$ be finite functions and $c\in{\rm supp}(f)$.
Assume  that 
there exists an ordinal
$d<c$ 
such that
$(d,c)\cap {\rm supp}(f)=(d,c)\cap {\rm supp}(g)=\emptyset$, 
$g_{d}=f_{d}$, 
$g(d)<f(d)\dot{+}\tilde{\theta}_{c-d}(f(c))\cdot\omega$,
and
$g<^{c}f(c)$.
Then
$Mh^{a}_{0}(g)\prec Mh^{a}_{0}(f)$ holds.
In particular if $\pi\in Mh^{a}_{0}(f)$ and
$SC_{\mK}(g)\subset\mathcal{H}_{a}(\pi)$, then
$\psi_{\pi}^{g}(a)<\pi$.
\ecor

\bprp\label{lem:psinucomparison}

Let $f,g:\mK\to\vphi_{\mK}(0)$.
If $f<^{0}_{lx}g$, then
$Mh^{a}_{0}(f)\prec Mh^{a}_{0}(g)$.

\eprp
\bprf
This is seen from Corollary \ref{cor:stepdownpi11}.
\eprf

\subsection{A Mostowski collapsing}
$OT(\Pi^{1}_{1})$ denotes a computable notation system of ordinals
with a constant $\mS$ for a stable ordinal,
collapsing functions $\psi_{\sig}^{g}(a)$ for finite functions $g$, where
${\rm supp}(g)=\{d\}$ for a $d<\mK=\mS^{+}$ and $g(d)<\veps_{\mK+1}$ if $\sig=\mS$.
Let $m(\alp)=g$ for $\alp=\psi_{\sig}^{g}(a)$ and $\sig<\mS$.
For $g\neq\emptyset$, $\alp=\psi_{\sig}^{g}(a)\in OT(\Pi^{1}_{1})$ only when
$g$ is obtained from $f=m(\sig)$ as follows, cf.\,Corollary \ref{cor:stepdownpi11}.
There are $c$ and $d$ such that $d<c\in {\rm supp}(f)$,
and $(d,c)\cap {\rm supp}(f)=\emptyset$.
Then 
$g_{d}=f_{d}$, $(d,c)\cap {\rm supp}(g)=\emptyset$
$g(d)<f(d)+\tilde{\theta}_{c-d}(f(c))\cdot\omega$, 
and $g<^{c}f(c)$.

In what follows, by ordinals we mean ordinal terms in $OT(\Pi^{1}_{1})$.
$\Psi_{\mS}$ denotes the set of ordinal terms $\psi_{\sig}^{f}(a)$ for some
$a,f$ and $\sig\in\Psi_{\mS}\cup\{\mS\}$.
Note that in $OT(\Pi^{1}_{1})$,
$\psi_{\sig}^{f}(a)\geq\mS$ only if $\sig=\mK=\mS^{+}$ and $f=\emptyset$.

We define a Mostowski collapsing $\alpha\mapsto\alpha[\rho/\mS]$,
which is needed to replace inference rules for stability by ones of reflections.
The domain of the collapsing $\alpha\mapsto\alpha[\rho/\mS]$ is a subset
$M_{\rho}$ of $E^{\mS}_{\rho}$.
For a reason of the restriction, see the beginning of subsection \ref{subsec:operatorcontpi11}.

\bdf
{\rm
For ordinal terms $\psi_{\sig}^{f}(a)\in \Psi_{\mS}\subset OT(\Pi^{1}_{1})$,
define
$m(\psi_{\sig}^{f}(a)):=f$ and $s(\psi_{\sig}^{f}(a)):=\max({\rm supp}(f))$.
Also
${\tt p}_{0}(\psi_{\sig}^{f}(a))={\tt p}_{0}(\sig)$ if $\sig<\mS$, and
${\tt p}_{0}(\psi_{\mS}^{f}(a))=a$.
}
\edf

\bdf
$M_{\rho}:=\calh_{b}(\rho)$ for $b={\tt p}_{0}(\rho)$ and $\rho\in\Psi_{\mS}$.
\edf

$\alp=\psi_{\sig}^{g}(a)\in OT(\Pi^{1}_{1})$ only when
$\{\sig,a\}\subset \calh_{a}(\alp)$ and
$SC_{\mK}(g)\subset M_{\alp}$.

$OT(\Pi^{1}_{1})$ is defined to be closed under
$\alpha\mapsto\alpha[\rho/\mS]$ for $\alp\in M_{\rho}$.
Specifically if $\{\alp,\rho\}\subset OT(\Pi^{1}_{1})$ with
$\alp\in M_{\rho}$ and $\rho\in\Psi_{\mS}$, then
$\alp[\rho/\mS]\in OT(\Pi^{1}_{1})$.

\bprp\label{prp:EK2pi11}
Let $\rho\in\Psi_{\mS}$.
\benu
\item\label{prp:EK2.1}
$\mathcal{H}_{\gamma}(M_{\rho})\subset M_{\rho}$ 
if $\gam\leq{\tt p}_{0}(\rho)$.

\item\label{prp:EK2.2}
$M_{\rho}\cap\mS=\rho$ and $\rho\not\in M_{\rho}$.

\item\label{prp:EK2.3}
If $\sig<\rho$ and ${\tt p}_{0}(\sig)\leq {\tt p}_{0}(\rho)$, then
$M_{\sig}\subset M_{\rho}$.
\eenu
\eprp

\bdf\label{df:Mostwskicollaps}
{\rm
Let $\alpha\in M_{\rho}$ with $\rho\in \Psi_{\mS}$.
We define an ordinal $\alpha[\rho/\mS]$ recursively as follows.
$\alpha[\rho/\mS]:=\alpha$ when $\alpha<\mathbb{S}$.
In what follows assume $\alpha\geq\mathbb{S}$.

$\mathbb{S}[\rho/\mS]:=\rho$.
$\mK[\rho/\mS]\equiv(\mathbb{S}^{+})[\rho/\mS]:=\rho^{+}$.
$\left(\psi_{\mK}(a)\right)[\rho/\mS]=\left(\psi_{\mathbb{S}^{+}}(a)\right)[\rho/\mS]=\psi_{\rho^{+}}(a[\rho/\mS])$.
The map commutes with $+$ and $\varphi$.
}
\edf

\blem\label{lem:Mostowskicollapspi11}
For $\rho\in\Psi_{\mS}$,
$\{\alpha[\rho/\mS]:\alpha\in M_{\rho}\}$ is a transitive collapse of $M_{\rho}$ in the sense that
$\beta<\alpha\Leftrightarrow\beta[\rho/\mS]<\alpha[\rho/\mS]$,
$\beta\in\mathcal{H}_{\alpha}(\gamma)\Leftrightarrow 
\beta[\rho/\mS]\in\mathcal{H}_{\alpha[\rho/\mS]}(\gamma[\rho/\mS]))
$ for $\gamma>\mathbb{S}$,
and
$OT(\Pi^{1}_{1})\cap\alpha[\rho/\mS]=\{\beta[\rho/\mS]:\beta\in M_{\rho}\cap\alpha\}$
for $\alpha,\beta,\gamma\in M_{\rho}$.
\elem

Let $\rho\leq\mathbb{S}$, and
$\iota$ an $RS$-term or an $RS$-formula such that
${\sf k}(\iota)\subset M_{\rho}$, where $M_{\mathbb{S}}=\mathbb{K}$.
Then
$\iota^{[\rho/\mathbb{S}]}$ denotes the result of replacing each unbounded quantifier
$Qx$ by $Qx\in L_{\mathbb{K}[\rho/\mS]}$,
and each ordinal term $\alpha\in {\sf k}(\iota)$
by $\alpha[\rho/\mS]$ for the Mostowski collapse in
Definition \ref{df:Mostwskicollaps}.

\bprp\label{prp:levelcollaps}
Let $\rho\in \Psi_{\mS}\cup\{\mathbb{S}\}$.

\begin{enumerate}
\item\label{prp:levelcollaps1}
Let $v$ be an $RS$-term with ${\sf k}(v)\subset M_{\rho}$,
and $\alpha=|v|$.
Then 
$v^{[\rho/\mathbb{S}]}$ is an $RS$-term of level $\alpha[\rho/\mS]$,
$\left| v^{[\rho/\mathbb{S}]} \right|=\alpha[\rho/\mS]$ and
${\sf k}(v^{[\rho/\mathbb{S}]})=\left({\sf k}(v)\right)^{[\rho/\mathbb{S}]}$.

\item\label{prp:levelcollaps2}
Let $\alpha\leq\mathbb{K}$ be such that $\alpha\in M_{\rho}$. Then
$
\left(Tm(\alpha)\right)^{[\rho/\mathbb{S}]}
:= 
\{v^{[\rho/\mathbb{S}]}: v\in Tm(\alpha), {\sf k}(v)\subset M_{\rho}\}
= Tm(\alpha[\rho/\mS])$.

\item\label{prp:levelcollaps3}
Assume $\mathcal{H}_{\gamma}(\rho)\cap\mathbb{S}\subset\rho$.
For an $RS$-formula $A$ with ${\sf k}(A)\subset \mathcal{H}_{\gamma}(\rho)$, 
$A^{[\rho/\mathbb{S}]}$ is an $RS$-formula such that
${\sf k}(A^{[\rho/\mathbb{S}]})\subset\{\alpha[\rho/\mS]: \alpha\in{\sf k}(A)\}\cup\{\mathbb{K}[\rho/\mS]\}$.
\end{enumerate}
\eprp

For each sentence $A$, either a disjunction is assigned as $A\simeq\bigvee(A_{\iota})_{\iota\in J}$, or
a conjunction is assigned as $A\simeq\bigwedge(A_{\iota})_{\iota\in J}$.
In the former case $A$ is said to be a \textit{$\bigvee$-formula}, and in the latter
$A$ is a \textit{$\bigwedge$-formula}.

\bdf\label{df:assigndc}
{\rm
Let
$[\rho]Tm(\alpha):=\{u\in Tm(\alpha) : {\sf k}(u)\subset M_{\rho}\}$.

}
\edf

\bprp\label{lem:assigncollaps}
Let $\rho\in \Psi_{\mS}\cup\{\mathbb{S}\}$.
For $RS$-formulas $A$, let
$A\simeq \bigvee(A_{\iota})_{\iota\in J}$ and assume ${\sf k}(A)\subset M_{\rho}$.
Then
$A^{[\rho/\mathbb{S}]}\simeq \bigvee\left((A_{\iota})^{[\rho/\mathbb{S}]}\right)_{\iota\in [\rho]J}$.
The case $A\simeq \bigwedge(A_{\iota})_{\iota\in J}$ is similar.
\eprp

\subsection{Operator controlled derivations for $\Pi^{1}_{1}$-reflection}\label{subsec:operatorcontpi11}

We define a derivability relation 
$(\mathcal{H}_{\gamma},\Theta;{\tt Q}_{\Pi})\vdash^{* a}_{c}\Gamma;\Pi^{[\cdot]}$
where 
${\tt Q}_{\Pi}$ is a finite set of ordinals in $\Psi_{\mS}$,
$c$ is a bound of ranks of the inference rules $({\rm stbl})$ and of ranks of cut formulas.
The relation depends on an ordinal $\gam_{0}$, and
should be written as $(\mathcal{H}_{\gamma},\Theta;{\tt Q}_{\Pi})\vdash^{*a}_{c,\gamma_{0}} \Gamma;\Pi^{[\cdot]}$.
However the ordinal $\gam_{0}$ will be fixed.
So let us omit it.

The r\^{o}le of the calculus $\vdash^{* a}_{c}$ is twofold:
first finite proof figures are embedded in the calculus, and
second the cut rank $c$ in $\vdash^{* a}_{c}$ is lowered to $\mK=\mS^{+}$.
In the next subsection \ref{subsec:operatorcontcap} the relation $\vdash^{* a}_{c}$
is embedded in another derivability relation $\vdash^{a}_{c,e,b_{1}}A^{(\rho)}$ with caps $\rho$.
In the latter calculus, cut ranks $c$ as well as the ranks of formulas to be reflected
are lowered to $\mS$, and the inferences for reflections are removed.
For this we need to distinguish formulas with smaller ranks$<\mS$ from higher ones.

As in
Lemma \ref{lem:lowerPi43b}, in eliminating of inferences for reflections,
{\small
\[
\infer[({\rm rfl}_{\rho})]{
\calh_{\gam}[\Tht]\vdash^{a}\Gam^{(\rho)}
}
{
\{
\mathcal{H}_{\gam}[\Tht]\vdash\Gamma^{(\rho)}, \lnot\del^{(\rho)}
\}_{\del\in\Del}
&
\{
\mathcal{H}_{\gam}[\Tht\cup \{\sig\}]\vdash\Gamma^{(\rho)}, 
\Del^{(\sig)}
\}_{\sig}
}
\]
}
 is rewritten to, cf.\,Recapping \ref{lem:recappingpi11}
{\small
\[
\hspace{-10mm}
\infer[({\rm rfl}_{\kap})]{\mathcal{H}_{\gam}[\Tht]\vdash\Gam^{(\kap)}}
{
\{\mathcal{H}_{\gam}[\Tht]\vdash\lnot\tht^{(\kap)},\Gam^{(\kap)}\}_{\tht\in\Gam}
&
\hspace{-15mm}
\infer[(cut)]{
\{\mathcal{H}_{\gam}[\Tht\cup \{\sig\}]\vdash
\Gamma^{(\kap)},\Gam^{(\sig)}\}_{\sig}
}
 {
 \{
 \infer*[\rho\leadsto\sig]{\mathcal{H}_{\gam}[\Tht]\vdash\Gamma^{(\sig)}, \lnot \del^{(\sig)}
 \}_{\del\in\Del}}{}
 &
 \infer*[\rho\leadsto\kap]{\mathcal{H}_{\gam}[\Tht\cup \{\sig\}]\vdash\Gamma^{(\kap)}, 
\Del^{(\sig)}}{}
 }
}
\]
}
where $\sig<\kap<\rho$.
In the rewriting, the inference $({\rm rfl}_{\rho})$ is replaced by $({\rm rfl}_{\kap})$
for a smaller $\kap<\rho$.
This means that $({\rm rfl}_{\rho})$ is replaced by $({\rm rfl}_{\sig})$ in the part $\rho\leadsto\sig$.
$\kap$ reflects $\Gam$ to some $\sig$, and
$\sig$ has to reflect $\Del$, where $\rk(\Del)>\rk(\Gam)$ is possible.
Therefore the termination of the whole process of removing is seen to be by induction on 
reflecting ordinals $\rho$, cf.\,Lemma \ref{lem:sum}.

The Mahlo degree $g=m(\kap)$ in $\kap=\psi_{\rho}^{g}(\alp)$
 is obtained by (an iteration of) a stepping-down
$(f,d,c)\mapsto g$, where $f=m(\rho)$, $d<c\in {\rm supp}(f)$,
$(d,c)\cap {\rm supp}(f)=\emptyset$,
$g_{d}=f_{d}$, $(d,c)\cap {\rm supp}(g)=\emptyset$,
$g(d)<f(d)+\tilde{\theta}_{c-d}(f(c))\cdot\omega$, 
and $g<^{c}f(c)$.
$g$ depends on $a$, $\rho$ and $\rk(\Gam^{(\rho)}):=\rk(\Gam)$.
In showing 
\[
SC_{\mK}(g)\subset\calh_{\alp}(\kap)
\]
$\rho$ and $\rk(\Gam^{(\rho)})$ are harmless since these relates to the given ordinal $\rho$,
while the ordinal $a$ causes trouble, since all of the reflecting ordinals $\rho,\ldots$ share
the ordinal depth $a$ of the derivation.
We need $a\in\calh_{\alp_{0}}(\rho)$ if $\rho=\psi_{\sig}^{f}(\alp_{0})$, and
$a\in\calh_{\bet}(\tau)$ if $\tau=\psi_{\lam}^{h}(\bet)$, and so forth.
This leads us to the set $M_{\rho}=\calh_{b}(\rho)$ for $b={\tt p}_{0}(\rho)$,
where $\displaystyle{\rho=\psi^{f}_{\ddots \psi_{\mS}(b)}(\alp_{0})}$,
and the condition (\ref{eq:controlderpi11ac}) that $a$ as well as ordinals occurring in the derivation
should be in $M_{\rho}$ for every reflecting ordinal $\rho$ occurring in derivations.
Note that $M_{\rho}=\calh_{b}(\rho)\subset\calh_{\alp_{0}}(\rho)$ by $b\leq\alp_{0}$, but
$E^{\mS}_{\rho}\not\subset\calh_{\alp_{0}}(\rho)$.
This is the reason why we restrict the domain of the Mostowski collapsing $\alp\mapsto \alp[\rho/\mS]$
to $\alp\in M_{\rho}\subsetneq E^{\mS}_{\rho}$.

${\tt Q}_{\Pi}$ in $(\mathcal{H}_{\gamma},\Theta;{\tt Q}_{\Pi})\vdash^{* a}_{c}\Gamma;\Pi^{[\cdot]}$, is the set of ordinals $\sig$ which is introduced in a right upper sequent 
$(\mathcal{H}_{\gamma},\Theta\cup\{\sig\};{\tt Q}_{\Pi}\cup\{\sig\}
)
\vdash^{* a_{0}}_{c}
\Gamma;\Pi^{[\cdot]}, \lnot B(u)^{[\sigma/\mathbb{S}]}$
of an inference $({\rm stbl})$
for stability occurring below $(\mathcal{H}_{\gamma},\Theta;{\tt Q}_{\Pi})\vdash^{* a}_{c}\Gamma;\Pi^{[\cdot]}$, while the set $\Pi^{[\cdot]}=\bigcup\{\Pi_{\sig}^{[\sig/\mS]}:\sig\in {\tt Q}_{\Pi}\}$
is the collection of formulas $\lnot B(u)^{[\sigma/\mathbb{S}]}$.
{\small
\[
\hspace{-10mm}
\infer[({\rm stbl})]{(\mathcal{H}_{\gamma},\Theta;{\tt Q}_{\Pi}
)
\vdash^{* a}_{c}
\Gamma;\Pi^{[\cdot]}}
{
(\mathcal{H}_{\gamma},\Theta;{\tt Q}_{\Pi}
)\vdash^{* a_{0}}_{c}
\Gamma, B(u);\Pi^{[\cdot]}
&
\{
(\mathcal{H}_{\gamma},\Theta\cup\{\sig\};{\tt Q}_{\Pi}\cup\{\sig\}
)
\vdash^{* a_{0}}_{c}
\Gamma;\Pi^{[\cdot]}, \lnot B(u)^{[\sigma/\mathbb{S}]}
\}_{\sig}
}
\]
}
These motivates the following Definitions \ref{df:QJpi11}, \ref{df:controldercollapspi11} and \ref{df:controldercollapscappi11}.

\bdf\label{df:QJpi11}
{\rm
Let ${\tt Q}\subset\Psi_{\mS}$ be a finite set of ordinals, and $A\simeq\bigvee(A_{\iota})_{\iota\in J}$.
Define $M_{{\tt Q}} :=  \bigcap_{\sig\in{\tt Q}} M_{\sig}$,
\[
[{\tt Q}]_{A}J  := [{\tt Q}]_{\lnot A}J :=  
\{\iota\in J: \rk(A_{\iota})\geq\mS \Rarw \sfk(\iota)\subset M_{{\tt Q}}\} 
\]
\[
\sfk^{\mS}(\Gam)   :=   \bigcup\{\sfk(A): A\in\Gam, \rk(A)\geq\mS\}
\]

}
\edf

\bdf\label{df:controldercollapspi11}
{\rm
Let $\Theta$ be a finite set of ordinals, 
$\gam\leq\gam_{0}$ and $a,c$ ordinals\footnote{In this subsection \ref{subsec:operatorcontpi11}
we can set $\gam=\mS$.}, and
${\tt Q}_{\Pi}\subset\Psi_{\mS}$ a finite set of ordinals such that
${\tt p}_{0}(\sig)\geq\gam_{0}$
for each $\sig\in{\tt Q}_{\Pi}$.
Let
$\Pi=\bigcup\{\Pi_{\sig}:\sig\in {\tt Q}_{\Pi}\}\subset\Delta_{0}(\mathbb{K})$ be a set of formulas
such that $\sfk(\Pi_{\sig})\subset M_{\sig}$ 
for each $\sig\in{\tt Q}_{\Pi}$, 
$\Pi^{[\cdot]}=\bigcup\{\Pi_{\sig}^{[\sig/\mS]}:\sig\in {\tt Q}_{\Pi}\}$,
$\Tht^{(\sig)}=\Tht\cap M_{\sig}$
and
$\Theta_{{\tt Q}_{\Pi}}=\Tht\cap M_{{\tt Q}_{\Pi}}$.

$(\mathcal{H}_{\gamma},\Theta; {\tt Q}_{\Pi})\vdash^{* a}_{c} \Gamma;\Pi^{[\cdot]}$ holds
for a set
$\Gamma$ of formulas
if 
\begin{equation}
\label{eq:controlderpi11}
{\sf k}(\Gam)\subset\mathcal{H}_{\gamma}[\Theta]
\spand
\fal\sig\in{\tt Q}_{\Pi}\left(
\sfk(\Pi_{\sig})\subset\calh_{\gam}[\Tht^{(\sig)}]
\right)
\end{equation}

\begin{equation}
\label{eq:controlderpi11ac}
\{\gam,a,c\}
\cup\sfk^{\mS}(\Gam)\cup\sfk^{\mS}(\Pi)
\subset\mathcal{H}_{\gamma}[\Theta_{{\tt Q}_{\Pi}}]
\end{equation}

and one of the following cases holds:

\begin{description}

\item[$(\bigvee)$]\footnote{The condition (\ref{eq:controlder1KP}), $|\iota|< a$ is absent in the inference $(\bigvee)$, cf.\,\textbf{Case 3} in Lemma \ref{lem:capping}.}
There exist 
$A\simeq\bigvee(A_{\iota})_{\iota\in J}$, an ordinal
$a(\iota)<a$ and an 
$\iota\in J$ such that
$A\in\Gamma$,
$(\mathcal{H}_{\gamma},\Theta;{\tt Q}_{\Pi})\vdash^{*a(\iota)}_{c}\Gamma,
A_{\iota};\Pi^{[\cdot]}$.

\item[$(\bigvee)^{[\cdot]}$]
There exist 
$A\equiv B^{[\sig/\mS]}\in\Pi^{[\cdot]}$,
$B\simeq\bigvee(B_{\iota})_{\iota\in J}$, an ordinal
$a(\iota)<a$ and an $\iota\in [\sig]J$ such that
$(\mathcal{H}_{\gamma},\Theta;{\tt Q}_{\Pi})\vdash^{*a(\iota)}_{c}\Gamma;\Pi^{[\cdot]},
A_{\iota}$ with $A_{\iota}\equiv B_{\iota}^{[\sig/\mS]}$.

\item[$(\bigwedge)$]
There exist 
$A\simeq\bigwedge(A_{\iota})_{ \iota\in J}$,  ordinals $a(\iota)<a$ such that
$A\in\Gamma$ and
$(\mathcal{H}_{\gamma},\Theta\cup\sfk(\iota);{\tt Q}_{\Pi}
)
\vdash^{*a(\iota)}_{c}\Gamma,
A_{\iota};\Pi^{[\cdot]}$
for each $\iota\in[{\tt Q}_{\Pi}]_{A} J$.

\item[$(\bigwedge)^{[\cdot]}$]
There exist 
$A\equiv B^{[\sig/\mS]}\in\Pi^{[\cdot]}$,
$B\simeq\bigwedge(B_{\iota})_{\iota\in J}$, ordinals
$a(\iota)<a$ such that 
$(\mathcal{H}_{\gamma},\Theta\cup\sfk(\iota);{\tt Q}_{\Pi})\vdash^{*a(\iota)}_{c}\Gamma;
A_{\iota},\Pi^{[\cdot]}$
for each $\iota\in [{\tt Q}_{\Pi}]_{B}J\cap[\sig]J$.

\item[$(cut)$]
There exist an ordinal $a_{0}<a$ and a formula $C$
such that 
$(\mathcal{H}_{\gamma},\Theta;{\tt Q}_{\Pi})\vdash^{*a_{0}}_{c}\Gamma,\lnot C;\Pi^{[\cdot]}$
and
$(\mathcal{H}_{\gamma},\Theta;{\tt Q}_{\Pi})\vdash^{*a_{0}}_{c}C,\Gamma;\Pi^{[\cdot]}$
with $\mbox{{\rm rk}}(C)<c$.

\item[$(\Sigma\mbox{{\rm -rfl}})$]

There exist ordinals
$a_{\ell}, a_{r}<a$ and a formula $C\in\Sigma(\pi)$ for a 
$\pi\in\{\Omega,\mK=\mS^{+}\}$
such that $c\geq \pi$, 
$(\mathcal{H}_{\gamma},\Theta;{\tt Q}_{\Pi}
)\vdash^{* a_{\ell}}_{c}\Gamma,C;\Pi^{[\cdot]}$
and
$(\mathcal{H}_{\gamma},\Theta;{\tt Q}_{\Pi}
)\vdash^{* a_{r}}_{c}
\lnot \exists x<\pi\,C^{(x,\pi)}, \Gamma;\Pi^{[\cdot]}$.

\item[({\rm stbl})]
There exist an ordinal $a_{0}<a$,
a $\bigwedge$-formula 
$B(0)\in\Delta_{0}(\mathbb{S})$,
and a $u\in Tm(\mathbb{K})$ 
for which the following hold:
$\mathbb{S}\leq {\rm rk}(B(u))<c$,
$(\mathcal{H}_{\gamma},\Theta;{\tt Q}_{\Pi}
)\vdash^{* a_{0}}_{c}
\Gamma, B(u);\Pi^{[\cdot]}$, 
and 
$(\mathcal{H}_{\gamma},\Theta\cup\{\sig\};{\tt Q}_{\Pi}\cup\{\sig\}
)
\vdash^{* a_{0}}_{c}
\Gamma;\Pi^{[\cdot]}, \lnot B(u)^{[\sigma/\mathbb{S}]}$
holds for every ordinal $\sigma\in\Psi_{\mS}$ such that
$\Tht\subset M_{\sigma}$.
{\small
\[
\hspace{-15mm}
\infer[({\rm stbl})]{(\mathcal{H}_{\gamma},\Theta;{\tt Q}_{\Pi}
)
\vdash^{* a}_{c}
\Gamma;\Pi^{[\cdot]}}
{
(\mathcal{H}_{\gamma},\Theta;{\tt Q}_{\Pi}
)\vdash^{* a_{0}}_{c}
\Gamma, B(u);\Pi^{[\cdot]}
&
\{
(\mathcal{H}_{\gamma},\Theta\cup\{\sig\};{\tt Q}_{\Pi}\cup\{\sig\}
)
\vdash^{* a_{0}}_{c}
\Gamma;\Pi^{[\cdot]}, \lnot B(u)^{[\sigma/\mathbb{S}]}
\}_{\Tht\subset M_{\sigma}}
}
\]
}
Note that $(\Tht\cup\{\sig\})_{{\tt Q}_{\Pi}\cup\{\sig\}}=\Tht_{{\tt Q}_{\Pi}}$ if $\Tht_{{\tt Q}_{\Pi}}\subset M_{\sigma}$.

\end{description}
}
\edf

\bprp\label{lem:tautologypi11}{\rm (Tautology)}
Let $\gam\in\calh_{\gam}[\sfk(A)]$ and $d=\mbox{{\rm rk}}(A)$.
\benu
\item\label{lem:tautologypi11.1}
$(\mathcal{H}_{\gamma},{\sf k}(A);\emptyset
)\vdash^{* 2d}_{0}
\lnot A, A; \emptyset$.

\item\label{lem:tautologypi11.2}
$(\mathcal{H}_{\gamma},{\sf k}(A)\cup\{\sig\};\{\sig\}
)\vdash^{* 2d}_{0}
\lnot A^{[\sig/\mS]};A^{[\sig/\mS]}$ if $\sfk(A)\subset M_{\sig}$ and $\gam\geq\mS$.
\eenu
\eprp
\bprf
Both are seen by induction on $d$.
Consider Proposition \ref{lem:tautologypi11}.\ref{lem:tautologypi11.2}.

We have $(\sfk(A)\cup\{\sig\})\cap M_{\sig}=\sfk(A)$ for (\ref{eq:controlderpi11}) and (\ref{eq:controlderpi11ac}),
and $\sfk(A^{[\sig/\mS]})\subset\calh_{\mS}((\sfk(A)\cap\mS)\cup\{\sig\})$ for (\ref{eq:controlderpi11}).
Note that 
$\sig\not\in\mathcal{H}_{\gamma}[{\sf k}(A)]$ since $\sig\not\in\sfk(A)\subset M_{\sig}$ and 
$\gam\leq\gam_{0}\leq{\tt p}_{0}(\sig)$, and 
$\rk(A^{[\sig/\mS]})\not\in\mathcal{H}_{\gamma}[(\sfk(A)\cup\{\sig\})\cap M_{\sig}]$.

Let $A\simeq\bigvee(A_{\iota})_{\iota\in J}$.
Then $A^{[\sig/\mS]}\simeq\bigvee(A_{\iota}^{[\sig/\mS]})_{\iota\in [\sig]J}$ by 
Proposition \ref{lem:assigncollaps} and
$\sfk(\iota^{[\sig/\mS]})\subset\calh_{\mS}[(\sfk(\iota)\cap\mS)\cup\{\sig\}]$.
Let $I=\{\iota^{[\sig/\mS]}: \iota\in[\sig]J\}$.
Then $A^{[\sig/\mS]}\simeq\bigvee(B_{\nu})_{\nu\in I}$ with $B_{\nu}\equiv A_{\iota}^{[\sig/\mS]}$ for $\nu=\iota^{[\sig/\mS]}$,
and
$[\{\sig\}]_{A^{[\sig/\mS]}} I=I$ 
by $\rk(A^{[\sig/\mS]})<\mS$.
For $d_{\iota}=\rk(A_{\iota})\in\calh_{\gam}[\sfk(A,\iota)]$ with $\iota\in[\sig]J=[\{\sig\}]_{A^{(\sig)}}J$ 
we obtain
\[
\infer[(\bigwedge)]{(\mathcal{H}_{\gamma},{\sf k}(A)\cup\{\sig\};\{\sig\})\vdash^{* 2d}_{0}
\lnot A^{[\sig/\mS]}; A^{[\sig/\mS]}}
{
\infer[(\bigvee)^{[\cdot]}]{(\mathcal{H}_{\gamma},{\sf k}(A,\iota)\cup\{\sig\};\{\sig\})\vdash^{* 2d_{\iota}+1}_{0}
\lnot A_{\iota}^{[\sig/\mS]}; A^{[\sig/\mS]}}
 {
  (\mathcal{H}_{\gamma},{\sf k}(A,\iota)\cup\{\sig\};\{\sig\})\vdash^{* 2d_{\iota}}_{0}
\lnot A_{\iota}^{[\sig/\mS]}; A_{\iota}^{[\sig/\mS]}
 }
}
\]
and
\[
\infer[(\bigwedge)^{[\cdot]}]{(\mathcal{H}_{\gamma},{\sf k}(A)\cup\{\sig\}; \{\sig\})\vdash^{* 2d}_{0}
A^{[\sig/\mS]}; \lnot A^{[\sig/\mS]}}
{
\infer[(\bigvee)]{(\mathcal{H}_{\gamma},{\sf k}(A)\cup\sfk(\iota)\cup\{\sig\}; \{\sig\})\vdash^{* 2d_{\iota}+1}_{0}
A^{[\sig/\mS]}; \lnot A_{\iota}^{[\sig/\mS]}}
 {
  (\mathcal{H}_{\gamma},{\sf k}(A)\cup\sfk(\iota)\cup\{\sig\}; \{\sig\})\vdash^{* 2d_{\iota}}_{0}
A_{\iota}^{[\sig/\mS]}; \lnot A_{\iota}^{[\sig/\mS]}
 }
}
\]
\eprf

\blem\label{th:embedpi11}{\rm (Embedding of Axioms)}
For each axiom $A$ in $S_{1}$, there is an $m<\omega$ such that  
 $(\mathcal{H}_{\mS},\emptyset;\emptyset)\vdash^{* \mathbb{K}\cdot 2}_{\mathbb{K}+m} A;$
holds for $\mK=\mS^{+}$.
\elem
\bprf
We show that the axiom $\exi x\, B(x,v) \land v\in L_{\mathbb{S}} \to \exi x\in L_{\mS}\, B(x,v)\,(B\in\Del_{0})$ follows 
by an inference $({\rm stbl})$.
In the proof let us omit the operator $\mathcal{H}_{\mS}$.
Let $B(0)\in\Delta_{0}(\mathbb{S})$ be a $\bigwedge$-formula 
and $u\in Tm(\mathbb{K})$.
We may assume that
$\mK>d={\rm rk}(B(u))\geq\mathbb{S}$.
Let ${\sf k}_{0}={\sf k}(B(0))$ and ${\sf k}_{u}={\sf k}(u)$.
Let ${\sf k}_{0}\cup{\sf k}_{u}\subset M_{\sigma}$.
Then for $\exists x\in L_{\mathbb{S}}B(x)\simeq\bigvee(B(v))_{v\in J}$, we obtain
$
u^{[\sigma/\mathbb{S}]}\in J=Tm(\mathbb{S})
$ by $\rk(\exists x\in L_{\mathbb{S}}B(x))=\mS$.
We have
$B(u^{[\sigma/\mathbb{S}]})\equiv B(u)^{[\sigma/\mathbb{S}]}$,
$\sfk_{u}^{[\sig/\mS]}=\sfk(u^{[\sig/\mS]})\subset\calh_{\mS}[\sfk(u)\cup\{\sig\}]$,
$(\sfk_{0}\cup\sfk_{u})_{\emptyset}=\sfk_{0}\cup\sfk_{u}$ and
$({\sf k}_{0}\cup\sfk_{u}\cup\{\sig\})\cap M_{\sig}={\sf k}_{0}\cup\sfk_{u}$.
{\small
\[
\hspace{-5mm}
\infer[(\bigwedge)]{ {\sf k}_{0};\vdash^{* \mathbb{K}+1}_{\mathbb{K}}
 \lnot \exists x\, B(x),\exists x\in L_{\mathbb{S}}B(x);}
 {
 \infer[({\rm stbl})]{{\sf k}_{0}\cup{\sf k}_{u}; \vdash^{* \mathbb{K}}_{\mathbb{K}}
    \lnot B(u),\exists x\in L_{\mathbb{S}}B(x);}
    {
     {\sf k}_{0}\cup{\sf k}_{u}; \vdash^{* 2d}_{0}\lnot B(u),B(u);
    &
    \{
    \infer[(\bigvee)]{{\sf k}_{0}\cup\sfk_{u}\cup\{\sig\};\{\sig\}\vdash^{* 2d+1}_{0}
     \exists x\in L_{\mathbb{S}}B(x);\lnot B(u)^{[\sigma/\mathbb{S}]}
     \}_{{\sf k}_{0}\cup{\sf k}_{u}\subset M_{\sigma}}
     }
     {
     {\sf k}_{0}\cup{\sf k}_{u}\cup\{\sig\};\{\sig\} \vdash^{* 2d}_{0}
      B(u^{[\sigma/\mathbb{S}]}); \lnot B(u)^{[\sigma/\mathbb{S}]}
     }
    }
 }
\]
}
\eprf

\bprp\label{lem:inversionpi11}{\rm (Inversion)}
Let
$A\simeq \bigwedge(A_{\iota})_{\iota\in J}$ with $A\in\Gam$, $\iota\in [{\tt Q}_{\Pi}]_{A}J$ and
$(\mathcal{H}_{\gam},\Theta;{\tt Q}_{\Pi})\vdash^{* a}_{c}\Gamma;\Pi^{[\cdot]}$.
Then
$(\mathcal{H}_{\gam},\Theta\cup\sfk(\iota);{\tt Q}_{\Pi}
)
\vdash^{* a}_{c}
\Gamma,A_{\iota};\Pi^{[\cdot]}$.
\eprp

\bprp\label{lem:prereduction}
Let $(\mathcal{H}_{\gamma},\Theta;{\tt Q}_{\Pi})\vdash^{* a}_{c}\Gamma;\Pi^{[\cdot]}$.
Assume 
$\Tht\subset M_{\sig}$.
Then
\\
$(\mathcal{H}_{\gamma},\Theta\cup\{\sig\};{\tt Q}_{\Pi}\cup\{\sig\})\vdash^{* a}_{c}\Gamma;\Pi^{[\cdot]}$.
\eprp
\bprf
By induction on $a$.
We obtain
$(\Theta\cup\{\sig\})_{{\tt Q}_{\Pi}\cup\{\sig\}}=\Tht_{{\tt Q}_{\Pi}}$ 
by the assumption.
In an inference $({\rm stbl})$,
the right upper sequents are restricted to $\tau$ such that
$\sig\in M_{\tau}$.
Also we need to prune some branches at $(\bigwedge)$ and $(\bigwedge)^{[\cdot]}$
since $[({\tt Q}_{\Pi}\cup\{\sig\})]_{A}J\subset [{\tt Q}_{\Pi}]_{A}J$.
\eprf

\bprp\label{prp:reductionpi11}{\rm (Reduction)}
Let $C\simeq\bigvee(C_{\iota})_{\iota\in J}$ and
$\mK=\mS^{+}\leq\mbox{{\rm rk}}(C)\leq c$.
Assume
$(\mathcal{H}_{\gamma},\Theta;{\tt Q}_{\Pi})\vdash^{* a}_{c}\Gamma,\lnot C;\Pi^{[\cdot]}$
and
$(\mathcal{H}_{\gamma},\Theta;{\tt Q}_{\Pi})\vdash^{* b}_{c}C,\Gamma;\Pi^{[\cdot]}$.

Then
$(\mathcal{H}_{\gamma},\Theta;{\tt Q}_{\Pi})\vdash^{* a+b}_{c}\Gamma;\Pi^{[\cdot]}$.
\eprp
\bprf
By induction on $b$ using
Inversion \ref{lem:inversionpi11} and Proposition \ref{lem:prereduction}.

Note that if $(\mathcal{H}_{\gamma},\Theta;{\tt Q}_{\Pi})\vdash^{* b(\iota)}_{c}C_{\iota},\Gamma;\Pi^{[\cdot]}$ for an $\iota\in J$ such that $\rk(C_{\iota})\geq\mK$,
we obtain 
$\sfk(C_{\iota})\subset\calh_{\gam}[\Tht_{{\tt Q}_{\Pi}(\mS)}]\subset M_{{\tt Q}_{\Pi}(\mS)}$ by (\ref{eq:controlderpi11ac}) and Proposition \ref{prp:EK2pi11}
with $\gam\leq\gam_{0}\leq {\tt p}_{0}(\sig)$ for $\sig\in{\tt Q}_{\Pi}$.
Hence 
$\iota\in[{\tt Q}_{\Pi}]_{C}J$ if $\sfk(\iota)\subset\sfk(C_{\iota})$.
\eprf

\bprp\label{lem:predcepi11*}{\rm (Cut-elimination)}
Assume
$(\mathcal{H}_{\gamma},\Theta;{\tt Q}_{\Pi})\vdash^{* a}_{c+1}\Gamma;\Pi^{[\cdot]}$
with $c\geq\mS^{+}=\mK$.
Then $(\mathcal{H}_{\gamma},\Theta;{\tt Q}_{\Pi})\vdash^{* \ome^{a}}_{c}\Gamma;\Pi^{[\cdot]}$.
\eprp
\bprf
This is seen by induction on $a$
using Reduction \ref{prp:reductionpi11}.
\eprf

\blem\label{lem:Kcollpase.1}{\rm (Collapsing)}
Let
$\Gamma\subset\Sigma$ be a set of formulas, and $\Pi\subset\Del_{0}(\mK)$.
Suppose
$\Theta\subset
\mathcal{H}_{\gamma}(\psi_{\mK}(\gamma))$
 and
$
(\mathcal{H}_{\gamma},\Theta;{\tt Q}_{\Pi}
)\vdash^{* a}_{\mK}\Gamma;\Pi^{[\cdot]}
$.
Let
$\beta=\psi_{\mK}(\hat{a})$ with $\hat{a}=\gamma+\omega^{a}$.
Then 
$(\mathcal{H}_{\hat{a}+1},\Theta;{\tt Q}_{\Pi})
\vdash^{* \beta}_{\beta}
\Gamma^{(\beta,\mK)};\Pi^{[\cdot]}$ holds.
\elem
\bprf
By induction on $a$ as in Theorem \ref{th:CollapsingthmKP}.
We have 
$\{\gamma, a\}\subset\mathcal{H}_{\gamma}[\Theta_{{\tt Q}_{\Pi}}]$ by (\ref{eq:controlderpi11ac}),
and
$\beta\in\mathcal{H}_{\hat{a}+1}[\Theta_{{\tt Q}_{\Pi}}]$.

When the last inference is a $({\rm stbl})$, let
$B(0)\in\Delta_{0}(\mathbb{S})$ be a $\bigwedge$-formula
and a term $u\in Tm(\mK)$ such that
$\mS\leq{\rm rk}(B(u))<\mathbb{K}$,
${\sf k}(B(u))\subset \mathcal{H}_{\gamma}[\Theta]$,
and
$
(\mathcal{H}_{\gamma},\Theta;{\tt Q}_{\Pi}
)\vdash^{* a_{0}}_{\mK}
\Gamma, B(u);\Pi^{[\cdot]}
$ for an ordinal $a_{0}\in\mathcal{H}_{\gamma}[\Theta_{{\tt Q}_{\Pi}}]\cap a$.
Then we obtain
${\rm rk}(B(u))<\beta$.

Consider the case when the last inference is a $(\Sigma\mbox{{\rm -rfl}})$ on $\mK$.
We have ordinals
$a_{\ell}, a_{r}<a$ and a formula $C\in\Sigma$ such that
$(\mathcal{H}_{\gamma},\Theta;{\tt Q}_{\Pi}
)\vdash^{* a_{\ell}}_{\mK}\Gamma,C;\Pi^{[\cdot]}$
and
$(\mathcal{H}_{\gamma},\Theta;{\tt Q}_{\Pi}
)\vdash^{* a_{r}}_{\mK}
\lnot \exists x\,C^{(x,\mK)} , \Gamma;\Pi^{[\cdot]}$.

Let $\beta_{\ell}=\psi_{\mK}(\widehat{a_{\ell}})\in\mathcal{H}_{\widehat{a_{\ell}}+1}[\Theta_{{\tt Q}_{\Pi}}]\cap\beta$ with 
$\widehat{a_{\ell}}=\gamma+\omega^{a_{\ell}}$.
IH yields
$(\mathcal{H}_{\hat{a}+1},\Theta;{\tt Q}_{\Pi}
)\vdash^{* \beta_{\ell}}_{\beta}\Gamma^{(\beta,\mK)},
C^{(\beta_{\ell},\mK)};\Pi^{[\cdot]}$.
On the other, Inversion \ref{lem:inversionpi11} yields
$(\mathcal{H}_{\widehat{a_{\ell}}+1},\Theta;{\tt Q}_{\Pi}
)\vdash^{* a_{r}}_{\mK}
\lnot C^{(\beta_{\ell},\mK)} , \Gamma;\Pi^{[\cdot]}$.
For $\beta_{r}=\psi_{\mK}(\widehat{a_{r}})\in\mathcal{H}_{\hat{a}+1}[\Theta_{{\tt Q}_{\Pi}}]\cap\beta$
with $\widehat{a_{r}}=\widehat{a_{\ell}}+\omega^{a_{r}}$, IH yields
$(\mathcal{H}_{\hat{a}+1},\Theta;{\tt Q}_{\Pi}
)\vdash^{* \beta_{r}}_{\beta}
\lnot C^{(\beta_{\ell},\mK)} , \Gamma^{(\beta,\mK)}; \Pi^{[\cdot]}$.
We obtain
$(\mathcal{H}_{\hat{a}+1},\Theta;{\tt Q}_{\Pi}
)\vdash^{* \beta}_{\beta}
\Gamma^{(\beta,\mK)};\Pi^{[\cdot]}$
by a $(cut)$.

Note that since $\Pi\subset\Delta_{0}(\mK)$, inferences $(\bigwedge)^{[\cdot]}$
are harmless for the condition
$\Theta\subset \calh_{\gam}(\psi_{\mK}(\gam))$.
\eprf

\subsection{Operator controlled derivations with caps}\label{subsec:operatorcontcap}

In this subsection we introduce
another
derivability relation $(\mathcal{H}_{\gamma},\Theta,{\tt Q})\vdash^{a}_{c,e,b_{1}} \Gamma$,
which depends again on an ordinal $\gam_{0}$, and
should be written as $(\mathcal{H}_{\gamma},\Theta,{\tt Q})\vdash^{a}_{c,e,\gamma_{0},b_{1}} \Gamma$.
However the ordinal $\gam_{0}$ will be fixed, and specified in the proof of Theorem \ref{th:Pi11}.
So let us omit it.

The inference rules $({\rm stbl})$ are replaced by inferences $({\rm rfl}(\rho,d,f,b_{1}))$
by putting a \textit{cap} $\rho$ on formulas 
 in Lemma \ref{lem:capping}. 
 In $(\mathcal{H}_{\gamma},\Theta,{\tt Q})\vdash^{a}_{c,e,b_{1}} \Gamma$, $c$ is a bound for cut ranks
 and $e$ a bound for ordinals $\rho$ in the inferences $({\rm rfl}(\rho,d,f,b_{1}))$ occurring in the derivation.
$b_{1}$ is a bound such that $s(\rho)=\max({\rm supp}(m(\rho)))\leq b_{1}$.
 Although the capped formula $A^{(\rho)}$ in Definition \ref{df:cap.1},
 is intended to denote the formula $A^{[\rho/\mathbb{S}]}$,
we need to distinguish it from $A^{[\rho/\mathbb{S}]}$.
Our main task is to eliminate inferences $({\rm rfl}(\rho,d,f))$ from a resulting derivation $\mathcal{D}_{1}$.
In Recapping \ref{lem:recappingpi11} the cap $\rho$ in inferences $({\rm rfl}(\rho,d,f,b_{1}))$
are replaced by another cap $\kappa<\rho$.
In this process new inferences $({\rm rfl}(\sigma,d_{1},f_{1},b_{1}))$ arise with $\sigma<\kappa$.
Iterating this process,
we arrive at a derivation $\mathcal{D}_{2}$ such that
$s(\rho)\leq\mathbb{S}$, i.e., ${\rm supp}(m(\rho))\subset\mS+1$.
Then caps play no r\^{o}le, i.e.,
$A^{(\rho)}$ is `equivalent' to $A$ for $A\in\Delta_{0}(\mathbb{S})$.
Finally inferences $({\rm rfl}(\rho,d,f,b_{1}))$ are removed from $\mathcal{D}_{2}$ by 
throwing up caps and replacing 
these by a series of $(cut)$'s,
cf.\,Lemma \ref{lem:sum}.

The ordinal, i.e., the threshold $\gamma_{0}$ will be specified 
in the end of this section.

\bdf\label{df:cap.1}
{\rm
By a \textit{capped formula} we mean a pair $(A,\rho)$ of $RS$-sentence $A$
and an ordinal 
$\rho<\mathbb{S}$ such that
${\sf k}(A)\subset M_{\rho}$.
Such a pair is denoted by $A^{(\rho)}$.
A \textit{sequent} is a finite set of capped 
formulas, denoted by
$\Gamma_{0}^{(\rho_{0})},\ldots,\Gamma_{n}^{(\rho_{n})}$,
where each formula in the set $\Gamma_{i}^{(\rho_{i})}$ puts on the cap 
$\rho_{i}\in\mathbb{S}$.
When we write $\Gamma^{(\rho)}$, we tacitly assume that
${\sf k}(\Gamma)\subset M_{\rho}$.
A capped formula $A^{(\rho)}$ is said to be a $\Sigma(\pi)$-formula if
$A\in\Sigma(\pi)$.
Let ${\sf k}(A^{(\rho)}):={\sf k}(A)$.
}
\edf

\bdf\label{df:special}
{\rm 
Let $f$ be a non-empty (and irreducible) finite function.
Then $f$ is said to be \textit{special} if there exists an ordinal $\alpha$
such that $f(c_{\max})=\alpha+\mathbb{K}$ for $c_{\max}=\max({\rm supp}(f))$.
For a special finite function $f$, $f^{\prime}$ denotes a finite function such that
${\rm supp}(f^{\prime})={\rm supp}(f)$,
$f^{\prime}(c)=f(c)$ for $c\neq c_{\max}$, and
$f^{\prime}(c_{\max})=\alpha$ with $f(c_{\max})=\alpha+\mathbb{K}$.
}
\edf
The ordinal $\mK$ in $f(c_{\max})=\alpha+\mathbb{K}$ is a `room'
to be replaced by a smaller ordinal, cf.\,Definition \ref{df:hstepdownpi11}.

\bdf\label{df:cap.2}
{\rm
A finite set ${\tt Q}\subset\Psi_{\mS}$
is said to be a \textit{finite family} for ordinals $\gamma_{0}$ and $b_{1}$
if 
$\rho\in\calh_{\gam_{0}+\mS}=\calh_{\gam_{0}+\mS}(0)$,
$m(\rho):\mathbb{K} \to\varphi_{\mathbb{K}}(0)$ is special such that 
$s(\rho)=\max({\rm supp}(m(\rho)))\leq b_{1}$
and 
${\tt p}_{0}(\rho)\geq\gam_{0}$
for each $\rho\in {\tt Q}$.
}
\edf

The resolvent class $H_{\rho}(f,b_{1},\gamma_{0},\Theta)$ in the following Definition \ref{df:resolventpi11}
is the set of ordinals $\sig<\rho$, which are candidates of substitutes for $\rho$ in the inference
$({\rm rfl}(\rho,d,f,b_{1}))$ for reflection.
Note that if ${\tt p}_{0}(\sig)\leq {\tt p}_{0}(\rho)$ and $\sig<\rho$, then 
$M_{\sig}\subset M_{\rho}=\calh_{{\tt p}_{0}(\rho)}(\rho)$.
Moreover if ${\tt p}_{0}(\sig)\geq\gam_{0}\geq\gam$ and $\Tht\subset M_{\sig}$, then
$\calh_{\gam}[\Tht]\subset M_{\sig}$ by  Proposition \ref{prp:EK2pi11}.

\bdf\label{df:resolventpi11}
{\rm
$H_{\rho}(f,b_{1},\gamma_{0},\Theta)$
 denotes the \textit{resolvent class} for 
finite functions $f$, 
ordinals $\rho,b_{1},\gamma_{0}$ 
and finite sets $\Theta$ of ordinals defined by
$\sig\in H_{\rho}(f,b_{1}, \gamma_{0},\Theta)$ iff
$\sig\in\calh_{\gam_{0}+\mS}\cap\rho$,
$SC_{\mK}(m(\sig))\subset\mathcal{H}_{\gamma_{0}}[\Theta]$,
$\Theta\subset M_{\sig}$,
${\tt p}_{0}(\sig)={\tt p}_{0}(\rho)\geq\gam_{0}$,
and $m(\sig)$ is special such that
$s(f)=\max({\rm supp}(f))\leq s(\sig)\leq b_{1}$ and 
$f^{\prime} \leq (m(\sig))^{\prime}$,
where
$f\leq g\Leftrightarrow
\forall i(f(i)\leq g(i))$.
}
\edf

We define a derivability relation
$(\calh_{\gam},\Tht,{\tt Q})\vdash^{a}_{c,e}\Gam$,
where $\mS\leq\gam\leq\gam_{0}$ is an ordinal, $\Tht$ a finite set of ordinals,
${\tt Q}$ a finite family for $\gam_{0}, b_{1}$, and 
$a,c<\mK=\mS^{+}$.
 $c$ a bound of cut ranks, 
 $e$ a bound of $\rho$ in inference rules $({\rm rfl}(\rho,d,f,b_{1}))$, and
$b_{1}$ a bound on $s(\rho)$.
The relation $\vdash^{a}_{c,e}$ depends on fixed ordinals $\gam_{0}$ and $b_{1}$.

For $d=\rk(A)<\mS$, it may be $\sfk(A)\cup\{d\}\not\subset M_{{\tt Q}}$.
Let us avoid deriving the tautology $\lnot A,A$ by a standard derivation
to show $\vdash^{2d}\lnot A,A$.

\bdf\label{df:controldercollapscappi11}
{\rm
Let 
$\Tht^{(\rho)}=\Tht\cap M_{\rho}$, $[{\tt Q}]_{A^{(\rho)}}J=[{\tt Q}]_{A}J\cap[\rho]J$,
$\mS\leq\gamma\leq\gamma_{0}$ and $e\in\calh_{\gam_{0}+\mS}(0)$.

$(\mathcal{H}_{\gamma},\Theta,{\tt Q})\vdash^{a}_{c,e ,\gam_{0},b_{1}} \Gamma$ holds
for a set
$\Gamma=\bigcup\{\Gamma_{\rho}^{(\rho)}:\rho\in {\tt Q}\}$ of formulas
if

\begin{equation}
\label{eq:controlderpi11cap}
\forall \rho\in{\tt Q}
\left(
{\sf k}(\Gam_{\rho})
\subset\mathcal{H}_{\gamma}[\Theta^{(\rho)}]
\right)
\eeqn

\begin{equation}
\label{eq:controlderpi11picap}
\{\gamma,a,c,b_{1}\}\cup\sfk^{\mS}(\Gam)
\subset\mathcal{H}_{\gamma}[
\Theta_{{\tt Q}}]
\end{equation}

and one of the following cases holds:

\begin{description}

\item[{\rm (Taut)}]
$\{\lnot A^{(\rho)},A^{(\rho)}\}\subset\Gam$ for
a $\rho\in{\tt Q}$ and a formula $A$ such that $\rk(A)<\mS$.

\item[$(\bigvee)$]
There exist 
$A\simeq\bigvee(A_{\iota})_{\iota\in J}$, a cap $\rho\in {\tt Q}$, an ordinal
$a_{\iota}<a$ and an $\iota\in [\rho]J$
such that 
$A^{(\rho)}\in\Gamma$ and
$(\mathcal{H}_{\gamma},\Theta,{\tt Q})\vdash^{a_{\iota}}_{c,e,\gam_{0},b_{1}}\Gamma,
\left(A_{\iota}\right)^{(\rho)}$.

Note that if $\rk(A_{\iota})\geq\mS$, then 
$\sfk(A_{\iota})\subset\calh_{\gam}[\Tht_{{\tt Q}}]\subset M_{{\tt Q}}$ by (\ref{eq:controlderpi11picap}).
Hence $\iota\in [{\tt Q}]_{A}J=\{\iota\in J: \rk(A_{\iota})\geq\mS \Rarw \sfk(\iota)\subset M_{{\tt Q}}\}$.

\item[$(\bigwedge)$]
There exist 
$A\simeq\bigwedge(A_{\iota})_{ \iota\in J}$, 
a cap $\rho\in {\tt Q}$, ordinals $a_{\iota}<a$ 
for each $\iota\in [{\tt Q}]_{A^{(\rho)}}J$ such that
$A^{(\rho)}\in\Gamma$ and
$(\mathcal{H}_{\gamma},\Theta\cup{\sf k}(\iota), {\tt Q}
)
\vdash^{a_{\iota}}_{c,e,\gam_{0},b_{1}}\Gamma,
\left(A_{\iota}\right)^{(\rho)}$.

Note that if $\rk(A_{\iota})\geq\mS$, then $\sfk(\iota)\subset M_{{\tt Q}}$
by $\iota\in [{\tt Q}]_{A^{(\rho)}}J$.
Hence $\sfk^{\mS}(A_{\iota})\subset\calh_{\gam}[(\Tht\cup\sfk(\iota))_{{\tt Q}}]$
for (\ref{eq:controlderpi11picap}), where $(\Tht\cup\sfk(\iota))_{{\tt Q}}=\Tht_{{\tt Q}}\cup\sfk(\iota)$.

\item[$(cut)$]
There exist a cap $\rho\in {\tt Q}$, an ordinal $a_{0}<a$ and a formula $C$
such that 
$(\mathcal{H}_{\gamma},\Theta,{\tt Q})\vdash^{a_{0}}_{c,e,\gam_{0},b_{1}}\Gamma,\lnot C^{(\rho)}$
and
$(\mathcal{H}_{\gamma},\Theta,{\tt Q})\vdash^{a_{0}}_{c,e,\gam_{0},b_{1}}C^{(\rho)},\Gamma$
with $\mbox{{\rm rk}}(C)<c$.

\item[$(\Sigma\mbox{{\rm -rfl}}(\Ome))$]
There exist a cap $\rho\in{\tt Q}$, ordinals
$a_{\ell}, a_{r}<a$, and an uncapped formula $C\in\Sigma(\Ome)$ 
such that 
$c\geq \Ome$,
$(\mathcal{H}_{\gamma},\Theta,{\tt Q}
)\vdash^{a_{\ell}}_{c,e,\gam_{0},b_{1}}\Gamma,C^{(\rho)}$
and
$(\mathcal{H}_{\gamma},\Theta,{\tt Q}
)\vdash^{a_{r}}_{c,e,\gam_{0},b_{1}}
\lnot \left(\exists x<\pi\,C^{(x,\Ome)}\right)^{(\rho)}, \Gamma$.

\item[$({\rm rfl}(\rho,d,f,b_{1}))$]
There exist a cap $\rho\in{\tt Q}$ such that
$\Tht\subset M_{\rho}$,
ordinals $d\in {\rm supp}(m(\rho))$, and $a_{0}<a$,
a special finite function $f$,
and a
finite set $\Delta$ of uncapped formulas enjoying the following conditions.

\begin{enumerate}

\item[(r0)]
$\rho< e$ if $s(\rho)>\mathbb{S}$.

\item[(r1)]
$\Delta\subset\bigvee(d):=\{\delta: \mbox{{\rm rk}}(\delta)<d, \delta \mbox{ {\rm is a}}
\bigvee\mbox{{\rm -formula}}\}\cup\{\delta: \mbox{{\rm rk}}(\delta)<\mS\}$.

\item[(r2)]
For the special finite function  $g=m(\rho)$, $s(f)\leq b_{1}$,
$SC_{\mK}(f, g)\subset\mathcal{H}_{\gamma_{0}}[\Theta^{(\rho)}]$ and 
$
f_{d}=g_{d} \,\&\, f^{d}<^{d} g^{\prime}(d)
$.

 \item[(r3)]
For each $\delta\in\Delta$,
 $(\mathcal{H}_{\gamma},\Theta,{\tt Q}
 )\vdash^{a_{0}}_{c,e,\gam_{0},b_{1}}\Gamma, \lnot\delta^{(\rho)}$.

\item[(r4)]
$
(\mathcal{H}_{\gamma},\Theta\cup\{\sig\},
{\tt Q}\cup\{\sigma\}
)\vdash^{a_{0}}_{c,e,\gam_{0},b_{1}}\Gamma, 
\Delta^{(\sigma)}$ holds
for every 
$\sigma\in H_{\rho}(f,b_{1},\gamma_{0},\Theta^{(\rho)})$.

\end{enumerate}
{\small
\[
\hspace{-20mm}
\infer[({\rm rfl}(\rho,d,f,b_{1}))]{(\mathcal{H}_{\gamma},\Theta,{\tt Q})\vdash^{a}_{c,e}\Gamma}
{
\{(\mathcal{H}_{\gamma},\Theta,{\tt Q}
 )\vdash^{a_{0}}_{c,e}\Gamma, \lnot\delta^{(\rho)}\}_{\del\in\Del}
 &
\{
(\mathcal{H}_{\gamma},\Theta\cup\{\sig\},
{\tt Q}\cup\{\sigma\}
)\vdash^{a_{0}}_{c,e}\Gamma, 
\Delta^{(\sigma)}
\}_{\sigma\in H_{\rho}(f,b_{1},\gamma_{0},\Theta^{(\rho)})}
}
\]
}
Note that $(\Tht\cup\{\sig\})_{{\tt Q}\cup\{\sig\}}=\Tht_{{\tt Q}\cup\{\sig\}}=\Tht_{{\tt Q}}$ by 
$\Tht^{(\rho)}\subset M_{\sig}$ and $\rho\in{\tt Q}$.
\end{description}
}
\edf

$\{e\}\cup{\tt Q}\subset\calh_{\gam}[\Tht]$ need not to hold.

Suppose $(\mathcal{H}_{\gamma},\Theta,{\tt Q})\vdash^{a}_{c,e} \Gamma$ holds with 
$A^{(\rho)}\in\Gam$ and $\rho\in{\tt Q}$.
By (\ref{eq:controlderpi11cap}) 
we have
${\sf k}(A)\subset\mathcal{H}_{\gamma}[\Theta^{(\rho)}]$.
We obtain $\sfk(A)\subset M_{\rho}$
by Proposition \ref{prp:EK2pi11}.

In this subsection the ordinals $\gam_{0}$ and $b_{1}$ will be fixed, and
we write $\vdash^{a}_{c,e}$ for $\vdash^{a}_{c,e,\gam_{0},b_{1}}$.

\bprp\label{lem:tautologycap}{\rm (Tautology)}
Let $\{\gam\}\cup\sfk^{\mS}(A)\subset\calh_{\gam}[\Tht_{{\tt Q}}]$ and
$\sigma\in{\tt Q}$, $\sfk(A)\subset\calh_{\gam}[\Tht^{(\sig)}]$.
Then
$(\mathcal{H}_{\gamma},\Tht,{\tt Q}
)\vdash^{2d}_{0,0}
\lnot A^{(\sigma)}, A^{(\sigma)}$ holds for $d=\max\{\mS,\mbox{{\rm rk}}(A)\}$.
\eprp
\bprf By induction on $d$.
Let $A\simeq\bigvee(A_{\iota})_{\iota\in J}$ with $\rk(A)\geq\mS$. 
For $\iota\in [{\tt Q}]_{A^{(\sig)}}J\subset [\sig]J$,
let $d_{\iota}=0$ if $\rk(A_{\iota})<\mS$.
Otherwise $d_{\iota}=\max\{\mS,\rk(A_{\iota})\}$.
In each case we have $d_{\iota}<d$.
IH yields 
\[
\infer[(\bigwedge)]{(\mathcal{H}_{\gamma},\Tht,{\tt Q})\vdash^{2d}_{0,0}
\lnot A^{(\sigma)}, A^{(\sigma)}}
{
 \infer[(\bigvee)]{(\mathcal{H}_{\gamma},\Tht\cup{\sf k}(\iota),{\tt Q})\vdash^{2d_{\iota}+1}_{0,0}\lnot A_{\iota}^{(\sigma)}, A^{(\sigma)}
 }
 {(\mathcal{H}_{\gamma},\Tht\cup{\sf k}(\iota),{\tt Q})\vdash^{2d_{\iota}}_{0,0}\lnot A_{\iota}^{(\sigma)}, A_{\iota}^{(\sigma)}
 }
}
\]
\eprf

\bprp\label{lem:inversionpi11cap}{\rm (Inversion)}
Let
$A\simeq \bigwedge(A_{\iota})_{\iota\in J}$ with $A^{(\rho)}\in\Gam$ and $\rk(A)\geq\mS$, 
$\iota\in [{\tt Q}]_{A^{(\rho)}}J$ 
with $\rho\in{\tt Q}$ and
$(\mathcal{H}_{\gam},\Theta,{\tt Q})\vdash^{a}_{c,e}\Gamma$.
Then
$(\mathcal{H}_{\gam},\Theta\cup\sfk(\iota),{\tt Q}
)
\vdash^{a}_{c,e}
\Gamma,A_{\iota}$.
\eprp

\bprp\label{lem:predceregcap}{\rm (Cut-elimination)}
Let $(\mathcal{H}_{\gamma},\Theta,{\tt Q})\vdash^{a}_{c+d,e}\Gamma$ with
$\calh_{\gam}[\Tht_{{\tt Q}}]\ni c\geq\mS$.
Then $(\mathcal{H}_{\gamma},\Theta,{\tt Q})\vdash^{\vphi_{d}(a)}_{c,e}\Gamma$.
\eprp
\bprf
By main induction on $d$ with subsidiary induction on $a$ using
an analogue to
Reduction \ref{prp:reductionpi11} with (\ref{eq:controlderpi11picap}).
Note that $\rk(C)\in\calh_{\gam}[\Tht_{{\tt Q}}]$ when
$\rk(C)\geq\mS$ and $(\mathcal{H}_{\gamma},\Theta,{\tt Q})\vdash^{a}_{c,e}\Gamma,C$.
\eprf

\blem\label{lem:capping}{\rm (Capping)}
Let
$\Gam\cup\Pi\subset\Del_{0}(\mK)$ with $\Pi=\bigcup\{\Pi_{\sig}: \sig\in{\tt Q}_{\Pi}\}$.
Suppose
$
(\mathcal{H}_{\gamma},\Theta;{\tt Q}_{\Pi}
)\vdash^{* a}_{c,\gam_{0}}\Gamma;\Pi^{[\cdot]}
$ for $a,c<\mK$ and $\Pi^{[\cdot]}=\bigcup\{\Pi_{\sig}^{[\sig/\mS]}: \sig\in{\tt Q}_{\Pi}\}$.
Let
$\rho=\psi_{\mathbb{S}}^{g}(\gamma_{1})$
be an ordinal such that 
${\tt Q}_{\Pi}\subset\rho$,
\beqn\label{eq:capping}
\Tht\subset M_{\rho}
\eeqn
and
$g=m(\rho)$ a special finite function such that
${\rm supp}(g)=\{c\}$ and $g(c)=\alpha_{0}+\mathbb{K}$, where
$\mathbb{K}(2a+1)\leq\alpha_{0}+\mathbb{K}\leq\gamma_{0}\leq\gam_{1}$
with 
$\{\gamma_{1},c,\alpha_{0}\}\subset\mathcal{H}_{\gam}[\Theta]\cap\calh_{\gam_{0}}$,
and ${\tt p}_{0}(\sig)\leq{\tt p}_{0}(\rho)=\gam_{1}$ for each $\sig\in{\tt Q}_{\Pi}$.
Let
$\widehat{\Gam}=\bigcup\{A^{(\rho)}: A\in\Gam\}$,
$\widehat{\Pi}=\bigcup\{\Pi_{\sigma}^{(\sigma)}: \sigma\in{\tt Q}_{\Pi}\}$ and
${\tt Q}= {\tt Q}_{\Pi}\cup\{\rho\}$.

Then 
$(\mathcal{H}_{\gam},\Theta_{\Pi}, {\tt Q})
\vdash^{a}_{c,\rho+1,\gam_{0},c}
\widehat{\Gamma},\widehat{\Pi}$ holds holds for
$\Tht_{\Pi}=\Tht\cup{\tt Q}_{\Pi}$.

\elem
\bprf
By induction on $a$. 
Let us 
write $\vdash^{a}_{c}$ for $\vdash^{a}_{c,\rho+1,\gam_{0},c}$ in the proof.
By assumptions we have
$\Tht \subset M_{\rho}$ and ${\tt Q}_{\Pi}\subset\rho$.
Hence $\Tht=\Tht^{(\rho)}$ and $\Tht_{{\tt Q}_{\Pi}}=\Tht_{{\tt Q}}$.
On the other hand we have
$\sfk(\Gam)\subset\calh_{\gam}[\Tht]$ and for $\sig\in{\tt Q}_{\Pi}$,
$\sfk(\Pi_{\sig})\subset\calh_{\gam}[\Tht^{(\sig)}]$ by (\ref{eq:controlderpi11}).
Therefore (\ref{eq:controlderpi11cap}) is enjoyed.
We have $\{\gam,a,c\}\subset\calh_{\gam}[\Tht_{{\tt Q}_{\Pi}}]$ 
by
(\ref{eq:controlderpi11ac}).
Hence (\ref{eq:controlderpi11picap}) is enjoyed.
Moreover we have 
$SC_{\mK}(g)\subset \calh_{\gam}[\Tht]\subset M_{\rho}$.
\\
\textbf{Case 1}. 
First consider the case when the last inference is a $({\rm stbl})$: 
{\small
\[
\hspace{-5mm}
\infer[({\rm stbl})]{(\mathcal{H}_{\gamma},\Theta;{\tt Q}_{\Pi}
)
\vdash^{* a}_{c}
\Gamma;\Pi^{[\cdot]}}
{
(\mathcal{H}_{\gamma},\Theta;{\tt Q}_{\Pi}
)\vdash^{* a_{0}}_{c}
\Gamma, B(u);\Pi^{[\cdot]}
&
\{
(\mathcal{H}_{\gamma},\Theta\cup\{\sig\};{\tt Q}_{\Pi}\cup\{\sig\}
)
\vdash^{* a_{0}}_{c}
\Gamma; \lnot B(u)^{[\sigma/\mathbb{S}]},\Pi^{[\cdot]}
\}_{\Tht\subset M_{\sigma}}
}
\]
}
Note that it may be  the formula $B(u)^{[\sigma/\mathbb{S}]}$ is in $\Gam$, 
cf.\,Embedding \ref{th:embedpi11}.
$\sig$ in $\Tht\cup\{\sig\}$ ensures us 
$\sfk(B(u)^{[\sigma/\mathbb{S}]})\subset\calh_{\gam}[\Tht\cup\{\sig\}]$ in (\ref{eq:controlderpi11}).
This explains the additional set ${\tt Q}_{\Pi}$ in 
$
(\mathcal{H}_{\gamma},\Theta_{\Pi},{\tt Q}
)\vdash^{a}_{c}
\widehat{\Gamma},\widehat{\Pi}
$, and the addition would be an obstacle to $a\in\Theta_{{\tt Q}}$ in (\ref{eq:controlderpi11picap}).

We have 
an ordinal $a_{0}<a$, 
a $\bigwedge$-formula
$B(0)\in\Delta_{0}(\mathbb{S})$, and
a term $u\in Tm(\mathbb{K})$
such that $\mS\leq{\rm rk}(B(u))<c$.
We have
$
(\mathcal{H}_{\gamma},\Theta;{\tt Q}_{\Pi}
)\vdash^{* a_{0}}_{c}
\Gamma, B(u);\Pi^{[\cdot]}
$.
$
(\mathcal{H}_{\gamma},\Theta_{\Pi},{\tt Q}
)\vdash^{a_{0}}_{c}
\widehat{\Gamma},\left(B(u)\right)^{(\rho)},\widehat{\Pi}
$ follows from IH.

On the other hand we have
$
(\mathcal{H}_{\gamma},\Theta\cup\{\sig\};{\tt Q}_{\Pi}\cup\{\sig\}
)\vdash^{* a_{0}}_{c}
\Gamma; \lnot B(u)^{[\sigma/\mathbb{S}]},\Pi^{[\cdot]}
$
for every ordinal $\sigma$ such that
$\Tht\subset M_{\sigma}$.

Let $h$ be a special finite function such that ${\rm supp}(h)=\{c\}$ and
$h(c)=\mathbb{K}(2a_{0}+1)$.
Then $h_{c}=g_{c}=\emptyset$ and $h^{c}<^{c}g^{\prime}(c)$
by $h(c)=\mathbb{K}(2a_{0}+1)<\mathbb{K}(2a)\leq\alpha_{0}=g^{\prime}(c)$.
Let $\sigma\in H_{\rho}(h,c,\gamma_{0},\Theta)$.
For example $\sig=\psi_{\rho}^{h}(\gam_{1}+\eta)$ with $\eta=\max(\{1\}\cup E_{\mS}(\Tht))$, where
$E_{\mS}(\Tht)=\bigcup_{\alp\in\Tht}E_{\mS}(\alp)$ with the set
$E_{\mS}(\alp)$ of subterms$<\mS$ of $\alp$.
We obtain $\Tht\subset \calh_{\gam_{1}}(\sig)=M_{\sig}$ by $\Tht\subset M_{\rho}$, and
$\{\gam_{1},c,a_{0}\}\subset\mathcal{H}_{\gam}[\Theta]\subset
\calh_{\gam_{1}}(\sig)$.

We have 
$\sfk^{\mS}(B(u))={\sf k}(B(u))\subset \mathcal{H}_{\gamma}[\Theta_{{\tt Q}}]\subset M_{\sigma}$ for (\ref{eq:controlderpi11picap}),
and
$
(\mathcal{H}_{\gamma},\Theta_{\Pi}\cup\{\sig\},
{\tt Q}\cup\{\sigma\}
)
\vdash^{a_{0}}_{c}
\widehat{\Gamma},\lnot B(u)^{(\sigma)},\widehat{\Pi}
$ follows from IH with $\sig\in M_{\rho}$.
Since this holds for every such $\sigma$,
we obtain
$
(\mathcal{H}_{\gamma},\Theta_{\Pi},{\tt Q})
\vdash^{a}_{c,\rho+1}\widehat{\Gamma},\widehat{\Pi}
$
by an inference $({\rm rfl}(\rho,c,h,c))$ with ${\rm rk}(B(u))<c\in {\rm supp}(m(\rho))$.
In the following figure let us omit the operator $\calh_{\gam}$.
{\small
\[
\infer[({\rm rfl}(\rho,c,h,c))]{(\Theta_{\Pi},{\tt Q})\vdash^{a}_{c}
\widehat{\Gamma},\widehat{\Pi}}
{
(\Theta_{\Pi},{\tt Q}
)\vdash^{a_{0}}_{c}
\widehat{\Gamma},B(u)^{(\rho)},\widehat{\Pi}
 &
\{
(\Theta_{\Pi}\cup\{\sig\},
{\tt Q}\cup\{\sigma\}
)\vdash^{a_{0}}_{c}
\widehat{\Gamma},\lnot B(u)^{(\sigma)},\widehat{\Pi}
\}_{\sig}
}
\]
}
\noindent
\textbf{Case 2}.  
Second the last inference introduces
a $\bigvee$-formula $A$.
\\
\textbf{Case 2.1}.
First let $A\in\Gamma$ be introduced by a $(\bigvee)$,
and $A\simeq\bigvee\left(A_{\iota}\right)_{\iota\in J}$.
There are an $\iota\in J$ an ordinal
 $a(\iota)<a$
such that
$(\mathcal{H}_{\gamma},\Theta;{\tt Q}_{\Pi})\vdash^{* a(\iota)}_{c}
\Gamma,A_{\iota};\Pi^{[\cdot]}$.
Let ${\sf k}(\iota)\subset{\sf k}(A_{\iota})$.
We obtain 
${\sf k}(\iota)\subset\mathcal{H}_{\gamma}[\Theta]\subset 
M_{\rho}$ by (\ref{eq:controlderpi11}), $\Theta\subset M_{\rho}$ and
$\gam\leq\gam_{0}\leq\gam_{1}$.
Hence $\iota\in[\rho]J$.
IH yields $(\mathcal{H}_{\gam},\Tht_{\Pi},{\tt Q})
\vdash^{a(\iota)}_{c}
\widehat{\Pi},\widehat{\Gamma},
\left(A_{\iota} \right)^{(\rho)}$.
$(\mathcal{H}_{\gam},\Tht_{\Pi},{\tt Q})
\vdash^{a}_{c}
\widehat{\Pi},\widehat{\Gamma}$ 
follows from a $(\bigvee)$.
\\
\textbf{Case 2.2}.
Second 
$A\equiv B^{[\sigma/\mathbb{S}]}\in\Pi^{[\cdot]}$ is introduced by a $(\bigvee)^{[\cdot]}$
with $B^{(\sigma)}\in\widehat{\Pi}$ and $\sigma\in {\tt Q}_{\Pi}$.
Let $B\simeq\bigvee\left(B_{\iota}\right)_{\iota\in J}$.
Then $A\simeq\bigvee\left(B_{\iota}^{[\sigma,\mathbb{S}]}\right)_{\iota\in [\sigma]J}$ by
Proposition \ref{lem:assigncollaps}.
There are an $\iota\in [\sig]J$ and an ordinal
 $a(\iota)<a$ 
such that
$(\mathcal{H}_{\gamma},\Theta;{\tt Q}_{\Pi})\vdash^{a(\iota)}_{c}\Gamma;
B_{\iota}^{[\sig/\mS]},\Pi^{[\cdot]}$ for $A_{\iota}\equiv B_{\iota}^{[\sigma/\mathbb{S}]}$.
IH yields $(\mathcal{H}_{\gam},\Tht_{\Pi},{\tt Q})
\vdash^{a(\iota)}_{c}
\widehat{\Pi},\widehat{\Gamma},\left(B_{\iota}\right)^{(\sigma)}$.
We obtain
$(\mathcal{H}_{\gam},\Tht_{\Pi},{\tt Q})
\vdash^{a}_{c}
\widehat{\Pi},\widehat{\Gamma}$ 
by a $(\bigvee)$.
\\
\textbf{Case 3}.
Third the last inference introduces a $\bigwedge$-formula $A$.
\\
\textbf{Case 3.1}.
First let $A\in\Gamma$ be introduced by a $(\bigwedge)$, and
$A\simeq\bigwedge\left(A_{\iota}\right)_{\iota\in J}$.
For every
$\iota\in [{\tt Q}_{\Pi}]_{A} J$ there exists an $a(\iota)<a$ 
such that
$(\mathcal{H}_{\gamma},\Theta\cup{\sf k}(\iota);{\tt Q}_{\Pi}
)
\vdash^{* a(\iota)}_{c}\Gamma, A_{\iota};\Pi^{[\cdot]}$.

IH yields 
$(\mathcal{H}_{\gam},\Tht_{\Pi}\cup\sfk(\iota),{\tt Q}
)
\vdash^{a(\iota)}_{c}
\widehat{\Pi},\widehat{\Gamma},
\left(A_{\iota}\right)^{(\rho)}$
for each $\iota\in [{\tt Q}]_{A^{(\rho)}}J\subset [{\tt Q}_{\Pi}]_{A}J$, where $\sfk(\iota)\subset M_{\rho}$.
We obtain
$(\mathcal{H}_{\gam},\Tht_{\Pi},{\tt Q})\vdash^{a}_{c}
\widehat{\Pi},\widehat{\Gamma}$ by a $(\bigwedge)$.
\\
\textbf{Case 3.2}.
Second 
$A\equiv B^{[\sigma/\mathbb{S}]}\in\Pi^{[\cdot]}$ is introduced by a $(\bigwedge)^{[\cdot]}$
with $B^{(\sigma)}\in\widehat{\Pi}$ and $\sigma\in {\tt Q}_{\Pi}$.
Let $B\simeq\bigwedge\left(B_{\iota}\right)_{\iota\in J}$
with $A\simeq\bigwedge\left(B_{\iota}^{[\sigma/\mathbb{S}]}\right)_{\iota\in [\sigma]J}$.
For each $\iota\in [{\tt Q}_{\Pi}]_{B}J\cap [\sig]J$
there is an ordinal
 $a(\iota)<a$ 
such that
$(\mathcal{H}_{\gamma},\Theta\cup\sfk(\iota);{\tt Q}_{\Pi}
)
\vdash^{* a(\iota)}_{c}
\Gamma; A_{\iota},\Pi^{[\cdot]}$ for $A_{\iota}\equiv B_{\iota}^{[\sigma/\mathbb{S}]}$.
IH yields $(\mathcal{H}_{\gam},\Tht_{\Pi}\cup\sfk(\iota),{\tt Q}
)
\vdash^{a(\iota)}_{c}
\widehat{\Pi},\widehat{\Gamma},\left(B_{\iota}\right)^{(\sigma)}$
for each $\iota\in[{\tt Q}]_{B^{(\sig)}}J\subset [{\tt Q}_{\Pi}]_{B}J\cap [\sig]J$, where
$\sfk(\iota)\subset M_{\sig}\subset M_{\rho}$.
$(\mathcal{H}_{\gam},\Tht_{\Pi},{\tt Q})
\vdash^{a}_{c}
\widehat{\Pi},\widehat{\Gamma}$ 
follows from a $(\bigwedge)$.

The other cases
$(cut)$ or $(\Sigma\mbox{{\rm -rfl}})$ on $\Ome$ are seen from IH.
\eprf

\subsection{Eliminations of inferences (rfl)}\label{subsec:elimrfl}

In this subsection, $({\rm rfl}(\rho,c,\gamma))$ 
are removed from operator controlled derivations of $\Sigma_{1}$-sentences 
$\theta^{L_{\Omega}}$ over $\Omega$.

\bdf\label{df:hstepdownpi11}
{\rm
For a special finite function $g$ and
ordinals $a<\mathbb{K}$, $b< c_{\max}=\max({\rm supp}(g))<\mathbb{K}$, 
let us define a special finite function $h=h^{b}(g;a)$ as follows.
$\max({\rm supp}(h))=b$, and
$h_{b}=g_{b}$.
To define $h(b)$,
let $\{b=b_{0}<b_{1}<\cdots<b_{n}=c_{\max}\}=\{b,c_{\max}\}\cup\left((b,c_{\max})\cap {\rm supp}(g)\right)$.
Define recursively ordinals $\alpha_{i}$ by
$\alpha_{n}=\alpha+a$ with $g(c_{\max})=\alpha+\mathbb{K}$.
$\alpha_{i}=g(b_{i})+\tilde{\theta}_{c_{i}}(\alpha_{i+1})$ for
$c_{i}=b_{i+1}-b_{i}$.
Finally put $h(b)=\alpha_{0}+\mathbb{K}$.

}
\edf

\bprp\label{prp:hstepdown}
Let $f$ and $g$ be special finite functions with $c_{\max}=\max({\rm supp}(g))$.

\begin{enumerate}

\item\label{prp:hstepdown.1}
Let $b<e<c_{\max}$ and $a_{0},a_{1}<a$.
Then
$h^{b}(h^{e}(g;a_{0});a_{1})\leq (h^{b}(g;a))^{\prime}$.

\item\label{prp:hstepdown.2}
Suppose
$f<^{d}g^{\prime}(d)$ for a $d\in{\rm supp}(g)$.
Let $b<d$.
Then $f_{b}=(h^{b}(g;a))_{b}$ and
$f<^{b}(h^{b}(g;a))^{\prime}(b)$.

\end{enumerate}
\eprp

Recall that $s(\rho)=\max({\rm supp}(m(\rho)))$.

\blem\label{lem:recappingpi11}{\rm (Recapping)}\\
Let $
(\mathcal{H}_{\gamma},\Theta,{\tt Q})
\vdash^{a}_{c_{1},e,\gamma_{0},b_{2}}
\Pi,\widehat{\Gamma}
$
for a finite family ${\tt Q}$ for $\gam_{0}, b_{2}$, 
${\tt Q}^{t}\subset{\tt Q}$, $\fal\rho\in{\tt Q}^{t}(s(\rho)>\mS)$
and
${\tt Q}^{f}={\tt Q}\setm{\tt Q}^{t}$, 
$\Gamma\cup\Pi\subset\Delta_{0}(\mathbb{K})$,
$\widehat{\Gam}=\bigcup\{\Gam_{\rho}^{(\rho)} :\rho\in{\tt Q}^{t}\}$,
where
each $\tht\in \Gam$ is either a $\bigvee$-formula or $\rk(\tht)<\mS$,
and
$\Pi$ a set of formulas such that $\tau\in{\tt Q}^{f}$ for every $A^{(\tau)}\in\Pi$.

Let $\max\{s(\rho):\rho\in{\tt Q}^{t}\}\leq b_{1}$.
For each $\rho\in{\tt Q}^{t}$, let $\mS\leq b^{(\rho)}\in\calh_{\gam}[\Tht^{(\rho)}]$ 
with $\rk(\Gamma_{\rho})<b^{(\rho)}<s(\rho)$, and
$\kappa(\rho)\in H_{\rho}(h^{b^{(\rho)}}(m(\rho);\ome(b_{1},a)),b_{2}, \gamma_{0},\Theta^{(\rho)})$
 with $\ome(b,a)=\ome^{\ome^{b}}a$.
Assume 
$b_{1}\in\calh_{\gam}[\Tht_{{\tt Q}}]$.

Then 
$
(\mathcal{H}_{\gamma},\Theta,{\tt Q}(\kappa)
)
\vdash^{\ome(b_{1},a)}_{c_{b_{1}},e^{\kap}, \gamma_{0},b_{2}}
\Pi,\widehat{\Gamma}_{\kappa}
$
holds,
where ${\tt Q}(\kappa)={\tt Q}^{f}\cup\{\kap(\rho):\rho\in{\tt Q}^{t}\}$,
$c_{b_{1}}=\max\{c_{1},b_{1}\}$,
$e^{\kap}=\max(\{\tau\in{\tt Q}^{f} : s(\tau)>\mS\}\cup\{\kap(\rho): \rho\in{\tt Q}^{t}\})+1$,
$\widehat{\Gamma}_{\kappa}=\bigcup\{\Gam_{\rho}^{(\kap(\rho))}:\rho\in{\tt Q}^{t}\}$.

$e^{\kap}<e$ holds
when
${\tt Q}^{t}=\{\rho\in{\tt Q}: s(\rho)>\mS\}\neq\emptyset$.

\elem
\bprf
We show the lemma
by main induction on $b_{1}$
with subsidiary induction on $a$.
The subscripts $\gam_{0}, b_{2}$ are omitted in the proof.
We obtain $\{\gam,b_{1},a,c_{1}\}\cup\sfk^{\mS}(\Pi,\Gam)\subset\calh_{\gam}[\Tht_{{\tt Q}}]$ by the assumption and (\ref{eq:controlderpi11picap}).
Then 
$\{\gam,\ome(b_{1},a),c_{b_{1}}\}\cup\sfk^{\mS}(\Pi,\Gam)\subset\calh_{\gam}[\Tht_{{\tt Q}(\kap)}]$ since $\Tht^{(\rho)}\subset M_{\kap(\rho)}$ for
each $\rho\in{\tt Q}^{t}$.
Hence (\ref{eq:controlderpi11picap}) is enjoyed in 
$
(\mathcal{H}_{\gamma},\Theta,{\tt Q}(\kappa)
)
\vdash^{\ome(b_{1},a)}_{c_{b_{1}}, e,\gamma_{0},b_{2}}
\Pi,\widehat{\Gamma}_{\kappa}
$.

Let $\rho\in{\tt Q}^{t}$.
We have 
$b^{(\rho)}\in\calh_{\gam}[\Tht^{(\rho)}]$, 
$SC_{\mK}(m(\rho))\subset \calh_{\gam_{0}}[\Tht^{(\rho)}]$ and
$\Theta^{(\rho)}\subset M_{\kappa(\rho)}$.
$SC_{\mK}(h^{b^{(\rho)}}(m(\rho);\ome(b_{1},a)))\subset \calh_{\gam_{0}}[\Tht^{(\rho)}]$
follows.
Moreover we have
$SC_{\mK}(m(\kap(\rho)))\subset\calh_{\gam_{0}}[\Tht^{(\rho)}]\subset M_{\kap(\rho)}$.

Consider the case when the last inference is a $({\rm rfl}(\rho,d,f,b_{2}))$ for a $\rho\in{\tt Q}$.
The case $\rho\in{\tt Q}^{f}$ is seen from SIH.
Assume $\rho\in{\tt Q}^{t}$. Let $b=b^{(\rho)}$,
$g=m(\rho)$, $b_{1}\geq s(\rho)\geq d\in {\rm supp}(g)$, $\kap=\kap(\rho)$,
$\Gam=\Gam_{\rho}$,
$\widehat{\Lam}=\bigcup_{\rho\neq\tau\in{\tt Q}^{t}}\{\Gam_{\tau}^{(\tau)}\}$, and
$\widehat{\Lam}_{\kap}=\bigcup_{\rho\neq\tau\in{\tt Q}^{t}}\{\Gam_{\tau}^{\kap(\tau)}\}$.
We have a sequent $\Delta\subset\bigvee(d)$ such that ${\rm rk}(\Delta)<d\leq s(\rho)\leq b_{1}$ 
and
$\sfk^{\mS}(\Del)\subset\mathcal{H}_{\gamma}[\Theta_{{\tt Q}}]\subset M_{{\tt Q}}$ by (\ref{eq:controlderpi11picap})
and $\sfk^{\mS}(\Del)\subset M_{{\tt Q}(\kap)}$ by $\Tht_{{\tt Q}}=\Tht_{{\tt Q}(\kap)}$.
There is an ordinal $a_{0}\in\calh_{\gam}[\Tht_{{\tt Q}}]\cap a$ such that 
$
(\mathcal{H}_{\gamma},\Theta,
{\tt Q}
)
\vdash^{a_{0}}_{c_{1},e}
\Pi,\widehat{\Lam},\Gamma^{(\rho)}, 
\lnot\delta^{(\rho)}
$
for each $\delta\in\Delta$.
For each $\delta\in\Delta\subset\bigvee(d)$ with $\rk(\del)\geq\mS$, we have 
$\delta\simeq \bigvee\left(\delta_{\iota}\right)_{\iota\in J}$.
Let $b_{0}=\max(\{\mathbb{S}\}\cup\{\mbox{{\rm rk}}(\delta) :\delta\in\Delta\})$.
Then 
$s(\rho)>b_{0}\in\mathcal{H}_{\gamma}[\Theta_{{\tt Q}}]$.
Inversion \ref{lem:inversionpi11cap} yields for $\rk(\del)\geq\mS$
\begin{equation}\label{eq:Case1api11}
(\mathcal{H}_{\gamma},\Theta\cup{\sf k}(\iota),
{\tt Q}
)
\vdash^{a_{0}}_{c_{1},e}
\Pi,\widehat{\Lam},\Gamma^{(\rho)},\lnot(\delta_{\iota})^{(\rho)}
\end{equation}
for each 
$\iota\in[{\tt Q}]_{\del^{(\rho)}}J$,
where $J\subset Tm(b_{0})$ and $\lnot\delta_{\iota}\in\bigvee(b_{0})$ by
$\mbox{{\rm rk}}(\delta_{\iota})<\mbox{{\rm rk}}(\delta)$.

On the other side for each $\sigma\in H_{\rho}(f,b_{2},\gamma_{0},\Theta^{(\rho)})$
\begin{equation}\label{eq:Case1bpi11}
(\mathcal{H}_{\gamma},\Theta\cup\{\sigma\},
{\tt Q}\cup\{\sigma\}
)
\vdash^{a_{0}}_{c_{1},e}
\Pi,\widehat{\Lam},\Gamma^{(\rho)},  \Delta^{(\sigma)}
\end{equation}
$f$ is a special finite function such that $s(f)\leq b_{2}$,
$f_{d}=g_{d}$, $f^{d}<^{d}g^{\prime}(d)$ and
$SC_{\mK}(f)\subset\mathcal{H}_{\gamma_{0}}[\Theta^{(\rho)}]$.
Let $({\tt Q}\cup\{\sig\})^{f}={\tt Q}^{f}\cup\{\sig\}$.
\\
\textbf{Case 1}. 
$b_{0}<b$:
Then ${\rm rk}(\Delta)<b$.
Let $\rk(\del)\geq\mS$.
From (\ref{eq:Case1api11}) 
we obtain by SIH with $b>b_{0}\geq\mathbb{S}$,
$(\mathcal{H}_{\gamma},\Theta\cup{\sf k}(\iota),
{\tt Q}(\kappa)
)
\vdash^{\ome(b_{1},a_{0})}_{c_{b_{1}},e^{\kap}}
\Pi,\widehat{\Lam}_{\kap},\Gamma^{(\kappa)},\lnot(\delta_{\iota})^{(\kappa)}$
for each $\iota\in [{\tt Q}(\kap)]_{\del^{(\kap)}}J\subset [{\tt Q}]_{\del^{(\rho)}}J$.
An inference $(\bigwedge)$ yields
\begin{equation}\label{eq:Case1cpi11}
(\mathcal{H}_{\gamma},\Theta,
{\tt Q}(\kappa)
)
\vdash^{\ome(b_{1},a_{0})+1}_{c_{b_{1}},e^{\kap}}
\Pi,\widehat{\Lam}_{\kap},\Gamma^{(\kappa)},\lnot\delta^{(\kappa)}
\end{equation}
Moreover SIH yields (\ref{eq:Case1cpi11}) for $\rk(\del)<\mS$.
Let 
$d_{1}=\min\{b,d\}$.
Then $\Delta\subset\bigvee(d_{1})$ by $b> b_{0}$.

We claim for the special finite function $h=h^{b}(g;\ome(b_{1},a))$ that
\begin{equation}\label{eq:zigpi11}
f_{d_{1}}=h_{d_{1}} \,\&\,
f^{d_{1}}<^{d_{1}}h^{\prime}(d_{1})
\end{equation}
If $d_{1}=d\leq b$, then $h_{d}=g_{d}$ and $g^{\prime}(d)=g(d)\leq h^{\prime}(d)$.
Proposition \ref{prp:idless} yields the claim.
If $d_{1}=b< d$, then Proposition \ref{prp:hstepdown}.\ref{prp:hstepdown.2}
yields the claim.

On the other hand, for each 
$\sigma\in H_{\kap}(f,b_{2}, \gamma_{0},\Theta^{(\rho)})\subset H_{\rho}(f,b_{2},\gamma_{0},\Theta^{(\rho)})$
we have by (\ref{eq:Case1bpi11}) and SIH,
\begin{equation}\label{eq:Case1dpi11}
(\mathcal{H}_{\gamma},\Theta\cup\{\sigma\},
{\tt Q}(\kappa)\cup\{\sigma\}
)
\vdash^{\ome(b_{1},a_{0})}_{c_{b_{1}},e^{\kap}}
\Pi,\widehat{\Lam}_{\kap},\Gamma^{(\kappa)},  \Delta^{(\sigma)}
\end{equation}
We have $\kap=\kap(\rho)<\kap(\rho)+1\leq e^{\kap}$ for (r0).
An inference $({\rm rfl}(\kappa,d_{1},f,b_{2}))$ with (\ref{eq:zigpi11}), (\ref{eq:Case1cpi11}) and (\ref{eq:Case1dpi11}) yields
$
(\mathcal{H}_{\gamma},\Theta,
{\tt Q}(\kappa)
)
\vdash^{\ome(b_{1},a)}_{c_{b_{1}},e^{\kap}}
\Pi,\widehat{\Lam}_{\kap},\Gamma^{(\kappa)}
$, where $d_{1}\in {\rm supp}(m(\kappa))$ and $\sfk^{\mS}(\Del)\subset\calh_{\gam}[\Tht_{{\tt Q}(\kap)}]$.
\\
\textbf{Case 2}. $b\leq b_{0}$:
When $b=b_{0}$, let $\tau=\kappa$.
When $b<b_{0}$, let
$\tau\in H_{\rho}(h,b_{2}, \gamma_{0},\Theta^{(\rho)})$
be such that
$\kappa<\tau$ and  $m(\tau)=h=h^{b_{0}}(g;a_{1})$
with $a_{1}=\ome(b_{1},a_{0})+1$.

Let $\sigma\in H_{\tau}(f,b_{2},\gamma_{0},\Theta^{(\rho)})$.
SIH with (\ref{eq:Case1bpi11}) and $b_{0}<s(\rho)$
yields 
\begin{equation}\label{eq:Case1epi11}
(\mathcal{H}_{\gamma},\Theta\cup\{\sig\},
{\tt Q}_{\tau}\cup\{\sigma\}
)
\vdash^{\ome(b_{1},a_{0})}_{c_{b_{1}},e^{\tau}}
\Delta^{(\sigma)},\Pi,\widehat{\Lam}_{\kap},\Gamma^{(\tau)}
\end{equation}
where
${\tt Q}_{\tau}={\tt Q}^{f}\cup\{\kap(\lam):\rho\neq\lam\in{\tt Q}^{t}\}\cup\{\tau\}$, and
$e^{\tau}=\max(\{\lam\in{\tt Q}^{f} : s(\lam)>\mS\}\cup\{\kap(\lam): \rho\neq\lam\in{\tt Q}^{t}\}\cup\{\tau\})+1$.
Let
$\sig\in R
:=\{\sigma\in  H_{\tau}(f,b_{2}, \gamma_{0},\Theta^{(\rho)}): 
(m(\sigma))^{\prime}\geq (h^{b_{0}}(g;\ome(b_{1},a_{0})))^{\prime}
\}
$.
We see $\sig\in H_{\rho}(h^{b_{0}}(g;\ome(b_{1},a_{0})),b_{2},\gamma_{0},\Theta^{(\rho)})$.
Moreover $\rk(\lnot\del_{\iota})<b_{0}$ if $\rk(\del)\geq\mS$, and
$\rk(\lnot\del)<b_{0}$ if $\rk(\del)<\mS\leq b_{0}$.

For each $\iota\in[{\tt Q}]_{\del^{(\rho)}}J$ and $\rk(\del)\geq\mS$, we obtain
$
(\mathcal{H}_{\gamma},\Theta\cup{\sf k}(\iota),
{\tt Q}_{\sigma}
)
\vdash^{\ome(b_{1},a_{0})}_{c_{b_{1}},e^{\sig}}
\Pi,\widehat{\Lam}_{\kap},\Gamma^{(\sigma)},\lnot(\delta_{\iota})^{(\sigma)}
$ by $\rk(\lnot\del_{\iota})<b_{0}$, SIH and (\ref{eq:Case1api11}), where ${\tt Q}_{\sigma}\cup\{\tau\}={\tt Q}_{\tau}\cup\{\sig\}$.
A $(\bigwedge)$ yields
$
(\mathcal{H}_{\gamma},\Theta,
{\tt Q}_{\sigma}
)
\vdash^{\ome(b_{1},a_{0})+1}_{c_{b_{1}},e^{\sig}}
\Pi,\widehat{\Lam}_{\kap},\Gamma^{(\sigma)},\lnot\delta^{(\sigma)}
$.
When $\rk(\del)<\mS$, this follows from SIH.
Also $M_{{\tt Q}_{\sig}}=M_{{\tt Q}_{\sig}\cup\{\tau\}}$ and
$e^{\sig}\leq e^{\tau}$ by $\tau>\sig$.
Therefore
\begin{equation}\label{eq:Case1fpi11}
(\mathcal{H}_{\gamma},\Theta,
{\tt Q}_{\tau}\cup\{\sigma\}
)
\vdash^{\ome(b_{1},a_{0})+1}_{c_{b_{1}},e^{\tau}}
\Pi,\widehat{\Lam}_{\kap},\Gamma^{(\sigma)},\lnot\delta^{(\sigma)}
\end{equation}
From (\ref{eq:Case1epi11})  and (\ref{eq:Case1fpi11}) 
by several $(cut)$'s of $\delta$ with 
${\rm rk}(\delta)<d\leq b_{1}\leq c_{b_{1}}$
we obtain 
for a $p<\omega$,
\begin{equation} \label{eq:Case1gpi11}
\forall \sigma\in R
\left[
(\mathcal{H}_{\gamma},\Theta\cup\{\sig\},
{\tt Q}_{\tau}\cup\{\sigma\}
)
\vdash^{\ome(b_{1},a_{0})+p}_{c_{b_{1}},e^{\tau}}
\Pi,\widehat{\Lam}_{\kap},\Gamma^{(\sigma)},\Gamma^{(\tau)}
\right]
\end{equation}
On the other hand we have $r=\max\{\mS,\rk(\Gam)\}\leq b<b_{1}$
and $\sfk^{\mS}(\Gam)\subset \calh_{\gam}[\Tht_{{\tt Q}}]=\calh_{\gam}[\Tht_{{\tt Q}_{\tau}}]\subset M_{{\tt Q}_{\tau}}$ by (\ref{eq:controlderpi11picap}), where
$\Tht_{{\tt Q}}=\Tht_{{\tt Q}_{\tau}}$ by $\Tht^{(\rho)}\subset M_{\tau}$.
Tautology \ref{lem:tautologycap} yields
for each $\theta\in\Gamma$
\begin{equation}\label{eq:Case1hpi11}
(\mathcal{H}_{\gamma},\Tht,{\tt Q}_{\tau}
)\vdash^{2r}_{0,0}
\Gamma^{(\tau)},\lnot \theta^{(\tau)}
\end{equation}
Let us define a finite function $h$ by
${\rm supp}(h)={\rm supp}(g_{b_{0}})\cup{\rm supp}(f^{b_{0}+1})\cup\{b_{0}\}$,
$h_{b_{0}}=g_{b_{0}}$ and $h^{b_{0}+1}=f^{b_{0}+1}$.
Let $(h^{b_{0}}(g;\ome(b_{1},a_{0})))(b_{0})=\alpha+\mathbb{K}$.
Then $h(b_{0})=\alpha$ if $f^{b_{0}+1}\neq\emptyset$.
Otherwise $h(b_{0})=\alpha+\mathbb{K}$.
We see that 
$R=H_{\tau}(h,\gamma_{0},\Theta^{(\rho)})$,
and $
h^{b_{0}}<^{b_{0}}
\left(m(\tau)\right)^{\prime}(b_{0})$.

By an inference 
$({\rm rfl}(\tau,b_{0},h,b_{2}))$ 
with 
its resolvent class 
$R=H_{\tau}(h,b_{2},\gamma_{0},\Theta^{(\rho)})$ and $\Gamma\subset\bigvee(b_{0})$ 
we conclude
from (\ref{eq:Case1hpi11})
and (\ref{eq:Case1gpi11})
for $\rk(\Gam)<b\leq b_{0}\leq s(\tau)$
\begin{equation}\label{eq:Case1ipi11}
(\mathcal{H}_{\gamma},\Theta,
{\tt Q}_{\tau}
)
\vdash^{a_{2}}_{c_{b_{1}},e^{\tau}}
\Pi,\widehat{\Lam}_{\kap},\Gamma^{(\tau)}
\end{equation}
where 
$a_{2}=\max\{2r,\ome(b_{1},a_{0})+p\}+1<\ome(b_{1},a)=\ome^{\ome^{b_{1}}} a$.
If $b_{0}=b$, we are done.
In what follows assume $b<b_{0}$.
We have $a_{1}<\ome(b_{1},a)$ and 
$\ome(b_{0},a_{2})=\ome^{\ome^{b_{0}}}a_{2}<\ome(b_{1},a)$ by 
$b_{0}<b_{1}$.
Moreover Proposition \ref{prp:hstepdown}.\ref{prp:hstepdown.1} 
for $m(\tau)=h^{b_{0}}(g;a_{1})$ yields
$(h^{b}(m(\tau);\ome(b_{0},a_{2})))^{\prime}=
(h^{b}(h^{b_{0}}(g;a_{1});\ome(b_{0},a_{2})))^{\prime}
\leq (h^{b}(g;\ome(b_{1},a)))^{\prime}$.

Let $({\tt Q}_{\tau})^{t}=\{\tau\}$ and $\kap(\tau)=\kap(\rho)=\kap$. Then 
$(e^{\tau})^{\kap}=\max(\{\lam\in({\tt Q}_{\tau})^{f} : s(\lam)>\mS\}\cup\{\kap\})+1=e^{\kap}$.
We have $\sfk^{\mS}(\Gam)\cup\{b_{0}\}\subset\calh_{\gam}[\Tht_{{\tt Q}_{\tau}}]$, 
$\rk(\Gam_{\rho})<b^{(\rho)}=b<b_{0}=s(\tau)<b_{1}$ for $\Gam=\Gam_{\rho}$ and
$b\in\calh_{\gam}[\Tht^{(\tau)}]$,
$\ome(b_{0},a_{2})<\ome(b_{1},a)$ and
$\max\{c_{b_{1}},b_{0}\}=c_{b_{1}}$.
Also
$\kappa\in H_{\rho}(h^{b}(g;\ome(b_{1},a)),b_{2}, \gam_{0},\Theta^{(\rho)})\cap\tau
\subset
H_{\tau}(h^{b}(m(\tau);\ome(b_{1},a_{2})),b_{2}, \gam_{0},\Theta^{(\rho)})$.
MIH with (\ref{eq:Case1ipi11}) yields
$
(\mathcal{H}_{\gamma},\Theta,
{\tt Q}(\kappa)
)
\vdash^{\ome(b_{1},a)}_{c_{b_{1}},e^{\kap}}
\Pi,\Gamma^{(\kappa)}
$.

Second consider the case when the last inference $(\bigvee)$ introduces a $\bigvee$-formula
$B$:
If $B\in\Pi$, SIH yields the lemma.
Assume that
$B\equiv A^{(\rho)}\in\Gamma_{\rho}^{(\rho)}$
with $A\simeq\bigvee\left(A_{\iota}\right)_{\iota\in J}$ and $\rho\in {\tt Q}$.
We may assume $\rho\in{\tt Q}^{t}$.
We have 
$
(\mathcal{H}_{\gamma},\Tht,{\tt Q}
)
\vdash^{a_{0}}_{c_{1},e}
\Pi,\widehat{\Gamma},\left(A_{\iota}\right)^{(\rho)}
$,
where $a_{0}<a$, $\iota\in [\rho]J$.
We claim that
$\iota\in [\kap(\rho)]J$.
We may assume $\sfk(\iota)\subset\sfk(A_{\iota})$.
We have $\sfk(A_{\iota})\subset\calh_{\gam}[\Tht^{(\rho)}]$ by (\ref{eq:controlderpi11cap}).
$\Tht^{(\rho)}\subset M_{\kap(\rho)}$ yields 
$\sfk(A_{\iota})\subset M_{\kap(\rho)}$.

Let
$A_{\iota}\simeq\bigwedge\left(B_{\nu}\right)_{\nu\in I}$
for $\bigvee$-formulas $B_{\nu}$, and assume $\rk(A_{\iota})\geq\mS$.
Inversion \ref{lem:inversionreg} yields for each $\nu\in[{\tt Q}]_{A_{\iota}^{(\rho)}} I$,
$
(\mathcal{H}_{\gamma},\Tht\cup{\sf k}(\nu),
{\tt Q}
)\vdash^{a_{0}}_{c_{1},e}
\Pi,\widehat{\Gamma},
\left(B_{\nu}\right)^{(\rho)}
$.

SIH yields for each $\nu\in[{\tt Q}(\kap)]_{A_{\iota}^{(\rho)}} I\subset[{\tt Q}]_{A_{\iota}^{(\rho)}} I$ that
$
(\mathcal{H}_{\gamma},\Tht\cup{\sf k}(\nu),
{\tt Q}(\kap)
)\vdash^{\ome(b_{1},a_{0})}_{c_{b_{1}},e^{\kap}}
\Pi,\widehat{\Gamma}_{\kap},
\left(B_{\nu}\right)^{(\kappa)}
$.
$
(\mathcal{H}_{\gamma},\Tht,
{\tt Q}(\kap)
)\vdash^{\ome(b_{1},a_{0})+1}_{c_{b_{1}},e^{\kap}}
\Pi,\widehat{\Gamma}_{\kap},
\left(A_{\iota}\right)^{(\kappa)}
$ follows from a $(\bigwedge)$.
An inference $(\bigvee)$ yields
$
(\mathcal{H}_{\gamma},\Tht,
{\tt Q}(\kap)
)\vdash^{\ome(b_{1},a)}_{c_{b_{1}},e^{\kap}}
\Pi,\widehat{\Gamma}_{\kap}
$.

Other cases are seen from SIH.
\eprf
\\

For $c\leq\mS$,
$(\calh_{\gam},\Tht)\vdash^{* a}_{c}\Gam$ denotes 
$(\calh_{\gam},\Tht;\emptyset)\vdash^{*a}_{c}\Gam;\emptyset$.
Since $\Tht_{\emptyset}=\Tht$, (\ref{eq:controlderpi11}) and (\ref{eq:controlderpi11ac}) 
amount to
(\ref{eq:controlderKP}) $\{\gam,a,c\}\cup\sfk(\Gam)\subset\mathcal{H}_{\gamma}[\Theta]$,
and there occurs no inferences $(\bigvee)^{[\cdot]}$, $(\bigwedge)^{[\cdot]}$ nor 
$({\rm stbl})$.
The inference $(\Sigma\mbox{{\rm -rfl}})$ is only on $\Ome$.
This means that $(\calh_{\gam},\Tht)\vdash^{* a}_{c}\Gam$ is equivalent to
$\calh_{\gam}[\Tht]\vdash^{a}_{c}\Gam$ in
Definition \ref{df:controlderreg}.

\blem\label{lem:sum}{\rm (Elimination of inferences (rfl))}\\
Let ${\tt Q}$ be a finite family for $\gam_{0}$ and $b_{1}\geq\mS$.
Let
 $\max ({\rm rk}(\Gamma))<\mS$,
$\widehat{\Gam}=\bigcup\{\Gam_{\rho}^{(\rho)} :\rho\in{\tt Q}\}$
and
$\Gam=\bigcup\{\Gam_{\rho} :\rho\in{\tt Q}\}$,
where $\sfk(\Gam_{\rho})\subset M_{\rho}$.
Suppose
$
(\mathcal{H}_{\gamma},\Theta,{\tt Q})
\vdash^{a}_{\mS,e,\gamma_{0},b_{1}}
\widehat{\Gamma}
$.

Then
$
(\mathcal{H}_{\gam_{1}},\Theta
)
\vdash^{* \tilde{a}}_{\mS}
\Gamma
$
holds for $\gam_{1}=\gam_{0}+\mS$, 
$\tilde{a}=\vphi_{e}(b_{1}+a)$.

\elem
\bprf 
By main induction on $e$  with subsidiary induction on $a$.
We have $\{e\}\cup{\tt Q}\subset\calh_{\gam_{1}}$ by 
Definitions \ref{df:controldercollapscappi11} and \ref{df:cap.2}, $b_{1}\in\calh_{\gam}[\Tht_{{\tt Q}}]$
by (\ref{eq:controlderpi11picap}), 
and
$\emptyset=\sfk^{\mS}(\Gam)\subset\calh_{\gam}[\Tht_{{\tt Q}}]$.
\\
\textbf{Case 1}.
First let $\{\lnot A^{(\sig)},A^{(\sig)}\}\subset\widehat{\Gam}$ with
$\rk(A)<\mS$ by (Taut).
Then
$(\calh_{0},\sfk(A))\vdash^{* \mS}_{0}\lnot A,A$ by Tautology \ref{lem:tautologypi11}.\ref{lem:tautologypi11.1} and
$
(\mathcal{H}_{\gam_{1}},\Theta
)
\vdash^{* \tilde{a}}_{\mS}
\Gamma
$ by $\tilde{a}>\mS$.
\\
\textbf{Case 2}.
Second consider the case when the last inference is a $({\rm rfl}(\rho,d,f,b_{1}))$ for a $\rho\in{\tt Q}$.
Let ${\tt Q}^{t}=\{\tau\in{\tt Q}: s(\tau)>\mS\}$,
${\tt Q}^{f}={\tt Q}\setm{\tt Q}^{t}$,
and
$\kappa(\tau)\in H_{\tau}(h^{\mS}(m(\tau);\ome(b,a)),b_{1},\gamma_{0},\Theta^{(\tau)})$
for each $\tau\in{\tt Q}^{t}$.
  Let $g=m(\rho)$, $s(\rho)\geq d\in {\rm supp}(g)$, $\kap=\kap(\rho)$ when $\rho\in{\tt Q}^{t}$,
  $\widehat{\Pi}=\bigcup_{\rho\neq\tau\in{\tt Q}^{f}}\Gam_{\tau}^{(\tau)}$,
$\widehat{\Lam}=\bigcup_{\rho\neq\tau\in{\tt Q}^{t}}\Gam_{\tau}^{(\tau)}$, and
$\widehat{\Lam}_{\kap}=\bigcup_{\rho\neq\tau\in{\tt Q}^{t}}\Gam_{\tau}^{\kap(\tau)}$.
We have a sequent $\Delta\subset\bigvee(d)$
and an ordinal $a_{0}<a$ such that 
${\rm rk}(\Delta)<d\leq s(\rho)$ and
$
(\mathcal{H}_{\gamma},\Theta,
{\tt Q}
)
\vdash^{a_{0}}_{\mS,e,\gam_{0},b_{1}}
\widehat{\Pi},\widehat{\Lam},\Gamma_{\rho}^{(\rho)}, 
\lnot\delta^{(\rho)}
$
for each $\delta\in\Delta$.
On the other hand we have 
$
(\mathcal{H}_{\gamma},\Theta\cup\{\sig\},
{\tt Q}\cup\{\sigma\}
)
\vdash^{a_{0}}_{\mS,e,\gam_{0},b_{1}}
\widehat{\Pi},\widehat{\Lam},\Gamma_{\rho}^{(\rho)},  \Delta^{(\sigma)}
$, where
$\sigma\in H_{\rho}(f,b_{1},\gamma_{0},\Theta^{(\rho)})$,
$f$ is a special finite function such that $s(f)\leq b_{1}$,
$f_{d}=g_{d}$, $f^{d}<^{d}g^{\prime}(d)$ and
$SC_{\mK}(f)\subset\mathcal{H}_{\gamma_{0}}[\Theta^{(\rho)}]$.
\\
\textbf{Case 2.1} $s(\rho)\leq\mS$:
We have $\rk(\Del)<d\leq s(\rho)\leq\mS$.
Let $\tilde{a}_{0}=\vphi_{e}(b_{1}+a_{0})$.
By SIH we obtain
$
(\mathcal{H}_{\gam_{1}},\Theta
)
\vdash^{* \tilde{a}_{0}}_{\mS}
\Pi,\Lam,\Gam_{\rho},
\lnot\delta
$
for each $\delta\in\Delta$, and
$
(\mathcal{H}_{\gam_{1}},\Theta\cup\{\sig\}
)
\vdash^{* \tilde{a}_{0}}_{\mS}
\Pi,\Lam,\Gam_{\rho}, \Delta
$, where
$\sig\in\calh_{\gam_{0}+\mS}\subset\calh_{\gam_{1}}[\Tht]$.
Several $(cut)$'s of $\rk(\del)<\mS$ yields
$
(\mathcal{H}_{\gam_{1}},\Theta
)
\vdash^{* \tilde{a}}_{\mS}
\Pi,\Lam,\Gam_{\rho}
$ for $\Gam=\Pi\cup\Lam\cup\Gam_{\rho}$.
 \\
\textbf{Case 2.2}. $s(\rho)>\mS$:
Then $\rho\in{\tt Q}^{t}\neq\emptyset$.
$
(\mathcal{H}_{\gamma},\Theta,
{\tt Q}(\kap)
)
\vdash^{\ome(b_{1},a)}_{b_{1},e^{\kap},\gam_{0}, b_{1}}
\widehat{\Pi},\widehat{\Lam}_{\kap},\Gamma_{\rho}^{(\kap)}
$ follows by Recapping \ref{lem:recappingpi11},
where $b_{1}\geq\mS$ and
$e^{\kap}<e$.
Cut-elimination \ref{lem:predceregcap} yields for
$a_{1}=\vphi_{b_{1}}(\ome(b_{1},a))$,
$
(\mathcal{H}_{\gamma},\Theta,
{\tt Q}(\kap)
)
\vdash^{a_{1}}_{\mS,e^{\kap},\gam_{0}, b_{1}}
\widehat{\Pi},\widehat{\Lam}_{\kap},\Gamma_{\rho}^{(\kap)}
$.
MIH then yields
$
(\mathcal{H}_{\gam_{1}},\Theta
)
\vdash^{* \tilde{a}_{1}}_{\mS}
\Gam
$,
where $\Gam=\Pi\cup\Lam\cup\Gam_{\rho}$ and
$\tilde{a}_{1}=\vphi_{e^{\kap}}(b_{1}+a_{1})<\vphi_{e}(b_{1}+a)=\tilde{a}$ by 
$e^{\kap}<e$
and $a,b_{1}<\tilde{a}$.
\\
\textbf{Case 3}.
The last inference is a $(\bigwedge)$:
We have $a(\iota)<a$, $A^{(\rho)}\in\widehat{\Gam}$ with $A\simeq\bigwedge(A_{\iota})_{\iota\in J}$,
and
$
(\mathcal{H}_{\gamma},\Tht\cup\sfk(\iota),{\tt Q})
\vdash^{a(\iota)}_{\mS,e,\gam_{0}, b_{1}}
\widehat{\Gam},(A_{\iota})^{(\rho)}
$ for each 
$\iota\in [{\tt Q}]_{A^{(\rho)}}J$.
Since $A\in\Del_{0}(\mS)$, we obtain 
$[{\tt Q}]_{A^{(\rho)}}J=[\rho]J=J$.
SIH yields
$
(\mathcal{H}_{\gam_{1}},\Tht)
\vdash^{* \tilde{a}(\iota)}_{\mS}
\Gam,A_{\iota}
$
for each $\iota\in J$, where $\tilde{a}(\iota)=\vphi_{e}(b_{1}+a(\iota))<\tilde{a}$.
A $(\bigwedge)$ yields
$
(\mathcal{H}_{\gam_{1}},\Tht)
\vdash^{* \tilde{a}}_{\mS}
\Gam
$.

Other cases are seen from SIH.
\eprf

\bprp\label{prp:collapsingOme}{\rm (Collapsing)}
Suppose 
$
\Tht\subset\calh_{\gam}(\psi_{\Ome}(\gam))
$,
$(\calh_{\gam},\Tht)\vdash^{* a}_{\Ome}\Gam$
and $\Gam\subset\Sig(\Ome)$. 
Then for $\hat{a}=\gam+\ome^{a}$ and $\bet=\psi_{\Ome}(\hat{a})$,
$(\calh_{\hat{a}+1},\Tht)\vdash^{* \bet}_{\bet}\Gam^{(\bet,\Ome)}$ holds.
\eprp

\bprp\label{lem:predcepi11Ome}{\rm (Cut-elimination)}
Suppose
$(\mathcal{H}_{\gamma},\Theta)\vdash^{* a}_{c+d}\Gamma$
with $c+d\leq\mS$ and $\lnot(c<\Ome\leq c+d)$.
Then $(\mathcal{H}_{\gamma},\Theta)\vdash^{* \tht_{d}(a)}_{c}\Gamma$.
\eprp

\bth\label{th:Pi11}
Assume $S_{1}\vdash\tht^{L_{\Ome}}$ for $\tht\in\Sig$.
Then there exists an $n<\ome$ such that
$L_{\alp}\models\tht$ for $\alp=\psi_{\Ome}(\ome_{n}(\mK+1))$ in $OT(\Pi_{1}^{1})$.
\end{theorem}
\bprf
Let
$S_{1}\vdash\theta^{L_{\Omega}}$ for a $\Sigma$-sentence $\theta$.
By Embedding \ref{th:embedpi11} pick an $m$ so that 
 $(\mathcal{H}_{\mS},\emptyset;\emptyset)
 \vdash_{\mathbb{K}+m}^{* \mathbb{K}\cdot 2+m}
 \theta^{L_{\Omega}};\emptyset$.
Cut-elimination \ref{lem:predcepi11*} yields 
$(\mathcal{H}_{\mS},\emptyset;\emptyset)
\vdash_{\mathbb{K}}^{* a}\theta^{L_{\Omega}}$
for $a=\omega_{m}(\mathbb{K}\cdot 2+m)<\omega_{m+1}(\mathbb{K}+1)$.
Now let
$\gamma_{0}=\omega_{m+2}(\mathbb{K}+1)$.
Let
$\beta=\psi_{\mathbb{K}}(\omega^{a})>\mS$, where
$\omega^{a}<\gamma_{0}=\omega_{m+2}(\mathbb{K}+1)$.
Collapsing \ref{lem:Kcollpase.1} yields
$(\mathcal{H}_{\ome^{a}+1},\emptyset;\emptyset)
\vdash_{\bet}^{* \bet}\theta^{L_{\Omega}};\emptyset$.

Let $\rho=\psi_{\mS}^{g}(\gam_{0})$ with $g=\{(\bet,\bet+\mK)\}$, where $\mK(\bet+1)=\bet+\mK$.
We obtain
$(\mathcal{H}_{\ome^{a}+1},\emptyset,\{\rho\})
\vdash_{\bet,\rho+1,\gam_{0},\bet}^{\bet}(\theta^{L_{\Omega}})^{(\rho)}$
by Capping \ref{lem:capping}.
Cut-elimination \ref{lem:predceregcap} yields 
$(\mathcal{H}_{\ome^{a}+1},\emptyset,\{\rho\})
\vdash_{\mS,\rho+1,\gam_{0},\bet}^{a_{1}}(\theta^{L_{\Omega}})^{(\rho)}$
for $a_{1}=\vphi_{\bet}(\bet)$.

We obtain
$
(\calh_{\gam_{1}},\emptyset)\vdash^{* a_{2}}_{\mS}\theta^{L_{\Omega}}$
by Lemma \ref{lem:sum},
where
$a_{2}=\vphi_{\rho+1}(\bet+a_{1})$ and $\gam_{1}=\gam_{0}+\mS$.
Cut-elimination \ref{lem:predcepi11Ome} yields
$(\calh_{\gam_{1}},\emptyset)\vdash^{* a_{3}}_{\Ome}\theta^{L_{\Omega}}$
for $a_{3}=\tht_{\mS}(a_{2})$.
Collapsing \ref{prp:collapsingOme} yields
$(\calh_{\gam_{1}+a_{3}+1},\emptyset)\vdash^{* \eta}_{\eta}\theta^{L_{\eta}}$
for $\eta=\psi_{\Ome}(\gam_{1}+a_{3})<\psi_{\Omega}(\omega_{m+3}(\mathbb{K}+1))$.
Cut-elimination \ref{lem:predcepi11Ome} yields
$(\calh_{\gam_{1}+a_{3}+1},\emptyset)\vdash^{* \tht_{\eta}(\eta)}_{0}\theta^{L_{\eta}}$.
We then see $L_{\eta}\models\theta$ by induction up to $\theta_{\eta}(\eta)$.
\eprf
\\

Actually the bound is shown to be tight.

\begin{theorem}\label{th:wf}{\rm \cite{singlewfprfpi11}}\\
${\sf KP}\ome+(M\prec_{\Sigma_{1}}V)$
proves the well-foundedness up to $\psi_{\Omega}(\ome_{n}(\mathbb{S}^{+}+1))$ for 
{\rm each} $n$.
\end{theorem}

${\sf KP}\ome+(M\prec_{\Sigma_{1}}V)$ proves an axiom of
$\Sig_{1}$-Separation with parameters from $M$.
$\exi b\left[b=\{x\in a: \vphi(x,c)\}=\{x\in a: M\models\vphi(x,c)\}\right]$, where
$c\in M$, $a\in M\cup\{M\}$ and $\vphi\in\Sig_{1}$.
However it is open for us whether the parameter-free $\Sig^{1}_{2}$-Comprehension Axiom
holds in ${\sf KP}\ome+(M\prec_{\Sigma_{1}}V)$.

\section{$\Pi_{1}$-Collection}\label{sec:pi1collection}

The axioms of the set theory ${\sf KP}\ome+\Pi_{1}\mbox{-Collection}+(V=L)$ consist of those of 
${\sf KP}\ome+(V=L)$ plus the axiom schema $\Pi_{1}\mbox{-Collection}$:
for each $\Pi_{1}$-formula $A(x,y)$ in the language of set theory, 
$\fal x\in a\exi yA(x,y)\to \exi b\fal x\in a\exi y\in b A(x,y)$.
It is easy to see that the second order arithmetic
$\Sig^{1}_{3}{\rm -DC+BI}$ is interpreted to
${\sf KP}\ome+\Pi_{1}\mbox{{\rm -Collection}}+(V=L)$ canonically.

Next we show that ${\sf KP}\ome+\Pi_{1}\mbox{-Collection}+(V=L)$ is contained in a set theory 
$S_{\mI}$.
The language of the theory $S_{\mI}$ is $\{\in,St,\Ome\}$ with a unary predicate constant $St$
and an individual constant $\Ome$.
$\Del_{0}(St)$ denotes the set of bounded formulas in the language $\{\in,St,\Ome\}$,
in which atomic formulas $St(t)$ may occur.
Similarly
$\Sig_{1}(St)$ the set of $\Sig_{1}$-formulas in the expanded language.
$St(\alp)$ is intended to denote the fact that $\alp$ is a stable ordinal, $L_{\alp}\prec_{\Sig_{1}}L$,
and $\Ome=\ome_{1}^{CK}$.
The axioms of $S_{\mI}$ are obtained from those
\footnote{In the axiom schemata $\Del_{0}$-Separation and $\Del_{0}$-Collection,
 $\Del_{0}$-formulas remain to mean a $\Del_{0}$-formula in which $St$ does not occur, 
 while the axiom of foundation may be applied to a formula in which $St$ may occur.} of ${\sf KP}\ome$ by adding the following axioms.
Let $ON$ denote the class of all ordinals.
For ordinals $\alp$, 
$\alp^{\dagger}$ denotes the least stable ordinal above $\alp$.
A \textit{successor stable ordinal} is an ordinal $\alp^{\dagger}$ for an $\alp$.
Note that the least stable ordinal $0^{\dagger}$ is a successor stable ordinal.

\benu
\item
$V=L$, and
the axioms for recursively regularity of $\Ome$.

\item
$\Del_{0}(St)\mbox{-collection}$:
\[
\fal x\in a\exi y\, \tht(x,y) \to \exi b
\fal x\in a \exi y\in b\, \tht(x,y)
\]
for each $\Del_{0}(St)$-formula $\tht$ in which the predicate $St$ may occurs.

\item
$L=\bigcup\{L_{\sig}: St(\sig)\}$, i.e.,
\beqn\label{eq:stbl0}
\fal \alp\in ON \exi \sig \left(\alp<\sig \land St(\sig) \right)
\eeqn.

\item 
For a successor stable ordinal $\sig<\mI$, $L_{\sig}\prec_{\Sig_{1}}L=L_{\mI}$:
\beqn\label{eq:sucstable}
SSt(\sig)\land \vphi(u) \land u\in L_{\sig} \to \vphi^{L_{\sig}}(u)
\eeqn
for each $\Sig_{1}$-formula $\vphi$ in the language of set theory, i.e., the constant $St$ does not occur in $\vphi$.

\eenu

\blem\label{lem:9setI}
$S_{\mI}$ is an extension of ${\sf KP}\ome+\Pi_{1}\mbox{{\rm -Collection}}+(V=L)$.
Namely $S_{\mI}$ proves $\Pi_{1}\mbox{{\rm -Collection}}$.
\elem
\bprf
Argue in $S_{\mI}$.
Let $A(x,y)$ be a $\Pi_{1}$-formula in the language of set theory.
We obtain by the axioms (\ref{eq:stbl0}) 
and (\ref{eq:sucstable})
\beqn\label{eq:9setI}
A(x,y)\lrarw\exi \bet(St(\bet^{\dagger}) \land x,y\in L_{\bet^{\dagger}} \land A^{L_{\bet^{\dagger}}}(x,y))
\eeqn
Assume $\fal x\in a\exi yA(x,y)$. 
Then we obtain
$\fal x\in a\exi y\exi \bet (St(\bet^{\dagger}) \land x,y\in L_{\bet^{\dagger}} \land A^{L_{\bet^{\dagger}}}(x,y))$ by (\ref{eq:9setI}).
 Since $St(\bet^{\dagger})\land x,y\in L_{\bet^{\dagger}} \land A^{L_{\bet^{\dagger}}}(x,y)$ is a $\Sig_{1}(St)$-formula,
pick a set $c$ such that 
$\fal x\in a\exi y\in c\exi \bet\in c(St(\bet^{\dagger})\land x,y\in L_{\bet^{\dagger}} \land A^{L_{\bet^{\dagger}}}(x,y))$
by $\Del_{0}(St)$-Collection.
Again by (\ref{eq:9setI}) we obtain $\fal x\in a\exi y\in c A(x,y)$.
\eprf
\\

Conversely in ${\sf KP}\ome+\Pi_{1}\mbox{{\rm -Collection}}+(V=L)$,
the predicate $St(\alp)$ is defined by a $\Pi_{1}$-formula $st(\alp)$ so that
(\ref{eq:sucstable}) is provable, and $\Del_{0}(St)$-collection follows from
$\Pi_{1}\mbox{{\rm -Collection}}$.

\blem\label{lem:9Sep}
${\sf KP}\ome+\Pi_{1}\mbox{{\rm -Collection}}$ proves each of
$\Sig_{1}\mbox{{\rm -Separation}}$,
$\Del_{2}\mbox{{\rm -Separation}}$ and
$\Sig_{2}\mbox{{\rm -Replacement}}$.
\elem
\bprf
We show that $\{x\in a:\vphi(x)\}$ exists as a set for a $\Sig_{1}$-formula $\vphi\equiv \exi y\tht(x,y)$ 
with a $\Del_{0}$ matrix $\tht$. 
We have by logic $\fal x\in a\exi y(\exi z\tht(x,z)\rarw \tht(x,y))$.
By $\Pi_{1}$-Collection pick a set $b$ so that
$\fal x\in a\exi y\in b(\vphi(x)\rarw \tht(x,y))$.
In other words,
$\{x\in a:\vphi(x)\}=\{x\in a: \exi y\in b\, \tht(x,y)\}$.
\eprf
\\

Let ${\rm Hull}_{\Sig_{1}}(\alp)$ denote the $\Sig_{1}$-Skolem hull ${\rm Hull}_{\Sig_{1}}(\alp)$ of an ordinal $\alp$.
${\rm Hull}_{\Sig_{1}}(\alp)$ is the collection of
$\Sig_{1}$-definable elements from parameters$<\alp$ in the universe.
Specifically
let $\{\vphi_{i}:i\in\ome\}$ denote an enumeration of $\Sig_{1}$-formulas.
Each is of the form $\vphi_{i}\equiv\exi y\theta_{i}(x,y;u)\, (\tht\in\Del_{0})$ with fixed variables $x,y,u$. Set for $b\in \alp$
\beqnarrs
r(i,b) & \simeq & \mbox{ {\rm the }} <_{L} \mbox{{\rm -least }} c\in L
\mbox{ {\rm such that} } L\models\theta_{i}((c)_{0},(c)_{1}; b)
\\
h(i,b) & \simeq & (r(i,b))_{0}
\\
\mbox{{\rm Hull}}_{\Sig_{1}}(\alp) & = & \{h(i,b)\in L :i\in\ome, b\in \alp\}
\eeqnarrs

The domain of the partial $\Del_{1}$-map $r$ is a $\Sig_{1}$-subset of $\ome\times\alp$, and
from Lemma \ref{lem:9Sep} ($\Sig_{1}\mbox{{\rm -Separation}}$) we see that
the domain exists as a set, and so does $\mbox{{\rm Hull}}_{\Sig_{1}}(\alp)$.
Therefore
its Mostowski collapse\footnote{The collapse coincides with $L_{\bet}$ for the least ordinal $\bet$ not in $\mbox{{\rm Hull}}_{\Sig_{1}}(\alp)$.}
ordinal $\bet\geq\alp$.
This shows (\ref{eq:stbl0}).
\\

Note that a limit of admissible ordinals need not to be admissible since
there exists a $\Pi_{3}^{-}$-formula $ad$ such that for any transitive set $x$,
$x$ is admissible iff $ad^{x}$ holds.
On the other side every limit $\kap$ of stable ordinals is stable:
for $c\in L_{\kap}$, pick a stable ordinal $\sig<\kap$ such that $c\in L_{\sig}$.
Then for $\Sig_{1}$-formula $A$,
$L\models A(c)\Rarw L_{\sig}\models A(c)\Rarw L_{\kap}\models A(c)$.

\subsection{Ordinals for $\Pi_{1}$-Collection}\label{sect:ordinalnotationpi1}

In this subsection up to subsection \ref{subsec:largecardinal}
we work in a set theory ${\sf ZFC}(St)$,
where $St$ is a unary predicate symbol. 
We assume that $St$ is an unbounded class of ordinals below $\mI$ such that
the least element $\mS_{0}$ of $St$ is larger than $\Ome$.
$\alp^{\dagger}$ denotes the least ordinal$>\alp$ in the class $St$ when $\alp<\mI$.
$\alp^{\dagger}:=\mI$ if $\alp\geq\mI$.
Then $\mS_{0}=\Ome^{\dagger}$.
Let $SSt:=\{\alp^{\dagger}:\alp\in ON\}$ and $LS=St\setminus SSt$.
For natural numbers $k$,
$\alp^{\dagger k}$ is defined recursively by
$\alp^{\dagger 0}=\alp$ and $\alp^{\dagger(k+1)}=(\alp^{\dagger k})^{\dagger}$.

$\varphi_{b}(\xi)$ denotes the binary Veblen function on 
$\mI^{+}=\omega_{\mathbb{I}+1}$ with $\varphi_{0}(\xi)=\omega^{\xi}$
Let $\Lam\leq\mI$ be a strongly critical number.
As in Definition \ref{df:Lam},
$\tilde{\varphi}_{b}(\xi):=\varphi_{b}(\mI\cdot \xi)$.
Let $b,\xi<\mI^{+}$.
$\theta_{b}(\xi)$ [$\tilde{\theta}_{b}(\xi)$] denotes
a $b$-th iterate of $\varphi_{0}(\xi)=\omega^{\xi}$ [of $\tilde{\varphi}_{0}(\xi)=\mI^{\xi}$], resp.

\bdf\label{df:Lam3pi1}
{\rm
A finite function $f:\mI\to\vphi_{\mI}(0)$ is said to be a \textit{finite function}
if 
$\fal i>0(a_{i}=1)$ and $a_{0}=1$ when $b_{0}>1$
in
$f(c)=_{NF}\tilde{\theta}_{b_{m}}(\xi_{m})\cdot a_{m}+\cdots+\tilde{\theta}_{b_{0}}(\xi_{0})\cdot a_{0}$
for any $c\in{\rm supp}(f)$.
Let
$SC_{\mI}(f):=\bigcup\{\{c\}\cup SC_{\mI}(f(c)): c\in {\rm supp}(f)\}$.

For a finite function $f$,
$c<\mI$ and $\xi<\vphi_{\mI}(0)$.
A relation $f<_{\mI}^{c}\xi$ is defined by induction on the
cardinality of the finite set $\{d\in {\rm supp}(f): d>c\}$ as in Definition \ref{df:Lam3}.\ref{df:Exp2.5}.

}
\end{definition}

\begin{definition}
{\rm
Let $A\subset\mathbb{I}$ be a set, and $\alpha\leq\mathbb{I}$ a limit ordinal.
\[
\alpha\in M(A) :\Lrarw A\cap\alpha \mbox{ is stationary in }
\alpha
\Lrarw \mbox{ every club subset of } \alpha \mbox{ meets } A.
\]

}
\end{definition}

Classes $\mathcal{H}_{a}(X)\subset\Gamma_{\mI+1}$,
$Mh^{a}_{c}(\xi)\subset(\mathbb{I}+1)$, and 
ordinals $\psi_{\kappa}^{f}(a)\leq\kappa$ are defined simultaneously as follows.

$\mathcal{H}_{a}(X)$ denotes the closure of
$\{0,\Omega,\mI\}\cup X$ under $+,\vphi$, 
$a\mapsto\psi_{\Ome}(a)$, $a\mapsto\psi_{\mI}(a)\in LS$,
$\alp\mapsto\alp^{\dagger}\in SSt$, and $(\pi,b,f)\mapsto \psi_{\pi}^{f}(b)$.

$\pi\in Mh^{a}_{c}(\xi)$ iff
 $\{a,c,\xi\}\subset\mathcal{H}_{a}(\pi)$ and
 the following condition is met for any 
 finite functions $f,g:\mI\to\vphi_{\mI}(0)$
 such that $f<_{\mI}^{c}\xi$
 \[
 SC_{\mI}(f,g) \subset\mathcal{H}_{a}(\pi) 
 \spand \pi\in Mh^{a}_{0}(g_{c})
 \Rarw \pi\in M(Mh^{a}_{0}(g_{c}*f^{c}))
\]
where
\beqnarrs
Mh^{a}_{c}(f) & := & \bigcap\{Mh^{a}_{d}(f(d)): d\in {\rm supp}(f^{c})\}
\\
& = &
\bigcap\{Mh^{a}_{d}(f(d)): c\leq d\in {\rm supp}(f)\}
\eeqnarrs

 Let $a,\pi$ ordinals and
 $f:\mI\to\vphi_{\mI}(0)$ a finite function.
Then $\psi_{\pi}^{f}(a)$ denotes the least ordinal $\kap<\pi$ such that
\beqn\label{eq:Psivec}
\kappa\in Mh^{a}_{0}(f) \spand
 \mathcal{H}_{a}(\kappa)
 \cap\pi\subset\kappa
 \spand
   \{\pi,a\}\cup SC_{\mI}(f)\subset\mathcal{H}_{a}(\kappa)
\eeqn
if such a $\kap$ exists. Otherwise set $\psi_{\pi}^{f}(a)=\pi$.

\beqn\label{eq:psiI}
\psi_{\mI}(a):=
\min(
\{\mI\}\cup\{\kap\in LS : \mathcal{H}_{a}(\kap)\cap\mI\subset\kap\}
)
\eeqn

For classes $A\subset\mathbb{I}$, let
$\alpha\in M^{a}_{c}(A)$ iff $\alpha\in A$ and
for any finite functions
$g:\mI\to\vphi_{\mI}(0)$
\begin{equation}\label{eq:Mca}
\alpha\in Mh_{0}^{a}(g_{c}) \spand SC_{\mI}(g_{c})\subset\mathcal{H}_{a}(\alpha) \Rarw
\alpha\in M\left( Mh_{0}^{a}(g_{c}) \cap A \right)
\end{equation}

\bprp\label{prp:definability}
Each of 
$x\in\mathcal{H}_{a}(y)$,
$x\in Mh^{a}_{c}(f)$ and $x=\psi^{f}_{\kappa}(a)$
is a $\Delta_{1}(St)$-predicate in ${\sf ZFC}(St)$.
\eprp

\subsection{A small large cardinal hypothesis}\label{subsec:largecardinal}
It is convenient for us to assume the existence of a small large cardinal in justification of 
the above definition.

Subtle cardinals are introduced by R. Jensen and K. Kunen.
It is shown in Lemma 2.7 of \cite{RathjenAFML2} that the set of shrewd cardinals
in $V_{\pi}$
is stationary in a subtle cardinal $\pi$.
From this fact we see that the set 
of shrewd limits of shrewd cardinals
in $V_{\pi}$ 
is also stationary in
 a subtle cardinal $\pi$, where
 for a shrewd cardinal $\kap$ in $V_{\pi}$,
 $\kap$ is a shrewd limit 
 iff $\kap$ is a limit of shrewd cardinals in $V_{\pi}$.

Let $C$ be a closed subset of $\pi$, and
$C_{0}\subset C$ be a subset defined by
$\kap\in C_{0}$ iff $\kap\in C$ and $\kap$ is a limit of shrewd cardinals.
Since the set of shrewd cardinals is stationary in $V_{\pi}$,
$C_{0}$ is a club subset of $\pi$.
Hence the exists a shrewd cardinal in $C_{0}$.

In this subsection we work in an extension $T$ of ${\sf ZFC}$ by adding
the axiom stating that
there exists a regular cardinal $\mI$ such that
the set $St$ of shrewd cardinals in $V_{\mI}$ is stationary in $\mI$.
In this subsection $\Ome$ denotes the least uncountable ordinal $\ome_{1}$, and
$LS$ denotes the set of 
shrewd limits in $V_{\mI}$.
The class $LS$ is stationary in $\mI$.
A \textit{successor shrewd cardinal} is a shrewd cardinal in $V_{\mI}$, not in $LS$.

\blem\label{lem:welldefinedness.2}
$\forall a[
\psi_{\mI}(a)<\mI]$.
\elem
\bprf
The set
$C=\{\kappa<\mI: 
\mathcal{H}_{a}(\kappa)\cap\mI\subset\kappa\}$
 is a club subset of the regular cardinal $\mI$.
This shows the existence of a $\kappa\in LS \cap C$, and hence
$\psi_{\mI}(a)<\mI$ by the definition (\ref{eq:psiI}).
\eprf

\blem\label{lem:welldefinedness.1}
Let $\mS$ be a shrewd cardinal, 
$a<\veps(\mI)$, $h:\mI\to\vphi_{\mI}(0)$ a finite function 
with
$\{a\}\cup SC_{\mI}(h)\subset\mathcal{H}_{a}(\mS)$. Then
$\mS\in Mh^{a}_{0}(h)\cap M(Mh^{a}_{0}(h))$.
\elem
\bprf
By induction on $\xi<\vphi_{\mI}(0)$ we show $\mS\in Mh^{a}_{c}(\xi)$
for $\{a,c,\xi\}\subset\mathcal{H}_{a}(\mS)$ as in Lemma \ref{lem:welldefinedness.1pi11}.
\eprf

\blem\label{lem:limitcollapse.1}
Let $\mS$ be 
a shrewd cardinal, $a$ an ordinal, and
$f:\mI\to\vphi_{\mI}(0)$ a finite function 
such that $\{a\}\cup SC_{\mI}(f)\subset\mathcal{H}_{a}(\mathbb{S})$.
Then
$\psi_{\mathbb{S}}^{f}(a)<\mathbb{S}$ holds.
\elem

\bcor\label{cor:stepdown}
Let $f,g:\mI\to\vphi_{\mI}(0)$
 be finite functions and $c\in{\rm supp}(f)$.
Assume  that 
there exists an ordinal
$d<c$ 
such that
$(d,c)\cap {\rm supp}(f)=(d,c)\cap {\rm supp}(g)=\emptyset$, 
$g_{d}=f_{d}$, 
$g(d)<f(d)+\tilde{\theta}_{c-d}(f(c);\mI)\cdot\omega$,
and
$g<_{\mI}^{c}f(c)$.

Then
$Mh^{a}_{0}(g)\prec Mh^{a}_{0}(f)$ holds.
In particular if $\pi\in Mh^{a}_{0}(f)$ and
$SC_{\mI}(g)\subset\mathcal{H}_{a}(\pi)$, then
$\psi_{\pi}^{g}(a)<\pi$.
\ecor
\bprf
This is seen as in Corollary \ref{cor:stepdownpi11}.
\eprf
\\

An \textit{irreducibility} of finite functions $f:\mI\to\vphi_{\mI}(0)$ is defined as in Definition \ref{df:irreducible},
and a lexicographic order $f<^{b}_{lx}g$ on finite functions $f,g$ as in Definition \ref{df:lx}.
Then $f<^{0}_{lx}g \Rarw Mh^{a}_{0}(f)\prec Mh^{a}_{0}(g)$ is seen 
as in Proposition \ref{lem:psinucomparison}.
\\

A computable notation system $OT(\mI)$ for $\Pi_{1}$-collection
is defined so as to be closed under 
Mostowski collapsings.
A new constructor $\mI[\cdot]$ is used to generate terms in $OT(\mI)$.
Note that 
there is no clause for constructing $\kap=\psi_{\mS}(a)$ from $a$
for $\mS\in LS$.

\bdf\label{df:prec}
{\rm
\benu
\item
$\{(\rho,\sig): \rho\prec\sig\}$ denotes the transitive closure 
of the relation $\{(\rho,\sig): \exi f,a(\rho=\psi_{\sig}^{f}(a))\}$.
Let $\rho\preceq\sig:\Lrarw \rho\prec\sig \lor \rho=\sig$.

\item
Let $\alp\prec\mS$ for an $\mS\in SSt$ and $b={\tt p}_{0}(\alp)$.
Then let
\[
M_{\alp}  :=  \mathcal{H}_{b}(\alp)
.\]

\item
For $\alp\in\Psi$ an ordinal ${\tt p}_{0}(\alp)$ is defined.

\benu
\item
Let $\alp\preceq\psi_{\mS}^{g}(b)$ for an $\mS\in SSt$. 
Then ${\tt p}_{0}(\alp)=b$.

\item
There exists an $\mS=\mT^{\dagger}\in SSt$ and a $\mT<\tau<\mS$ such that
$\alp\prec\tau^{\dagger k}$ for a $k>0$.
Let $\rho\prec\mS$ be such that $\alp=\bet[\rho/\mS]$ for a 
$\bet\in M_{\rho}$.
Let 
${\tt p}_{0}(\alp)={\tt p}_{0}(\bet)$.

\item
${\tt p}_{0}(\alp)=0$ otherwise.
\eenu

\eenu
}
\edf

$\alp=\psi_{\mS}^{f}(a)\in OT(\mI)$ only if

\beqn\label{eq:notationsystem.112.4}
SC_{\mI}(f)\subset \mathcal{H}_{a}(SC_{\mI}(a))
\eeqn
where 
$a={\tt p}_{0}(\alp)$.

Let
$\{\pi,a,d\}\subset  OT(\mI)$ with
$\pi\prec\mS\in SSt$, 
$m(\pi)=f$,
 $d<c\in \supp(f)$,
and $(d,c)\cap \supp(f)=\emptyset$.

When $g\neq\emptyset$, let
$g$
be an irreducible finite function 
 such that 
$SC_{\mI}(g)\subset OT(\mI)$,
$g_{d}=f_{d}$, $(d,c)\cap \supp(g)=\emptyset$,
$g(d)<f(d)+\tilde{\theta}_{c-d}(f(c); \mI)\cdot\omega$, 
and $g<_{\mI}^{c}f(c)$.

Then
$\alp=\psi_{\pi}^{g}(a)\in OT(\mI)$  only if 

\beqn\label{eq:notationsystem.111}
SC_{\mI}(g)\subset M_{\alp}
\eeqn

The Mostowski collapsing $\alp\mapsto\alp[\rho/\mS]\, (\alp\in M_{\rho})$ is defined as follows.
$(\mS)[\rho/\mS]:=\rho$, $(\mS^{\dagger})[\rho/\mS]:=\rho^{\dagger}$, and
$(\mathbb{I})[\rho/\mS]:=\mathbb{I}[\rho]$.
$(\tau^{\dagger})[\rho/\mS]=(\tau[\rho/\mS])^{\dagger}$,
where $\mS<\tau^{\dagger}$.
$(\mI[\tau])[\rho/\mS]=\mI[\tau[\rho/\mS]]$, where $\mI[\tau]\neq\mI$.

A relation $\alp<\bet$ for $\alp,\bet\in OT(\mI)$ is defined so that
$\psi_{\kap}^{f}(a)<\kap$ and 
$\rho<\psi_{\rho^{\dagger}}^{g}(b)<\rho^{\dagger}<\tau=\psi_{\mI[\rho]}(c)<
\psi_{\tau^{\dagger}}^{h}(d)<\tau^{\dagger}<\mI[\rho]$
for every $\kap,\rho$, $a,b,c,d$ and $f,g,h$.

\bprp\label{prp:jumpover}

There is no $\psi_{\sig}^{f}(a)\in OT(\mI)$ such that
 $\rho<\psi_{\sig}^{f}(a)\leq\rho^{\dagger}<\sig$.
\eprp

\blem\label{lem:Mostowskicollaps}
For $\rho\prec\mS$ and $\mS\in SSt$,
$\{\alpha[\rho/\mS]:\alpha\in M_{\rho}\}$ is a transitive collapse of $M_{\rho}$ as in Lemma \ref{lem:Mostowskicollapspi11}.
\elem

\subsection{Operator controlled derivations for $\Pi_{1}$-Collection}\label{subsec:operatorcont}
We consider $RS$-formulas in a language with a unary predicate $St(a)$, where
$a=L_{\kap}$ for a stable ordinal $\kap$.
Specifically
$St(a):\simeq \bigvee( (\fal x\in \iota(x\in a))\land (\fal x\in a(x\in \iota)))_{\iota\in J}$ 
with $J=\{L_{\kap}: \kap\in St\cap(|a|+1)\}$ for
$St\subset OT(\mI)$.

\bdf\label{df:QJ}
{\rm
A \textit{finite family} 
is a finite function ${\tt Q}\subset\coprod_{\mS}\Psi_{\mS}$ such that
its domain $dom({\tt Q})$ is a finite set of successor stable ordinals,
and ${\tt Q}(\mS)$ is a finite set of ordinals in $\Psi_{\mS}$ for each $\mS\in dom({\tt Q})$.
Let
${\tt Q}(\mT)=\emptyset$ for $\mT\not\in dom({\tt Q})$ and
$\bigcup{\tt Q}=\bigcup_{\mS\in dom({\tt Q})}{\tt Q}(\mS)$.
Define
$M_{{\tt Q}(\mS)}=\bigcap_{\sig\in{\tt Q}(\mS)} M_{\sig}$.

For $A\simeq\bigvee(A_{\iota})_{\iota\in J}$ and $\iota\in J$
\[
\iota\in [{\tt Q}]_{A}J = [{\tt Q}]_{\lnot A}J   :\Lrarw
\fal\mU\in dom({\tt Q})
\left(
\rk(A_{\iota})\geq\mU \Rarw \sfk(\iota)\subset M_{{\tt Q}(\mU)}
\right)
\]
}
\edf

We define a derivability relation 
$(\mathcal{H}_{\gamma},\Theta;{\tt Q}_{\Pi})\vdash^{* a}_{c}\Gamma;\Pi^{\cdot]}$
where 
$c$ is a bound of ranks of the inference rules $({\rm stbl})$ and of ranks of cut formulas.
The relation depends on an ordinal $\gam_{0}$, and
should be written as $(\mathcal{H}_{\gamma},\Theta;{\tt Q}_{\Pi})\vdash^{*a}_{c,\gamma_{0}} \Gamma;\Pi^{\cdot]}$.
However the ordinal $\gam_{0}$ will be fixed.
So let us omit it.

\bdf\label{df:controldercollaps}
{\rm
Let $\Theta$ a finite set of ordinals, $a,c$ ordinals, and 
${\tt Q}_{\Pi}$ a finite family such that 
$\gam_{0}\leq{\tt p}_{0}(\sig)$
for each $(\mS,\sig)\in{\tt Q}_{\Pi}$.
Let 
$\Pi=\bigcup_{\sigma\in\bigcup{\tt Q}_{\Pi}}\Pi_{\sigma}\subset\Delta_{0}(\mI)$ be a set of formulas such that
${\sf k}(\Pi_{\sigma})\subset M_{\sig}$
for each $(\mS,\sigma)\in{\tt Q}_{\Pi}$.
Let 
$\Pi^{[\cdot]}=\bigcup_{\sigma\in\bigcup{\tt Q}_{\Pi}}\Pi_{\sigma}^{[\sigma/\mathbb{S}]}$
and $\Theta_{{\tt Q}_{\Pi}(\mS)}=\Tht\cap M_{{\tt Q}_{\Pi}(\mS)}$.

$(\mathcal{H}_{\gamma},\Theta;{\tt Q}_{\Pi})\vdash^{* a}_{c} \Gamma;\Pi^{[\cdot]}$ holds
for a set
$\Gamma$ of formulas
if $\gam\leq\gam_{0}$

\begin{equation}
\label{eq:controlder}
{\sf k}(\Gam)\subset\mathcal{H}_{\gamma}[\Theta]\spand
\fal\sig\in\bigcup{\tt Q}_{\Pi}\left(\sfk(\Pi_{\sig})\subset\calh_{\gam}[\Tht^{(\sig)}]\right)
\end{equation}

\begin{equation}
\label{eq:controlderac}
\fal\mS\in dom({\tt Q}_{\Pi})
\left(
\{\gam,a,c,\gam_{0}\}
\cup\sfk^{\mS}(\Gam, \Pi)
\subset\mathcal{H}_{\gamma}[\Theta_{{\tt Q}_{\Pi}(\mS)}]
\right)
\footnote{(\ref{eq:controlderac}) means $\{\gam,a,c,\gam_{0}\}\subset\mathcal{H}_{\gamma}[\Theta]$ when 
$dom({\tt Q}_{\Pi})=\emptyset$.}
\end{equation}

\begin{equation}
\label{eq:controlderS}
\fal\{\mU\leq\mS\}\subset  dom({\tt Q}_{\Pi})\left(
\mS\in \mathcal{H}_{\gamma}[\Theta_{{\tt Q}_{\Pi}(\mU)}]
\right)
\end{equation}

and one of the following cases holds:

\begin{description}

\item[$(\bigvee)$]\footnote{The condition $|\iota|< a$ is absent in the inference $(\bigvee)$.}
There exist 
$A\simeq\bigvee(A_{\iota})_{\iota\in J}$, an ordinal
$a(\iota)<a$ and an $\iota\in J$ such that 
$A\in\Gamma$ and
$(\mathcal{H}_{\gamma},\Theta;{\tt Q}_{\Pi})\vdash^{* a(\iota)}_{c}\Gamma,
A_{\iota};\Pi^{[\cdot]}$.

\item[$(\bigvee)^{[\cdot]}$]
There exist $\sigma\in\bigcup{\tt Q}_{\Pi}$,
$A\simeq\bigvee(A_{\iota})_{\iota\in J}$, an ordinal
$a(\iota)<a$ and an $\iota\in [\sigma] J$ such that 
$A^{[\sigma/\mathbb{S}]}\in\Pi^{[\cdot]}$,
$(\mathcal{H}_{\gamma},\Theta;{\tt Q}_{\Pi})\vdash^{* a(\iota)}_{c}\Gamma;
\left(A_{\iota}\right)^{[\sigma/\mathbb{S}]},\Pi^{[\cdot]}$.

\item[$(\bigwedge)$]
There exist 
$A\simeq\bigwedge(A_{\iota})_{ \iota\in J}$,  ordinals $a(\iota)<a$ such that
$A\in\Gamma$ and
for each $\iota\in [{\tt Q}_{\Pi}]_{A} J$,
$(\mathcal{H}_{\gamma},\Theta\cup{\sf k}(\iota);{\tt Q}_{\Pi}
)
\vdash^{* a(\iota)}_{c}\Gamma,
A_{\iota};\Pi^{[\cdot]}$.

\item[$(\bigwedge)^{[\cdot]}$]
There exist $\sigma\in\bigcup{\tt Q}_{\Pi}$,
$A\simeq\bigwedge(A_{\iota})_{ \iota\in J}$, 
ordinals $a(\iota)<a$ such that
$A^{[\sigma/\mathbb{S}]}\in\Pi^{[\cdot]}$, and
$(\mathcal{H}_{\gamma},\Theta\cup{\sf k}(\iota);{\tt Q}_{\Pi}
)
\vdash^{* a(\iota)}_{c}\Gamma;\Pi^{[\cdot]},\left(A_{\iota}\right)^{[\sigma/\mathbb{S}]}$
for each $\iota\in [{\tt Q}_{\Pi}]_{A}J\cap[\sig]J$.

\item[$(cut)$]
There exist an ordinal $a_{0}<a$ and a formula $C$
such that 
$(\mathcal{H}_{\gamma},\Theta;{\tt Q}_{\Pi})\vdash^{* a_{0}}_{c}\Gamma,\lnot C;\Pi^{[\cdot]}$
and
$(\mathcal{H}_{\gamma},\Theta;{\tt Q}_{\Pi})\vdash^{* a_{0}}_{c}C,\Gamma;\Pi^{[\cdot]}$
with $\mbox{{\rm rk}}(C)<c$.

\item[$(\Sigma(St)\mbox{{\rm -rfl}})$]
There exist ordinals
$a_{\ell}, a_{r}<a$ and a formula $C\in\Sigma(St)$ 
such that $c\geq\mI$, 
$(\mathcal{H}_{\gamma},\Theta;{\tt Q}_{\Pi}
)\vdash^{* a_{\ell}}_{c}\Gamma,C;\Pi^{[\cdot]}$
and
$(\mathcal{H}_{\gamma},\Theta;{\tt Q}_{\Pi}
)\vdash^{* a_{r}}_{c}
\lnot \exists x\,C^{(x,\mI)}, \Gamma;\Pi^{[\cdot]}$.

\item[$(\Sigma(\Ome)\mbox{{\rm -rfl}})$]
There exist ordinals
$a_{\ell}, a_{r}<a$ and a formula $C\in\Sigma(\Ome)$ 
such that $c\geq\Ome$, 
$(\mathcal{H}_{\gamma},\Theta;{\tt Q}_{\Pi}
)\vdash^{* a_{\ell}}_{c}\Gamma,C;\Pi^{[\cdot]}$
and
$(\mathcal{H}_{\gamma},\Theta;{\tt Q}_{\Pi}
)\vdash^{* a_{r}}_{c}
\lnot \exists x<\Ome\,C^{(x,\Ome)}, \Gamma;\Pi^{[\cdot]}$.

\item[$({\rm stbl}(\mS))$]
There exist an ordinal $a_{0}<a$,
a successor stable ordinal $\mS$,
a $\bigwedge$-formula 
$B(0)\in\Del_{0}(\mS)$ 
and a $u\in Tm(\mI)$ 
for which the following hold:
\beqn\label{eq:stblS0}
\mS\in\calh_{\gam}[\Tht_{{\tt Q}_{\Pi}(\mS)}] \spand
\fal\mU\in dom({\tt Q}_{\Pi})\cap\mS
\left(
\mS\in\calh_{\gam}[\Tht_{{\tt Q}_{\Pi}(\mU)}]
\right)
\eeqn
$
\mS\leq\rk(B(u))<c
$,
$(\mathcal{H}_{\gamma},\Theta;{\tt Q}_{\Pi}
)\vdash^{* a_{0}}_{c}
\Gamma, B(u);\Pi^{[\cdot]}$, and 
$(\mathcal{H}_{\gamma},\Theta\cup\{\sig\};
{\tt Q}_{\Pi}\cup\{(\mS,\sigma)\}
)
\vdash^{* a_{0}}_{c}
\Gamma; \lnot B(u)^{[\sigma/\mathbb{S}]},\Pi^{[\cdot]}$
holds for every ordinal $\sigma\in\Psi_{\mS}$ such that
${\tt p}_{0}(\sig)\geq\gam_{0}$ and
\beqn\label{eq:stblS}
\Theta\cup\{\mS\}\subset M_{\sig} 
\eeqn
where
$dom({\tt Q}_{\Pi}\cup\{(\mS,\sigma)\})=dom({\tt Q}_{\Pi})\cup\{\mS\}$,
and
$\left({\tt Q}_{\Pi}\cup\{(\mS,\sig)\}\right)(\mS)={\tt Q}_{\Pi}(\mS)\cup\{\sig\}$.
{\small
\[
\hspace{-10mm}
\infer{(\mathcal{H}_{\gamma},\Theta;{\tt Q}_{\Pi})\vdash^{* a}_{c} \Gamma;\Pi^{[\cdot]}}
{
(\mathcal{H}_{\gamma},\Theta;{\tt Q}_{\Pi}
)\vdash^{* a_{0}}_{c}
\Gamma, B(u);\Pi^{[\cdot]}
&
\{
(\mathcal{H}_{\gamma},\Theta\cup\{\sig\};
{\tt Q}_{\Pi}\cup\{(\mS,\sigma)\}
)
\vdash^{* a_{0}}_{c}
\Gamma; \lnot B(u)^{[\sigma/\mathbb{S}]},\Pi^{[\cdot]}
\}_{\sig}
}
\]
}
Assume (\ref{eq:stblS0}) and (\ref{eq:stblS}). Then
$(\Tht\cup\{\sig\})_{({\tt Q}_{\Pi}\cup\{(\mS,\sig)\})(\mS)}=\Tht_{{\tt Q}_{\Pi}(\mS)}$, and
$(\Tht\cup\{\sig\})_{({\tt Q}_{\Pi}\cup\{(\mS,\sig)\})(\mU)}=
(\Tht\cup\{\sig\})_{{\tt Q}_{\Pi}(\mU)}\supset\Tht_{{\tt Q}_{\Pi}(\mU)}$
for $\mU\in dom({\tt Q}_{\Pi})\cap\mS$.

\end{description}
}
\edf

\blem\label{lem:tautology*}{\rm (Tautology)}
Let $\gam\in\calh_{\gam}[\sfk(A)]$ and $d=\mbox{{\rm rk}}(A)$.
\benu
\item\label{lem:tautology*.1}
$(\mathcal{H}_{\gamma},{\sf k}(A);\emptyset
)\vdash^{* 2d}_{0}
\lnot A, A; \emptyset$.

\item\label{lem:tautology*.2}
$(\mathcal{H}_{\gamma},{\sf k}(A)\cup\{\mS,\sig\};\{(\mS,\sig)\}
)\vdash^{* 2d}_{0}
\lnot A^{[\sig/\mS]};A^{[\sig/\mS]}$ if $\sfk(A)\cup\{\mS\}\subset M_{\sig}$ and $\gam\geq\mS$.
\eenu

\elem
\bprf
Each is seen by induction on $d=\mbox{{\rm rk}}(A)$.
For example consider the lemma \ref{lem:tautology*}.\ref{lem:tautology*.2}.
We have $\rk(A^{[\sig/\mS]})<\mS$ and
$(\sfk(A)\cup\{\mS,\sig\})\cap M_{\sig}=\sfk(A)\cup\{\mS\}$ for (\ref{eq:controlderac}) and (\ref{eq:controlderS}),
and $\sfk(A^{[\sig/\mS]})\subset\calh_{\mS}((\sfk(A)\cap\mS)\cup\{\sig\})$ for (\ref{eq:controlder}).
\eprf

\blem\label{th:embedreg}{\rm (Embedding of Axioms)}
For each axiom $A$ in $S_{\mI}$ 
there is an $m<\omega$ such that
 $(\mathcal{H}_{\mI},\emptyset;\emptyset)
 \vdash^{* \mI\cdot 2}_{\mI+m} A;\emptyset$
holds.
\elem
\bprf
Let us suppress the operator $\calh_{\mI}$.
We show first that the axiom (\ref{eq:sucstable}),
$SSt(\sig)\land \vphi(u) \land u\in L_{\sig} \to \vphi^{L_{\sig}}(u)$
by an inference $({\rm stbl}(\mathbb{S}))$ for successor stable ordinals $\mathbb{S}<\mI$.
Let $B(0)\in\Delta_{0}(\mathbb{S})$ be a $\bigwedge$-formula,
and $u\in Tm(\mI)$.

We may assume that
$\mI>d={\rm rk}(B(u))\geq\mathbb{S}$.
Let ${\sf k}_{0}={\sf k}(B(0))$ and ${\sf k}_{u}={\sf k}(u)$.
Then $\sfk(B(0))\subset\calh_{0}(\sfk_{0})$.
Let $\sigma\in\Psi_{\mS}$ be an ordinal such that 
${\sf k}_{0}\cup{\sf k}_{u}\cup\{\mS\}\subset M_{\sig}$ and 
$\gam_{0}\leq{\tt p}_{0}(\sig)$. 

{\footnotesize
\[
\hspace{-5mm}
\infer[(\bigwedge)]{ {\sf k}_{0}\cup\{\mS\};\vdash^{* \mathbb{I}+1}_{\mathbb{I}}
 \lnot \exists x\, B(x),\exists x\in L_{\mathbb{S}}B(x);}
 {
 \infer[({\rm stbl}(\mS))]{{\sf k}_{0}\cup{\sf k}_{u}\cup\{\mS\}; \vdash^{* \mathbb{I}}_{\mathbb{I}}
    \lnot B(u),\exists x\in L_{\mathbb{S}}B(x);}
    {
     {\sf k}_{0}\cup{\sf k}_{u};  \vdash^{* 2d}_{0}\lnot B(u),B(u);
    &
    \{
    \infer[(\bigvee)]{{\sf k}_{0}\cup\sfk_{u}\cup\{\mS,\sig\};\{(\mS,\sig)\}\vdash^{* 2d+1}_{0}
     \exists x\in L_{\mathbb{S}}B(x);\lnot B(u)^{[\sigma/\mathbb{S}]}
     \}_{\sig}
     }
     {
     {\sf k}_{0}\cup{\sf k}_{u}\cup\{\mS,\sig\};\{(\mS,\sig)\} \vdash^{* 2d}_{0}
      B(u^{[\sigma/\mathbb{S}]}); \lnot B(u)^{[\sigma/\mathbb{S}]}
     }
    }
 }
\]
}

Therefore
$(\mathcal{H}_{\mI}, \emptyset;\emptyset)\vdash^{* \mI+\ome}_{\mI}
\fal\mS,v
\left[
SSt(\mS)\land
A(v) \land v\in L_{\mathbb{S}} \to A^{(\mathbb{S},\mI)}(v)
\right];\emptyset$, where
\\
$SSt(\alp):\Lrarw
\left(St(\alp)\land \exi\bet<\alp
\fal \gam<\alp(St(\gam)\to \gam\leq\bet)]\right)$.

Next we show the axiom (\ref{eq:stbl0}).
Let $\alp$ be an ordinal such that $\alp<\mI$.
We obtain for $\alp<\alp^{\dagger}<\mI$ with $d_{0}=\rk(\alp<\alp^{\dagger})$
and $\alp^{\dagger}\leq d_{1}=\rk(St(\alp^{\dagger}))<d_{2}=\ome(\alp^{\dagger}+1)$
with $\alp^{\dagger}\in\mathcal{H}_{0}[\{\alp\}]$
\[
\infer[(\bigwedge)]{
\emptyset;\emptyset \vdash^{* \mI}_{0}
\fal \alp\in ON \exi \sig \left(\alp<\sig \land St(\sig) \right);\emptyset
}
{
 \infer[(\bigvee)]{
\{\alp\};\emptyset \vdash^{* d_{2}+1}_{0} \exi \sig \left(\alp<\sig \land St(\sig)\right);\emptyset
 }
 {
  \infer[(\bigwedge)]{
 \{\alp\};\emptyset \vdash^{* d_{2}}_{0} 
  \alp<\alp^{\dagger} \land St(\alp^{\dagger});\emptyset
  }
  {
\{\alp\};\emptyset \vdash^{* d_{0}}_{0}\alp<\alp^{\dagger};\emptyset
  &
  \{\alp\};\emptyset \vdash^{* 2d_{1}}_{0} St(\alp^{\dagger});\emptyset
  }
 }
}
\]
\eprf

\blem\label{lem:predcepi1*}{\rm (Cut-elimination)}
Assume
$(\mathcal{H}_{\gamma},\Theta;{\tt Q}_{\Pi})\vdash^{* a}_{c+1}\Gamma;\Pi^{[\cdot]}$
with $c\geq\mI$.
Then $(\mathcal{H}_{\gamma},\Theta;{\tt Q}_{\Pi})\vdash^{* \ome^{a}}_{c}\Gamma;\Pi^{[\cdot]}$.
\elem
\bprf
Use the fact:
if $(\mathcal{H}_{\gamma},\Theta;{\tt Q}_{\Pi})\vdash^{* a}_{c}\Gamma;\Pi^{[\cdot]}$
and $\Tht\cup\{\mS\}\subset M_{\sig}$, then
$(\mathcal{H}_{\gamma},\Theta\cup\{\sig\};{\tt Q}_{\Pi}\cup\{(\mS,\sig)\})\vdash^{* a}_{c}\Gamma;\Pi^{[\cdot]}$.
\eprf

\blem\label{lem:Kcollpase.1}{\rm (Collapsing)}
Let
$\Gamma\subset\Sigma(St)$ be a set of formulas.
Suppose
$\Theta\subset
\mathcal{H}_{\gamma}(\psi_{\mI}(\gamma))$
 and
$
(\mathcal{H}_{\gamma},\Theta;{\tt Q}_{\Pi}
)\vdash^{* a}_{\mI}\Gamma;\Pi^{[\cdot]}
$.
Let
$\beta=\psi_{\mI}(\hat{a})$ with $\hat{a}=\gamma+\omega^{a}$.
Then 
$(\mathcal{H}_{\hat{a}+1},\Theta;{\tt Q}_{\Pi})
\vdash^{* \beta}_{\beta}
\Gamma^{(\beta,\mI)};\Pi^{[\cdot]}$ holds.
\elem
\bprf
By induction on $a$.
We have 
$\{\gamma, a\}\subset\mathcal{H}_{\gamma}[\Theta_{{\tt Q}_{\Pi}(\mS)}]$ by (\ref{eq:controlderac}),
and
$\beta\in\mathcal{H}_{\hat{a}+1}[\Theta_{{\tt Q}_{\Pi}(\mS)}]$ for $\mS\in dom({\tt Q}_{\Pi})$

When the last inference is a $({\rm stbl}(\mS))$, let
$B(0)\in\Delta_{0}(\mathbb{S})$ be a $\bigwedge$-formula
and a term $u\in Tm(\mI)$ such that
$\mS\leq{\rm rk}(B(u))<\mathbb{I}$,
${\sf k}(B(u))\subset \mathcal{H}_{\gamma}[\Theta]$,
and
$
(\mathcal{H}_{\gamma},\Theta;{\tt Q}_{\Pi}
)\vdash^{* a_{0}}_{\mI}
\Gamma, B(u);\Pi^{[\cdot]}
$ for an ordinal $a_{0}\in\mathcal{H}_{\gamma}[\Theta_{{\tt Q}_{\Pi}}]\cap a$.
Then we obtain
$\mS\leq{\rm rk}(B(u))<\beta$.
\eprf

\subsection{Operator controlled derivations with caps}

Let $
(\mathcal{H}_{\gamma},\Theta;{\tt Q}_{\Pi}
)\vdash^{* a}_{\mK}\Gamma;\Pi^{[\cdot]}
$ in the calculus for $\Pi^{1}_{1}$-reflection in subsection \ref{subsec:operatorcontpi11}.
In Capping \ref{lem:capping}, each formula $A\in\Gam$ puts on a cap
$\rho$ such that ${\tt Q}_{\Pi}\subset \rho$ and
(\ref{eq:capping}), $\Tht\subset M_{\rho}$.
(\ref{eq:capping}) is needed in \textbf{Case 3.1} of the proof.
Namely 
when $\Gam\ni A\simeq\bigvee\left(A_{\iota}\right)_{\iota\in J}$ is introduced by a $(\bigvee)$
such that
$(\mathcal{H}_{\gamma},\Theta;{\tt Q}_{\Pi})\vdash^{* a(\iota)}_{\mK}
\Gamma,A_{\iota};\Pi^{[\cdot]}$,
we need $\iota\in[\rho]J$, i.e., $\sfk(\iota)\subset M_{\rho}$, which follows from
${\sf k}(A_{\iota})\subset\mathcal{H}_{\gamma}[\Theta]\subset 
M_{\rho}$ by (\ref{eq:controlderpi11}) and $\Theta\subset M_{\rho}$.

We are concerned here with several stable ordinals $\mS,\mT,\ldots$.
It is convenient for us to regard \textit{uncapped formulas} $A$ as capped formulas
$A^{({\tt u})}$ with its cap ${\tt u}$.
Let $M_{{\tt u}}=OT(\mI)$.

In Capping \ref{lem:Mostowskicollpasecap} $\Gam$ is classified into 
$\Gam=\Gam_{{\tt u}}\cup\bigcup_{\mS\in dom({\tt Q}_{\Pi})}\Gam_{\mS}$.
$\Gam_{\mS}$ is the set of formulas $B(u)$ in inferences for the stability of a successor stable
ordinal $\mS$.
{\small
\[
\hspace{-10mm}
\infer{(\mathcal{H}_{\gamma},\Theta;{\tt Q}_{\Pi})\vdash^{* a}_{c} \Gamma;\Pi^{[\cdot]}}
{
(\mathcal{H}_{\gamma},\Theta;{\tt Q}_{\Pi}\cup\{\mS\}
)\vdash^{* a_{0}}_{c}
\Gamma, B(u);\Pi^{[\cdot]}
&
\{
(\mathcal{H}_{\gamma},\Theta\cup\{\sig\};
{\tt Q}_{\Pi}\cup\{(\mS,\sigma)\}
)
\vdash^{* a_{0}}_{c}
\Gamma; \lnot B(u)^{[\sigma/\mathbb{S}]},\Pi^{[\cdot]}
\}_{\sig}
}
\]
}
Each formula $A\in\Gam_{\mS}$ puts on a cap $\rho_{\mS}$ for the stable ordinal $\mS$.
Then (\ref{eq:capping}) runs $\Tht\subset M_{\rho_{\mS}}$ for every $\mS\in dom({\tt Q}_{\Pi})$.
This means $\Tht\subset M_{\partial{\tt Q}}:=\bigcap_{\kap\in\partial{\tt Q}}M_{\kap}$, where
\[
\partial{\tt Q}=\{\max({\tt Q}(\mS)) : \mS\in dom({\tt Q}), {\tt Q}(\mS)\neq\emptyset\}.
\]
Ordinals occurring in derivations are restricted to the set $M_{\partial{\tt Q}}$.
\\

In section \ref{sec:pi11} for $\Pi^{1}_{1}$-reflection, an ordinal $\gam_{0}$ is a threshold, which means that
every ordinal occurring in derivations is in $\calh_{\gam_{0}}(0)$ and
the subscript $\gam\leq\gam_{0}$ in $\calh_{\gam}$, while
each $\rho\in{\tt Q}$ exceeds $\gam_{0}$ in such a way that
${\tt p}_{0}(\rho)\geq\gam_{0}$.
This ensures us that $\calh_{\gam}(M_{\rho})\subset M_{\rho}$.
In the end, inferences $({\rm rfl}(\rho,d,f))$ are removed in Lemma \ref{lem:sum}
by moving outside $\calh_{\gam_{0}}(0)$.
Specifically ${\tt Q}\subset\calh_{\gam_{0}+\mS}(0)$.

Now we have several (successor) stable ordinals $\mS,\mT,\ldots\in dom({\tt Q})$.
Inferences $({\rm stbl}(\mS))$ and their children $({\rm rfl}_{\mS}(\rho,d,f))$ are eliminated
first for bigger $\mS>\mT$, and then smaller ones $({\rm stbl}(\mT))$.
Therefore we need assignment $dom({\tt Q})\ni\mS\mapsto\gam_{\mS}^{{\tt Q}}$ for 
thresholds so that
$\gam_{\mS}^{{\tt Q}}<\gam_{\mT}^{{\tt Q}}$ if $\mS>\mT$.
This is done by gapping, i.e., a gap $\mI\cdot 2^{a}$ between $\gam_{\mS}^{{\tt Q}}$ and 
$\gam_{\mT}^{{\tt Q}}$
in advance, when 
$
(\mathcal{H}_{\gamma},\Theta;{\tt Q}_{\Pi}
)\vdash^{* a}_{c}\Gamma;\Pi^{[\cdot]}
$
is embedded to
$(\mathcal{H}_{\gamma},\Theta_{\Pi}, {\tt Q})
\vdash^{a}_{c,c,\gam_{0}}
\widehat{\Gamma},\widehat{\Pi}$, cf.\,Capping \ref{lem:Mostowskicollpasecap}.

\bdf\label{df:caphat}
{\rm
A triple $({\tt Q},\gam^{{\tt Q}}, e^{{\tt Q}})$ 
is said to be a \textit{finite family for ordinals $\gamma_{0}$ and $b_{1}$}
if ${\tt Q}$ is a finite family in the sense of Definition \ref{df:QJ} and the following conditions are met:
 \benu
 
  \item
$\gam^{{\tt Q}}$ is a map $dom({\tt Q})\ni\mS\mapsto\gam_{\mS}^{{\tt Q}}$ such that
$\gam_{0}+\mI^{2}>\gam^{{\tt Q}}_{\mS}\geq \gam_{0}$,
$\gam^{{\tt Q}}_{\mS}\geq \gam^{{\tt Q}}_{\mT}+\mI$ for $\{\mS<\mT\}\subset dom({\tt Q})$
and
$\mS\in\calh_{\gam^{{\tt Q}}_{\mS}+\mI}$ for $\mS\in dom({\tt Q})$.

 ${\tt Q}$ is said to have \textit{gaps} $\eta$ if
$\gam^{{\tt Q}}_{\mS}\geq \gam^{{\tt Q}}_{\mT}+\mI\cdot\eta$ holds
for $\{\mS<\mT\}\subset dom({\tt Q})$, and
$\gam^{{\tt Q}}_{\mS}\geq \gam_{0}+\mI\cdot\eta$ for $\mS\in dom({\tt Q})$.

  \item
For each $\rho\in{\tt Q}(\mS)$, 
$m(\rho):\mI \to\vphi_{\mI}(0)$ is special, $s(\rho)\leq b_{1}$,
$\rho\in\calh_{\gam_{\mS}^{{\tt Q}}+\mI}(0)$,
and $\gam_{\mS}^{{\tt Q}}\leq{\tt p}_{0}(\rho)$.

\item
$e^{{\tt Q}}$ assigns an ordinal $e_{\mS}^{{\tt Q}}\in\calh_{\gam^{{\tt Q}}_{\mS}+\mI}\cap(\mS+1)$ to each $\mS\in dom({\tt Q})$ such that
\beqn\label{eq:familye}
\max(\{0\}\cup\{\rho\in{\tt Q}(\mS): s(\rho)>\mS\})< e^{{\tt Q}}_{\mS}
\eeqn
Let $e^{{\tt Q}}_{\mS}=\mS$ when $\mS\not\in dom({\tt Q})$.

\eenu

}
\edf

\bdf
{\rm
For a finite family ${\tt Q}$, 
and for $A\simeq\bigvee(A_{\iota})_{\iota\in J}$
\[
[{\tt Q}]_{A^{(\rho)}}J=[{\tt Q}]_{\lnot A^{(\rho)}}J=[{\tt Q}]_{A}J\cap [\partial{\tt Q}]J\cap [\rho]J
\]
where $[{\tt u}]J=J$ and
\[
[\partial{\tt Q}]J=\bigcap_{\kap\in\partial{\tt Q}} [\kap]J
.\]
}
\edf

\bdf
{\rm
\benu

\item
For a finite family ${\tt Q}$, let
$
\partial{\tt Q}=\{\max({\tt Q}(\mS)) : \mS\in dom({\tt Q}), {\tt Q}(\mS)\neq\emptyset\}
$ and
$M_{\partial{\tt Q}}=\bigcap_{\kap\in\partial{\tt Q}}M_{\kap}$.

\item
\[
[{\tt Q}]_{A^{(\rho)}}J=[{\tt Q}]_{\lnot A^{(\rho)}}J=[{\tt Q}]_{A}J\cap [\partial{\tt Q}]J\cap [\rho]J
\]
where 
$[{\tt u}]J=J$ and
$[\partial{\tt Q}]J=\bigcap_{\kap\in\partial{\tt Q}} [\kap]J$.

\eenu
}
\edf

\bdf\label{df:resolvent}
{\rm
$H_{\rho}^{{\tt Q}}(f,b_{1},\gamma,\Theta)$
 denotes the \textit{resolvent class} for ${\tt Q}$, $\rho$,
special functions $f$,
ordinals $b_{1},\gamma$, and
finite sets $\Theta$ of ordinals defined as follows:
$\sig\in H_{\rho}^{{\tt Q}}(f,\gamma,\Theta)$ iff 
$\sig\in\calh_{\gam+\mI}(0)\cap\rho\cap M_{\partial{\tt Q}}$,
$SC_{\mI}(m(\sig))\subset\calh_{\gam}[\Tht]$,
$\Tht
\subset M_{\sig}$,
$\gam\leq{\tt p}_{0}(\sig)\leq{\tt p}_{0}(\rho)$ and
$m(\sig)$ is special such that 
$s(f)\leq s(m(\sig))\leq b_{1}$,
 $f^{\prime} \leq (m(\sig))^{\prime}$,
where $\sig, \rho\prec\mS$ and
$f\leq g\Leftrightarrow
\forall i(f(i)\leq g(i))$.
}
\edf

We define another derivability  relation $(\mathcal{H}_{\gamma},\Theta,{\tt Q}
)\vdash^{a}_{c,\xi,\gam_{0},b_{1}}\Gamma$, 
where
$c$ is a bound of ranks of cut formulas,
and $\xi$ a bound of ordinals $\mS$ in the inference rules $({\rm rfl}_{\mS}(\rho,d,f,b_{1}))$.

\bdf\label{df:controldercollaps}
{\rm
Let 
$\Tht^{(\rho)}=\Tht\cap M_{\rho}$ and
$\Tht_{\partial{\tt Q}}=\Tht\cap M_{\partial{\tt Q}}$.
Let 
$a,b,c,\xi<\mI$, a finite set $\Tht\subset\mI$, 
and
${\tt Q}$ be a finite family for 
$\gamma_{0}, b_{1}$ such that
$dom({\tt Q})\subset(\xi+1)$.

$(\mathcal{H}_{\gamma},\Tht,{\tt Q})\vdash^{a}_{c,\xi,\gamma_{0},b_{1}} \Gamma$ holds
for a sequent
$\Gamma=\bigcup\{\Gam_{\rho}^{(\rho)}:\rho\in\{{\tt u}\}\cup\bigcup{\tt Q}\}$ 
if $\gam\leq\gam_{0}$

\begin{equation}
\label{eq:controldercap}
\forall \rho\in\{{\tt u}\}\cup\bigcup{\tt Q}
\left(
{\sf k}(\Gam_{\rho})
\subset\mathcal{H}_{\gamma}[\Theta^{(\rho)}]
\cap\mathcal{H}_{\gamma}[\Theta_{\partial{\tt Q}}]
\right)
\eeqn

\begin{equation}
\label{eq:controldercapace}
\fal\mS\in dom({\tt Q})
\left(
\{\gamma,a,c,\xi,\gam_{0},b_{1}\}
\subset\mathcal{H}_{\gamma}[
\Tht_{{\tt Q}(\mS)}
]
\cap\mathcal{H}_{\gamma}[\Theta_{\partial{\tt Q}}]
\right)
\footnote{(\ref{eq:controldercapace}) means $\{\gam,a,c,\xi,\gam_{0}\}\subset\mathcal{H}_{\gamma}[\Theta]$ when 
$dom({\tt Q})=\emptyset$.}
\end{equation}

\begin{equation}
\label{eq:controlderScapgam}
\fal\{\mU\leq\mS\}\subset dom({\tt Q})
\left(
\{\mS,\gam^{{\tt Q}}_{\mS}\}\subset
\mathcal{H}_{\gamma}[
\Tht_{{\tt Q}(\mU)}
]
\cap\mathcal{H}_{\gamma}[\Theta_{\partial{\tt Q}}]
\right)
\eeqn

\begin{equation}
\label{eq:controlderScap}
\fal\rho\in \{{\tt u}\}\cup\bigcup{\tt Q}\fal\mS\in dom({\tt Q})
\left(
\sfk^{\mS}(\Gam_{\rho})\subset\mathcal{H}_{\gamma}\left[
\Theta_{{\tt Q}(\mS)}
\right]
\cap\mathcal{H}_{\gamma}[\Theta_{\partial{\tt Q}}]
\right)
\end{equation}

\begin{equation}
\label{eq:controlderrhocap}
\fal(\mS,\rho)\in {\tt Q}
\left(
SC_{\mI}(m(\rho))\subset\mathcal{H}_{\gamma_{\mS}^{{\tt Q}}}
\left[
\Tht^{(\rho)}\cup\{\mS\}\cup\Theta_{\partial{\tt Q}}
\right]
\right)
\end{equation}

and one of the following cases holds:

\begin{description}

\item[{\rm (Taut)}]
$\{\lnot A^{(\rho)},A^{(\rho)}\}\subset\Gam$ for
a $\rho\in\{{\tt u}\}\cup\bigcup{\tt Q}$ and a formula $A$ such that $\rk(A)<\mS\leq\xi$
for some successor stable ordinal $\mS$.

If $\rk(A)<\mS$, then $(\mathcal{H}_{\gamma},\Tht,{\tt Q}
)\vdash^{0}_{0,\mS,\gam_{0},b_{1}}
\lnot A^{(\sigma)}, A^{(\sigma)}$ by (Taut) provided that (\ref{eq:controldercapace}) and (\ref{eq:controlderScap}) are met.

\item[$(\bigvee)$]
There exist 
$A\simeq\bigvee(A_{\iota})_{\iota\in J}$, a cap $\rho\in \{{\tt u}\}\cup\bigcup{\tt Q}$, an ordinal
$a(\iota)<a$ and an $\iota\in [\rho]J\cap[\partial{\tt Q}]J$ such that 
$A^{(\rho)}\in\Gamma$ and
$(\mathcal{H}_{\gamma},\Theta,{\tt Q})\vdash^{a(\iota)}_{c,\xi,\gam_{0},b_{1}}\Gamma,
\left(A_{\iota}\right)^{(\rho)}$.

\item[$(\bigwedge)$]
There exist 
$A\simeq\bigwedge(A_{\iota})_{ \iota\in J}$, 
a cap $\rho\in \{{\tt u}\}\cup\bigcup{\tt Q}$, ordinals $a(\iota)<a$ 
for each $\iota\in [{\tt Q}]_{A^{(\rho)}}J$ such that
$A^{(\rho)}\in\Gamma$ and
$(\mathcal{H}_{\gamma},\Theta\cup{\sf k}(\iota), {\tt Q}
)
\vdash^{a(\iota)}_{c,\xi,\gam_{0},b_{1}}\Gamma,
\left(A_{\iota}\right)^{(\rho)}$.

\item[$(cut)$]
There exist a cap $\rho\in \{{\tt u}\}\cup\bigcup{\tt Q}$, ordinals $a_{0}<a$ and a formula $C$
such that 
$\mbox{{\rm rk}}(C)<c$,
$(\mathcal{H}_{\gamma},\Theta,{\tt Q})\vdash^{a_{0}}_{c,\xi,\gam_{0},b_{1}}\Gamma,\lnot C^{(\rho)}$
and
$(\mathcal{H}_{\gamma},\Theta,{\tt Q})\vdash^{a_{0}}_{c,\xi,\gam_{0}, b_{1}}C^{(\rho)},\Gamma$.

\item[$(\Sigma(\Ome)\mbox{{\rm -rfl}})$]

There exist ordinals
$a_{\ell}, a_{r}<a$ and an uncapped formula $C\in\Sigma(\Ome)$
such that 
$c\geq \Ome$,
$(\mathcal{H}_{\gamma},\Theta,{\tt Q}
)\vdash^{a_{\ell}}_{c,\xi,\gam_{0},b_{1}}\Gamma,C$
and
$(\mathcal{H}_{\gamma},\Theta,{\tt Q}
)\vdash^{a_{r}}_{c,\xi,\gam_{0},b_{1}}
\lnot \exists x<\pi\,C^{(x,\Ome)}, \Gamma$.

\item[$({\rm rfl}_{\mS}(\rho,d,f,b_{1}))$]
There exist a successor stable ordinal $\mS\leq\xi$ 
and an ordinal $\rho\prec\mS$ such that
\beqn\label{eq:rflrho}
\Tht_{{\tt Q}(\mS)}\cup\{\mS\}\cup\Tht_{\partial{\tt Q}}  \subset M_{\rho}
\eeqn
and
$
\rho\in{\tt Q}(\mS)
$ if $\mS\in dom({\tt Q})$.
Let ${\tt R}={\tt Q}$ if $\mS\in dom({\tt Q})$.
Otherwise ${\tt R}={\tt Q}\cup\{(\mS,\rho)\}$, where
${\tt Q}\cup\{(\mS,\rho)\}$ is a finite family for $\gam_{0}$ 
extending ${\tt Q}$ such that
$dom({\tt R})=dom({\tt Q})\cup\{\mS\}$,
${\tt R}(\mS)={\tt Q}(\mS)\cup\{\rho\}$,
$e^{{\tt R}}_{\mT}=e^{{\tt Q}}_{\mT}$ for $\mS\neq\mT\in dom({\tt Q})$,
$\gam^{{\tt Q}}_{\mT}\geq\gam^{{\tt R}}_{\mS}+\mI$ for every $\mS>\mT\in dom({\tt Q})$
and $\gam^{{\tt R}}_{\mS}\geq\gam_{0}+\mI$.

Also there exist
an ordinal
$d\in {\rm supp}(m(\rho))$,
a special function $f$,
an ordinal $a_{0}<a$,  
and a 
finite set $\Delta$ of uncapped formulas enjoying the following conditions.

\begin{enumerate}

\item[(r0)]
$\rho< e^{{\tt R}}_{\mS}$
if $s(\rho)=\max({\rm supp}(m(\rho)))>\mS$.

\item[(r1)]
$\Delta\subset\bigvee_{\mS}\left( 
d
\right):=
\{\delta: \mbox{{\rm rk}}(\delta)<
d, \delta \mbox{ is a } \bigvee\mbox{-formula}\}\cup\{\delta: \mbox{{\rm rk}}(\delta)<\mS\}$.

\item[(r2)]
For $g=m(\rho)$, $s(f)\leq b_{1}$,
$SC_{\mI}(f)\cup SC_{\mI}(g)\subset\mathcal{H}_{\gamma_{\mS}^{{\tt R}}}[\Theta^{(\rho)}]$ and 
$
f_{d}=g_{d} \,\&\, f^{d}<_{\mI}^{d} g^{\prime}(d)$.

 \item[(r3)]
For each $\delta\in\Delta$,
 $(\mathcal{H}_{\gamma},\Tht,{\tt R}
 )\vdash^{a_{0}}_{c,\xi,\gam_{0}}\Gamma, \lnot\delta^{(\rho)}$.

\item[(r4)]
Let $\gam^{{\tt R}\cup\{(\mS,\sig)\}}=\gam^{{\tt R}}$,
$e^{{\tt R}\cup\{(\mS,\sig)\}}=e^{{\tt R}}$ and
$\sigma\in H^{{\tt R}}_{\rho}(f,b_{1}, \gam^{{\tt R}}_{\mS}, \Tht^{(\rho)}\cup\Tht_{\partial{\tt Q}})$.
Then
$
(\mathcal{H}_{\gamma},\Tht\cup\{\sig\},
{\tt R}\cup\{(\mS,\sig)\}
)\vdash^{a_{0}}_{c,\xi,\gam_{0}}\Gamma, 
\Delta^{(\sig)}$ holds.

In particular
$\sig<e^{{\tt R}}_{\mS}$ if $s(\sig)>\mS$
by (\ref{eq:familye}).

\end{enumerate}

\end{description}

}
\edf

Note that 
$\bigcup{\tt Q}\subset\mathcal{H}_{\gam}[\Tht]$
need not to hold.
Moreover $(\Tht\cup\{\sig\})_{({\tt R}(\mS)\cup\{\sig\})}=\Tht_{{\tt R}(\mS)}=\Tht_{{\tt Q}(\mS)}$ 
and $\Tht_{\partial{\tt R}}=\Tht_{\partial{\tt Q}}$
by 
$\Tht^{(\rho)}\subset M_{\sig}$ and (\ref{eq:rflrho}).

In this subsection the ordinals $\gam_{0}$ and $b_{1}$ will be fixed, and
we write $\vdash^{a}_{c,\xi}$ for $\vdash^{a}_{c,\xi,\gam_{0}, b_{1}}$.

\blem\label{lem:tautology}{\rm (Tautology)}
Let $\{\gam,\gam_{0},\mS\}\cup\sfk^{\mT}(A)\subset\calh_{\gam}[\Tht_{{\tt Q}(\mT)}]\cap
\calh_{\gam}[\Tht_{\partial{\tt Q}}]$ for every 
$\mT\in dom({\tt Q})\subset(\mS+1)$, 
$\sigma\in\{{\tt u}\}\cup\bigcup{\tt Q}$ and $\sfk(A)\subset M_{\sig}$.
Then
$(\mathcal{H}_{\gamma},\Tht,{\tt Q}
)\vdash^{2d}_{0,\mS}
\lnot A^{(\sigma)}, A^{(\sigma)}$ holds for $d=\max\{\mS,\mbox{{\rm rk}}(A)\}$.

\elem

\blem\label{lem:inversionreg}{\rm (Inversion)}
Let  $A\simeq \bigwedge(A_{\iota})_{\iota\in J}$ and
$(\mathcal{H},\Theta,{\tt Q})\vdash^{a}_{c,\xi}\Gamma$ with $A^{(\rho)}\in\Gamma$
and there is no $\mS\in SSt$ such that $\rk(A)<\mS\leq\xi$.
Then for any $\iota\in [{\tt Q}]_{A^{(\rho)}}J$, 
$(\mathcal{H},\Theta\cup{\sf k}(\iota), {\tt Q}
)
\vdash^{a}_{c,\xi}
\Gamma,\left(A_{\iota}\right)^{(\rho)}$.
\elem
\bprf
We need to assume that  there is no $\mS\in SSt$ such that $\rk(A)<\mS\leq\xi$
due to (Taut).
\eprf

\blem\label{prp:reductionpi1}{\rm (Reduction)}
Let $C\simeq\bigvee(C_{\iota})_{\iota\in J}$ and
$\Ome\leq\mbox{{\rm rk}}(C)\leq c$.
Assume
$(\mathcal{H}_{\gamma},\Theta, {\tt Q}\vdash^{a}_{c,\xi}\Gamma,\lnot C^{(\tau)}$
and
$(\mathcal{H}_{\gamma},\Theta,{\tt Q})\vdash^{b}_{c,\xi}C^{(\tau)},\Gamma$
with $SSt\cap(c,\xi]=\emptyset$.

Then
$(\mathcal{H}_{\gamma},\Theta,{\tt Q})\vdash^{a+b}_{c,\xi}\Gamma$.
\elem

\blem\label{lem:predcereg}{\rm (Cut-elimination)}
If 
$(\mathcal{H}_{\gamma},\Theta,{\tt Q})\vdash^{a}_{c+c_{1},\xi}\Gamma$
with $\Ome\leq c<\mI$, $\fal\mS\in dom({\tt Q})(c\in\calh_{\gam}[\Tht_{{\tt Q}(\mS)}]\cap
\calh_{\gam}[\Tht_{\partial{\tt Q}}])$ and $SSt\cap(c,\xi]=\emptyset$,
then $(\mathcal{H}_{\gamma},\Theta,{\tt Q})\vdash^{\vphi_{c_{1}}(a)}_{c,\xi}\Gamma$.
\elem

\blem\label{lem:Omegacollpase}{\rm (Collapsing)}
Let
$\Gamma\subset\Sigma(\Ome)$ 
be a sets of uncapped formulas.
Suppose
$\Theta\subset
\mathcal{H}_{\gamma}(\psi_{\Ome}(\gamma))$
 and
$
(\mathcal{H}_{\gamma},\Theta,\emptyset
)\vdash^{a}_{\Ome,0}\Gamma
$.
Let
$\beta=\psi_{\Ome}(\hat{a})$ with $\hat{a}=\gamma+\omega^{a}<\gamma_{0}$.
Then 
$(\mathcal{H}_{\hat{a}+1},\Theta, \emptyset)
\vdash^{\beta}_{\beta,0}
\Gamma^{(\beta,\Ome)}$ holds.

\elem

\subsection{Eliminations of stable ordinals}
\blem\label{lem:Mostowskicollpasecap}{\rm (Capping)}
Let
$\Gamma\cup\Pi\subset\Delta_{0}(\mathbb{I})$ be a set of uncapped formulas.
Suppose
$
(\mathcal{H}_{\gamma},\Theta;{\tt Q}_{\Pi}
)\vdash^{* a}_{c,\gam_{0}}\Gamma;\Pi^{[\cdot]}
$, where
$a,c<\mathbb{I}$, $dom({\tt Q}_{\Pi})\subset c$,
 $\Gam=\Gam_{{\tt u}}\cup\bigcup_{\mS\in dom({\tt Q}_{\Pi})}\Gam_{\mS}$,
 $\Pi^{[\cdot]}=\bigcup_{(\mS,\sigma)\in {\tt Q}_{\Pi}}\Pi_{\sigma}^{[\sigma/\mathbb{S}]}$.

For each $\mS\in dom({\tt Q}_{\Pi})$, 
let 
$\rho_{\mS}=\psi_{\mathbb{S}}^{g_{\mS}}(\delta_{\mS})$
 be an ordinal with an ordinal $\del_{\mS}\in\calh_{\gam}[\Tht]$ and
 a special finite function 
$g_{\mS}=m(\rho_{\mS}):\mI\to\vphi_{\mI}(0)$ such that
${\rm supp}(g_{\mS})=\{c\}$ with $g_{\mS}(c)=\alpha_{\mS}+\mathbb{I}$,
$\mathbb{I}(2a+1)\leq\alpha_{\mS}+\mathbb{I}$,
$SC_{\mI}(g_{\mS})=SC_{\mI}(c,\alp_{\mS})\subset 
\mathcal{H}_{0}(SC_{\mI}(\delta_{\mS}))\cap\mathcal{H}_{\gam}[\Tht]$,
cf.\,(\ref{eq:notationsystem.112.4}) and (\ref{eq:controlderrhocap}).
Also let $\widehat{\Pi}=\bigcup_{(\mS,\sig)\in{\tt Q}_{\Pi}}\Pi_{\sig}^{(\sig)}$, 
$\widehat{\Gam}=\Gam_{{\tt u}}^{({\tt u})}\cup\bigcup_{\mS\in dom({\tt Q}_{\Pi})}\Gam_{\mS}^{(\rho_{\mS})}$.

Let ${\tt Q}$ be a finite family for $\gamma_{0}\geq\gam$ such that
${\tt Q}(\mS)={\tt Q}_{\Pi}(\mS)\cup\{\rho_{\mS}\}$ for $\mS\in dom({\tt Q}_{\Pi})=dom({\tt Q})$,
$\rho_{\mS}
\in\mathcal{H}_{\gamma_{\mS}^{{\tt Q}}+\mI}(0)$ for $\mS\in dom({\tt Q})$,
and
$\alpha_{\mS}+\mathbb{I}\leq\gamma_{\mS}^{{\tt Q}}\leq\delta_{\mS}<\gamma_{\mS}^{{\tt Q}}+\mI$.
Also $e^{{\tt Q}}_{\mS}=\rho_{\mS}+1$.

Assume 
$\fal\mS\in dom({\tt Q}_{\Pi})(\gam_{\mS}^{{\tt Q}}\in\calh_{\gam}[\Tht])$,
${\tt Q}_{\Pi}(\mS)\subset\rho_{\mS}$,
$\Tht\cup\{\mS\}\subset M_{\rho_{\mS}}$,
${\tt p}_{0}(\sig)\leq{\tt p}_{0}(\rho_{\mS})=\del_{\mS}$ 
and $SC_{\mI}(m(\sig))\subset\calh_{\gam_{\mS}^{{\tt Q}}}[\Tht\cup\{\mS\}]$
for each $(\mS,\sig)\in{\tt Q}_{\Pi}$, 
$\fal\{\mU<\mS\}\subset dom({\tt Q}_{\Pi})(\rho_{\mS}\in M_{\rho_{\mU}})$,
and ${\tt Q}$ has gaps $2^{a}$.

Then
$(\mathcal{H}_{\gamma},\Theta_{\Pi}, {\tt Q})
\vdash^{a}_{c,c,\gam_{0},c}
\widehat{\Gamma},\widehat{\Pi}$ holds for
$\Tht_{\Pi}=\Tht\cup\bigcup{\tt Q}_{\Pi}$.

\elem

\brem\label{rem:Mostowskicollpasecap}
{\rm 
When $\alp_{\mS}=\mI(2a)$ and $\Tht=\emptyset$, 
$\delta_{\mS}<\gam_{\mS}^{{\tt Q}}+\mI$ denotes the natural sum
$\gam_{\mS}^{{\tt Q}}\#a\#c$.
Then
$\Tht\cup\{\mS\}\subset M_{\rho_{\mS}}$ and
$\{a,c\}\subset \mathcal{H}_{0}(SC_{\mI}(\del_{\mS}))$.
Hence
 (\ref{eq:notationsystem.112.4}) is enjoyed for $\rho_{\mS}$.
 Namely 
$SC_{\mI}(g_{\mS})=\{c,\alp_{\mS}+\mI\}\subset 
\mathcal{H}_{0}(SC_{\mI}(\del_{\mS}))\subset\mathcal{H}_{\del_{\mS}}(SC_{\mI}(\del_{\mS}))$
holds.

Let $\mU\in dom({\tt Q}_{\Pi})\cap\mS$.
We have $\{\gam_{0},\mS,a,c\}\subset\calh_{\gam}[\Tht_{{\tt Q}_{\Pi}(\mU)}]$
by (\ref{eq:controlderac}).
We intend to be $\gam_{\mS}^{{\tt Q}}=\gam_{0}+\mI\cdot 2^{a}\cdot n$
for $n=\#\{\mT\in dom({\tt Q}): \mT\geq \mS\}$.
Then
$\{\mS,a,c,\gam_{\mS}^{{\tt Q}}\}\subset\calh_{\gam}[\Tht_{{\tt Q}_{\Pi}(\mU)}]\cap
\calh_{\gam}[\Tht_{\partial{\tt Q}}]$
for (\ref{eq:controldercapace}) and (\ref{eq:controlderScapgam}).

On the other hand we have ${\tt Q}_{\Pi}(\mS)\subset\rho_{\mS}$, and
$\rho_{\mS}=\max({\tt Q}(\mS))$, i.e., $\partial{\tt Q}=\{\rho_{\mS}: \mS\in dom({\tt Q}_{\Pi})\}$.
Also $\{\mS,\del_{\mS}\}\cup SC_{\mI}(g_{\mS})\subset
\calh_{0}(\{\mS,a,c,\gam_{\mS}^{{\tt Q}}\}\cup \Tht)
\subset M_{\rho_{\mU}}=\calh_{\del_{\mU}}(\rho_{\mU})$ for $\mU\leq\mS$.
Therefore $\rho_{\mS}\in M_{\rho_{\mU}}$ for $\mU<\mS$ by 
$\del_{\mS},\gam_{\mS}^{{\tt Q}}+\mI\leq\gam_{\mU}^{{\tt Q}}$.
Moreover  $\rho_{\mS}\in M_{\rho_{\mU}}$ for $\mU>\mS$ since
$\rho_{\mS}<\mS<\rho_{\mU}$.
}
\erem
\textbf{Proof} of Lemma \ref{lem:Mostowskicollpasecap}.
This is seen
by induction on $a$ as in Capping \ref{lem:capping}.
Let us 
write $\vdash^{a}_{c}$ for $\vdash^{a}_{c,c,\gam_{0},c}$ in the proof.
By assumptions we have
${\tt Q}_{\Pi}(\mS)\subset\rho_{\mS}$ and
$\Tht\subset M_{\rho_{\mS}}$.
Hence $\Tht=\Tht^{(\rho_{\mS})}=\Tht_{\partial{\tt Q}}$ and
$\Tht_{{\tt Q}_{\Pi}(\mS)}=\Tht_{{\tt Q}(\mS)}$.
On the other hand we have
$\sfk(\Gam)\subset\calh_{\gam}[\Tht]$ and for
$\sig\in\bigcup{\tt Q}_{\Pi}$, $\sfk(\Pi_{\sig})\subset\calh_{\gam}[\Tht^{(\sig)}]$
by (\ref{eq:controlder}).
Therefore (\ref{eq:controldercap}) and (\ref{eq:controlderScap}) are enjoyed.
We have $\{\gam,a,c,\gam_{0},\gam_{\mS}^{{\tt Q}},\mS\}\subset\calh_{\gam}[\Tht_{{\tt Q}_{\Pi}(\mU)}]
\cap\calh_{\gam}[\Tht_{\partial{\tt Q}}]$ 
for
every $\{\mU\leq\mS\}\subset dom({\tt Q})=dom({\tt Q}_{\Pi})$ by the assumption, 
(\ref{eq:controlderac}) and (\ref{eq:controlderS}).
Hence (\ref{eq:controldercapace}) and (\ref{eq:controlderScapgam}) are enjoyed.
Moreover for (\ref{eq:controlderrhocap}) we have 
$SC_{\mI}(m(\rho_{\mS}))\subset\calh_{\gam}[\Tht]$ and $\gam\leq\gam_{\mS}^{{\tt Q}}$.
\\
\textbf{Case 1}.
First consider the case when the last inference is a $({\rm stbl}(\mS))$: 
We have a successor stable ordinal $\mS$,
an ordinal $a_{0}<a$,
a $\bigwedge$-formula
$B(0)\in\Delta_{0}(\mathbb{S})$, and
a term $u\in Tm(\mI)$ with $\mS\leq\rk(B(u))<c$.

For every ordinal $\sigma$ such that
$\Theta\cup\{\mS\}\subset M_{\sig}$ and ${\tt p}_{0}(\sig)\geq\gam_{0}$
{\small
\[
\infer{(\calh_{\gam},\Tht;{\tt Q}_{\Pi})\vdash^{* a}_{c}\Gam;\Pi^{[\cdot]}}
{
(\mathcal{H}_{\gamma},\Tht, {\tt Q}_{\Pi})\vdash^{* a_{0}}_{c}
\Gamma, B(u);\Pi^{[\cdot]}
&
(\mathcal{H}_{\gamma},\Tht\cup\{\mS,\sig\};{\tt Q}_{\Pi}\cup\{(\mS,\sig)\}
)\vdash^{* a_{0}}_{c}
\Gamma;\lnot B(u)^{[\sigma/\mS]}, \Pi^{[\cdot]}
}
\]
}

Let $h$ be a special finite function such that ${\rm supp}(h)=\{c\}$ and
$h(c)=\mathbb{I}(2a_{0}+1)$.
Then $h_{c}=(g_{\mS})_{c}=\emptyset$ and $h^{c}<_{\mI}^{c}(g_{\mS})^{\prime}(c)$
by $h(c)=\mathbb{I}(2a_{0}+1)<\mathbb{I}(2a)\leq\alpha_{0}=(g_{\mS})^{\prime}(c)$.
Let ${\tt R}={\tt Q}\cup\{(\mS,\rho_{\mS})\}$
and $\sigma\in H^{{\tt R}}_{\rho_{\mS}}(h,c,\gam^{{\tt R}}_{\mS},\Tht^{(\rho_{\mS})}\cup\{\mS\}\cup\Tht_{\partial{\tt Q}})$, 
where
$\Tht^{(\rho_{\mS})}\cup\Tht_{\partial{\tt Q}}=\Tht$.

For example let
$\sig=\psi_{\rho_{\mS}}^{h}(\del_{\mS}+\eta)$ with $\eta=\max(\{1\}\cup E_{\mS}(\Tht))$.
We obtain $\Tht\cup\{\mS\}\subset \calh_{\del_{\mS}}(\sig)=M_{\sig}$ by $\Tht\cup\{\mS\}\subset M_{\rho}$, and
$\{\del_{\mS},a_{0},c\}\subset\mathcal{H}_{\gam}[\Theta]$.
Let $\rho_{\mU}\in\partial{{\tt R}}$. We claim that $\sig\in M_{\rho_{\mU}}$.
If $\mU\geq\mS$, then $\sig<\rho_{\mU}$.
Let $\mU<\mS$. Then we have $\rho_{\mS}\in M_{\rho_{\mU}}$ by the assumption,
and $\sig\in M_{\rho_{\mU}}$ follows from 
$\{c,a_{0},\del_{\mS}\}\cup\Tht\subset\calh_{\gam}[\Tht]\subset\calh_{\del_{\mU}}(\rho_{\mU})$
and $\del_{\mS}+\eta<\gam_{\mS}^{{\tt Q}}+\mI\leq\gam_{\mU}^{{\tt Q}}\leq\del_{\mU}$.
Therefore $\sigma\in H^{{\tt R}}_{\rho_{\mS}}(h,c,\gam^{{\tt R}}_{\mS},\Tht^{(\rho_{\mS})}\cup\{\mS\}\cup
\Tht_{\partial{\tt Q}})$.

Since ${\tt Q}$ is assumed to have gaps $2^{a}$,
we may assume that
${\tt R}\cup\{(\mS,\sigma)\}$ as well as ${\tt R}$
has gaps $2^{a_{0}}$.

IH yields
$
(\mathcal{H}_{\gam},\Tht_{\Pi},{\tt R})\vdash^{a_{0}}_{c}
\widehat{\Gamma},B(u)^{(\rho_{\mS})},\widehat{\Pi}
$, and for $u^{[\sigma/\mS]}\in Tm(\mS)$ and 
$B(u^{[\sigma,\mS]})\equiv B(u)^{[\sigma,\mathbb{S}]}$,
$(\mathcal{H}_{\gamma},\Tht_{\Pi}\cup\{\mS,\sig\},{\tt R}\cup\{(\mS,\sigma)\}
)\vdash^{a_{0}}_{c}
\widehat{\Gamma},\lnot B(u)^{(\sigma)},\widehat{\Pi}
$ follows, where $\rho_{\mS}>\sig\in M_{\rho_{\mS}}$ and we have by (\ref{eq:controlderS}),
$\sfk(B(u))\subset\calh_{\gam}[\Tht_{{\tt Q}_{\Pi}(\mT)}]\cap\calh_{\gam}[\Tht_{\partial{\tt Q}}]$
if $\rk(B(u))\geq\mT$.
Hence $\sfk(B(u))\subset\calh_{\gam}[\Tht_{{\tt R}(\mT)}]\cap\calh_{\gam}[\Tht_{\partial{\tt Q}}]$ by $\Tht_{{\tt R}(\mT)}=\Tht_{{\tt Q}_{\Pi}(\mT)}$
for (\ref{eq:controlderS}).
Moreover
we have 
$\mS\in\calh_{\gam}[\Tht_{{\tt Q}(\mT)}]$ for every $\mT<c$,
$\Tht_{{\tt Q}_{\Pi}(\mS)}\cup\Tht_{\partial{\tt Q}}\subset M_{\rho_{\mS}}$ for (\ref{eq:rflrho}),
$\rho_{\mS}<e^{{\tt Q}}_{\mS}$
 for (r0),
$\rk(B(u))<c$ and $s(\rho_{\mS})\leq c$ for (r1).

We obtain by an inference $({\rm rfl}_{\mS}(\rho_{\mS},c,h,c))$
{\small
\[
\infer
{
(\mathcal{H}_{\gam},\Tht_{\Pi},{\tt Q})\vdash^{a}_{c}
\widehat{\Gamma},\widehat{\Pi}
}
{
(\mathcal{H}_{\gam},\Tht_{\Pi},{\tt R})\vdash^{a_{0}}_{c}
\widehat{\Gamma},B(u)^{(\rho_{\mS})},\widehat{\Pi}
&
(\mathcal{H}_{\gamma},\Tht_{\Pi}\cup\{\mS,\sig\}, {\tt R}\cup\{(\mS,\sigma)\}
)\vdash^{a_{0}}_{c}
\widehat{\Gamma},\lnot B(u)^{(\sigma)},\widehat{\Pi}
}
\]
}
in the right upper sequents
$\sigma$ ranges over the resolvent class 
$\sigma\in H^{{\tt R}}_{\rho_{\mS}}(h,c,\gam^{{\tt R}}_{\mS},\Tht^{(\rho_{\mS})}\cup\{\mS\}\cup
\Tht_{\partial{\tt Q}})$.
\\
\textbf{Case 2}. 
When the last inference is a 
$(cut)$:
There exist $a_{0}<a$ and $C$
such that $\mbox{{\rm rk}}(C)<c$,
$(\mathcal{H}_{\gamma},\Theta;{\tt Q}_{\Pi})
\vdash^{* a_{0}}_{c}
\Gamma,\lnot C;\Pi^{[\cdot]}$
and
$(\mathcal{H}_{\gamma},\Theta;{\tt Q}_{\Pi})
\vdash^{* a_{0}}_{c}\Gamma,C;\Pi^{[\cdot]}$.
IH followed by a $(cut)$ with an uncapped cut formula $C^{({\tt u})}$ yields the lemma.
\\
\textbf{Case 3}.  
Third the last inference introduces
a $\bigvee$-formula $A$ in $\Gam$.
Let $A\simeq\bigvee\left(A_{\iota}\right)_{\iota\in J}$.
Then $A^{(\rho_{\mS})}\in\Gam_{\mS}^{(\rho_{\mS})}$.
There are an $\iota\in J$, an ordinal
 $a(\iota)<a$
such that
$(\mathcal{H}_{\gamma},\Theta;{\tt Q}_{\Pi})\vdash^{* a(\iota)}_{c}
\Gamma,A_{\iota};\Pi^{[\cdot]}$.
We can assume ${\sf k}(\iota)\subset{\sf k}(A_{\iota})$, and claim that 
$\iota\in [\partial{\tt Q}]J$ with $\rho_{\mS}\in\partial{\tt Q}$.
We obtain
$\sfk(\iota)\subset\calh_{\gam}[\Tht_{\partial{\tt Q}}] \subset M_{\partial{\tt Q}}$
by (\ref{eq:controlder}) for $\Tht_{\partial{\tt Q}}=\Tht$
and
$\gam\leq\gam_{0}\leq\gam^{{\tt Q}}_{\mS}\leq\del_{\mS}\leq{\tt p}_{0}(\rho_{\mS})$.

IH yields $(\mathcal{H}_{\gamma},\Theta,{\tt Q})
\vdash^{a(\iota)}_{c}
\widehat{\Gamma},
(A_{\iota})^{(\rho_{\mS})},\widehat{\Pi}$.
$(\mathcal{H}_{\gamma},\Theta,{\tt Q})
\vdash^{a}_{c}
\widehat{\Gamma},\widehat{\Pi}$ 
follows from a $(\bigvee)$.

Other cases are seen from IH as in Capping \ref{lem:capping}.
\eprf

\blem\label{lem:recapping}{\rm (Recapping)}\\
Let $\mS$ be a successor stable ordinal,
$
(\mathcal{H}_{\gamma},\Theta,{\tt Q})
\vdash^{a}_{c_{1},\mS,\gam_{0},b_{2}}
\Pi,\widehat{\Gamma}
$ with a finite family ${\tt Q}$ for $\gam_{0}, b_{2}$, $\Gamma\cup\Pi\subset\Delta_{0}(\mathbb{I})$,
and $\widehat{\Gam}=\bigcup\{\Gam_{\rho}^{(\rho)} :\rho\in{\tt Q}^{t}(\mS)\}$, where
each $\tht\in \widehat{\Gam}$ is either a $\bigvee$-formula or $\rk(\tht)<\mS$,
${\tt Q}^{t}\subset{\tt Q}$ 
such that ${\tt Q}^{t}(\mS)\subset{\tt Q}(\mS)$ with $dom({\tt Q}^{t})\subset\{\mS\}$
and $\fal\rho\in{\tt Q}^{t}(\mS)(s(\rho)>\mS)$, and
${\tt Q}^{f}$ is a family such that ${\tt Q}^{f}(\mS)={\tt Q}(\mS)\setm{\tt Q}^{t}(\mS)$ and
${\tt Q}^{f}(\mT)={\tt Q}(\mT)$ for $\mT\neq\mS$.
$\Pi$ is a set of formulas such that $\tau\in\{{\tt u}\}\cup\bigcup{\tt Q}^{f}$ for every $A^{(\tau)}\in\Pi$.

Let $\max\{s(\rho):\rho\in{\tt Q}^{t}(\mS)\}\leq b_{1}$ and $\ome(b,a)=\ome^{\ome^{b}}a$.
For each $\rho\in{\tt Q}^{t}(\mS)$, let $\mS\leq b^{(\rho)}\in\calh_{\gam}[\Tht^{(\rho)}]
\cap\calh_{\gam}[\Tht_{\partial{\tt Q}}]$ 
with $\rk(\Gamma_{\rho})<b^{(\rho)}<s(\rho)$,
and 
$\kap(\rho)$ be ordinals such that
$\kappa(\rho)\in H^{{\tt Q}}_{\rho}(h^{b^{(\rho)}}(m(\rho);\ome(b_{1},a)),b_{2},\gamma_{\mS}^{{\tt Q}},\Theta^{(\rho)}\cup\{\mS\}\cup \Tht_{\partial{\tt Q}})$.
Assume 
$\fal\mT\leq\mS(b_{1}\in\calh_{\gam}[\Tht_{{\tt Q}(\mT)}]\cap\calh_{\gam}[\Tht_{\partial{\tt Q}}])$.

Then 
$
(\mathcal{H}_{\gamma},\Theta,{\tt Q}^{\kappa}
)
\vdash^{\ome(b_{1},a)}_{c_{b_{1}},\mS,\gam_{0} ,b_{2} }
\Pi,\widehat{\Gamma}_{\kappa}
$
holds,
where
$\widehat{\Gamma}_{\kappa}=\bigcup\{\Gam_{\rho}^{(\kap(\rho))}:\rho\in{\tt Q}^{t}(\mS)\}$,
$c_{b_{1}}=\max\{c_{1},b_{1}\}$, 
${\tt Q}^{\kappa}={\tt Q}^{f}\cup\{(\mS,\kap(\rho)): \rho\in{\tt Q}^{t}(\mS)\}$, 
$\gam_{\mT}^{{\tt Q}^{\kap}}=\gam_{\mT}^{{\tt Q}}$,
$e^{{\tt Q}^{\kap}}_{\mT}=e^{{\tt Q}}_{\mT}$ for $\mT\neq\mS$ and
$e^{{\tt Q}^{\kap}}_{\mS}=\max(\{\tau\in{\tt Q}^{f}(\mS): s(\tau)>\mS\}\cup\{\kap(\rho): \rho\in{\tt Q}^{t}(\mS)\})+1$.

$e^{{\tt Q}^{\kap}}_{\mS}<e^{{\tt Q}}_{\mS}$ holds
when
${\tt Q}^{t}=\{(\mS,\rho)\in{\tt Q}: s(\rho)>\mS\}\neq\emptyset$.

\elem
\bprf
This is shown by main induction on $b_{1}$
with subsidiary induction on $a$ as in Recapping \ref{lem:recappingpi11}.
\eprf

\blem\label{lem:main.2}{\rm (Elimination of one stable ordinal)}\\
Let $\mS=\mT^{\dagger}$ be a successor stable ordinal and
$
(\mathcal{H}_{\gamma},\Theta,{\tt Q})
\vdash^{a}_{\mS,\mS,\gam_{0},b_{1}}
\Pi,\widehat{\Gam}
$ with a finite family ${\tt Q}$ for $\gam_{0}$ and
$b_{1}\geq\mS$, $\Pi\subset\Delta_{0}(\mathbb{I})$,
$\Gam\subset\Del_{0}(\mS)$,
$\widehat{\Gam}=\bigcup\{\Gam_{\rho}^{(\rho)}: \rho\in{\tt Q}(\mS)\}$,
and 
${\tt Q}^{t}=\{(\mS,\tau)\in{\tt Q}: s(\tau)>\mS\}$,
${\tt Q}^{f}={\tt Q}\setm{\tt Q}^{t}$.
$\Pi$ is a set of formulas such that 
for each $A^{(\tau)}\in\Pi$, 
$\tau\in\{{\tt u}\}\cup\bigcup_{\mU<\mS}{\tt Q}(\mU)$.

Let $\tilde{a}=\vphi_{b_{1}+e^{{\tt Q}}_{\mS}}(a)$,
${\tt Q}_{1}={\tt Q}\restrict\mS=\{(\mT,\rho)\in{\tt Q}: \mT<\mS\}$ and
$\gam_{1}=\gam_{\mS}^{{\tt Q}}+\mI<\gam_{0}+\mI^{2}$.

Then 
${\tt Q}_{1}$ is a finite family for $\gam_{1},b_{1}$ and
$
(\mathcal{H}_{\gam_{1}},\Theta, {\tt Q}_{1}
)
\vdash^{\tilde{a}}_{\mT,\mT,\gam_{1},b_{1}}
\Pi,\Gam^{({\tt u})}
$
holds for 
$\Gam^{({\tt u})}=\bigcup\{\Gam_{\rho}^{({\tt u})}: \rho\in{\tt Q}(\mS)\}$.

\elem
\bprf
This is seen
by main induction on $e^{{\tt Q}}_{\mS}$ with subsidiary induction on $a$ as in Lemma \ref{lem:sum}.
When $\mS\in dom({\tt Q})$, we have ${\tt Q}(\mS)\subset\calh_{\gam_{1}}$ and
$e^{{\tt Q}}_{\mS}\in\calh_{\gam_{1}}$ for $\gam_{1}=\gam_{\mS}^{{\tt Q}}+\mI$
by Definition \ref{df:caphat}.
${\tt Q}_{1}$ is a finite family for $\gam_{1}, b_{1}$.
Then $\gam_{1}\in\calh_{\gam}[\Tht_{{\tt Q}_{1}(\mT)}]\cap\calh_{\gam}[\Tht_{\partial{\tt Q}}]$ for every $\mT\in dom({\tt Q}_{1})$
by (\ref{eq:controldercapace}).

First assume ${\tt Q}^{t}(\mS)\neq\emptyset$.
For each $\rho\in{\tt Q}^{t}(\mS)$, let
 $\kap(\rho)$ be an ordinal such that
$\kappa(\rho)\in H^{{\tt Q}}_{\rho}(h^{\mS}(m(\rho);\ome(b_{1},a)),b_{1}, \gamma_{\mS}^{{\tt Q}},\Theta^{(\rho)}
\cup\{\mS\}\cup\Tht_{\partial{\tt Q}})$
 with $\ome(b,a)=\ome^{\ome^{b}}a$.
We obtain $(\calh_{\gam},\Tht,{\tt Q}^{\kap})\vdash^{\ome(b_{1},a)}_{b_{1},\mS,\gam_{0},b_{1}}\Pi,\widehat{\Gam}_{\kap}$
by Recapping \ref{lem:recapping}.
Cut-elimination \ref{lem:predcereg} with $SSt\cap(\mS,\mS]=\emptyset$
yields for $a_{1}=\vphi_{b_{1}}(\ome(b_{1},a))$,
 $(\calh_{\gam},\Tht,{\tt Q}^{\kap})\vdash^{a_{1}}_{\mS,\mS,\gam_{0},b_{1}}\Pi,\widehat{\Gam}_{\kap}$,
where 
$e^{{\tt Q}^{\kap}}_{\mS}=\max\{\kap(\rho):\rho\in{\tt Q}^{t}(\mS)\}+1<e^{{\tt Q}}_{\mS}$.
MIH yields
 $(\calh_{\gam_{1}},\Tht,{\tt Q}_{1})\vdash^{\tilde{a}_{1}}_{\mT,\mT,\gam_{1},b_{1}}\Pi,\Gam^{({\tt u})}$, where
$\tilde{a}_{1}=\vphi_{b_{1}+e^{{\tt Q}^{\kap}}_{\mS}}(a_{1})<
\vphi_{b_{1}+e^{{\tt Q}}_{\mS}}(a)$ and $\gam_{1}=\gam_{\mS}^{{\tt Q}}+\mI$.

In what follows assume ${\tt Q}^{t}(\mS)=\emptyset$. 
\\
\textbf{Case 1}.
First let $\{\lnot A^{(\sig)},A^{(\sig)}\}\subset\Pi\cup\widehat{\Gam}$ with $\sig\in\{{\tt u}\}\cup\bigcup{\tt Q}$ and
$d=\rk(A)<\mS$ by (Taut).
If $d<\mT$, then 
$(\mathcal{H}_{\gam_{1}},\Theta, {\tt Q}_{1})\vdash^{\tilde{a}}_{\mT,\mT,\gam_{1},b_{1}}\Pi,\Gamma^{({\tt u})}$
by (Taut).

Let $\mT\leq d<\mS$.
Then
$(\calh_{\gam_{1}},\Tht,{\tt Q}_{1})\vdash^{2d}_{0,\mT,\gam_{1},b_{1}}\Pi,\Gamma^{({\tt u})}$ by Tautology \ref{lem:tautology} and
$(\calh_{\gam_{1}},\Tht,{\tt Q}_{1})\vdash^{\tilde{a}}_{0,\mT,\gam_{1},b_{1}}\Pi,\Gamma^{({\tt u})}$
by $\tilde{a}>\mS>d$.
\\
\textbf{Case 2}.
Second consider the case when the last inference is a $({\rm rfl}_{\mU}(\rho,d,f,b_{1}))$.
If $\mU\leq\mT$, then SIH followed by a $({\rm rfl}_{\mU}(\rho,d,f,b_{1}))$ yields the lemma.
Let $\mU=\mS$.

Let 
$g=m(\rho)$ and $s(\rho)\geq d\in {\rm supp}(g)$.
Let ${\tt R}={\tt Q}\cup\{(\mS,\rho)\}$ and $\gam_{1}=\gam_{\mS}^{{\tt R}}+\mI$.
We have a sequent $\Delta\subset\bigvee_{\mS}(d)$
and an ordinal $a_{0}<a$ such that 
${\rm rk}(\Delta)<d\leq s(\rho)$ and
$
(\mathcal{H}_{\gamma},\Theta,
{\tt R}
)
\vdash^{a_{0}}_{\mS,\mS,\gam_{0},b_{1}}
\Pi,\widehat{\Gam},
\lnot\delta^{(\rho)}
$
for each $\delta\in\Delta$.
On the other hand we have 
$
(\mathcal{H}_{\gamma},\Theta\cup\{\sig\},
{\tt R}\cup\{(\mS,\sigma)\}
)
\vdash^{a_{0}}_{\mS,\mS,\gam_{0},b_{1}}
\Pi, \widehat{\Gam},\Delta^{(\sigma)}
$, where
$\sigma\in H^{{\tt Q}}_{\rho}(f,b_{1},\gamma_{\mS}^{{\tt R}},\Theta^{(\rho)}\cup\{\mS\}\cup\Tht_{\partial{\tt Q}})$,
$f$ is a special finite function such that $s(f)\leq b_{1}$,
$f_{d}=g_{d}$, $f^{d}<^{d}g^{\prime}(d)$ and
$SC_{\mI}(f)\subset\mathcal{H}_{\gamma_{\mS}^{{\tt R}}}[\Theta^{(\rho)}]$.
\\
\textbf{Case 2.1}. $s(\rho)\leq\mS$: Then $\Del\subset\Del_{0}(\mS)$.
Let $\tilde{a}_{0}=\vphi_{b_{1}+e^{{\tt R}}_{\mS}}(a_{0})$.
SIH yields $
(\mathcal{H}_{\gamma_{1}},\Theta,
{\tt Q}_{1}
)
\vdash^{\tilde{a}_{0}}_{\mT,\mT,\gam_{1},b_{1}}
\Pi,\Gam^{({\tt u})},
\lnot\delta^{({\tt u})}
$
for each $\delta\in\Delta$, and
$
(\mathcal{H}_{\gamma_{1}},\Theta\cup\{\sig\},
{\tt Q}_{1}
)
\vdash^{\tilde{a}_{0}}_{\mT,\mT,\gam_{1},b_{1}}
\Pi, \Gam^{({\tt u})},\Delta^{({\tt u})}
$ for $\sig\in\calh_{\gam_{\mS}^{{\tt R}}+\mI}=\calh_{\gam_{1}}$.
We obtain
$
(\mathcal{H}_{\gamma_{1}},\Theta,
{\tt Q}_{1}
)
\vdash^{\tilde{a}_{0}+p}_{\mS,\mT,\gam_{1},b_{1}}
\Pi,\Gam^{({\tt u})}
$ by several $(cut)$'s for a $p<\ome$.
Cut-elimination \ref{lem:predcereg} with $SSt\cap(\mT,\mT]=\emptyset$ yields
$
(\mathcal{H}_{\gamma_{1}},\Theta,
{\tt Q}_{1}
)
\vdash^{\vphi_{\mS}(\tilde{a}_{0}+p)}_{\mT,\mT,\gam_{1},b_{1}}
\Pi,\Gam^{({\tt u})}
$,
where $\vphi_{\mS}(\tilde{a}_{0}+p)<\tilde{a}=\vphi_{b_{1}+e^{{\tt Q}}_{\mS}}(a)$ by $b_{1}+e_{\mS}^{{\tt Q}}>\mS$.
\\
\textbf{Case 2.2}. $s(\rho)>\mS$: Then $\mS\not\in dom({\tt Q})$ and $\Gam=\emptyset$.
We have
$
(\mathcal{H}_{\gamma},\Theta,{\tt R})
\vdash^{a}_{\mS,\mS,\gam_{0},b_{1}}
\Pi
$.
Let ${\tt R}^{t}=\{(\mS,\rho)\}$.
Recapping \ref{lem:recapping} yields
$
(\mathcal{H}_{\gamma},\Theta,{\tt R}^{\kap})
\vdash^{\ome(b_{1},a)}_{\mS,\mS,\gam_{0},b_{1}}
\Pi
$ and $e^{{\tt R}^{\kap}}_{\mS}=\kap+1<\rho<e^{{\tt R}}_{\mS}$.
MIH yields
$
(\mathcal{H}_{\gamma_{1}},\Theta,{\tt Q}_{1})
\vdash^{a_{1}}_{\mT,\mT,\gam_{1},b_{1}}
\Pi
$ with
$a_{1}=\vphi_{b_{1}+e^{{\tt R}^{\kap}}_{\mS}}(\ome(b_{1},a))<\vphi_{b_{1}+e^{{\tt Q}}_{\mS}}(a)=\tilde{a}$ by 
$e^{{\tt R}^{\kap}}_{\mS}<\mS=e^{{\tt Q}}_{\mS}$.
\\
\textbf{Case 3}.
The last inference is a $(\bigwedge)$:
We have $a(\iota)<a$, $A^{(\rho)}\in\widehat{\Gam}$ and for each 
$\iota\in [{\tt Q}]_{A^{(\rho)}}J$ with
$A\simeq\bigwedge(A_{\iota})_{\iota\in J}$, we have
$
(\mathcal{H}_{\gamma},\Tht\cup\sfk(\iota),{\tt Q})
\vdash^{a(\iota)}_{\mS,\mS,\gam_{0},b_{1}}
\Pi,\widehat{\Gam},(A_{\iota})^{(\rho)}
$.
Since $A\in\Del_{0}(\mS)$, we obtain $\sfk(A)\subset\calh_{\gam}[\Tht^{(\rho)}]\cap\mS\subset
M_{\rho}\cap\mS=\rho$ for $\rho\in{\tt Q}(\mS)$.
This means 
$A\in\Del_{0}(\rho)$, and
$[\rho]J=J$. Hence
$[{\tt Q}]_{A^{(\rho)}}J=[{\tt Q}_{1}]_{A^{({\tt u})}}J$.
SIH yields
$
(\mathcal{H}_{\gam_{1}},\Tht\cup\sfk(\iota),{\tt Q}_{1})
\vdash^{\tilde{a}(\iota)}_{\mT,\mT,\gam_{1},b_{1}}
\Pi,\Gam^{({\tt u})},(A_{\iota})^{({\tt u})}
$
for each $\iota\in [{\tt Q}]_{A}J$, where $\tilde{a}(\iota)=\vphi_{b_{1}+e^{{\tt Q}}_{\mS}}(b+a(\iota))<\tilde{a}$.
A $(\bigwedge)$ yields
$
(\mathcal{H}_{\gam_{1}},\Tht,{\tt Q}_{1})
\vdash^{\tilde{a}}_{\mT,\mT,\gam_{1},b_{1}}
\Pi,\Gam^{({\tt u})}
$.

Other cases are seen from SIH.
\eprf

\bdf\label{df:srank}
{\rm
We define the \textit{S-rank} ${\rm srk}(A^{(\rho)})$ of a capped formula $A^{(\rho)}$ as follows.
Let ${\rm srk}(A^{({\tt u})})=0$, and
${\rm srk}(A^{(\rho)})=\mS$ for $\rho\prec\mS\in SSt$.

${\rm srk}(\Gam)=\max\{{\rm srk}(A^{(\rho)}): A^{(\rho)}\in\Gam\}$.
}
\edf

\blem\label{lem:main.1}{\rm (Elimination of stable ordinals)}\\
Suppose
$
(\mathcal{H}_{\gamma},\Tht,{\tt Q})
\vdash^{a}_{\xi,\xi,\gam_{0},b_{1}}
\Gamma
$ and
${\rm srk}(\Gam)\leq\mS<\xi\leq b_{1}<\mI$, where $\mS$ is either a stable ordinal or $\mS=\Ome$ such that
$\fal\mU\in dom({\tt Q}_{\mS})(\mS\in\calh_{\gam}[\Tht_{{\tt Q}(\mU)}]\cap\calh_{\gam}[\Tht_{\partial{\tt Q}}])$ for ${\tt Q}_{\mS}={\tt Q}\restrict\mS$.

Then there exists an ordinal $\gam_{0}\leq\gam_{\mS}<\gam_{0}+\mI^{2}$ such that
${\tt Q}_{\mS}$ is a finite family for $\gam_{\mS},b_{1}$ and
$
(\mathcal{H}_{\gam_{\mS}},\Tht,{\tt Q}_{\mS})
\vdash^{f(\xi,a)}_{\mS,\mS,\gam_{\mS},b_{1}}
\Gamma
$ holds for $f(\xi,a)=\vphi_{b_{1}+\xi+1}(a)$.
\elem
\bprf
By main induction on $\xi$ with subsidiary induction on $a$.
(\ref{eq:controldercapace}) in \\
$
(\mathcal{H}_{\gam_{\mS}},\Tht,{\tt Q}_{\mS})
\vdash^{f(\xi,a)}_{\mS,\mS,\gam_{\mS},b_{1}}
\Gamma
$
follows from (\ref{eq:controldercapace}) and (\ref{eq:controlderScapgam}) in $
(\mathcal{H}_{\gamma},\Tht,{\tt Q})
\vdash^{a}_{\xi,\xi,\gam_{0},b_{1}}
\Gamma
$.
\\
\textbf{Case 1}.
Consider the case when the last inference is a $({\rm rfl}_{\mT}(\rho,d,f,b_{1}))$
for a $\mT=\mU^{\dagger}\leq\xi$.
If $\mT\leq\mS$, then SIH yields the lemma.
Let $\mS<\mT\in dom({\tt R})$ for ${\tt R}={\tt Q}\cup\{(\mT,\rho)\}$.
We have $\fal\mU\in dom({\tt Q}_{\mT})(\mT\in\calh_{\gam}[\Tht_{{\tt Q}(\mU)}]\cap\calh_{\gam}[\Tht_{\partial{\tt Q}}])$ by (\ref{eq:controldercapace}).
Let $\Delta$ be a finite set of sentences such that
 $(\mathcal{H}_{\gamma},\Tht,{\tt R})\vdash^{a_{0}}_{\xi,\xi,\gam_{0},b_{1}}
\Gam, \lnot\delta^{(\rho)}$
for each $\delta\in\Delta$, and
$
(\mathcal{H}_{\gamma},\Tht,
{\tt R}\cup\{(\mT,\sigma)\}
)\vdash^{a_{0}}_{\xi,\xi,\gam_{0},b_{1}}
\Gam,
\Delta^{(\sigma)}
$
for each $\sigma\in H^{{\tt Q}}_{\rho}(f,b_{1}\gamma_{\mT}^{{\tt R}},\Tht^{(\rho)}\cup\{\mT\}\cup\Tht_{\partial{\tt Q}})$, and 
$a_{0}<a$.
We have ${\rm srk}(\delta^{(\rho)})={\rm srk}(\Del^{(\sig)})=\mT$.
By SIH there exists a $\gam_{\mT}<\gam_{0}+\mI^{2}$ such that for $a_{1}=f(\xi,a_{0})=\vphi_{b_{1}+\xi+1}(a_{0})$, 
 $(\mathcal{H}_{\gam_{\mT}},\Tht,{\tt Q}_{\mT})\vdash^{a_{1}}_{\mT,\mT,\gam_{\mT},b_{1}}
\Gam, \lnot\delta^{(\rho)}$
for each $\delta\in\Delta$, and
$
(\mathcal{H}_{\gam_{\mT}},\Tht,
{\tt Q}_{\mT}\cup\{(\mT,\sigma)\}
)\vdash^{a_{1}}_{\mT,\mT,\gam_{\mT},b_{1}}
\Gam,
\Delta^{(\sigma)}
$.
$({\rm rfl}_{\mT}(\rho,d,f,b_{1}))$ yields
$
(\mathcal{H}_{\gam_{\mT}},\Tht,
{\tt Q}_{\mT}
)\vdash^{a_{2}}_{\mT,\mT,\gam_{\mT},b_{1}}
\Gam
$ for $a_{2}=a_{1}+1$.

On the other hand we have
${\rm srk}(\Gam)\leq\mS<\mT=\mU^{\dagger}\leq\xi$.
By Lemma \ref{lem:main.2} pick a $\gam_{\mU}<\gam_{\mT}+\mI^{2}=\gam_{0}+\mI^{2}$ such that
$
(\mathcal{H}_{\gam_{\mU}},\Tht,
{\tt Q}_{\mU}
)\vdash^{a_{3}}_{\mU,\mU,\gam_{\mU},b_{1}}
\Gam
$, where
$a_{3}=\vphi_{b_{1}+e^{{\tt Q}_{1}}_{\mT}}(a_{2})=\vphi_{b_{1}+e^{{\tt Q}_{1}}_{\mT}}(f(\xi,a_{0})+1)<
\vphi_{b_{1}+\xi+1}(a)=f(\xi,a)$ by $e^{{\tt Q}_{1}}_{\mT}\leq\mT\leq\xi$.
If $\mS=\mU$, then we are done.
Let $\mS<\mU$ with $\mU<\xi$.
Then by MIH pick a $\gam_{\mS}$ such that
$
(\mathcal{H}_{\gam_{\mS}},\Tht,
{\tt Q}_{\mS}
)\vdash^{a_{4}}_{\mS,\mS,\gam_{\mS},b_{1}}
\Gam
$ for
$a_{4}=f(\mU,a_{3})=\vphi_{b_{1}+\mU+1}(a_{3})<\vphi_{b_{1}+\xi+1}(a)=f(\xi,a)$
by $\mU<\xi$.
\\
\textbf{Case 2}. Next consider the case when the last inference is a $(cut)$ of a cut formula
$C^{(\sig)}$ wth $\rk(C)<\xi$ and $\mT={\rm srk}(C^{(\sig)})\leq\xi$.
We have an ordinal $a_{0}<a$ such that 
$
(\mathcal{H}_{\gamma},\Tht,{\tt Q})
\vdash^{a_{0}}_{\xi,\xi,\gam_{0},b_{1}}
\Gamma,\lnot C^{(\sig)}
$
and
$
(\mathcal{H}_{\gamma},\Tht,{\tt Q})
\vdash^{a_{0}}_{\xi,\xi,\gam_{0},b_{1}}
C^{(\sig)},\Gamma
$.

Let $\mU=\max\{\mS,\mT\}$.
First assume $\mU<\xi$.
By SIH pick a $\gam_{\mU}$ such that
$
(\mathcal{H}_{\gamma_{\mU}},\Tht,{\tt Q}_{\mU})
\vdash^{a_{1}}_{\mU,\mU,\gam_{\mU},b_{1}}
\Gamma,\lnot C^{(\sig)}
$
and
$
(\mathcal{H}_{\gamma_{\mU}},\Tht,{\tt Q}_{\mU})
\vdash^{a_{1}}_{\mU,\mU,\gam_{\mU},b_{1}}
C^{(\sig)},\Gamma
$, where $a_{1}=f(\xi,a_{0})=\vphi_{b_{1}+\xi+1}(a_{0})$.
A $(cut)$ yields
$
(\mathcal{H}_{\gamma_{\mU}},\Tht,{\tt Q}_{\mU})
\vdash^{a_{1}+1}_{\xi,\mU,\gam_{\mU},b_{1}}
\Gamma
$.
Cut-elimination \ref{lem:predcereg} with $SSt\cap(\mU,\mU]=\emptyset$ yields
$
(\mathcal{H}_{\gamma_{\mU}},\Tht,{\tt Q}_{\mU})
\vdash^{a_{2}}_{\mU,\mU,\gam_{\mU},b_{1}}
\Gamma
$, where
$a_{2}=\vphi_{\xi}(a_{1}+1)<\vphi_{b_{1}+\xi+1}(a)=f(\xi,a)$ by $\xi<b_{1}+\xi+1$.
If $\mU=\mS$, then we are done.
Let $\mU=\mT>\mS$. By MIH with $\mU<\xi$ we obtain
$
(\mathcal{H}_{\gamma_{\mS}},\Tht,{\tt Q}_{\mS})
\vdash^{a_{3}}_{\mS,\mS,\gam_{\mS},b_{1}}
\Gamma
$
for a $\gam_{\mS}$,
where $a_{3}=f(\mU,a_{2})=\vphi_{b_{1}+\mU+1}(a_{2})<\vphi_{b_{1}+\xi+1}(a)=f(\xi,a)$
by $\mU<\xi$.

Second let $\mT=\mU=\xi=\mW^{\dagger}>\mS$. Then $C\in\Del_{0}(\mT)$.
By Lemma \ref{lem:main.2} pick a $\gam_{\mW}$ such that
$
(\mathcal{H}_{\gam_{\mW}},\Tht,{\tt Q}_{\mW})
\vdash^{\tilde{a}_{0}}_{\mW,\mW,\gam_{\mW},b_{1}}
\Gamma,\lnot C^{({\tt u})}
$
and
$
(\mathcal{H}_{\gam_{\mW}},\Tht,{\tt Q}_{\mW})
\vdash^{\tilde{a}_{0}}_{\mW,\mW,\gam_{\mW}, b_{1}}
C^{({\tt u})},\Gamma
$, where $\tilde{a}_{0}=\vphi_{b_{1}+e^{{\tt Q}}_{\mT}}(a_{0})$.
A $(cut)$ yields
$
(\mathcal{H}_{\gam_{\mW}},\Tht,{\tt Q}_{\mW})
\vdash^{\tilde{a}_{0}+1}_{\mT,\mW,\gam_{\mW},b_{1}}
\Gamma
$, and we obtain
$
(\mathcal{H}_{\gam_{\mW}},\Tht,{\tt Q}_{\mW})
\vdash^{a_{4}}_{\mW,\mW,\gam_{\mW},b_{1}}
\Gamma
$
by Cut-elimination \ref{lem:predcereg}, where
$a_{4}=\vphi_{\mT}(\tilde{a}_{0}+1)$ and $SSt\cap(\mW,\mW]=\emptyset$.
By MIH pick a $\gam_{\mS}$ such that
$
(\mathcal{H}_{\gam_{\mS}},\Tht,{\tt Q}_{\mS})
\vdash^{a_{5}}_{\mS,\mS,\gam_{\mS},b_{1}}
\Gamma
$
for $\mW<\xi$
and $a_{5}=f(\mW,a_{4})=\vphi_{b_{1}+\mW+1}(a_{4})<\vphi_{b_{1}+\xi+1}(a)$
by $\mW<\xi$, $\mT=\xi<b_{1}+\xi+1$, $e^{{\tt Q}}_{\mT}\leq\mT=\xi<\xi+1$ and $a_{0}<a$.
\\
\textbf{Case 3}. There exists an $A$ such that $\{\lnot A^{(\rho)},A^{(\rho)}\}\subset\Gam$
with ${\rm srk}(A^{(\rho)})\leq\mS$
and $d=\rk(A)<\mT\leq\xi$ for a $\mT\in SSt$ by (Taut).
We may assume $d\geq\mS$.
Then
$(\calh_{\gam},\Tht,{\tt Q}_{\mS})\vdash^{2d}_{0,\mS,\gam_{0},b_{1}}\Gamma$ by Tautology \ref{lem:tautology} and the lemma follows from $d<\xi<f(\xi,a)$.

Other cases are seen from SIH.
\eprf

\begin{theorem}\label{thm:2}
Suppose ${\sf KP}\ome+\Pi_{1}\mbox{{\rm -Collection}}+(V=L)\vdash\tht^{L_{\Ome}}$
for a $\Sig_{1}$-sentence $\tht$.
Then 
$L_{\psi_{\Ome}(\veps_{\mI+1})}\models\tht$ holds.
\end{theorem}
\bprf
Let 
$S_{\mI}\vdash\theta^{L_{\Omega}}$ for a $\Sigma$-sentence $\theta$.
By Embedding \ref{th:embedreg} pick an $m>0$ so that 
 $(\mathcal{H}_{\mI},\emptyset; \emptyset)
 \vdash_{\mI+m}^{* \mI\cdot 2+m}
\theta^{L_{\Omega}}$.
Cut-elimination \ref{lem:predcepi1*} yields 
$(\mathcal{H}_{\mI},\emptyset;\emptyset)
\vdash_{\mI}^{* a}\theta^{L_{\Omega}}$
for $a=\omega_{m}(\mI\cdot 2+m)<\omega_{m+1}(\mI+1)$.
Then Collapsing \ref{lem:Kcollpase.1} yields
$(\mathcal{H}_{\hat{a}+1},\emptyset;\emptyset)
\vdash_{\bet}^{* \bet}\theta^{L_{\Omega}}$
for $\bet=\psi_{\mI}(\hat{a})\in LS$ with 
$\hat{a}=\ome^{\mI+a}=\omega_{m+1}(\mI\cdot 2+m)>\bet$.
Capping \ref{lem:Mostowskicollpasecap} then yields
$(\mathcal{H}_{\hat{a}+1},\emptyset,\emptyset)
\vdash_{\bet,\bet,\gam_{0},\bet}^{\bet}\theta^{L_{\Omega}}$
where $\gam_{0}=\hat{a}+1$ and $\theta^{L_{\Omega}}\equiv(\theta^{L_{\Omega}})^{({\tt u})}$.

Let 
$\alp=\vphi_{\bet\cdot 2+1}(\bet)$.
By Lemma \ref{lem:main.1} we obtain
$(\mathcal{H}_{\gam_{\Ome}},\emptyset,\emptyset)
\vdash_{\Ome,\Ome,\gam_{\Ome},\bet}^{\alp}
\theta^{L_{\Omega}}$ for a $\gam_{\Ome}<\gam_{0}+\mI^{2}$.
This means
$(\mathcal{H}_{\gam_{\Ome}},\emptyset,\emptyset)
\vdash_{\Ome,0,\gam_{\Ome},\bet}^{\alp}
\theta^{L_{\Omega}}$.
$(\mathcal{H}_{\gam_{\Ome}+\alp+1},\emptyset,\emptyset)
\vdash_{\del,0,\gam_{\Ome},\bet}^{\del}\theta^{L_{\del}}$ follows from Collapsing \ref{lem:Omegacollpase} 
for
$\del=\psi_{\Ome}(\gam_{\Ome}+\alp)$ with $\ome^{\alp}=\alp$.
Cut-elimination \ref{lem:predcereg} yields 
$(\mathcal{H}_{\gam_{\Ome}+\alp+1},\emptyset,\emptyset)
\vdash_{0,0,\gam_{\Ome},\bet}^{\vphi_{\del}(\del)}\theta^{L_{\del}}$.
We see that $\theta^{L_{\del}}$ is true by induction up to $\vphi_{\del}(\del)$,
where $\del<\psi_{\Ome}(\ome_{m+2}(\mI+1))<\psi_{\Ome}(\veps_{\mI+1})$.

\subsection{Well-foundedness proof in $\Sig^{1}_{3}\mbox{{\rm -DC+BI}}$}

\begin{theorem}\label{th:wf}{\rm \cite{pi1collection}}\\
$\Sig^{1}_{3}\mbox{{\rm -DC+BI}}\vdash Wo[\alp]$
for {\rm each} $\alp<\psi_{\Ome}(\veps_{\mI+1})$.
\end{theorem}

To prove Theorem \ref{th:wf}, let us introduce $1$-distinguished sets $D_{1}[X]$, which is obtained from
Definition \ref{df:3wfdtg32}.\ref{df:3wfdtg.832} of distinguished sets $D[X]$,
first by replacing the next regular $\alp^{+}$ by the next stable $\alp^{\dagger}$, and
second by changing the well-founded part $W(\calc^{\alp}(X))$ to the maximal
distinguished set $\calw_{1}^{\alp}(X)=\bigcup\{P: D_{0}^{\alp}[P;X]\}$ relative to $\alp$ and $X$,
where $P\cap\alp=X\cap\alp$ if $D_{0}^{\alp}[P;X]$ and $\alp$ is stable.
We see that $\calw=\bigcup\{X: D_{1}[X]\}$ is the maximal $1$-distinguished and $\Sig^{1}_{3}$-class.

In this subsection let us sketch a part of a well-foundeness proof in $\Sig^{1}_{3}\mbox{{\rm -DC+BI}}$
by pinpointing the lemma for which we need $\Sig^{1}_{3}\mbox{{\rm -DC}}$.
\\

An ordinal term $\sig$ in $OT(\mI)$ is said to be \textit{regular} if
$\psi_{\sig}^{f}(a)$ is in $OT(\mI)$ for some $f$ and $a$.
$Reg$ denotes the set of regular terms.
In this section we need the next regular ordinal above an ordinal $\alp$
in defining distinguished sets.
Although it is customarily denoted by $\alp^{+}$, it is hard to discriminate $\alp^{+}$ from
the next stable ordinal $\alp^{\dagger}$.
Therefore let us write for $\alp<\mI$,
$\alp^{+^{1}}=\min\{\sig\in SSt: \sig>\alp\}$
for the next stable ordinal $\alp^{\dagger}$, and
$\alp^{+^{0}}=\min\{\sig\in Reg: \sig>\alp\}$ for the next regular ordinal $\alp^{+}$.
Let $\alp^{+^{1}}:=\alp^{+^{0}}:=\infty$ if $\alp\geq\mI$.
Let 
$\alp^{-^{1}}:=\max\{\sig\in St_{\mI}\cup\{0\}: \sig\leq\alp\}$ when $\alp<\mI$, and
$\alp^{-^{1}}:=\mI$ if $\alp\geq\mI$.
Since $SSt \subset Reg$, we obtain $\alp^{+^{0}}\leq\alp^{+^{1}}$ and
$\bet^{+^{0}}<\sig$ if $\bet<\sig\in St$ since each $\sig\in St$ is a limit of regular ordinals.

\begin{definition}\label{df:CX}
{\rm 
$\mathcal{C}^{\alpha}(X)$ is the closure of 
$\{0,\Omega,\mathbb{I}\}\cup(X\cap\alpha)$ under $+,\vphi$,
$\{\sigma,\beta\}\cup SC_{\mI}(f)\mapsto \psi_{\sigma}^{f}(\beta)$ for $\sigma>\alpha$, and
$\rho\mapsto\mI[\rho],\rho^{\dagger}$ for $\mI[\rho],\rho^{\dagger}\geq\alp$ in $OT(\mI)$.
}
\end{definition}

\begin{definition} 
{\rm 
For $P,X\subset OT(\mI)$ and
$\gam\in OT(\mI)\cap\mI$, let
\beqnarr
W_{0}^{\alp}(P) & := & W(\mathcal{C}^{\alpha}(P))
\nonumber
\\
D_{0}^{\gam}[P;X] & :\Lrarw &
P\cap\gam^{-^{1}}=X\cap\gam^{-^{1}} \spand Wo[X\cap\gam^{-^{1}}] \spand
\label{eq:distinguishedclass}
\\
&&
\forall\alpha
\left(
\gam^{-^{1}}\leq\alpha\leq P \to W_{0}^{\alp}(P)\cap\alpha^{+^{0}}= P\cap\alpha^{+^{0}}
\right)
\nonumber
\\
\mathcal{W}_{1}^{\gam}(X) & := & \bigcup\{P\subset OT(\mI) :D_{0}^{\gam}[P;X]\}
\nonumber
\\
D_{1}[X] & :\Lrarw  &
Wo[X]  \spand
\forall\gam
\left(
\gam\leq X\to \mathcal{W}_{1}^{\gam}(X)\cap\gam^{+^{1}}= X\cap\gam^{+^{1}}
\right)
\label{df:3wfdtg.832.2}
\\
\mathcal{W}_{2} & := & \bigcup\{X\subset OT(\mI) :D_{1}[X]\}
\nonumber
\eeqnarr
A set $P$ is said to be a $0$-\textit{distinguished set} for $\gam$ and $X$ if $D_{0}^{\gam}[P;X]$,
and a set $X$ is a $1$-\textit{distinguished set} if $D_{1}[X]$.

}
\end{definition}

Observe that in $\Sig^{1}_{2}\mbox{{\rm -AC}}$,
$W^{\alp}_{0}(P)$ is $\Pi^{1}_{1}$, 
$D^{\gam}_{0}[P;X]$ is $\Del^{1}_{2}$, $\mathcal{W}_{1}^{\gam}(X)$ is $\Sig^{1}_{2}$, and
$D_{1}[X]$ is $\Del^{1}_{3}$.
Hence $\mathcal{W}_{2}$ is a $\Sig^{1}_{3}$-class.

Let $\alpha\in P$ for a $0$-distinguished set $P$ for $\gam<\mI$ and $X$. 
If $\alp<\gam^{-^{1}}$, then $\alp\in X$ with $Wo[X]$.
Otherwise
$W(\mathcal{C}^{\alpha}(P))\cap\alpha^{+^{0}}=W^{\alp}_{0}(P)\cap\alp^{+^{0}}= P\cap\alpha^{+^{0}}$ with $\alp<\alp^{+^{0}}$.
Hence $P$ is a well order.

\blem\label{lem:3wf6}$(\Sigma^{1}_{2}\mbox{{\rm -CA}})$\\
Suppose $Wo[X\cap\gam^{-^{1}}]$. Then
$\mathcal{W}^{\gam}_{1}(X)$ is the maximal $0$-distinguished {\rm set}
for $\gam$ and $X$, i.e.,
$D^{\gam}_{0}[\mathcal{W}^{\gam}_{1}(X);X]$ and 
$\exists Y(Y=\mathcal{W}^{\gam}_{1}(X))$.

\elem
\bprf
This is seen as in Proposition \ref{lem:3wf6max}.
\eprf

\blem\label{lem:8wf2}
\benu
\item\label{lem:8wf2.1}
Let $X$ and $Y$ be $1$-distinguished sets.

Then $\gam\leq X\spand \gam\leq Y
\Rarw 
X\cap\gam^{+^{1}}=Y\cap\gam^{+^{1}}$.

\item\label{lem:8wf2.4}
$\calw_{2}$ is the $1$-maximal distinguished class, i.e.,
$D_{1}[\calw_{2}]$.

\item\label{lem:8wf2.5}
For a family $\{Y_{j}\}_{j\in J}$ of $1$-distinguished sets,
the union $Y=\bigcup_{j\in J} Y_{j}$ is also a $1$-distinguished set.

\eenu
\elem

\blem\label{lem:6.34}
\benu
\item\label{lem:6.34.2}
$\mathcal{C}^{\mI}(\mathcal{W}_{2})\cap\mI=\mathcal{W}_{2}\cap\mI=
W(\mathcal{C}^{\mI}(\mathcal{W}_{2}))\cap\mI$.

\item\label{lem:6.34.3}
{\rm (BI)}
For {\rm each} $n<\ome$, 
$TI[\mathcal{C}^{\mI}(\mathcal{W}_{2})
\cap\ome_{n}(\mI+1)
]$, i.e.,
for each class $\mathcal{X}$,
$Prg[\mathcal{C}^{\mI}(\mathcal{W}_{2}),\mathcal{X}] \to 
\mathcal{C}^{\mI}(\mathcal{W}_{2})
\cap\ome_{n}(\mI+1)
\subset\mathcal{X}$.

\item\label{lem:6.34.4}
For {\rm each} $n<\ome$, 
$\mathcal{C}^{\mI}(\mathcal{W}_{2})
\cap\ome_{n}(\mI+1)
\subset 
W(\mathcal{C}^{\mI}(\mathcal{W}_{2}))$.
In particular 
$\{\mI,\ome_{n}(\mI+1)\}\subset W(\mathcal{C}^{\mI}(\mathcal{W}_{2}))$.

\eenu
\elem

As in Definition \ref{df:calg},
$\mathcal{G}^{X}:=\{\alpha\in OT(\mI) :\alpha\in \mathcal{C}^{\alpha}(X)
\spand \mathcal{C}^{\alpha}(X)\cap\alpha\subset X\}$.

\blem\label{lem:6.21pi1}$(\Sigma^{1}_{2}\mbox{{\rm -CA}})$\\
Suppose $D_{1}[Y]$ and $\alp\in\mathcal{G}^{Y}$.
Let 
$X=\mathcal{W}_{1}^{\alp}(Y)\cap\alp^{+^{1}}$.
Assume that one of the following conditions (\ref{eq:6.21.55}) and (\ref{eq:6.21.56})
is fulfilled.
Then $\alp\in X$ and $D_{1}[X]$.
In particular $\alp\in\mathcal{W}_{2}$ holds.
Moreover if $\alp^{-^{1}}\leq Y$, then $\alp\in Y$ holds.

\beqnarr
&&
\fal\bet\left(
Y\cap\alp^{+^{1}}<\bet \spand \bet^{+^{0}}<\alp^{+^{0}} \to
W_{0}^{\bet}(Y)\cap\bet^{+^{0}}\subset Y
\right)
\label{eq:6.21.55}
\\
&&
\fal\bet\geq\alp^{-^{1}}
\left(Y\cap\alp^{+^{1}}<\bet \spand \bet^{+^{0}}<\alp^{+^{0}} \to
W_{0}^{\bet}(Y)\cap\bet^{+^{0}}\subset Y
\right)
\nonumber
\\
&&
\spand \fal\bet<\alp^{-^{1}}\exi\gam(\bet<\gam^{+^{1}} \spand \gam^{-^{1}}\leq Y)
\label{eq:6.21.56}
\eeqnarr

\elem
\bprf
This is seen as in Lemma \ref{lem:6.21} by showing that
$D_{0}^{\alp}[P;Y]$, $\alp\in X$ and $D_{1}[X]$
for
$P=W_{0}^{\alp}(Y)\cap\alp^{+^{0}}=W(\mathcal{C}^{\alp}(Y))\cap\alp^{+^{0}}$.
\eprf

\blem\label{cor:6.21} 
Assume $D_{1}[Y]$, $\mI>\mS\in Y\cap(St\cup\{0\})$ and 
$\{0,\Ome\}\subset Y$.
Then $\mS^{+^{1}}=\mS^{\dagger}\in\mathcal{W}_{2}$.
\elem
\bprf
Since the condition (\ref{eq:6.21.56}) in Lemma \ref{lem:6.21}
is fulfilled with $(\mS^{+^{1}})^{-^{0}}=(\mS^{+^{1}})^{-^{1}}=\mS^{+^{1}}$ and 
$\mS^{-^{1}}=\mS$,
it suffices to show that $\mS^{+^{1}}\in\mathcal{G}^{Y}$.
Let $\alp=\mS^{+^{1}}$.
$\alp\in\mathcal{C}^{\alp}(Y)$ follows from
$\mS\in Y\cap\alp$.
Moreover
$\gam\in\mathcal{C}^{\alp}(Y)\cap\alp \Rarw \gam\in Y$ is seen by induction on 
$\ell\gam$ using the assumption $\{0,\Ome\}\subset Y$.
Therefore $\alp\in\mathcal{G}^{Y}$.
\eprf

\blem\label{lem:6.43}$(\Sig^{1}_{3}\mbox{{\rm -DC}})$\\
If 
$\alp\in\mathcal{G}^{\mathcal{W}_{2}}$, 
then there exists a $1$-distinguished {\rm set}
$Z$ such that $\{0,\Ome\}\subset Z$, 
$\alp\in\mathcal{G}^{Z}$ and
$\fal\mS\in Z\cap(St\cup\{\Ome\})[\mS^{\dagger}\in Z]$.
\elem
\bprf
Let $\alp\in\mathcal{G}^{\mathcal{W}_{2}}$.
We have $\alp\in\mathcal{C}^{\alp}(\mathcal{W}_{2})$.
Pick a $1$-distinguished set $X_{0}$ such that $\alp\in\mathcal{C}^{\alp}(X_{0})$.
We can assume $\{0,\Ome\}\subset X_{0}$.
On the other hand we have
$\mathcal{C}^{\alp}(\mathcal{W}_{2})\cap\alp\subset\mathcal{W}_{2}$
and
$\fal\mS\in \mathcal{W}_{2}\cap(St_{\mI}\cup\{\Ome\})[\mS^{\dagger}\in \mathcal{W}_{2}]$
by Lemma \ref{cor:6.21}.
We obtain 
\beqnarrs
&&
\fal n\fal X\exi Y
\{ 
D_{1}[X] \to D_{1}[Y]
\\
& \land &
\fal\bet\in OT(\mI)
\left(
\ell\bet\leq n \land \bet\in\mathcal{C}^{\alp}(X)\cap\alp \to
\bet\in Y
\right)
\\
& \land &
\fal \mS\in (St\cup\{\Ome\})
\left( 
\ell \mS\leq n \land \mS\in X \to \mS^{\dagger}\in Y
\right)
\}
\eeqnarrs

Since $D_{1}[X]$ is $\Del^{1}_{3}$, 
$\Sig^{1}_{3}\mbox{-DC}$ yields a set $Z$ such that $Z_{0}=X_{0}$ and
\beqnarrs
&&
\fal n
\{ 
D_{1}[Z_{n}] \to D_{1}[Z_{n+1}]
\\
& \land &
\fal\bet\in OT(\mI)
\left(
\ell\bet\leq n \land \bet\in\mathcal{C}^{\alp}(Z_{n})\cap\alp \to
\bet\in Z_{n+1}
\right)
\\
& \land &
\fal \mS\in (St\cup\{\Ome\})
\left( 
\ell\mS\leq n \land \mS\in Z_{n} \to \mS^{\dagger}\in Z_{n+1}
\right)
\}
\eeqnarrs
Let $Z=\bigcup_{n}Z_{n}$.
We see by induction on $n$ that $D_{1}[Z_{n}]$ for every $n$.
Lemma \ref{lem:8wf2}.\ref{lem:8wf2.5} yields $D_{1}[Z]$.
Let $\bet\in\mathcal{C}^{\alp}(Z)\cap\alp$.
Pick an $n$ such that $\bet\in\mathcal{C}^{\alp}(Z_{n})$ and $\ell\bet\leq n$.
We obtain $\bet\in Z_{n+1}\subset Z$.
Therefore $\alp\in\mathcal{G}^{Z}$.
Furthermore let $\mS\in Z\cap(St\cup\{\Ome\})$.
Pick an $n$ such that $\mS\in Z_{n}$ and $\ell\mS\leq n$.
We obtain $\mS^{\dagger}\in Z_{n+1}\subset Z$.
\eprf

\brem
{\rm Lemma \ref{lem:6.43} is a $\Sig^{1}_{4}$-statement, which is proved in
$\Sig^{1}_{3}\mbox{-DC}$.
Alternatively we could prove the lemma in $\Sig^{1}_{3}\mbox{-AC}$
if we assign fundamental sequences to limit ordinals as in \cite{J2}.
}
\erem

\end{document}